\documentclass{memo-l}

\usepackage[alphabetic]{amsrefs}
\usepackage{mathrsfs}
\usepackage{tikz-cd}
\usepackage[all]{xy}
\newdir{ >}{{}*!/-5pt/@{>}} 
\SelectTips{cm}{}

\DeclareMathOperator{\cok}{coker}
\DeclareMathOperator*{\colim}{colim}
\DeclareMathOperator{\cone}{cone}
\DeclareMathOperator{\Ext}{Ext}
\DeclareMathOperator{\hatExt}{\widehat{Ext}}
\DeclareMathOperator{\Fun}{Fun}
\DeclareMathOperator*{\hocolim}{hocolim}
\DeclareMathOperator*{\holim}{holim}
\DeclareMathOperator{\Hom}{Hom}
\DeclareMathOperator{\gm}{gm}
\DeclareMathOperator{\hm}{hm}
\DeclareMathOperator{\id}{id}
\DeclareMathOperator{\im}{im}
\DeclareMathOperator{\In}{in}
\DeclareMathOperator{\Mod}{Mod}
\DeclareMathOperator{\pr}{pr}
\DeclareMathOperator{\res}{res}
\DeclareMathOperator*{\Rlim}{Rlim}
\DeclareMathOperator{\Tel}{Tel}
\DeclareMathOperator{\THH}{THH}
\DeclareMathOperator{\TP}{TP}

\newcommand{\ev}{\mathrm{ev}}
\newcommand{\ind}{\mathrm{Ind}^{\Gamma}_k}
\newcommand{\coind}{\mathrm{Coind}^{\Gamma}_k}

\newcommand{\wEG}{\widetilde{EG}}

\newcommand{\bC}{\mathbb{C}}
\newcommand{\bF}{\mathbb{F}}
\newcommand{\bH}{\mathbb{H}}
\newcommand{\bN}{\mathbb{N}}
\newcommand{\bR}{\mathbb{R}}
\newcommand{\bS}{\mathbb{S}}
\newcommand{\bT}{\mathbb{T}}
\newcommand{\bU}{\mathbb{U}}
\newcommand{\bZ}{\mathbb{Z}}

\newcommand{\cA}{\mathcal{A}}
\newcommand{\cI}{\mathcal{I}}

\newcommand{\sE}{\mathscr{E}}
\newcommand{\sU}{\mathscr{U}}

\newcommand{\ds}{\displaystyle}
\newcommand{\op}{\mathrm{op}}
\newcommand{\och}{\quad \text{and} \quad}
\newcommand{\qqandqq}{\qquad\text{and}\qquad}
\newcommand{\longto}{\longrightarrow}
\renewcommand{\:}{\colon}

\newcommand{\bsm}{\left[\begin{smallmatrix}}
\newcommand{\esm}{\end{smallmatrix}\right]}

\newtheorem{theorem}{Theorem}[chapter]
\newtheorem{lemma}[theorem]{Lemma}
\newtheorem{proposition}[theorem]{Proposition}
\newtheorem{corollary}[theorem]{Corollary}

\theoremstyle{definition}
\newtheorem{definition}[theorem]{Definition}
\newtheorem{notation}[theorem]{Notation}
\newtheorem{construction}[theorem]{Construction}
\newtheorem{example}[theorem]{Example}

\theoremstyle{remark}
\newtheorem{remark}[theorem]{Remark}

\numberwithin{section}{chapter}
\numberwithin{equation}{chapter}

\makeindex

\begin{document}

\frontmatter

\title[A multiplicative Tate spectral sequence]{A Multiplicative Tate Spectral Sequence \\ for Compact Lie Group Actions}


\author{Alice Hedenlund}
\address{Department of Mathematics, University of Oslo, Norway}
\curraddr{}
\email{alice.hedenlund@gmail.com}
\thanks{}

\author{John Rognes}
\address{Department of Mathematics, University of Oslo, Norway}
\curraddr{}
\email{rognes@math.uio.no}
\thanks{}

\date{}




\begin{abstract}
Given a compact Lie group~$G$ and a commutative orthogonal ring
spectrum~$R$ such that $R[G]_* = \pi_*(R \wedge G_+)$ is finitely
generated and projective over $\pi_*(R)$, we construct a multiplicative
$G$-Tate spectral sequence for each $R$-module $X$ in orthogonal
$G$-spectra, with $E^2$-page given by the Hopf algebra Tate cohomology of
$R[G]_*$ with coefficients in $\pi_*(X)$.  Under mild hypotheses, such
as $X$ being bounded below and the derived page $RE^\infty$ vanishing,
this spectral sequence converges strongly to the homotopy $\pi_*(X^{tG})$
of the $G$-Tate construction $X^{tG} = [\wEG \wedge F(EG_+, X)]^G$.
\end{abstract}

\maketitle

\tableofcontents


\mainmatter
%
%
%

\chapter{Introduction}

This memoir grew out of an attempt to spell out the details for the Tate spectral sequence for the circle group $\bT$. The construction of a multiplicative Tate spectral sequence for finite groups has been around for a while now: the first construction, due to Greenlees--May, can be found in~\cite{GM95}, and another one, due to Hesselholt--Madsen, which makes the multiplicative properties of the spectral sequence more transparent, can be found in~\cite{HM03}. However, while multiplicativity of the $\bT$-Tate spectral sequences has been used in computations, the authors of this memoir have found references discussing the details for how such a spectral sequence is constructed surprisingly lacking. We hope that this memoir will fill that gap in the literature.  

The authors' motivation for considering the $\bT$-Tate spectral sequence  comes from the study of topological Hochschild homology and its refinements, such as topological cyclic homology. Given an $E_1$-ring spectrum $B$, the topological Hoch\-schild homology\index{topological Hochschild homology} $\THH(B)$, first defined in the unpublished manuscript~\cite{Bok85}, is a genuine~$\bT$-equivariant spectrum. The study of the Tate construction on this spectrum using the entire circle action goes back to~\cite{BM93} and~\cite{AR02}, and was put in the spotlight by Hesselholt in~\cite{Hes18} under the name of periodic topological cyclic homology\index{periodic topological cyclic homology}:
\[
\TP(B) = \THH(B)^{t\bT}.
\]
Recently, Bhatt--Morrow--Scholze showed that there is a tight connection between periodic topological cyclic homology and crystalline cohomology\index{crystalline cohomology}~\cite{BMS2}. 

\section*{Background and aim}

Classically, Tate cohomology is a way to combine group homology and
group cohomology into a single multiplicative cohomology theory, and was first introduced by Tate in his study of class field theory~\cite{tate}. We sketch the main ideas involved following~\cite{CE56}*{Section XII.3} and~\cite{Brown}. Given a finite group~$G$, the main observation of Tate cohomology is this: if we dualise a projective resolution of $\bZ$ as a trivial module over~$\bZ[G]$, we end up with a `coresolution' of $\bZ$ by projective~$\bZ[G]$-modules. This `coresolution' $\Hom_{\bZ}(P_*,\bZ)$ can be spliced with the original projective resolution~$P_*$, and we so obtain a bi-infinite resolution~$\hat{P}_*$ of $\bZ$ called a complete resolution\index{complete resolution}. Tate cohomology\index{Tate cohomology!of a finite group} of $G$ with coefficients in a~$G$-module~$M$ is defined as 
\[
\hat{H}^n(G,M) = H^n(\Hom_{G}(\hat{P}_*,M)) \,.
\]
The Tate construction in the category of $G$-spectra can be seen as a generalisation of Tate cohomology in the context of higher algebra. Given a compact Lie group $G$ and orthogonal $G$-spectrum $X$, we define the $G$-homotopy orbits and~$G$-homotopy fixed points of $X$ as 
\[
X_{hG} = EG_+ \wedge_G X \och X^{hG} = F(EG_+,X)^G \,,
\]
respectively\index{homotopy orbits}\index{homotopy fixed points}. Here $EG$ denotes a free contractible $G$-space. These can be regarded as generalisations of group homology and group cohomology. Indeed, if $G$ is a finite group and $X = HM$ is the Eilenberg--Mac Lane spectrum on the $G$-module~$M$, then the homotopy groups of the $G$-homotopy orbits and $G$-homotopy fixed points of $HM$ recover group homology and group cohomology of $G$ with coefficients in~$M$, respectively. Following Greenlees~\cite{Gre87,GM95}, we define the $G$-Tate construction\index{Tate construction} on $X$ as the $G$-fixed point spectrum
\[
X^{tG} = \left(\wEG \wedge F(EG_+, X)\right)^G
\]
with respect to the diagonal $G$-action. Here, $\wEG$ denotes the mapping cone of the collapse map $c \: EG_+ \to S^0$. This is a generalisation of Tate cohomology in the sense that the homotopy groups of the Tate construction on $HM$ for a $G$-module~$M$ recover the Tate cohomology groups of the finite group $G$ with coefficients in~$M$.

One important property of the Tate construction is that it is multiplicative\index{Tate construction!multiplicative} in the sense that any pairing $X \wedge Y \to Z$ of orthogonal $G$-spectra gives rise to a pairing $X^{tG} \wedge Y^{tG} \longto Z^{tG}$ of their Tate constructions. This relies on the existence of $G$-maps $EG_+ \to EG_+ \wedge EG_+$ and $\wEG \wedge \wEG \to \wEG$. It is well-known that the diagonal map $EG_+ \to EG_+ \wedge EG_+$ induces a pairing  
\[
X^{hG} \wedge Y^{hG} \longto Z^{hG} \,,
\]
making the $G$-homotopy fixed points construction a lax symmetric monoidal functor. The inclusion $S^0 \to \wEG$ and the canonical identifications $S^0 \wedge \wEG \cong \wEG \cong \wEG \wedge S^0$ induce a natural map
\[
X^{hG} \longto X^{tG}
\]
and pairings $X^{hG} \wedge Y^{tG} \to Z^{tG}$ and $X^{tG} \wedge Y^{hG} \to Z^{tG}$. There is a $G$-map $N \: \wEG \wedge \wEG \to \wEG$ extending the canonical identifications, and any two such extensions are homotopic. Any choice of extension then induces a pairing
\[
X^{tG} \wedge Y^{tG} \longto Z^{tG} 
\]
compatible with the above-mentioned map and pairings\footnote{Work by Nikolaus--Scholze shows that this multiplicative structure is actually unique, in a homotopy theoretical sense; see~\cite{NS18}*{Theorem I.3.1}. This will not be important for our work, though.}.  In general, the extension~$N$ will only be commutative and associative up to (coherent) homotopy, so~$X \mapsto X^{tG}$ is not a lax symmetric monoidal functor to the category of orthogonal spectra, but only satisfies a homotopy coherent version of this property, which could be made precise using operad actions.  For our purposes it suffices to note that it is lax symmetric monoidal as a functor to the stable homotopy category.

Given an orthogonal $G$-spectrum $X$, the aim of the present memoir is to construct a $G$-Tate spectral sequence
\[
\hat E^r_{s,t}(X)
	\Longrightarrow \pi_{s+t}(X^{tG}) \,,
\]
with an algebraically specified $E^2$-page, converging, in some suitable sense, to the homotopy groups of the $G$-Tate construction on
$X$. Moreover, we would like this spectral sequence to be multiplicative,
in the sense that a pairing $X \wedge Y \to Z$ of orthogonal $G$-spectra
should induce a pairing 
\[
(\hat E^r(X) , \hat E^r(Y)) \longto \hat E^r(Z)
\]
of $G$-Tate spectral sequences. Finally, we want the pairing
of $E^\infty$-pages to be compatible with the pairing 
\[
\pi_*(X^{tG}) \otimes \pi_*(Y^{tG}) \longto \pi_*(Z^{tG})
\]
of abutments. In particular, if $X$ is an orthogonal $G$-ring spectrum, then the $G$-Tate spectral sequence of~$X$ should be an algebra spectral sequence converging multiplicatively to $\pi_*(X^{tG})$. As already mentioned, how to construct such spectral sequences is well-known in the situation of $G$ being a finite group. Our goal is to generalise this to higher dimensional compact Lie groups.

\section*{Main results}

Let us start by describing roughly, without going into too much detail, what we will do in this memoir. We will carry out the construction of multiplicative and conditionally convergent Tate spectral sequences for compact Lie groups $G$ such that~$\bS[G]_* = \pi_*(\bS[G])$ is finitely generated projective as a module over $\bS_* = \pi_*(\bS)$. Here $\bS$ denotes the sphere spectrum and 
\[
\bS[G] = \bS \wedge G_+
\]
is the unreduced suspension spectrum of~$G$. Under these assumptions, $\bS[G]_*$ is a finitely generated projective and cocommutative Hopf algebra over~$\bS_*$, and we will show that we have access to a multiplicative~$G$-Tate spectral sequence with $E^2$-page given by the complete $\Ext$-groups\index{complete $\Ext$}
\[
\hat E^2_{s,*}(X) = \widehat{\Ext}^{-s}_{\bS[G]_*}(\bS_*, \pi_*(X))
\]
of $\bS_*$ over $\bS[G]_*$ with coefficients in the $\bS[G]_*$-module $\pi_*(X)$. The multiplicative structure in complete $\Ext$ is given by a graded commutative and associative cup product, and this will serve as a substitute for the failure of $X \mapsto X^{tG}$ to be lax symmetric monoidal.  This spectral sequence will be strongly convergent under mild hypotheses, such as for instance in the case when the derived $E^\infty$-page $RE^\infty$ vanishes and the spectrum $X$ is bounded below.

We note that this generality includes the case where $G = \bT$ is the circle group, our main interest, but does not cover cases such as $G = SO(3)$. We therefore broaden our scope by considering a commutative\footnote{For somewhat technical reasons, it is not sufficient for
us to assume that $R$ is homotopy commutative.  We analyse
the product in the filtered $R$-module $G$-spectrum
\[
\wEG \wedge F(EG_+, R \wedge X) \cong L \wedge_R M \,,
\]
with $L = R \wedge \wEG$ and $M = F(EG_+, R \wedge X)$, as
a composition
\[
L \wedge_R M \wedge_R L \wedge_R M
	\overset{(23)}\longto L \wedge_R L \wedge_R M \wedge_R M
	\overset{\phi\wedge\psi}\longto L \wedge_R M
\]
for filtered products $\phi \: L \wedge_R L \to L$ and $\psi \: M
\wedge_R M \to M$. Homotopy commutativity is not sufficient to ensure
that the twist map $\tau \: M \wedge_R L \to L \wedge_R M$ implicit in
the definition of~$(23)$ is an $R$-$R$-bimodule map.} `ground' orthogonal ring spectrum~$R$ and a compact Lie group $G$ such that $R[G]_* = \pi_*(R[G])$ is finitely generated and projective over $R_* = \pi_*(R)$, where 
\[
R[G] = R \wedge G_+ \,.
\]
This then includes cases such as $R = \bS[1/2]$ and $R = H\bF_2$,
with $G = SO(3)$.
(We have not classified the pairs $(R, G)$ for which this condition
holds, but it is easy to see that if $H_*(G; Z)$ is finitely generated
and free over a PID $Z$ with a ring homomorphism to $\pi_0(R)$, and
the Atiyah--Hirzebruch spectral sequence $E^2_{*,*} = H_*(G; R_*)
\Longrightarrow R[G]_*$ collapses at~$E^2$, then $R[G]_*$ is finitely
generated and free over~$R_*$.  This includes all cases where $G$ is
topologically a product of spheres.)
Given an $R$-module $X$ in orthogonal $G$-spectra we shall construct a multiplicative $G$-Tate spectral sequence
\[
\hat E^2_{s,*}(X) = \widehat{\Ext}^{-s}_{R[G]_*}(R_*, \pi_*(X))
	\Longrightarrow \pi_{s+*}(X^{tG})
\]
where the $E^2$-page is now given as complete $\Ext$ of $R_*$ over $R[G]_*$ with coefficients in $\pi_*(X)$. This will be strongly convergent under the same conditions as before.

\subsection*{Tate cohomology of Hopf algebras}
\index{Tate cohomology!of a Hopf algebra|(}
In Chapter~\ref{sec:TateforHopf} we develop a theory of Tate cohomology of a finitely generated and projective Hopf algebra $\Gamma$ over a (possibly graded) commutative ring $k$, with the aim being to algebraically describe the $E^2$-page of a suitable Tate spectral sequence. Our approach will be different from the complete resolution approach, and we instead rely on the so-called Tate complex\index{Tate complex}. Given a projective $\Gamma$-resolution $P_*$ of $k$, we will denote the mapping cone of the augmentation map $\epsilon \: P_* \to k$ as $\widetilde{P}_*$. The Tate complex\index{Tate complex} of a $\Gamma$-module $M$, first defined in~\cite{Gre95}, is the~$\Gamma$-chain complex 
\[
\hm_*(M) = \widetilde{P}_* \otimes_k \Hom_k(P_*,M)
\]  
where $\Gamma$ acts diagonally on the tensor product and by conjugation on $\Hom(P_*,M)$. In the aforementioned paper, the author shows that in the classical case, meaning $k=\bZ$ and $\Gamma = \bZ[G]$ for a finite group $G$, there is a zigzag of maps
\[
\xymatrix{
\widetilde{P}_* \otimes_k \Hom_k(P_*,M) \ar[r] & \widetilde{P}_* \otimes_k \Hom_k(\hat{P}_*,M) & \ar[l] \Hom_k(\hat{P}_*,M)
} 
\]
which become quasi-isomorphisms after taking~$G$-invariants. The conclusion is that Tate cohomology can also be computed as the (co)homology groups of the~$G$-invariants of the Tate complex. Recall that $\hat{P}_*$ denoted a complete resolution. We show that a similar result holds true in our setting: under the assumption that $\Gamma$ is a finitely generated and projective Hopf algebra over~$k$, the homology of the~$\Gamma$-invariants of $\hm_*(M)$, which we can reasonably refer to as the Tate cohomology of~$\Gamma$ with coefficients in~$M$, is isomorphic to the complete $\Ext$ of $k$ over $\Gamma$ with coefficients in~$M$. 

\begin{theorem}\index{complete $\Ext$}
If $\Gamma$ is a finitely generated projective and cocommutative Hopf algebra over $k$, then
\[
\hatExt_{\Gamma}^n(k,M) \cong H_{-n}(\Hom_\Gamma(k,\hm_*(M)))  \,.
\]
\end{theorem}

The above result, which in the text corresponds to Theorem~\ref{prop:hatExtdeltafunctor} and Remark~\ref{rmk:CK}, relies crucially on a result by Pareigis which exhibits the $k$-dual of a Hopf algebra~$\Gamma$ as an induced $\Gamma$-module.

\begin{theorem}[Pareigis]
Let $\Gamma$ be a finitely generated projective Hopf algebra over~$k$.  Then there is an
isomorphism
\[
\Hom_k(\Gamma,k) \cong \ind P(\Hom_k(\Gamma,k))
\]
of right $\Gamma$-modules, where $P(\Hom_k(\Gamma,k))$ is a finitely generated projective $k$-module of constant rank~$1$, given as the primitives for the right $\Gamma$-coaction on $\Hom_k(\Gamma,k)$.
\end{theorem}

Our main reason for working primarily with Tate complexes, as opposed to complete resolutions, has to do with multiplicative structures\index{cup product!on complete $\Ext$}. Recall that the cup product\index{cup product!on $\Ext$} 
\[
\smile \: \Ext^*_\Gamma(k,M) \otimes_k \Ext^*_\Gamma(k,N) \longto \Ext^*_\Gamma(k,M \otimes_k N)
\]
relies on the existence of a $\Gamma$-linear chain map $\Psi \: P_* \to P_* \otimes_k P_*$ covering the identity map $\id: k \to k \otimes_k k$. Such a chain map exists and is unique up to chain homotopy, by elementary homological algebra. One can extend this cup product to a product on Hopf algebra Tate cohomology by the existence of a $\Gamma$-linear chain map $\Phi \: \widetilde{P}_* \otimes_k \widetilde{P}_* \to \widetilde{P}_*$ extending the fold map $\widetilde{P}_* \oplus_k \widetilde{P}_* \to \widetilde{P}_*$. For $\Gamma$-modules $M$ and~$N$ the composite pairing
\begin{align*}
\widetilde{P}_* \otimes_k \Hom_k(P_*, M)
	\otimes_k \widetilde{P}_* \otimes_k &\Hom_k(P_*, N) \\
        & \overset{1 \otimes \tau \otimes 1}\longto
\widetilde{P}_* \otimes_k \widetilde{P}_*
	\otimes_k \Hom_k(P_*,M) \otimes \Hom_k(P_*, N) \\
	    & \overset{1 \otimes 1 \otimes \alpha}\longto \widetilde{P}_* \otimes_k \widetilde{P}_*
	\otimes_k \Hom_k(P_* \otimes_k P_*,M \otimes_k N) \\
        & \overset{\Phi \otimes \Psi^*}\longto
\widetilde{P}_* \otimes_k \Hom_k(P_*, M \otimes_k N)
\end{align*}
is $\Gamma$-linear, and it induces an associative, unital, and graded commutative pairing 
\[
\smile \: \hatExt^*_\Gamma(k,M) \otimes_k \hatExt^*_\Gamma(k,N) \longto \hatExt^*_\Gamma(k,M \otimes_k N)
\]
after passing to homology, which we refer to as the cup product on Tate cohomology. This extends the cup product on ordinary $\Ext$, in a suitable sense. See Proposition~\ref{prop:cupproductcompatible}.

Finally, in Section~\ref{subsec:computation}, we do a full computation of the Tate cohomology, together with the cup product, of the Hopf algebra
\[
\Gamma = k[s]/(s^2 = \eta s) \,, \quad |s|=1 \,,
\]
where $s$ is a primitive element and $k$ is a graded commutative ring with an element~$\eta$ in internal degree $1$ satisfying $2\eta = 0$. This has relevance in the situation $G=\bT$, which is our main case of interest. Indeed, we have 
\[
\pi_*(\bS[\bT]) \cong \pi_*(\bS)[s]/(s^2 = \eta s)
\]
where $\eta$ is the image of the complex Hopf map in $\pi_1(\bS) \cong \bZ/2$. See Proposition~\ref{prop:s2isetas}. The conclusion of the computation is the following theorem, which in the text is Theorem~\ref{thm:completeExtcup} and Remark~\ref{rmk:DGA}.

\begin{theorem}
Tate cohomology of $\Gamma = k[s]/(s^2 = \eta s)$ with coefficients in the~$\Gamma$-module $M$ is isomorphic to the homology of the differential graded $\Gamma$-module 
\[
M[t,t^{-1}]
\]
with differential  
\[
d(m) = tms \och d(t) = t^2 \eta \,,
\]
where $m$ is an element of $M$ and $t$ has homological degree~$-1$, internal degree~$|t| = -1$ and total degree~$\|t\| = -2$. If $\mu \: M \otimes N \to L$ is a pairing of $\Gamma$-modules, then the cup product 
\[
\smile \: \hatExt^{c_1}_\Gamma(k, M) \otimes \hatExt^{c_2}_{\Gamma}(k,N) \longto \hatExt^{c_1+c_2}_\Gamma(k,M \otimes N) \longto \hatExt^{c_1+c_2}_\Gamma(k,L)
\]
is precisely the one induced by the obvious pairing 
\[
M[t,t^{-1}] \otimes N[t,t^{-1}] \longto L[t,t^{-1}]
\]
on homology. 
\end{theorem}
\index{Tate cohomology!of a Hopf algebra|)}

\subsection*{Sequences of spectra and spectral sequences}

The main difficulty of the memoir lies in verifying that there is a construction of the Tate spectral sequence that is multiplicative. To deal with multiplicative structures on spectral sequences we have decided to employ Cartan--Eilenberg systems. These are mathematical gadgets, first introduced in~\cite{CE56}, which determine a spectral sequence. For us, the advantage is that there is a useful notion of pairings of Cartan--Eilenberg systems, and that one can prove that a pairing of Cartan--Eilenberg systems gives rise to a pairing of the associated spectral sequences. Our contribution is a detailed and explicit proof that a pairing of sequences of orthogonal $G$-spectra gives rise to a pairing of Cartan--Eilenberg systems\index{pairing!of Cartan--Eilenberg systems}. Here, \emph{sequence}\index{sequence} simply means a sequential diagram
\[
\cdots \longto X_{i-1} \longto X_i \longto X_{i+1} \longto \cdots 
\]
of maps of orthogonal $G$-spectra, and \emph{pairing}\index{pairing!of sequences} $\phi \: (X_\star, Y_\star) \to Z_\star$ refers to a collection of $G$-maps
\[
\phi_{i,j} \: X_i \wedge Y_j \longto Z_{i+j}
\]
for all integers~$i$ and~$j$, making the squares
\[
\begin{tikzcd}[column sep = large]
X_{i-1} \wedge Y_j \arrow[r,"\phi_{i-1,j}"] \arrow[d] & Z_{i+j-1} \arrow[d]
	& X_i \wedge Y_{j-1} \arrow[l,"\phi_{i,j-1}"'] \ar[d] \\
X_i \wedge Y_j \arrow[r,"\phi_{i,j}"] & Z_{i,j}
	& X_i \wedge Y_j \arrow[l,"\phi_{i,j}"']
\end{tikzcd}
\]
commute strictly. It is well-known that a sequence of orthogonal $G$-spectra gives rise to an unrolled exact couple on equivariant homotopy groups, which in turn gives rise to a spectral sequence. That a pairing of sequences gives rise of a pairing of the corresponding spectral sequences can also reasonably be regarded as folklore, but as the authors feel that an explicit reference for this is not available at the time of writing, we have decided to give a complete proof of this fact.

For homotopical control in the proofs, some sort of `cofibrant replacement' of the sequence $X_\star$ is needed. In this memoir we have chosen to use the classical telescope construction\index{telescope construction} to deal with these sorts of issues. See Section~\ref{subsec:filt}. Our main reason for this is that these `cofibrant replacements' behave well with respect to monoidal properties. This allows us to always approximate a sequence $X_\star$ with an equivalent sequence $T_\star(X)$ in a way that will make our analysis of multiplicative structures more manageable. 

The main result of Chapter 4 of the memoir is the following, which in the text corresponds to Theorem~\ref{thm:mainmultII}.

\begin{theorem}\index{spectral sequence!associated to a sequence}\index{pairing!of spectral sequences}
A pairing $\phi \: (X_\star , Y_{\star}) \to Z_{\star}$ of sequences of orthogonal $G$-spectra gives rise to a pairing $\phi \: (E^*(X_\star),E^*(Y_\star)) \to E^*(Z_\star)$. Explicitly, we have access to a collection of homomorphisms
\[
\phi^r \: E^r(X_\star) \otimes E^r(Y_\star) \longto E^r(Z_\star)
\]
for all $r \geq 1$, such that:
\begin{enumerate}
\item The Leibniz rule\index{Leibniz rule}
\[
d^r \phi^r = \phi^r(d^r \otimes 1) + \phi^r(1 \otimes d^r)
\]
holds as an equality of homomorphisms
$E^r_i(X_\star) \otimes E^r_j(Y_\star) \longto E^r_{i+j-r}(Z_\star)$
for all $i,j \in \bZ$ and $r\ge1$.
\item The diagram 
\[
\begin{tikzcd}
E^{r+1}(X_\star) \otimes E^{r+1}(Y_\star) \arrow[r,"\phi^{r+1}"] \arrow[d] & E^{r+1}(Z_\star) \arrow[d,"\cong"] \\
H(E^{r}(X_\star) \otimes E^{r}(Y_\star)) \arrow[r,"H(\phi^r)"] & H(E^{r}(Z_\star))
\end{tikzcd}
\]
commutes for all $r \geq 1$.
\end{enumerate} 

Moreover, the induced pairing $\phi_*$ on filtered abutments is compatible with the pairing $\phi^\infty$ of $E^\infty$-pages in the sense of Proposition~\ref{prop:multabutment}. Explicitly, the diagram
\[
\xymatrix{
\ds \frac{F_i A_\infty(X_\star)}{F_{i-1} A_\infty(X_\star)}
\otimes \frac{F_j A_\infty(Y_\star)}{F_{j-1} A_\infty(Y_\star)}
\ar[r]^-{\bar\phi_*} \ar[d]_-{\beta \otimes \beta}
& \ds \frac{F_{i+j} A_\infty(Z_\star)}{F_{i+j-1} A_\infty(Z_\star)}
	\ar@{ >->}[d]^-{\beta} \\
E^\infty_i(X_\star) \otimes E^\infty_j(Y_\star) \ar[r]^-{\phi^\infty}
& E^\infty_{i+j}(Z_\star)
}
\]
commutes, for all $i,j \in \bZ$. Here the abutments are given as 
\begin{align*}
A_\infty(X_\star) &\cong \pi^G_* \Tel(X_\star) \\
A_\infty(Y_\star) &\cong \pi^G_* \Tel(Y_\star) \\
A_\infty(Z_\star) &\cong \pi^G_* \Tel(Z_\star)
\end{align*}
with filtrations by the images
\begin{align*}
F_i A_\infty(X_\star) &= \im(\pi^G_*(X_i) \longto A_\infty(X_\star)) \\
F_j A_\infty(Y_\star) &= \im(\pi^G_*(Y_j) \longto A_\infty(Y_\star)) \\
F_k A_\infty(Z_\star) &= \im(\pi^G_*(Z_k) \longto A_\infty(Z_\star)) \,,
\end{align*}
respectively.
\end{theorem}

\subsection*{The $G$-Tate spectral sequence}

Given an $R$-module $X$ in orthogonal $G$-spectra, there are a number of ways of constructing Tate spectral sequences additively; as mentioned, the difficulty lies in establishing multiplicative properties of the constructions. The standard way of constructing a Tate spectral sequence seems to be by filtering the Tate construction 
\[
X^{tG} = \left(\wEG \wedge F(EG_+, X)\right)^G \simeq \left((R \wedge \wEG) \wedge_R F_R(R \wedge EG_+, X)\right)^G
\]
by filtering $\wEG$, in some suitable sense, dualising this filtration, and splicing, in analogy with the construction of complete resolutions by dualising and splicing projective resolutions\index{Greenlees filtration}. This is far from ideal if one aims to prove any multiplicative properties of the Tate spectral sequence. We will instead prove multiplicativity of the Tate spectral sequence using a construction along the lines of~\cite{HM03}. In this construction, we filter $F(EG_+,X)$ and $\widetilde{EG}$ separately, and totalise to get a filtration on the Tate construction.
In the key case $R = \bS$ and $G = \bT$, essentially the same construction was
considered by Blumberg and Mandell in their preprint~\cite{BM17}*{Section~3}. See Remark~\ref{rmk:Blumberg-Mandell}.

In more detail, we proceed as follows in Chapter 6. We start by giving the free~$G$-space~$EG$ the simplicial skeletal filtration $F_\star EG$ coming from the construction of~$EG$ using the simplicial bar construction. This induces a filtration 
\[
E_\star = R \wedge F_\star EG_+
\]
on $R \wedge EG_+$, which in turn induces a filtration
\[
M_\star(X) = F_R(E/E_{-\star-1}, X)
\]
on $F_R(R \wedge EG_+,X)$, and a filtration 
\[
\widetilde{E}_\star = \cone(E_{\star -1} \longto R)
\]
on $R \wedge \widetilde{EG}$. The convolution filtration 
\[
HM_\star(X) = (\widetilde{E} \wedge T(M(X)))_\star = \colim_{i + j \leq \star} \widetilde{E}_i \wedge_R T(M(X))_j
\]
is referred to as the \emph{Hesselholt--Madsen filtration}\index{Hesselholt--Madsen filtration}. For homotopical control we have `cofibrantly replaced' the filtration $M_\star(X)$ with its telescopic approximation~$T_\star(M(X))$. Under our projectivity assumptions, we show that the $E^1$-page of the spectral sequence arising from the Hesselholt--Madsen filtration is given by 
\[
\hat{E}^1_{c,*} \cong \Hom_{R[G]_*}(R_*,\hm_c(\pi_*(X))) \,,
\]
so that the $E^2$-page is given as the Hopf algebra Tate cohomology groups 
\[
\hat{E}^2_{c,*} \cong \hatExt^{-c}_{R[G]_*}(R_*,\pi_*(X)) \,,
\]
as defined in Chapter 2. See Proposition~\ref{prop:TateE1isTatecx} and Theorem~\ref{thm:TateE2isTateExt}. We note that the Hesselholt--Madsen~$G$-Tate spectral sequence is not obviously conditionally convergent, so for convergence issues we need to do some additional work.
(In the key case $R = \bS$, $G = \bT$, Blumberg and Mandell establish
conditional convergence in~\cite{BM17}*{Lemma~3.16}.)\index{Hesselholt--Madsen--Tate spectral sequence!conditional convergence}

The existence of a multiplicative structure on the Hesselholt--Madsen $G$-Tate spectral sequence relies on the existence of filtration-preserving maps 
\[
EG_+ \longto EG_+ \wedge EG_+ \och \widetilde{EG} \wedge \widetilde{EG} \longto \widetilde{EG} \,.
\]
The first is known to exist, and we prove by obstruction theory that the second one exists under the assumption that $R[G]_*$ is projective over $R_*$. See Proposition~\ref{prop:nablaextendstoN}. This guarantees that a pairing $X \wedge_R Y \to Z$ of~$R$-modules in orthogonal spectra induces a pairing
\[
(HM_\star(X) , HM_\star(Y)) \longto HM_\star(Z)
\]
of the corresponding Hesselholt--Madsen filtrations. The work done in Chapter 4 then guarantees that the $G$-Tate spectral sequence constructed from the Hesselholt--Madsen filtration has a multiplicative structure. Moreover, we show that the multiplicative structure on the $E^2$-page agrees with the one given by cup product on Tate cohomology. See Theorem~\ref{thm:GTatespseqismult} and Theorem~\ref{thm:E2Tatepairingiscup}.

To settle questions about convergence we compare the Hesselholt--Madsen filtration to another possible filtration of the Tate construction. The filtration we are referring to is the filtration $GM_\star(X)$ given in each degree as  
\[
GM_k(X) = \begin{cases}
\widetilde{E}_k \wedge_R T_0(M(X)) & \text{for $k\ge0$,} \\
\widetilde{E}_0 \wedge_R T_k(M(X)) & \text{for $k\le0$.}
\end{cases}
\]
Here, the structure maps $GM_{k-1}(X) \to GM_k(X)$ for $k\ge1$ are
induced by the maps~$\widetilde{E}_{k-1} \to \widetilde{E}_k$ in the filtration $\widetilde{E}_\star$, while the maps for $k\le 0$ are those of~$T_\star(M(X))$. This filtration is referred to as the \emph{Greenlees--May filtration}\index{Greenlees--May filtration}. It is straight-forward to show that the spectral sequence arising from the Greenlees--May filtration is conditionally convergent; see Lemma~\ref{lem:checkEXcondconv}. Moreover, in Lemma~\ref{lem:alphamapfilt} we show that there is a map of filtrations 
\[
\alpha \: GM_\star(X) \longto HM_\star(X) \,,
\]
which induces an isomorphism of spectral sequences from the $E^2$-page and on. See Proposition~\ref{prop:GMHMfiltrationiso}. We can then deduce convergence results for the Hesselholt--Madsen $G$-Tate spectral sequence in certain favourable situations, such as in the case when the spectrum $X$ is bounded below and the derived $E^\infty$-page $RE^\infty$ vanishing. In particular, we have the following result, which in the text corresponds to Theorem~\ref{thm:strongimpliesstrong}. 

\begin{theorem}\index{Hesselholt--Madsen--Tate spectral sequence!strong convergence}
If the Greenlees--May $G$-Tate spectral sequence for $X$ is strongly convergent, then so is the Hesselholt--Madsen $G$-Tate spectral sequence for $X$. 
\end{theorem}

\section*{Organisation}

Let us discuss the various chapters contained in this memoir, and how they relate to one another. 

\begin{description}
\item[Chapter 2] In this chapter we develop a theory of Tate cohomology for finitely generated projective Hopf algebras, with a view toward being able to satisfactorily describe the $E^2$-page of a $G$-Tate spectral sequence for compact Lie groups. 
\item[Chapter 3] In this chapter we do a quick review of orthogonal $G$-spectra. Most of this chapter can be regarded as well-known to people working in genuine equivariant stable homotopy theory. However, we want to highlight Proposition~\ref{prop:omegaXiso}, for which we have not found a reference, and which will be important in later parts of the memoir.
\item[Chapter 4] In this chapter we discuss sequences of orthogonal $G$-spectra, Cartan--Eilenberg systems, and spectral sequences, with a special focus on multiplicative structures. This chapter may well be read separately from the rest of the memoir, possibly in addition to Section 3.1, which contains a quick recap on orthogonal $G$-spectra. We hope it can be of use as a reference for multiplicative structures on spectral sequences coming from sequences of spectra.
\item[Chapter 5] In this chapter we discuss the $G$-homotopy fixed point spectral sequence for an orthogonal $G$-spectrum. This is meant as a warm-up to the~$G$-Tate spectral sequence, but can absolutely be read in its own right.  
\item[Chapter 6] In this chapter we discuss various constructions of the $G$-Tate spectral sequence of an orthogonal $G$-spectrum. The reader who only cares for the $\bT$-Tate spectral sequence will find a summary of the relevant results at the very end of the memoir, in Section~\ref{subsec:TTatesummary}.
\end{description}

\section*{Acknowledgements}

We thank the anonymous referee for constructive comments on an earlier version of this memoir.

%
%
%

\chapter{Tate Cohomology for Hopf Algebras}
\label{sec:TateforHopf}

The algebraic objects that we are led to work with when constructing the Tate spectral sequence are Hopf algebras and chain complexes of modules over these. The topological context we will discuss later in the memoir allows for Hopf algebras over fairly complicated rings, which forces us to work in the generality of  Hopf algebras over arbitrary, possibly graded, commutative rings. We give a brief account of this in Section 2.1 and Section 2.2, referring e.g.~to~\cite{Mac95}*{Chapter~VI} for a fuller discussion. We go on to give a suitable definition of Tate cohomology of Hopf algebras via the so-called Tate complex in Section 2.3. In Section 2.4 we relate this definition to the ordinary definition of Tate cohomology in terms of complete resolutions. In particular, we show in Theorem~\ref{prop:hatExtdeltafunctor} that our definition agrees with what is traditionally referred to as Tate cohomology or complete $\Ext$, in the case when our Hopf algebra $\Gamma$ is finitely generated and projective over its base ring $k$. The crucial point that allows us to do this is a result of Pareigis, which in particular forces the~$k$-dual of $\Gamma$ to be finitely generated and projective over $\Gamma$, under the same hypotheses. We discuss the multiplicative structure of Tate cohomology in Section 2.5, and finish with an explicit computation in Section 2.6.

\section{Modules over Hopf algebras}
\label{sec:modulesHopf}
\index{module, over Hopf algebra|(}

Let $k$ be a graded commutative ring, where we mean commutative in the graded sense. All unlabelled tensors and homs are to be taken over $k$. We denote the closed symmetric monoidal category of right~$k$-modules by $\Mod(k)$. Note that such modules are implicitly graded, and that morphisms of such modules, which we will refer to as $k$-linear homomorphisms, are degree-preserving.

\begin{definition}
\label{def:Hopfalg}
A \emph{Hopf algebra $\Gamma$ over $k$}\index{Hopf algebra} is a $k$-module equipped with five $k$-linear homomorphisms: \emph{multiplication} $\phi \: \Gamma \otimes \Gamma \to \Gamma$, \emph{comultiplication} $\psi \: \Gamma \to \Gamma \otimes \Gamma$, \emph{unit} $\eta \: k \to \Gamma$, \emph{counit} $\epsilon \: \Gamma \to k$, and \emph{antipode} $\chi \: \Gamma \to \Gamma$. These are subject to the following conditions:
\begin{enumerate}
\item Multiplication and unit provide $\Gamma$ with the structure of a $k$-algebra.
\item Comultiplication and counit provide $\Gamma$ with the structure of a $k$-coalgebra.
\item Comultiplication and counit are $k$-algebra morphisms, or equivalently, multiplication and unit are $k$-coalgebra morphisms.
\item The antipode satisfies the formulae $\phi(1 \otimes \chi)\psi = \eta \epsilon = \phi(\chi \otimes 1)\psi$.
\end{enumerate}
\end{definition}

We say that a Hopf algebra is \emph{cocommutative}\index{Hopf algebra!cocommutative} if the comultiplication satisfies~$\tau \psi = \psi$, where $\tau$ denotes the twist in~$\Mod(k)$. We are going to assume that all Hopf algebras we work with are cocommutative in this memoir.

A module over a Hopf algebra is just a module over the underlying $k$-algebra. For a right $\Gamma$-module~$M$ we denote the right action by $\rho_M \: M \otimes \Gamma \to M$. We denote the category of right $\Gamma$-modules by $\Mod(\Gamma)$. This is a closed symmetric\footnote{Symmetry uses that $\Gamma$ is cocommutative.} monoidal category if we endow the category with the tensor products and internal homs over~$k$ together with appropriate $\Gamma$-actions on these objects. Here let $M$,~$N$, and $L$ be~$\Gamma$-modules. The tensor product $M \otimes N$ is endowed with the diagonal~$\Gamma$-action\index{diagonal action}\index{tensor product!of modules over Hopf algebra}. This is the composition
\[
\xymatrix{
M \otimes N \otimes \Gamma \ar[r]^-{1 \otimes 1 \otimes \psi} & M \otimes N \otimes \Gamma \otimes \Gamma \ar[r]^-{1 \otimes \tau \otimes 1} & M \otimes \Gamma \otimes N \otimes \Gamma \ar[r]^-{\rho_M \otimes \rho_N} & M \otimes N 
} \,.
\]
The unit of the tensor product is $k$ regarded as a trivial $\Gamma$-module via the counit:
\[
\xymatrix{
k \otimes \Gamma \ar[r]^-{1 \otimes \epsilon} & k \otimes k = k 
} \,.
\]
The internal Hom\index{internal Hom!of modules over Hopf algebra} $\Hom(N,L)$ becomes a $\Gamma$-module by giving it the conjugate~$\Gamma$-action. This is the~$\Gamma$-action that needs to be on the internal Hom to make sure that $\Hom(N,-)$ is right adjoint to $(-) \otimes N \: \Mod(\Gamma) \to \Mod(\Gamma)$. In other words, the characterising feature of the conjugate $\Gamma$-action\index{conjugate action} is that it is the~$\Gamma$-action on $\Hom(N,L)$ that makes the counit $\Hom(N,L) \otimes N \to L$ and the unit $M \to \Hom(M \otimes N, N)$ into~$\Gamma$-linear maps. Explicitly, the conjugate action is adjoint to the composition
\begin{align*}
\Hom(N,L) \otimes \Gamma \otimes N &\xrightarrow{1 \otimes \tau} \Hom(N,L) \otimes N \otimes \Gamma \xrightarrow{1 \otimes 1 \otimes \psi} \Hom(N,L) \otimes N \otimes \Gamma \otimes \Gamma \\ & \xrightarrow{1 \otimes 1 \otimes \chi \otimes 1} \Hom(N,L) \otimes N \otimes \Gamma \otimes \Gamma \\ &\xrightarrow{1 \otimes \rho_N \otimes 1} \Hom(N,L) \otimes N \otimes \Gamma\xrightarrow{\ev \otimes 1} L \otimes \Gamma \xrightarrow{\rho_L} L \,.
\end{align*}
These actions on tensor and hom-objects ensure that the forgetful functor\index{forgetful functor}
\[
U \: \Mod(\Gamma) \to \Mod(k)
\]
is strict closed monoidal. 

\begin{lemma}
\label{lem:proj}
Let $M$ and $N$ be $\Gamma$-modules, where we assume that $M$ is projective over $\Gamma$ and $N$ is projective over $k$. Then $M \otimes N$ is projective over $\Gamma$.
\begin{proof}
By the tensor-hom adjunction we have a natural isomorphism
\[
\Hom_\Gamma(M \otimes N,-) \cong \Hom_\Gamma(M , \Hom(N,-))
\]
of functors. Since $N$ is projective over $k$ the functor~$\Hom(N,-)$ is exact, and since~$M$ is projective over~$\Gamma$ the functor $\Hom_\Gamma(M,-)$ is exact. The left hand functor $\Hom_\Gamma(M \otimes N , -)$ is then also exact, being naturally isomorphic to the composition of two exact functors. This is equivalent to the assertion that~$M \otimes N$ is projective over $\Gamma$. 
\end{proof}
\end{lemma}

The forgetful functor~$U$ admits a left adjoint 
\[
\ind \: \Mod(k) \longto \Mod(\Gamma) \,,
\]
which we refer to as \emph{induction}\index{induction}. This functor sends a $k$-module $C$ to $C \otimes \Gamma$ with the~$\Gamma$-action given by
\[
C \otimes \Gamma \otimes \Gamma \xrightarrow{1 \otimes \phi} C \otimes \Gamma \,.
\] 
The forgetful functor~$U$ also admits a right adjoint
\[
\coind \: \Mod(k) \longto \Mod(\Gamma) \,,
\]
which is referred to as \emph{coinduction}\index{coinduction}. This functor sends a $k$-module $C$ to $\Hom(\Gamma,C)$ with $\Gamma$-action given as the adjoint of 
\begin{align*}
\Hom(\Gamma,C) \otimes \Gamma \otimes \Gamma \xrightarrow{1 \otimes \tau} \Hom(\Gamma,C) \otimes \Gamma \otimes \Gamma &\xrightarrow{1 \otimes 1 \otimes \chi} \Hom(\Gamma,C) \otimes \Gamma \otimes \Gamma \\ & \xrightarrow{1 \otimes \phi} \Hom(\Gamma,C) \otimes \Gamma \xrightarrow{\ev} C \,.
\end{align*}

The fact that the forgetful functor is strict monoidal makes sure that induction and coinduction interact with the forgetful functor in various useful ways. In~\cite{LMS86}*{Section 2.4} the following formulae, in the context of equivariant stable homotopy theory, are called \emph{untwisting isomorphisms}\index{untwisting isomorphism}, and we will refer to them as such also in this memoir.  

\begin{proposition}
\label{prop:inductiontwisting}
Let $M$ be a $\Gamma$-module and let $C$ be a $k$-module. There are natural~$\Gamma$-module isomorphisms:
\begin{enumerate}
\item 
\[
\ind(C \otimes U(M)) \cong \ind(C) \otimes M
\]
\item 
\[
\Hom(M,\coind(C)) \cong \coind(\Hom(U(M),C))
\]
\item 
\[
\Hom(\ind(C),M) \cong \coind(\Hom(C,U(M))) \,.
\]
\end{enumerate}
\begin{proof}
The result follows formally from the Yoneda lemma together with the fact that $\Mod(\Gamma) \to \Mod(k)$ is strict closed monoidal. We show the first isomorphism and leave the others for the reader, as they are proven in a similar manner.

Consider the functor corepresented by $\ind(C \otimes U(M))$. By adjunctions we have natural isomorphisms 
\begin{align*}
\Hom_{\Gamma}(\ind(C \otimes U(M)),-) &\cong \Hom(C \otimes U(M),U(-)) \\ &\cong \Hom(C,\Hom(U(M),U(-))) \,.
\end{align*}
Since the forgetful functor is strict closed monoidal we have the identity
\[
\Hom(C,\Hom(U(M),U(-))) = \Hom(C,U(\Hom(M,-)))
\]
and by adjunctions again
\begin{align*}
\Hom(C,U(\Hom(M,-))) &\cong \Hom_\Gamma(\ind(C),\Hom(M,-)) \\ &\cong \Hom_\Gamma(\ind(C) \otimes M,-) \,.
\end{align*}
The Yoneda lemma now asserts that we have a natural isomorphism, as wanted.
\end{proof}
\end{proposition}

\begin{corollary}
\label{cor:untwist}
Let $M$ be a $\Gamma$-module. There are natural isomorphisms
\[
\ind(U(M)) \cong \Gamma \otimes M \och \coind(U(M)) \cong \Hom(\Gamma,M)
\]
where the $\Gamma$-actions on the right hand sides are the ordinary diagonal and conjugate actions, respectively.
\begin{proof}
Use that $\Gamma = \ind(k)$.
\end{proof}
\end{corollary}

We will also deal a lot with functional duals of modules over Hopf algebras, so let us now recall this story.

\begin{definition}
For each $\Gamma$-module $M$ let
\[
DM = \Hom(M, k)
\]
be its \emph{functional dual}\index{functional dual!of module over Hopf algebra}. This is a $\Gamma$-module by using the usual conjugate $\Gamma$-action.
\end{definition}

Note that the evaluation pairing\index{evaluation pairing} $\ev \: \Hom(N,L) \otimes N \to L$ gives rise to a natural~$\Gamma$-linear pairing
\[
\alpha \: \Hom(N,L) \otimes \Hom(N',L') \longto \Hom(N \otimes N' , L \otimes L')
\] 
adjoint to the composition
\begin{align*}
\Hom(N,L) \otimes \Hom(N',L') \otimes N \otimes N' &\xrightarrow{1 \otimes \tau \otimes 1} \Hom(N,L) \otimes N \otimes \Hom(N',L') \otimes N' \\ &\xrightarrow{\ev \otimes \ev}  L \otimes L'\,.
\end{align*}
In the case $N' = L = k$ this specialises to a natural $\Gamma$-linear homomorphism
\[
\nu \: DN \otimes L' \longto \Hom(N,L') \,.
\]
This map is an isomorphism when $N$ is finitely generated and projective over~$k$. 

So far we have only discussed $\Gamma$-modules, but we can also talk about (right) comodules over $\Gamma$, by which we mean comodules over the underlying coalgebra structure of $\Gamma$. If $\Gamma$ is finitely generated projective over $k$ then we can endow its functional dual $D\Gamma$ with such a $\Gamma$-coaction. Moreover, this $\Gamma$-coaction is compatible with the $\Gamma$-action in a suitable way. See~\cite{Par71}*{Prop.~2}. This allows us to conclude the following.

\begin{theorem}[\cite{Par71}*{Lem.~2, Prop.~3}]
\label{thm:pareigis}
Let $\Gamma$ be a finitely generated projective Hopf algebra over~$k$.  Then there is an
isomorphism
\[
D\Gamma \cong \ind P(D\Gamma)
\]
of right $\Gamma$-modules, where $P(D\Gamma)$ is a finitely generated
projective $k$-module of constant rank~$1$, given as the primitives for the right
$\Gamma$-coaction on $D\Gamma$.
\end{theorem}

Note in particular that a direct consequence of $\Gamma$ being finitely generated and projective over $k$ is that $D\Gamma$ is itself finitely generated and projective over $\Gamma$. This result will be crucial in our treatment of Tate cohomology of Hopf algebras. For now, let us simply note that the result implies that we have a `Wirthm\"uller isomorphism'.

\begin{corollary}
\label{cor:wirthmuller}
Let $\Gamma$ be a finitely generated and projective Hopf algebra over~$k$ and let $C$ be a $k$-module. There is a natural isomorphism 
\[
\ind(P(D\Gamma) \otimes C) \cong \coind(C)
\]
of $\Gamma$-modules. 
\begin{proof}
We can assume that $C$ is obtained from a $\Gamma$-module $M$ by forgetting the~$\Gamma$-action, as in $C = UM$. By the first untwisting isomorphism of Proposition~\ref{prop:inductiontwisting} we have  
\[
\ind(P(D\Gamma) \otimes U(M)) \cong \ind(P(D\Gamma)) \otimes M \,.
\]
By Pareigis' result it follows that 
\[
\ind(P(D\Gamma)) \otimes M \cong D\Gamma \otimes M \,.
\]
Finally, by $\Gamma$ being finitely generated projective over $k$ and untwisting, more specifically Corollary~\ref{cor:untwist}, we have 
\[
D\Gamma \otimes M \cong \Hom(\Gamma,M) \cong \coind(U(M)) \,.
\]
\end{proof}
\end{corollary}

Since $P(D\Gamma)$ is tensor-invertible over $k$ the above tells us that induced modules are the same things as coinduced modules. We also note that that duals of finitely generated projectives over $\Gamma$ are themselves finitely generated projective over $\Gamma$. 

\begin{corollary}
\label{cor:Pareigis}
Let $\Gamma$ be a finitely generated projective Hopf algebra over~$k$ and let~$M$ be a finitely generated projective $\Gamma$-module. Then its dual $DM$ is also finitely generated projective over $\Gamma$. 
\begin{proof}
Since $M$ is finitely generated, we can find a short exact sequence
\[
\xymatrix{
0 \ar[r] & \ker(r) \ar[r] & \bigoplus_{i \in I} \Sigma^{n_i} \Gamma \ar[r]^-{r} & M \ar[r] & 0 
}
\]
of $\Gamma$-modules, where $I$ is a finite indexing set. Since $M$ is projective, this short exact sequence splits, so that we can find a $\Gamma$-linear map $u \: M \to \bigoplus_{i \in I} \Sigma^{n_i} \Gamma$ such that $ru = \id_M$. Consider the $k$-dual picture. Applying $\Hom(-,k)$ gives us a long exact sequence of $\Gamma$-modules:
\[
\xymatrix{
\cdots &\ar[l] \Ext_k^1(M,k) & \ar[l] D \ker(r) & \ar[l] D\left(\bigoplus_{i \in I} \Sigma^{n_i}\Gamma \right) & \ar[l]_-{r^*} DM & \ar[l] 0 \,.
}
\]
Since $M$ is projective over $\Gamma$, which is in turn projective over $k$, it follows that~$M$ is projective over $k$, from which we conclude that $\Ext_k^1(M,k) \cong 0$. We are left with a short exact sequence, which is also split since $u^* r^* = \id_{DM}$. We conclude that $DM$ is finitely generated projective over $\Gamma$ since it is a retract of the finitely generated projective $\Gamma$-module
\[
 D\left(\bigoplus_{i \in I} \Sigma^{n_i}\Gamma \right) \cong \bigoplus_{i \in I} \Sigma^{-n_i} D\Gamma \cong \bigoplus_{i \in I} \Sigma^{-n_i} \ind P(D\Gamma) \,.
\] 
\end{proof}
\end{corollary}

\index{module, over Hopf algebra|)}

\section{Chain complexes of $\Gamma$-modules}
\label{subsec:chaincomplex}

In this section we give the conventions for chain complexes of~$\Gamma$-modules. A lot is standard: the category $\Mod(\Gamma)$ is an abelian category and what we mean by chain complexes of $\Gamma$-modules is nothing more than the ordinary category of chain complexes in this abelian category. However, we want to make a point of clarifying certain subtle points, especially related to grading and signs. 

\begin{definition} 
A \emph{chain complex $X_*$ of $\Gamma$-modules}\index{chain complex} is a family $(X_n)_{n \in \bZ}$ of $\Gamma$-modules together with morphisms of $\Gamma$-modules $\partial \: X_n \to X_{n-1}$, called \emph{boundaries}\index{boundary, in chain complex}, such that $\partial^2 = 0$. A \emph{chain map}\index{chain map} $f \: X_* \to Y_*$ is a family of $\Gamma$-module homomorphisms~$f_n \: X_n \to Y_n$ that commute with the boundaries.
\end{definition}

The $\Gamma$-module $X_n$ in the chain complex is of course implicitly graded:
\[
X_n = \bigoplus_{\ell \in \bZ} X_{n,\ell},
\]
and if we want to emphasize the bigrading we will write $X_{*,*}$ for the complex. We remark on the following point: as the boundaries $\partial \: X_n \to X_{n-1}$ are morphisms of $\Gamma$-modules, they preserve the $\Gamma$-module grading in the sense that they are given as a direct sum of maps $\partial: X_{n,\ell} \rightarrow X_{n-1,\ell}$. We use the following terminology for the different degrees.

\begin{definition}
Let $X_{*}$ be a chain complex of $\Gamma$-modules. If $x$ is an element in~$X_{n,\ell}$ we say that $x$ has \emph{homological degree}\index{degree!homological} $n$, \emph{internal degree}\index{degree!internal} $|x| = \ell$, and \emph{total degree}~$\|x\| = n + \ell$\index{degree!total}.
\end{definition}

The category of chain complexes of $\Gamma$-modules is closed symmetric monoidal. If~$X_*$, $Y_*$, and $Z_*$ are chain complexes of $\Gamma$-modules then the tensor product~$X_* \otimes Y_*$ of the complexes\index{tensor product!of chain complexes} $X_*$ and $Y_*$ is (either the evident bicomplex or) the complex given in degree~$n$ as 
\[
(X \otimes Y)_n = \bigoplus_{i+j=n} X_i \otimes Y_j, \quad \partial(x \otimes y) = \partial(x) \otimes y + (-1)^{||x||} x \otimes \partial(y) \,.
\]
The unit for the tensor product is $k$ concentrated in homological degree 0. Note that the twist isomorphism is given as 
\[
\tau \: X_* \otimes Y_* \longto Y_* \otimes X_* \,, \quad x \otimes y \mapsto (-1)^{\|x\|\|y\|} y \otimes x \,.
\]
The Hom complex~$\Hom(Y_*,Z_*)$ of $Y_*$ and $Z_*$\index{internal Hom! of chain complexes} is (either the evident
bicomplex or) the complex given in degree~$n$ as 
\[
\Hom(Y,Z)_n = \prod_{i+j=n} \Hom(Y_{-i},Z_j), \quad (\partial f)(x) = \partial(f(x)) - (-1)^{||f||} f(\partial(x)) \,.
\]
We will in particular be interested in the case when $Z_\ast$ is a $\Gamma$-module $M$, regarded as a chain complex concentrated in homological degree 0. In this case, we will often denote the differential in the resulting function complex as $\partial^* = \Hom(\partial,1)$. Explicitly, we have 
\[
\Hom(Y,M)_n = \Hom(Y_{-n},M)\,, \quad (\partial^* f)(x) = -(-1)^{||f||} f(\partial(x)) \,.
\]
As before, the evaluation pairing\index{evaluation pairing} $\ev \: \Hom(Y_*,Z_*) \otimes Y_* \to Z_*$ gives rise to a natural~$\Gamma$-chain map 
\[
\alpha \: \Hom(Y_*,Z_*) \otimes \Hom(Y'_*,Z'_*) \longto \Hom(Y_* \otimes Y'_* , Z_* \otimes Z'_*)
\] 
adjoint to the composition
\begin{align*}
\Hom(Y_*,Z_*) \otimes \Hom(Y'_*,Z'_*) \otimes Y_* &\otimes Y'_* \\ &\xrightarrow{1 \otimes \tau \otimes 1} \Hom(Y_*,Z_*) \otimes Y_* \otimes \Hom(Y'_*,Z'_*) \otimes Y'_* \\ &\xrightarrow{\ev \otimes \ev}  Z_* \otimes Z'_* \,.
\end{align*}
Note that this introduces a sign in the formula for $\alpha$ coming from the twist $\tau$. Explicitly,
\begin{align}
\label{eq:sign}
\alpha(f \otimes g)(x \otimes y) = (-1)^{\|g\|\|x\|} f(x) \otimes g(y) \,.
\end{align}

We now turn to suspensions and mapping cones. These are determined by specifying a `circle chain complex' and an `interval chain complex'. 

\begin{notation}
The interval object\index{interval object} is the chain complex $I_\ast$ given as 
\[
0 \longto k\{i_1\} \overset{\partial}\longto k\{i_0\} \longto 0 \,, \quad \partial(i_1) = i_0 \,.
\]
Both of the generators $i_0$ and $i_1$ are regarded as a having internal degree 0 and the subscripts indicate the homological degrees. 
\end{notation}

\begin{notation}
The circle object\index{circle object} is the chain complex~$C_\ast$ given as
\[
0 \longto k\{c_1\} \longto 0 \,,
\]
again with $c_1$ regarded as having internal degree 0 and with the subscript indicating the homological degree. 
\end{notation}

The convention we will use in this memoir is that chain complexes are suspended on the left. In more precise terms, the suspension\index{suspension!of a chain complex} of a chain complex $X_*$ is the chain complex $X[1]_\ast = C_* \otimes X_*$. From the definition of the symmetric monoidal structure and the appropriate identifications we get  
\[
X[1]_n \cong X_{n-1} \,, \quad \partial_{X[1]}(x) = - \partial_X(x) \,.
\]   

\begin{definition}
\label{def:mappingcone}
The \emph{mapping cone}\index{mapping cone!of chain map} of a chain map $f \: X_\ast \to Y_\ast$ is the chain complex $\cone(f)_\ast$ given as the pushout in the diagram
\[
\xymatrix{
X_\ast \ar[r]^-{f} \ar[d]_-{i_0} & Y_\ast \ar[d] \\
I_\ast \otimes X_\ast \ar[r] & \cone(f)_\ast
}
\]
where $i_0(x) = i_0 \otimes x$.
\end{definition}

Explicitly, we have 
\[
\cone(f)_n \cong X_{n-1} \oplus Y_{n} \,, \quad \partial(x,y) = (-\partial(x),\partial(y) + f(x)) \,.
\]
We have a short exact sequence of chain complexes 
\[
\xymatrix{
0 \ar[r] & Y_* \ar[r] & \cone(f)_* \ar[r] & X[1]_* \ar[r] & 0
}
\]
where the first map is $y \mapsto (0,y)$ and the second one is $(x,y) \mapsto x$. We leave it to the reader to convince themself that these are indeed chain maps.

\section{Tate complexes}

Consider a projective $\Gamma$-resolution $\epsilon \: P_\ast \to k$ of the trivial~$\Gamma$-module $k$. We will denote the mapping cone of the map~$\epsilon$ as~$\widetilde{P}_\ast$. With the conventions from before we hence have
\[
\widetilde{P}_n \cong \begin{cases}
k & \text{if $n = 0$} \\
P_{n-1} & \text{otherwise}
\end{cases}
\]
with boundary $\tilde\partial \: \widetilde P_n \to \widetilde P_{n-1}$ given as 
\[
\tilde \partial(x) = \begin{cases}
-\partial(x) & \text{if $n \geq 2$} \\
\epsilon(x) & \text{if $n=1$.}
\end{cases}
\]
Let us use the notation $i \: k \to \widetilde{P}_*$ for the inclusion. We now define the so-called Tate complex.

\begin{definition} \label{def:Tatecx}
For each $\Gamma$-module $M$ let
\[
\hm_*(M) = \widetilde{P}_* \otimes \Hom(P_*, M)
\]
be the \emph{Tate complex}\index{Tate complex} of~\cite{Gre95}*{\S3}.
\end{definition}

Complexes of this type arise from a filtration of the Tate construction on a $G$-spectrum that we call the Hesselholt--Madsen filtration, which is adapted from~\cite{HM03}, and this explains the notation `$\hm$'. See Section~\ref{subsec:HMfilt}. Explicitly, the Tate complex is given in each homological degree by 
\[
\hm_n(M) = \bigoplus_{i+j=n} \widetilde{P}_i \otimes \Hom(P_{-j}, M)
\]
with boundary given as 
\[
\partial_{\hm}(x \otimes f) = \tilde{\partial}(x) \otimes f + (-1)^{||x||} x \otimes \partial^*(f) \,.
\]

\begin{definition}
\label{def:completeExt}
For an integer $n$ let
\[
\hatExt_\Gamma^{n}(k, M) = H_{-n}(\Hom_\Gamma(k, \hm_*(M)))
\]
be the $k$-module given by the $(-n)$th homology of the chain complex
\[
\Hom_\Gamma(k, \hm_*(M)) = \Hom_\Gamma(k,\widetilde{P}_* \otimes \Hom(P_*, M)) \,.
\]
We call this the \emph{$n$th Tate cohomology group of $\Gamma$ with coefficients in the $\Gamma$-module~$M$}\index{Tate cohomology!of a Hopf algebra}\index{complete $\Ext$}.
\end{definition}

To be able to compare this definition to the standard definition of Tate cohomology in terms of complete resolutions, it is convenient to introduce an alternative, quasi-isomorphic chain complex.

\begin{definition} \label{def:gmcx}
For each $\Gamma$-module $M$, let $\gm_\ast(M)$ be the pushout in the diagram
\[
\begin{tikzcd}
M \arrow[r,"\epsilon^*"] \arrow[d,"i \otimes 1"]& \Hom(P_*,M) \arrow[d] \\
\widetilde{P}_\ast \otimes M \arrow[r] & \gm_*(M) \,.
\end{tikzcd}
\]
Here the top horizontal morphism is the map
\[
\epsilon^* = \Hom(\epsilon,1) \: M \cong \Hom(k,M) \longto \Hom(P_*,M)
\]
contravariantly induced by the augmentation, and the left hand vertical morphism is the map 
\[
i \otimes 1 \: M \cong k \otimes M \longto \widetilde{P}_* \otimes M
\]
induced by the inclusion of $k$ into the mapping cone $\widetilde{P}_* = \cone(\epsilon)$.
\end{definition}

Complexes of this type arise from a filtration of the Tate construction on a $G$-spectrum that we call the Greenlees--May filtration, which is adapted from~\cite{GM95}, and this explains the notation `$\gm$'. See Section~\ref{subsec:GMfilt}. 

\begin{proposition} \label{prop:gmMexplicit}
Explicitly, the complex $\gm_*(M)$ is given in each homological degree as  
\[
\gm_n(M) \cong \begin{cases} \widetilde{P}_n \otimes M & \text{if $n \geq 1$} \\ \Hom(P_{-n},M) & \text{if $n \leq 0$} \end{cases}
\]
and under these identifications the boundary $\partial_{\gm} \: \gm_n(M) \to \gm_{n-1}(M)$ is given as
\[
\partial_{\gm} = \begin{cases}
\tilde{\partial}\otimes 1 & \text{if $n \geq 2$} \\
\partial^* & \text{if $n \leq 0$} \\
\widetilde{P}_1 \otimes M \xrightarrow{\epsilon \otimes 1} M \xrightarrow{\epsilon^*} \Hom(P_0,M) & \text{if $n=1$.}
\end{cases}
\]
\begin{proof}
The only non-trivial case happens when the homological degree is $n=0$. In this case we have a pushout square
\[
\begin{tikzcd}
M \arrow[r,"\epsilon^*"] \arrow[d,"i \otimes 1"] & \Hom(P_0,M) \arrow[d] \\
\widetilde{P}_0 \otimes M \arrow[r] & \gm_0(M) \,.
\end{tikzcd}
\]
Since $i \otimes 1$ is an isomorphism in this homological degree it follows that so is the map $\Hom(P_0,M) \to \gm_0(M)$.  

It is straight-forward to see that the boundary $\partial_{\gm} \: \gm_n(M) \to \gm_{n-1}(M)$ is given by $\tilde{\partial} \otimes 1$ and $\partial^*$ when $n \geq 2$ and $n\leq 0$, respectively. For the remaining boundary, note that the element $1 \otimes m$ in~$\widetilde{P}_0 \otimes M$ is identified with the element $f \: y \mapsto m \epsilon(y)$ in $\Hom(P_0,M)$ when both are viewed as elements of the pushout~$\gm_0(M)$. The boundary wants to take the element $x \otimes m$ in $\widetilde{P}_1 \otimes M \cong \gm_1(M)$ to~$\epsilon(x) \otimes m$ in~$\widetilde{P}_0 \otimes M$. This is identified with the map $y \mapsto (-1)^{|m||x|}m \epsilon(x)\epsilon(y)$ in $\Hom(P_0,M)$. Schematically, we are taking the composite
\[
\widetilde{P}_1 \otimes M
	\overset{\epsilon \otimes 1}\longto
\widetilde{P}_0 \otimes M \cong M \cong \Hom(k, M)
	\overset{\epsilon^*}\longto \Hom(P_0, M) \,. \qedhere
\]
\end{proof}

\end{proposition}

Visually, $\gm_*(M)$ is the complex
\[
\begin{tikzcd}[column sep = small]
\cdots \arrow[r] & \widetilde{P}_2 \otimes M \arrow[r,"\widetilde{\partial} \otimes 1"] & \widetilde{P}_1 \otimes M \arrow[rr] \arrow[rd,"\epsilon \otimes 1"'] & & \Hom(P_0,M) \arrow[r,"\partial^*"] & \Hom(P_1,M) \arrow[r] & \cdots \,. \\
& & & M \arrow[ur,"\epsilon^*"'] & & & 
\end{tikzcd}
\]
We will often refer to the complex $\gm_*(M)$ as being obtained by `splicing' $\widetilde{P}_* \otimes M$ and $\Hom(P_*,M)$ together. 

Note that the universal property of the pushout ensures that we have an induced~$\Gamma$-chain map $\theta \: \gm_\ast(M) \to \hm_*(M)$ in the commutative diagram
\[
\begin{tikzcd}
M \arrow[r,"\epsilon^*"] \arrow[d,"i \otimes 1"] & \Hom(P_*,M) \arrow[d] \arrow[rdd , bend left, "i \otimes 1"] & \\
\widetilde{P}_* \otimes M \arrow[r] \arrow[rrd, bend right , "1 \otimes \epsilon^*"] & \gm_*(M) \arrow[rd, "\theta"] & \\
& & \hm_*(M) \,.
\end{tikzcd}
\]
Here the `bendy' map 
\[
1 \otimes \epsilon^* \: \widetilde{P}_\ast \otimes M \cong \widetilde{P}_\ast \otimes \Hom(k,M) \longto \widetilde{P}_* \otimes \Hom(P_\ast,M)
\]
is again the map contravariantly induced by the augmentation, and the other `bendy' map 
\[
i \otimes 1 \: \Hom(P_*,M) \cong k \otimes \Hom(P_*,M) \longto \widetilde{P}_\ast \otimes \Hom(P_*,M)
\]
is induced by the inclusion $i \: k \to \widetilde{P}_*$.

\begin{proposition} \label{prop:gmtohmisqi}
The $k$-linear chain map
\[
\Hom(1,\theta) \: \Hom_\Gamma(k, \gm_*(M))
	\longto \Hom_\Gamma(k, \hm_*(M))
\]
is a quasi-isomorphism, inducing isomorphisms
\[
H_n(\Hom_\Gamma(k, \gm_*(M)))
	\overset{\cong}\longto \hatExt_\Gamma^{-n}(k, M)
\]
for all integers $n$.
\end{proposition}

\begin{proof}
We compatibly filter $\gm_*(M)$ and $\hm_*(M)$, setting
\[
F_s \gm_k(M) = \begin{cases}
	0 & \text{for $k > s$}, \\
	\gm_k(M) & \text{for $k \le s$}
\end{cases}
\]
and
\[
F_s \hm_k(M) = \bigoplus_{\substack{i+j=k \\ i \le s}}
	\widetilde{P}_i \otimes \Hom(P_{-j}, M) \,.
\]
We obtain a vertical map of short exact sequences
\[
\xymatrix{
0 \ar[r]
	& F_{s-1} \gm_*(M) \ar[r] \ar[d]_-{\theta_{s-1}}
	& F_s \gm_*(M) \ar[r] \ar[d]_-{\theta_s}
	& \ds \frac{F_s \gm_*(M)}{F_{s-1} \gm_*(M)} \ar[r] \ar[d]_-{\bar\theta_s}
	& 0 \\
0 \ar[r]
	& F_{s-1} \hm_*(M) \ar[r]
	& F_s \hm_*(M) \ar[r]
	& \ds \frac{F_s \hm_*(M)}{F_{s-1} \hm_*(M)} \ar[r]
	& 0 \rlap{\,.}
}
\]
Each horizontal short exact sequence is degree-wise split as an
extension of $\Gamma$-modules, hence remains short exact
after applying $\Hom_\Gamma(k, -)$.

For $s=0$, the map $F_0 \gm_*(M) \to F_0 \hm_*(M)$ is an isomorphism.
We claim for each $s\ge1$ that the map of filtration
subquotients
\[
\bar\theta_s \: \frac{F_s \gm_*(M)}{F_{s-1} \gm_*(M)}
	\longto
\frac{F_s \hm_*(M)}{F_{s-1} \hm_*(M)}
\]
induces a quasi-isomorphism $\Hom(1, \bar\theta_s)$ after applying
$\Hom_\Gamma(k, -)$. It follows by induction that $\Hom(1, \theta_s)$
is a quasi-isomorphism for each $s\ge0$.  Passing to colimits over~$s$ it follows that $\Hom(1, \theta)$ is a quasi-isomorphism.

It remains to prove the claim.  We can rewrite $\bar\theta_s$
for $s\ge1$ as
\[
1 \otimes \epsilon^* \: \widetilde{P}_s \otimes M
	\longto \widetilde{P}_s \otimes \Hom(P_*, M) \,.
\]
Here $\widetilde{P}_s$ is $\Gamma$-projective, so it suffices to prove
that
\[
\Hom_\Gamma(1, 1 \otimes \epsilon^*) \: \Hom_\Gamma(k, L \otimes M)
\longto \Hom_\Gamma(k, L \otimes \Hom(P_*, M))
\]
is a quasi-isomorphism for any projective $\Gamma$-module $L$.
By preservation of quasi-isomorphisms under passage to retracts, we may
assume that $L$ is free. Since the functor~$\Hom_\Gamma(k, -)$ commutes with direct sums, it suffices to consider the case $L = \Gamma$.  Using the Hopf
algebra structure of $\Gamma$, there is a natural untwisting isomorphism
\[
\Gamma \otimes N \cong \ind(N)
\]
for any $\Gamma$-module~$N$, where $\Gamma$ acts diagonally on the left
hand side and we use the induced $\Gamma$-action on the right hand side. See Corollary~\ref{cor:untwist}. The augmentation $\epsilon \: P_* \to k$ admits a $k$-linear chain
homotopy inverse.  Hence $\epsilon^* \: M \to \Hom(P_*, M)$ also
admits such a chain homotopy inverse, and
\[
\ind(\epsilon^*) \: \ind(M)
	\longto \ind(\Hom(P_*, M))
\]
admits a $\Gamma$-linear chain homotopy inverse. By naturality
of the untwisting isomorphism,
\[
1 \otimes \epsilon^* \: \Gamma \otimes M 
	\longto \Gamma \otimes \Hom(P_*, M) 
\]
admits a $\Gamma$-module chain homotopy inverse, and therefore
induces a $k$-module chain homotopy equivalence after applying
$\Hom_\Gamma(k,-)$. This proves the claim that
$\Hom_\Gamma(1, 1 \otimes \epsilon^*)$ is a quasi-isomorphism.
\end{proof}

\begin{corollary} \label{cor:ExtvshatExt}
The inclusion $\Hom(P_*, M) \to \gm_*(M)$ induces an isomorphism
\[
\gamma \: \Ext_\Gamma^n(k, M) \longto \hatExt_\Gamma^n(k, M)
\]
for each $n \ge 1$, and a surjection for $n = 0$.
\end{corollary}
 
\section{Complete resolutions}
\label{sec:completeres}

In this section, we make the standing assumption that the cocommutative Hopf algebra $\Gamma$ is finitely generated and projective over~$k$. Let us now relate the complex~$\gm_*(M)$ to the complete resolutions often used when defining Tate cohomology. See~\cite{CK97} for a standard treatment of this topic.

\begin{definition}
\label{def:completeresolution}
\index{complete resolution}
Let $\hat{P}_\ast$ be the pullback in the diagram
\[
\xymatrix{
\hat{P}_\ast \ar[r] \ar[d] & P_\ast \ar[d]^-{\epsilon} \\
D\widetilde{P}_\ast \ar[r]^-{i^*} & k
} 
\]
where $D\widetilde{P}_\ast = \Hom(\widetilde{P}_\ast, k)$.
\end{definition}

\begin{proposition}
Explicitly, the chain complex $\hat{P}_\ast$ is given in each homological degree as
\[
\hat P_n \cong \begin{cases}
P_n & \text{for $n\ge0$,} \\
D(\widetilde{P}_{-n}) & \text{for $n < 0$}
\end{cases}
\]
with boundary given as 
\[
\hat\partial_n = \begin{cases}
\partial_n & \text{for $n>0$,} \\
\epsilon^* \circ \epsilon & \text{for $n=0$,} \\
D(\tilde\partial_{1-n}) & \text{for $n<0$}
\end{cases}
\]
under these identifications.
\begin{proof}
The only non-trivial case is when we are dealing with something involving homological degree $0$. Since $D\widetilde{P}_0 \to k$ is an isomorphism, it follows that the projection $\hat{P}_0 \to P_0$ is one, as well. This shows that the chain complex is given in each homological degree as asserted. The only thing left to prove is that the boundaries are given as claimed. To do so, note that the inverse to the projection is the map~$P_0 \to \hat{P}_0$ given by
\[
P_0 \longto \hat{P}_0 = P_0 \times_k D\widetilde{P}_0 \,, \quad x \mapsto (x,\epsilon(x)) \,.
\] 
It is clear that the boundary $\hat{\partial} \: \hat{P}_{n} \to \hat{P}_{n-1}$ is given by $\partial_n$ and $D(\tilde \partial_{1-n})$ when $n > 0$ and $n < 0$, respectively. When $n=0$, we are looking at the boundary
\[
P_0 \times_k D\widetilde{P}_0 \longto D\widetilde{P}_1 \,, \quad (x,f) \mapsto \epsilon^*(f)
\]
which under the identifications made above corresponds to the composition
\[
P_0 \longto P_0 \times_k D\widetilde{P}_0 \longto D\widetilde{P}_1 \,, \quad x \mapsto (\epsilon^* \circ \epsilon)(x) \,. \qedhere
\]
\end{proof}
\end{proposition}

Diagrammatically, we can visualise $\hat{P}_*$ as the `spliced' complex
\[
\xymatrix@C-1pc{
\cdots \ar[rr]
	&& \hat P_1 \ar[rr]
	&& \hat P_0 \ar[rr] \ar[dr]_-{\epsilon}
	&& \hat P_{-1} \ar[rr]
	&& \hat P_{-2} \ar[rr]
	&& \cdots \\
&& && & k \ar[ur]_-{\epsilon^*}
} \,.
\]
We will show that if $P_*$ is assumed to be a projective resolution of finite type\footnote{Recall that a chain complex $C_*$ of projective modules is said to be of finite type\index{finite type} if it is finitely generated in each homological degree.}, then this is a complete resolution. See Remark~\ref{rmk:CK}. First we need a lemma.

\begin{lemma}
\label{lem:pullbackchainhomotopy}
Let
\[
\begin{tikzcd}
Q_* \arrow[r,"g'"] \arrow[d,"f'"] & A_* \arrow[d,"f"] \\
C_* \arrow[r,"g"] & B_*
\end{tikzcd}
\]
be a pullback diagram of chain complexes. Assume that there is some chain map $\phi \: B_* \to A_*$ such that $f \phi = \id_{B_*}$ and $\phi f \simeq \id_{A_*}$ witnessed by a chain homotopy $H \: A_n \to A_{n+1}$ satisfying $f H = 0$. Then there is a chain map $\phi' \: C_* \to Q_*$ such that $f'\phi' = \id_{C_*}$ and $\phi' f' \simeq \id_{Q_*}$.
\begin{proof}
Consider the diagram
\[
\begin{tikzcd}
C_*
\arrow[drr, bend left, "\phi g"]
\arrow[ddr, bend right, "\id_{C_*}"']
\arrow[dr, dotted , "\phi'"] & & \\
& Q_* \arrow[r, "g'"] \arrow[d, "f'"]
& A_* \arrow[d, "f"] \\
& C_* \arrow[r, "g"]
& B_*
\end{tikzcd}
\]
in which we have an induced chain map $\phi' \: C_* \to Q_*$ by the universal property of a pullback. This shows that $f'\phi' = \id_{C_*}$. Let us show that this constitutes a left homotopy inverse, as well.

Let $H_n \: A_n \to A_{n+1}$ be the chain homotopy between $\phi f$ and $\id_{A_*}$. That is
\[
\id_{A_*} - \phi f = \partial H_n + H_{n-1} \partial \,.
\]
We want to use this data to build a chain homotopy between $\id_{Q_*}$ and $\phi' f'$. To do this, consider the diagram
\[
\begin{tikzcd}
A_{n}
\arrow[drr, bend left, "H_n"]
\arrow[ddr, bend right, "0"']
\arrow[dr, dotted , "h"] & & \\
& Q_{n+1} \arrow[r, "g'"] \arrow[d, "f'"]
& A_{n+1} \arrow[d, "f"] \\
& C_{n+1} \arrow[r, "g"]
& B_{n+1}
\end{tikzcd}
\]
which commutes since $f H_n = 0$, by assumption. Again, by the universal property of a pullback we have induced maps $h \: A_n \to Q_{n+1}$. Let us set
\[
H'_n = h g' \: Q_n \longto Q_{n+1} \,.
\]
We claim that these maps constitute a chain homotopy between $\id_{Q_*}$ and $\phi' f'$. That is, we claim that they satisfy
\[
\id_{Q_*} - \phi' f' = \partial H'_n + H_{n-1}' \partial \,.
\]
To show this, we appeal to the uniqueness of maps induced from pullbacks. Consider the diagram
\[
\begin{tikzcd}
Q_n
\arrow[drr, bend left, "\partial H_n g' + H_{n-1} \partial g'"]
\arrow[ddr, bend right, "0"']
\arrow[dr, dotted , "?"] & & \\
& Q_n \arrow[r, "g'"] \arrow[d, "f'"]
& A_n \arrow[d, "f"] \\
& C_n \arrow[r, "g"]
& B_n \,.
\end{tikzcd}
\]
This diagram commutes, since 
\begin{align*}
f\partial H_n g' + f H_{n-1} \partial g' &= \partial f H_n g' + f H_{n-1} \partial g' \\ &= 0 \,, 
\end{align*}
so we do indeed have a unique induced map in the diagram. We claim that the question-mark in the diagram can be filled by both $\id_{Q_*} - \phi' f'$ and $\partial H'_n + H_{n-1}' \partial$, so they must agree by uniqueness of the induced map. Checking this claim is straight-forward. The checks
\begin{align*}
g'(\id_{Q_*} - \phi'f') &= g' - g'\phi'f' \\ &= g' - \phi g f' \\&= g' - \phi f g' \\ &= (\id_{A_*} - \phi f)g' \\ &= (\partial H_n + H_{n-1} \partial)g' \,, 
\end{align*}
and 
\begin{align*}
f'(\id_{Q_*} - \phi'f') &= f' - f'\phi'f' \\ &= f' - f' \\&= 0  
\end{align*}
show that the map $\id_{Q_*} - \phi' f'$ fits into the diagram. The checks 
\begin{align*}
g'(\partial H_n' + H_{n-1}' \partial) &= g'\partial h g' + g' h g' \partial \\ &= \partial g' h g' + g' h \partial g' \\&= \partial H_n g' + H_{n-1} \partial g' \\ &= (\partial H_n + H_{n-1} \partial) g'
\end{align*}
and 
\begin{align*}
f'(\partial H_n' + H_{n-1}' \partial) &= f'\partial H_n' + f' H_{n-1}' \partial \\ &= \partial' f' H_{n}' + f' H_{n-1}' \partial \\&=   \partial' f' hg' + f' hg' \partial \\ &= 0
\end{align*}
show that the map $\partial H_n' + H_{n-1}' \partial$ also fits into the diagram, which concludes the proof.
\end{proof}
\end{lemma}

\begin{proposition}
\label{prop:complete}
Assume that $P_*$ is of finite type over $\Gamma$. Then $\hat{P}_*$ is an acyclic complex of projective~$\Gamma$-modules such that $\Hom_\Gamma(\hat{P}_*,Q)$ is acyclic for every coinduced $\Gamma$-module $Q$.
\begin{proof}
Since $P_n$ is finitely generated and projective over $\Gamma$ in each homological degree $n$, it follows that~$\hat{P}_n$ must be finitely generated and projective over $\Gamma$, as well, by Corollary~\ref{cor:Pareigis}. 

To show that $\hat{P}_*$ is acyclic, we show that it is $k$-linearly contractible. Since $\epsilon \: P_* \to k$ is a chain homotopy equivalence, we can find a homotopy inverse $k \to P_*$. In this case we can pick $\eta \: k \to P_*$, so that $\epsilon \eta = \id_k$ on the nose. Since $k$ is concentrated in homological degree $0$ we know that the chain homotopy~$\eta \epsilon \simeq \id_{P_*}$ is zero after post-composition with $\epsilon$, so that Lemma~\ref{lem:pullbackchainhomotopy} applies. This shows that the map $\hat{P}_* \to D\widetilde{P}_*$ is a chain homotopy equivalence. Since~$\widetilde{P}_*$ is chain contractible, we conclude that so is its dual~$D\widetilde{P}_*$ and hence also~$\hat{P}_*$. 

If $Q = \coind(C)$ for some $k$-module~$C$, then 
\[
\Hom_\Gamma(\hat{P}_*,\coind(C)) \cong \Hom(\hat{P}_*,C) \,.
\]
Since $\hat{P}_*$ is $k$-linearly contractible it follows that $\Hom(\hat{P}_*,C)$ is contractible, and therefore acyclic. 
\end{proof}
\end{proposition}

Let $M$ be a $\Gamma$-module. The chain map $\hat{P}_\ast \to P_\ast$ induces a chain map 
\[
\Hom(P_*,M) \longto \Hom(\hat{P}_*,M)
\]
which is an isomorphism in homological degrees $\ast \leq 0$. In addition to this map, we also have a chain map composition
\[
\widetilde{P}_\ast \otimes M \longto DD\widetilde{P}_\ast \otimes M \overset{\nu}\longto \Hom(D\widetilde{P}_*,M) \longto \Hom(\hat{P}_*,M) \,.
\]
In particular, this is an isomorphism in homological degrees $\ast \geq 1$ under the assumption that $\widetilde{P}_*$ is of finite type over $k$. Note that the chain maps described above fit into the commutative diagram
\[
\begin{tikzcd}
M \arrow[r,"\epsilon^*"] \arrow[d,"i \otimes 1"] & \Hom(P_*,M) \arrow[d] \arrow[rdd , bend left] & \\
\widetilde{P}_* \otimes M \arrow[r] \arrow[rrd, bend right] & \gm_*(M) \arrow[rd, "\beta"] & \\
& & \Hom(\hat{P}_*,M) \,,
\end{tikzcd}
\]
so that we have an induced chain map $\beta$ by the universal property of $\gm_\ast(M)$. 

\begin{proposition} \label{prop:gmishomhatP}
Suppose that $\widetilde{P}_*$ is of finite type over~$k$. Then the map
\[
\beta \: \gm_*(M) \overset{\cong}\longto \Hom(\hat P_*, M) \,
\]
is a natural isomorphism of $\Gamma$-chain complexes.
\begin{proof}
The assertion is clear in homological degrees $n > 0$ and $n < 0$. If $n=0$ we are looking at the diagram
\[
\begin{tikzcd}
M \arrow[r] \arrow[d,"\cong"] & \Hom(P_0,M) \arrow[d,"\cong"] \arrow[rdd , bend left , "\cong"] & \\
\widetilde{P}_0 \otimes M \arrow[r] \arrow[rrd, bend right] & \gm_0(M) \arrow[rd, "\beta"] & \\
& & \Hom(\hat{P}_0,M) \,,
\end{tikzcd}
\]
in which we have marked the obvious isomorphisms. It is then clear that $\beta$ is an isomorphism, as well.
\end{proof}
\end{proposition}

\begin{corollary} \label{cor:hatExtfromhatP}
Suppose that the projective $\Gamma$-resolution $P_\ast$ is of finite type over~$k$. Then there is a natural isomorphism
\[
\hatExt_\Gamma^{n}(k, M)
	\cong H^n(\Hom_\Gamma(\hat P_*, M))
\]
for all $n \in \bZ$.
\begin{proof}
Combine Proposition~\ref{prop:gmtohmisqi} and Proposition~\ref{prop:gmishomhatP}.
\end{proof}
\end{corollary}

It turns out that, under the assumption that $\Gamma$ is finitely generated projective over $k$, we can always construct the projective resolution $P_*$ so that it is of finite type over $\Gamma$. It is then necessarily also of finite type over $k$. This can be done via the bar construction, which we now review. See (\cite{May72}*{\S 9, \S 10, \S 11} and~\cite{GM74}*{App. A}) for more details.

\begin{construction}[The bar construction]
\label{const:barconst}
\index{bar construction}
Let $\Gamma$ be a $k$-algebra, $M$ a right $\Gamma$-module, and $N$ a left~$\Gamma$-module. We form a simplicial object $B_\bullet(M,\Gamma,N) \: \Delta^{\op} \rightarrow \Mod(k)$ as follows. In simplicial degree~$q$ we let 
\[
B_q(M,\Gamma,N) = M \otimes \Gamma^{\otimes q} \otimes N .
\]
It is customary to write 
\[
m[\gamma_1 | \cdots | \gamma_q]n = m \otimes \gamma_1 \otimes \cdots \otimes \gamma_q \otimes n 
\]
for an element in the $q$th simplicial degree; hence the terminology \emph{bar} construction. In this notation, the face maps are given as 
\[
d_i(m[\gamma_1 | \cdots | \gamma_q]n) = \begin{cases}
m \gamma_1 [ \gamma_2 | \cdots | \gamma_q] n & i=0 \\
m [ \gamma_1 | \cdots |  \gamma_i \gamma_{i+1} |  \cdots | \gamma_q ] n & 0 < i < q \\
m [ \gamma_1 | \cdots | \gamma_{q-1} ] \gamma_q n & i=q
\end{cases}
\]
and the degeneracy maps are given by
\[
s_i(m[\gamma_1 | \cdots | \gamma_q]n) = m [ \gamma_1 | \cdots | \gamma_i | 1 | \gamma_{i+1} | \cdots | \gamma_q] n \, .
\]
\end{construction}

The simplicial $k$-module $B_\bullet(M,\Gamma,N)$ can be turned into a non-negative $k$-complex in essentially two ways. 
\begin{itemize}
\item The most straight-forward way to turn $B_\bullet(M,\Gamma,N)$ into a $\Gamma$-chain complex $B_*(M,\Gamma,N)$ is by taking the $\Gamma$-module in homological degree $n$ to be equal to the $n$-simplices of $B_\bullet(M,\Gamma,N)$ and to let the boundary in the chain complex be the alternating sum of the face maps:
\[
B_n = B_n(M,\Gamma,N) \och \partial = \sum_{i=0}^{n} (-1)^i d_i \: B_n \longto B_{n-1} \,.
\]
This is referred to as the \emph{bar complex}\index{bar complex}.
\item To get a smaller but quasi-isomorphic chain complex that is more convenient for computations, we can turn $B_\bullet(M,\Gamma,N)$ into a chain complex $NB_* = NB_*(M,\Gamma,N)$ by quotienting out by the degenerate simplices. Explicitly, in homological degree $n$ we have 
\begin{align*}
NB_n &= B_n / (s_0 B_{n-1} + \cdots + s_{n-1} B_{n-1}) \\
&\cong M \otimes \overline{\Gamma}^{\otimes n} \otimes N \,,
\end{align*}
where
\[
\overline{\Gamma} = \cok(\eta) \cong \ker(\epsilon) \,.\footnote{This isomorphism follows from $\epsilon \eta = \id_k$.}
\]
The boundary $\partial \: NB_n \to NB_{n-1}$ is given by the same formula as before, which makes sense because 
\[
\partial(s_0 B_{n-1} + \cdots + s_{n-1} B_{n-1}) = 0 \,.
\]
We refer to $(NB_*, \partial)$ as the \emph{normalised bar complex}\index{normalised bar complex}.
\end{itemize}

There is a natural $\Gamma$-action on the simplicial $k$-module $B_\bullet(M,\Gamma,\Gamma)$ arising from viewing $N = \Gamma$ as a $\Gamma$-$\Gamma$-bimodule. Explicitly, in each simplicial degree we have the right $\Gamma$-action $B_q(M,\Gamma,\Gamma) \otimes \Gamma \to B_q(M,\Gamma,\Gamma)$ given by
\[
m[\gamma_1 | \cdots | \gamma_q]\gamma_{q+1} \otimes \gamma \mapsto m[\gamma_1 | \cdots | \gamma_q]\gamma_{q+1}\gamma
\]
and this $\Gamma$-action commutes with the simplicial structure maps of $B_\bullet(M,\Gamma,\Gamma)$, so that $B_\bullet(M,\Gamma,\Gamma)$ extends to a simplicial $\Gamma$-module. It is a standard exercise in simplicial homotopy theory to check that $B_\bullet(M,\Gamma,\Gamma)$ is simplicially homotopy equivalent to $M$ viewed as a constant simplicial $\Gamma$-module. As a consequence, the complexes 
\[
B_*(M,\Gamma,\Gamma) \och NB_*(M,\Gamma,\Gamma)
\]
are resolutions of $M$ as a $\Gamma$-module. We refer to these as the \emph{bar resolution}\index{bar resolution} and \emph{normalised bar resolution}\index{normalised bar resolution} of $M$ as a $\Gamma$-module, respectively. See \cite{May72}*{Prop. 9.9} and~\cite{GM74}*{Lem. A.8}.

\begin{proposition}
\label{prop:barcomplexfinitetype}
Assume that $\Gamma$ is finitely generated and projective over $k$. If~$M$ is finitely generated projective over $k$, then the bar resolution $B_*(M,\Gamma,\Gamma)$ and the normalised bar resolution $NB_*(M,\Gamma,\Gamma)$ are~$\Gamma$-projective resolutions of $M$ of finite type.
\begin{proof}
Since $B_\bullet(M,\Gamma,\Gamma)$ is simplicially homotopy equivalent to $M$, the bar resolution is a resolution of~$M$. It is finitely generated, and projective in each degree by an application of Lemma~\ref{lem:proj}.  The proof in the normalised case is very similar.
\end{proof}
\end{proposition}

\begin{theorem} \label{prop:hatExtdeltafunctor}
When $\Gamma$ is finitely generated and projective over~$k$, each
short exact sequence
\[
0 \longto M' \longto M \longto M'' \longto 0
\]
of $\Gamma$-modules induces a long exact sequence
\[
\dots \longto \hatExt_\Gamma^n(k, M')
\longto \hatExt_\Gamma^n(k, M)
\longto \hatExt_\Gamma^n(k, M'')
\overset{\delta}\longto \hatExt_\Gamma^{n+1}(k, M')
	\longto \dots \,.
\]
Furthermore, if $M$ is an induced or coinduced $\Gamma$-module\footnote{Due to the assumption that $\Gamma$ is finitely generated projective induced modules are coinduced by Corollary~\ref{cor:wirthmuller}, and vice versa, so these are actually equivalent conditions.},
then $\hatExt_\Gamma^n(k, M) = 0$ for all $n \in \bZ$.
\end{theorem}

\begin{proof}
If $\Gamma$ is finitely generated projective over $k$, then the bar complex $B_\ast(k,\Gamma,\Gamma)$ constitutes a projective~$\Gamma$-resolution of $k$ of finite type, so that Proposition~\ref{prop:complete} applies. The long exact sequence is then induced by the short exact sequence
\[
0 \longto \Hom_\Gamma(\hat P_*, M')
	\longto \Hom_\Gamma(\hat P_*, M)
	\longto \Hom_\Gamma(\hat P_*, M'') \longto 0
\]
of $k$-module chain complexes. Here we are using Corollary~\ref{cor:hatExtfromhatP} to identify the terms in the long exact homology sequence. That Tate cohomology vanishes on induced/coinduced modules is a direct consequence of Proposition~\ref{prop:complete}. 
\end{proof}

\begin{remark}
\label{rmk:CK}\index{complete $\Ext$}\index{complete resolution}
Note that if $\Gamma$ is finitely generated projective, the chain complex~$\hat{P}_*$ constructed from $P_* = B_*(k,\Gamma,\Gamma)$ is indeed a complete $\Gamma$-resolution of $k$ in the sense of~\cite{CK97}*{Definition 1.1}. This uses Proposition~\ref{prop:complete} and the fact that all projective~$\Gamma$-modules are retracts of induced $\Gamma$-modules, which are coinduced~$\Gamma$-modules by Corollary~\ref{cor:wirthmuller}. In particular, the results of~\cite{CK97} apply and we can conclude that our $\hatExt_\Gamma(k,-)$ agrees with what is traditionally referred to as `complete Ext', in this case.
\end{remark}

\begin{remark}
In this approach to Tate cohomology\index{Tate cohomology!of a Hopf algebra} and complete $\Ext$ the Hopf algebra\index{Hopf algebra} structure of~$\Gamma$, rather that just its
algebra structure, enters in two ways: First, it is needed for
Pareigis' theorem (Theorem~\ref{thm:pareigis}), which is used in
Definition~\ref{def:completeresolution} to ensure that the spliced
complex $\hat{P}_*$ consists of projective $\Gamma$-modules, as required
for a complete resolution.  Second, the Hopf algebra diagonal and
conjugation are used in Definitions~\ref{def:Tatecx} and~\ref{def:gmcx}
to make sense of the $\Gamma$-module structure of the Tate complex
$\hm_*(M) = \widetilde{P}_* \otimes \Hom(P_*,M)$, as well as for its
variant~$\gm_*(M)$.
\end{remark}

\section{Multiplicative structure of Tate cohomology}
\label{sec:multTate}\index{cup product!on complete $\Ext$|(}

We will now define a suitable pairing on Hopf algebra Tate cohomology. As before, we will assume that $\Gamma$ is finitely generated and projective over $k$, so that this theory coincides with complete $\Ext$.

\begin{proposition}
\label{prop:Pidentityextension}
There is a unique, up to chain homotopy, $\Gamma$-linear chain map $\Psi \: P_\ast \to P_\ast \otimes P_\ast$ covering the identity map $\id \: k \to k \otimes k$. 
\begin{proof}
We first note that the chain complex $P_\ast \otimes P_\ast$, with diagonal $\Gamma$-action, is a~$\Gamma$-resolution of $k \otimes k = k$. To see this, use the spectral sequence associated to~$P_* \otimes P_*$ viewed as a double complex. This converges strongly to the homology of $P_* \otimes P_*$. The first page of the spectral sequence is given by
\[
E^1_{s,t} \cong H_t(P_s \otimes P_*) \cong P_s \otimes H_t(P_*) \,,
\]
since $P_s$ is projective over $\Gamma$, and hence over $k$. This is zero unless $t = 0$, where it is $P_s$. The~$d^1$-differential is induced by the horizontal differential in the double complex, so that the $E^2$-page is $k$ concentrated at the origin. 

Classical homological algebra then asserts that there is a unique (up to chain homotopy) $\Gamma$-linear chain map as asserted; see~\cite{Mac95}*{Chapter III  Thm. 6.1}.
\end{proof}
\end{proposition}

The $\Gamma$-linear chain map $\Psi \: P_\ast \to P_\ast \otimes P_\ast$ described above induces a product on~$\Ext_\Gamma^*(k,-)$ via the pairing 
\begin{align*}
\Hom_\Gamma(P_\ast,M) \otimes \Hom_\Gamma(P_\ast,N) &\overset{\alpha}\longto \Hom_\Gamma(P_* \otimes P_*,M \otimes N) \\ &\overset{\Psi^*}\longto \Hom_\Gamma(P_\ast,M \otimes N)
\end{align*}
of $k$-module complexes. By cocommutativity of $\Gamma$ and uniqueness (up to chain homotopy) of $\Psi$, we have that $\Psi \simeq \tau \circ \Psi$. Passing to homology, this gives us an associative, unital, and graded commutative multiplication
\[
\smile \: \Ext^*_\Gamma(k,M) \otimes \Ext^*_\Gamma(k,N) \longto \Ext^*_{\Gamma}(k,M \otimes N)
\]
that we will refer to as the \emph{cup product}\index{cup product!on $\Ext$}. In particular, $\Ext^*_{\Gamma}(k,k)$ is a $k$-algebra, and $\Ext_\Gamma^*(k,M)$ is an $\Ext_\Gamma^*(k,k)$-module for each $\Gamma$-module $M$. If $M$ is a $\Gamma$-module algebra, then $\Ext^*_\Gamma(k,M)$ is an $\Ext^*_\Gamma(k,k)$-algebra. 

We proceed to define the cup product in Tate cohomology for Hopf algebras. For this, we need a unique (up to chain homotopy) $\Gamma$-linear extension of the fold map\index{fold map|(} in the category of chain complexes of~$\Gamma$-modules under $k$. Explicitly, the fold map~$\nabla$ is the induced map in the commutative diagram
\[
\begin{tikzcd}
k \otimes k \arrow[r,"i \otimes 1"] \arrow[d,"1 \otimes i"] & \widetilde{P}_\ast \arrow[d] \arrow[rdd , bend left , "\id"] & \\
\widetilde{P}_* \arrow[r] \arrow[rrd, bend right , "\id"] & \widetilde{P}_* \oplus_k \widetilde{P}_* \arrow[rd, dashed, "\nabla"] & \\
& & \widetilde{P}_* \,,
\end{tikzcd}
\]
where the inner square is a pushout diagram. Let us start with a more general result.\index{fold map|)}

\begin{lemma}
\label{lem:relativeextension}
Let $A_*$, $B_*$ and $C_*$ be chain complexes of $\Gamma$-modules, where we assume that $C_*$ is non-negative and exact. Let $i \: A_* \to B_*$ be an injective chain map and assume that $Q_* = \cok(i)$ is projective over $\Gamma$ in each homological degree. Then, for each chain map $f \: A_* \to C_*$ there is a chain map $g \: B_* \to C_*$ such that~$gi = f$. Moreover, this chain map is unique up to a chain homotopy that is zero on the image of $i$.
\begin{proof}
Consider the diagram
\[
\begin{tikzcd}
0 \arrow[r] & A_* \arrow[r,"i"] \arrow[d, "f"'] & B_* \arrow[r,"r"] \arrow[ld, dashed, "g"]& Q_* \arrow[r] & 0 \\
& C_* \,, & & & 
\end{tikzcd}
\]
where the top sequence is short exact in each homological degree. Since $C_*$ is non-negative, we must have $g_n = 0$ for $n < 0$. To construct the rest of the chain map we proceed by induction. Assume inductively that we have constructed $g_m$ satisfying
\[
g_m i_m = f_m \och g_{m-1} \partial = \partial g_m
\]
for all $m < n$. Since 
\[
0 = g_{n-2} \partial^2 = \partial g_{n-1} \partial
\]
and $C_*$ is exact we know that $g_{n-1} \partial$ lands in $\partial(C_n)$. Consider the diagram 
\[
\begin{tikzcd}
0 \arrow[r] & A_n \arrow[r,"i_n"] \arrow[d, "f_n"] & B_n \arrow[r,"r_n"] \arrow[ld, dashed, "g_n"] \arrow[d,"g_{n-1} \partial"] & Q_n \arrow[r] & 0 \\
& C_n \arrow[r] & \partial(C_n) & 
\end{tikzcd}
\]
in which we want to find a dashed map $g_n \: B_n \to C_n$ that makes both triangles commute. Since $Q_n$ is projective, the short exact sequence at the top of the diagram splits, and we can find $s_n \: Q_n \to B_n$ and $t_n \: B_n \to A_n$ such that
\[
i_n t_n + s_n r_n = \id_{B_n} \,.
\]
Moreover, we can find a map $h_n \: Q_n \to C_n$ such that~$g_{n-1}\partial s_n = \partial h_n$. We define $g_n \: B_n \to C_n$ by setting  
\[
g_n = f_n t_n + h_n r_n \,.
\]
This map satisfies
\[
g_n i_n = f_n t_n i_n + h_n r_n i_n = f_n + 0 = f_n
\]
and 
\begin{align*}
\partial g_n &= \partial f_n t_n + \partial h_n r_n \\ &= f_{n-1} \partial t_n + \partial h_n r_n \\ &= g_{n-1} i_{n-1} \partial t_n + g_{n-1} \partial s_n r_n \\ &= g_{n-1}\partial(i_{n} t_n + s_n r_n) \\ &= g_{n-1} \partial \,,
\end{align*}
which concludes the construction of $g$.

We now show that the map $g \: B_* \to C_*$ is unique up to a chain homotopy that is zero on $i(A_*)$. Let $g' \: B_* \to C_*$ be another chain map satisfying $f = g'i$. We want to show that we can find a chain homotopy between $k = g -g'$ and~$0$ that is zero on the image of $i$. That is, we want to find a collection of maps $H_n \: B_n \to C_{n+1}$ such that
\[
k_n = H_{n-1} \partial + \partial H_n \och H_n i_n = 0
\]  
for all $n$. Again, we use induction. Since $C_*$ is non-negative we must have $H_{n} = 0$ for $n < 0$. Assume inductively that we have constructed $H_m \: B_m \to C_{m+1}$ satisfying
\[
k_m = H_{m-1} \partial + \partial H_m \och H_mi_m  = 0
\]
for all $m < n$. Consider the map $k_n - H_{n-1} \partial \: B_n \to C_n$. Since 
\begin{align*}
\partial(k_n - H_{n-1} \partial) &= \partial k_n - \partial H_{n-1} \partial \\ &= \partial k_n - (k_{n-1} - H_{n-2} \partial) \partial \\ &= \partial k_n - k_{n-1} \partial \\ &= 0
\end{align*}
and $C_*$ is exact we know that $k_n - H_{n-1}\partial$ lands in $\partial(C_{n+1})$. Consider the diagram
\[
\begin{tikzcd}
& Q_n \arrow[d,"k_n s_n - H_{n-1} \partial s_n"] \arrow[dl,dashed,"\beta_n"']& \\
C_{n+1} \arrow[r] & \partial(C_{n+1}) \arrow[r] & 0
\end{tikzcd}
\]
in which we have a map $\beta_n$ since $Q_n$ is projective. We define
\[
H_n = \beta_n r_n \: B_n \longto C_{n+1} \,,
\]
which vanishes on the image of $i_n$ since $r_n i_n = 0$. Furthermore, we have 
\begin{align*}
H_{n-1} \partial + \partial H_n &= H_{n-1} \partial + \partial \beta_n r_n \\ &= H_{n-1} \partial + k_n s_n r_n - H_{n-1} \partial s_n r_n \\ & = H_{n-1} \partial + k_n - H_{n-1} \partial \\ &= k_n
\end{align*}
where the penultimate equality sign follows from the fact that $k_n$ and $H_{n-1}$ vanish on the image of $i_n$ and $i_{n-1}$, respectively, so that $k_n = k_n(i_n t_n + s_n r_n) = k_n s_n r_n$ and 
\begin{align*}
H_{n-1} \partial &= H_{n-1} \partial (i_n t_n + s_n r_n) \\
&= H_{n-1} \partial i_n t_n + H_{n-1} \partial s_n r_n \\
&= H_{n-1} i_n \partial t_n + H_{n-1} \partial s_n r_n \\
&= H_{n-1} \partial s_n r_n \,. \qedhere
\end{align*}
\end{proof}
\end{lemma}

\begin{proposition}
\label{prop:tildePfoldextension}
There is a unique, up to chain homotopy, $\Gamma$-linear chain map $\Phi \: \widetilde{P}_* \otimes \widetilde{P}_* \to \widetilde{P}_*$ that makes the diagram 
\[
\begin{tikzcd}
\widetilde{P}_* \oplus_k \widetilde{P}_* \arrow[rr, "\nabla"] \arrow[rd] & & \widetilde{P}_* \\
& \widetilde{P}_\ast \otimes \widetilde{P}_\ast \arrow[ru , "\Phi" , dashed] & 
\end{tikzcd}
\]
commute. 
\begin{proof}
This is an application of Lemma~\ref{lem:relativeextension}. The diagram we are considering is 
\[
\begin{tikzcd}
0 \arrow[r] & \widetilde{P}_* \oplus_k \widetilde{P}_* \arrow[r,"i"] \arrow[d, "\nabla"'] & \widetilde{P}_* \otimes \widetilde{P}_* \arrow[r,"r"] \arrow[ld, dashed, "\Phi"] & Q_* \arrow[r] & 0 \\
& \widetilde{P}_* \,. & & & 
\end{tikzcd}
\]
We only need to check that the cokernel of the map $i$ is projective over $\Gamma$ in each homological degree. By construction, the cokernel is the total cokernel of the commutative diagram
\[
\begin{tikzcd}
k \otimes k \arrow[r,"i \otimes 1"] \arrow[d,"1 \otimes i"'] & \widetilde P_* \otimes k \arrow[d,"1 \otimes i"] \\
k \otimes \widetilde{P}_* \arrow[r,"i \otimes 1"] & \widetilde{P}_* \otimes \widetilde{P}_*  \,.
\end{tikzcd}
\]
This can be calculated by computing the cokernels of the two horizontal maps followed by the cokernel of induced vertical map. Explicitly:
\begin{align*}
\cok(i) &\cong \cok(\cok(i \otimes 1) \longto \cok(i \otimes 1)) \\ & \cong \cok(1 \otimes i :P[1]_* \otimes k \longto P[1]_* \otimes \widetilde{P}_*) \\ &\cong P[1]_* \otimes P[1]_* \,.
\end{align*}
In particular, we note that the cokernel is a complex of projective~$\Gamma$-modules.
\end{proof}
\end{proposition}

We can now define a pairing on $\hm_\ast(-)$ using $\Phi$ and $\Psi$.  For $\Gamma$-modules $M$ and~$N$ the composite pairing
\begin{align*}
\widetilde{P}_* \otimes \Hom(P_*, M)
	\otimes \widetilde{P}_* \otimes &\Hom(P_*, N) \\
        & \overset{1 \otimes \tau \otimes 1}\longto
\widetilde{P}_* \otimes \widetilde{P}_*
	\otimes \Hom(P_*,M) \otimes \Hom(P_*, N) \\
	    & \overset{1 \otimes 1 \otimes \alpha}\longto \widetilde{P}_* \otimes \widetilde{P}_*
	\otimes \Hom(P_* \otimes P_*,M \otimes N) \\
        & \overset{\Phi \otimes \Psi^*}\longto
\widetilde{P}_* \otimes \Hom(P_*, M \otimes N)
\end{align*}
is $\Gamma$-linear, so it induces a pairing 
\[
\Hom_\Gamma(k,\hm_\ast(M)) \otimes \Hom_\Gamma(k,\hm_*(N)) \longto \Hom_\Gamma(k,\hm_\ast(M \otimes N))
\]
of $k$-module complexes. Note that the uniqueness of $\Phi$ up to chain homotopy guarantees that $\Phi \circ \tau \simeq \Phi$, and we have already observed that $\tau \circ \Psi \simeq \Psi$, which ensures that we get an associative, unital, and graded commutative pairing
\[
\smile \: \hatExt^*_\Gamma(k,M) \otimes \hatExt^*_\Gamma(k,N) \longto \hatExt^*_\Gamma(k,M \otimes N)
\]
after passing to homology. The inclusion $\Hom(P_*,M) \to \hm_*(M)$ provides us with a map
\[
\Ext_\Gamma^*(k,M) \longto \hatExt^*_\Gamma(k,M) \,.
\]
This map is compatible with the multiplicative structures we have defined above, in the sense of the following proposition.

\begin{proposition}
\label{prop:cupproductcompatible}
The two diagrams
\[
\xymatrix{
\Ext^{b_1}_\Gamma(k,M) \otimes \Ext^{b_2}_\Gamma(k,N) \ar[r]^-{\smile} \ar[d]
& \Ext^{b_1+b_2}_\Gamma(k,M \otimes N) \ar[d] \\
\Ext^{b_1}_\Gamma(k,M) \otimes \hatExt^{b_2}_\Gamma(k,N) \ar[r] \ar[d]
& \hatExt^{b_1 + b_2}_\Gamma(k,M \otimes N) \ar@{=}[d] \\
\hatExt^{b_1}_\Gamma(k,M) \otimes \hatExt^{b_2}_\Gamma(k, N) \ar[r]^-{\smile} & \hatExt^{b_1 + b_2}_\Gamma(k,M \otimes N)
}
\]
and
\[
\xymatrix{
\Ext^{b_1}_\Gamma(k,M) \otimes \Ext^{b_2}_\Gamma(k,N) \ar[r]^-{\smile} \ar[d]
& \Ext^{b_1+b_2}_\Gamma(k,M \otimes N) \ar[d] \\
\hatExt^{b_1}_\Gamma(k,M) \otimes \Ext^{b_2}_\Gamma(k,N) \ar[r] \ar[d]
& \hatExt^{b_1 + b_2}_\Gamma(k,M \otimes N) \ar@{=}[d] \\
\hatExt^{b_1}_\Gamma(k,M) \otimes \hatExt^{b_2}_\Gamma(k, N) \ar[r]^-{\smile} & \hatExt^{b_1 + b_2}_\Gamma(k,M \otimes N)
}
\]
commute.  In particular, it follows that  $\hatExt_\Gamma^*(k,k)$ is an $\Ext^*_\Gamma(k,k)$-algebra, and $\Ext^*_\Gamma(k,M) \to \hatExt^*_\Gamma(k,M)$ is an $\Ext^*_\Gamma(k,k)$-module homomorphism.  If $M$ is a~$\Gamma$-module algebra, then $\Ext^*_\Gamma(k,M) \to \hatExt^*_\Gamma(k, M)$ is an $\Ext^*_\Gamma(k,k)$-algebra homomorphism.
\begin{proof}
This follows from the commutative diagrams 
\[
\begin{tikzcd}
k \otimes k \arrow[r,equal] \ar[d]
& k \arrow[d] \\
k \otimes \widetilde{P}_\ast \arrow[r,equal] \ar[d]
& \widetilde{P}_\ast \arrow[d,equal] \\
\widetilde{P}_\ast \otimes \widetilde{P}_\ast \arrow[r,"\Phi"] & \widetilde{P}_*
\end{tikzcd} \quad \text{and} \quad 
\begin{tikzcd}
k \otimes k \arrow[r,equal] \ar[d]
& k \arrow[d] \\
\widetilde{P}_\ast \otimes k \arrow[r,equal] \ar[d]
& \widetilde{P}_\ast \arrow[d,equal] \\
\widetilde{P}_\ast \otimes \widetilde{P}_\ast \arrow[r,"\Phi"] & \widetilde{P}_* \,.
\end{tikzcd}
\]
\end{proof}
\end{proposition}
\index{cup product!on complete $\Ext$|)}

\section{Computation}
\label{subsec:computation}

In this section we look at a sample computation of the Tate cohomology of a Hopf algebra. Let $k$ be a graded commutative ring with an element~$\eta$ in degree 1 such that $2 \eta = 0$. We will consider the Hopf algebra
\[
\Gamma = k[s]/(s^2 = \eta s) \,, \quad |s|=1 \,.
\]
Here $s$ is a primitive element, so that comultiplication is given by $\psi = s \otimes 1 + 1 \otimes s$, counit by $\epsilon(s) =0$, and antipode by $\chi(s) = -s$. To clarify: our goal is to compute $\hatExt^*_\Gamma(k,M)$ where $M$ is a $\Gamma$-module. This situation naturally appears when we consider the Tate construction on a spectrum~$X$ with an action of the circle $\bT$. In this situation, we will have 
\[
\Gamma = \pi_*(\bS[\bT]) \quad \quad k = \pi_*(\bS) \quad \quad M = \pi_*(X) \,.
\]
See Proposition~\ref{prop:s2isetas} and Chapter~\ref{sec:Tatespseq}. 

A projective resolution $P_*$ of $k$ as a trivial $\Gamma$-module is 
\[
\xymatrix{
\cdots \ar[r] & \Gamma\{p_3\} \ar[r]^{\partial_3} & \Gamma\{p_2\} \ar[r]^-{\partial_2} & \Gamma\{p_1\} \ar[r]^-{\partial_1} & \Gamma\{p_0\}  \ar[r] & 0
}
\]
with the internal degree of the generator $p_b$ being $|p_b| = b$ and the total degree being~$\|p_b\|=2b$. As (right) $k$-modules we have $P_b = \Gamma\{p_b\} = k\{p_b,p_b s\}$ where~$\|p_b s\|=2b+1$. The boundary of the complex is given by
\[
\partial_{b+1}(p_{b+1}) = \begin{cases}
p_b s & \text{$b \geq 0$ even } \\
p_b (s+\eta) & \text{$b \geq 1$ odd }
\end{cases}
\] 
and the augmentation $\epsilon \: P_* \to k$ is given by $\epsilon(p_0) = 1$. 

By definition, the mapping cone $\widetilde{P}_\ast = \cone(P_* \to k)$ is isomorphic to the complex
\[
\xymatrix{
\cdots \ar[r] & \Gamma\{\tilde p_3\} \ar[r]^-{\tilde{\partial}_3} & \Gamma\{\tilde p_2\} \ar[r]^-{\tilde{\partial}_2} & \Gamma\{\tilde p_1\} \ar[r]^-{\tilde{\partial}_1} & k\{\tilde p_0\} \ar[r] & 0
}
\]
where the internal degrees of the generators are $|\tilde p_0| = 0$ and $|\tilde p_a| = a-1$ for $a \geq 1$, and the total degrees are $\|\tilde p_0\| = 0$ and $\|\tilde p_a\| = 2a-1$ for~$a \geq 1$. As before, we can also write this as a complex of (right)~$k$-modules $\Gamma\{\tilde p_a\}=k\{\tilde p_a, \tilde p_a s\}$ for $a \geq 1$, where~$|\tilde p_a s| = a$ and $\|\tilde p_a s\| = 2a$. The boundary is given by
\[
\tilde{\partial}_a(\tilde{p}_{a}) = \begin{cases}
\tilde{p}_0 & a=1 \\
-\tilde{p}_{a-1}s & a \geq 2 \text{ even } \\
-\tilde{p}_{a-1}(s + \eta) & a \geq 3 \text{ odd.}
\end{cases}
\] 
The chain complex we want to consider is the Tate complex $\hm_*(M)$, or rather, its~$\Gamma$-invariants. Recall that the Tate complex $\hm_*(M)$ is given in each homological degree by 
\[
\hm_c(M) = \bigoplus_{a+b = c} \widetilde{P}_a \otimes \Hom(P_{-b},M)
\]
with boundary given as 
\[
\partial_{\hm}(x,f) = \tilde{\partial}(x) \otimes f + (-1)^{\|x\|} x \otimes \partial^*(f) \,.
\]
When calculating the second term we also remember that 
\[
(\partial^*f)(v) = -(-1)^{\|f\|}f(\partial(v))
\]
for an element $f \in \Hom(P_*,M)$ since $M$ is concentrated in homological degree~$0$. 

It will also be useful to consider the Tate complex in its bicomplex version. In this case, we will write $(\hm_{*,*}(M),\partial^h , \partial^v)$ where
\[
\hm_{a,b}(M) = \widetilde{P}_a \otimes \Hom(P_{-b},M)
\]
and the horizontal and vertical boundaries are the first and the second term in the formula for the boundary in $\hm_*$, respectively. The total complex of this bicomplex is equal to the Tate complex, by definition. Moreover, let us write $(U_{*,*},\partial^h,\partial^v)$ for the restriction of the bicomplex to the $\Gamma$-invariants 
\[
U_{a,b} = \Hom_\Gamma(k,\widetilde{P}_a \otimes \Hom(P_{-b},M)) \,.
\]
We refer to the total complex of this bicomplex as $(U_*,\partial^h + \partial^v)$; it is isomorphic to $\Hom_\Gamma(k,\hm_*(M))$. Let us introduce some notation for the different elements in the bicomplex $U_{*,*}$ to keep our computations from becoming too messy.

\begin{notation}
Let $x$ be an element of $M$ and write  
\[
f_{b} \cdot x = \tilde{p}_0 \otimes \binom{p_b \mapsto x}{p_b s \mapsto xs}
\]
for an element in $\Hom_{\Gamma}(k,\widetilde{P_0} \otimes \Hom(P_b,M))$ for $b \geq 0$.
\end{notation}

\begin{notation}
Let $y$ be an element of $M$ and write
\[
g_{a,b} \cdot y = 
\tilde p_a \otimes \binom{p_b \mapsto y}{p_b s \mapsto ys}
	\ +\  \tilde p_a s \otimes \binom{p_b \mapsto 0}{p_b s \mapsto (-1)^{|y|}y}
\]
for an element in $\Hom_{\Gamma}(k,\widetilde{P_a} \otimes \Hom(P_b,M))$ for $a \geq 1$ and $b \geq 0$.
\end{notation}

\begin{notation}
Let $z$ be an element of $M$ and write 
\[
h_{a,b} \cdot z = 
\tilde{p}_a s \otimes \binom{p_b \mapsto z}{p_bs \mapsto z (s+\eta)} 
\]
for an element in $\Hom_{\Gamma}(k,\widetilde{P_a} \otimes \Hom(P_b,M))$ for $a \geq 1$ and $b \geq 0$. 
\end{notation}

It is a straight-forward computation to check that these are indeed $\Gamma$-invariant elements, in the sense that 
\[
(g_{a,b} \cdot y) \cdot s = 0 \,,
\]
and analogously for $f_b \cdot x$ and $h_{a,b} \cdot z$, using the following lemma.

\begin{lemma}
\label{lem:primconj}
The (right) conjugate action of $s$ on an element $f$ of $\Hom(M,N)$ is given by
\[
(fs)(m) = (-1)^{\|m\|}(f(m)s - f(ms)) \,.
\]
\begin{proof}
Recall that the characterising property of the conjugate action is that it is the action on a function object $\Hom(M,N)$ such that the evaluation pairing $\ev \: \Hom(M,N) \otimes M \to N$ is $\Gamma$-linear. Explicitly, the $\Gamma$-action on $\Hom(M,N)$ makes the diagram 
\[
\xymatrix{
\Hom(M,N) \otimes M \otimes \Gamma \ar[r]^-{\ev \otimes 1} \ar[d]_-{\rho_{\Hom(M,N) \otimes M}} & N \otimes \Gamma \ar[d]^-{\rho_N} \\
\Hom(M,N) \otimes M \ar[r]^-{\ev} & N
}
\] 
commute. The top composition sends a generic element to 
\[
f \otimes m \otimes s \mapsto f(m) \otimes s \mapsto f(m)s \,,
\]
while the bottom composition sends the generic element to 
\[
f \otimes m \otimes s \mapsto f \otimes ms + (-1)^{\|m\|} fs \otimes m \mapsto f(ms) + (-1)^{\|m\|}(fs)(m) \,.
\]
These must agree, which necessarily gives us the assertion.
\end{proof}
\end{lemma}

Furthermore, these form an `$M$-basis' of the $\Gamma$-invariants of the Tate complex in the sense of the following proposition.

\begin{proposition}
\label{prop:Mbasis}
Let $b \geq 0$ and $a \geq 1$. There are $k$-module isomorphisms 
\begin{align*}
\Sigma^{-b} M &\overset{\cong}\longto \Hom_\Gamma(k,\widetilde{P}_0 \otimes \Hom(P_b,M)) \\ x &\mapsto f_b \cdot x
\end{align*}
and 
\begin{align*}
\Sigma^{a-b-1} M \oplus \Sigma^{a-b} M &\overset{\cong}\longto \Hom_\Gamma(k,\widetilde{P}_a \otimes \Hom(P_b,M)) \\ (y,z) &\mapsto g_{a,b} \cdot y + h_{a,b} \cdot z \,.
\end{align*}
\begin{proof}
The maps are clearly injective, so we only need to show that they are surjective. 

A general element in $\widetilde{P}_0 \otimes \Hom(P_b,M)$ is on the form 
\[
\tilde{p}_0 \otimes \binom{p_b \mapsto x}{p_b s \mapsto y} \,.
\]
By Lemma~\ref{lem:primconj}, the right action of the primitive element $s$ on such an element is 
\[
\tilde{p}_0 \otimes \binom{p_b \mapsto xs-y}{p_b s \mapsto y\eta - ys} \,.
\]
For our original element to be $\Gamma$-invariant this must be zero, which gives us $y = xs$. In other words, a $\Gamma$-invariant element in $\widetilde{P}_0 \otimes \Hom(P_b,M)$ can be written $f_b \cdot x$, where we let $x$ range throughout~$M$. The grading suspension appearing in the isomorphism makes sure that this is actually a map of graded~$k$-modules. Indeed, the internal degree of our element is 
\[
|f_b \cdot x| = |\tilde{p}_0| + |x| - |p_b| = |x| - b \,.
\]

A general element in $\widetilde{P}_a \otimes \Hom(P_b,M)$ is on the form
\[
\tilde p_a \otimes \binom{p_b \mapsto x}{p_b s \mapsto y}
	\ +\  \tilde p_a s \otimes \binom{p_b \mapsto z}{p_b s \mapsto w} \,.
	\]
We assume that this is a homogeneous element, so that $|y| = |x| + 1$, $|z| = |x|-1$, and $|w| = |x|$. Letting the primitive element $s$ act on this element from the right we obtain
\begin{align*}
& (-1)^{|x|} \tilde p_a s \otimes \binom{p_b \mapsto x}{p_b s \mapsto y} +
\tilde p_a \otimes \binom{p_b \mapsto xs-y}{p_b s \mapsto y\eta - ys}
	\ \\ &-\  (-1)^{|x|} \tilde p_a s \eta \otimes \binom{p_b \mapsto z}{p_b s \mapsto w}
	\ +\  \tilde p_a s \otimes \binom{p_b \mapsto zs-w}{p_b s \mapsto w\eta - ws} \,.
\end{align*}
For our element to be $\Gamma$-invariant we want this to add up to zero. In other words, we need to solve the following system of equations
\[
\begin{cases}
xs - y &= 0 \\
y\eta - ys &= 0 \\
(-1)^{|x|} x + z \eta + z s - w &= 0 \\
(-1)^{|x|} y + w \eta + w \eta - w s &= 0 \,.
\end{cases} 
\]
It is straight-forward to check that the solutions are given by the two independent equations 
\[
\begin{cases}
y &= xs \\
w &= (-1)^{|x|} x + z \eta + z s , \\
\end{cases}
\]
which tells us that a $\Gamma$-invariant element can be written
\[
g_{a,b} \cdot x + h_{a,b} \cdot z
\]
where we are free to vary $x$ and $z$ in $M$. The suspensions in the source of the~$k$-isomorphism are again there to make sure that the grading is preserved by the isomorphism. Indeed, 
\[
|g_{a,b} \cdot x| = |\tilde{p}_a| + |x| - |p_b| = a - 1 + |x| -b  
\]
and 
\[
|h_{a,b} \cdot z| = |\tilde{p}_a s| + |z| - |p_b| = a + |z| - b \,.
\qedhere
\]
\end{proof}
\end{proposition}

We now need to figure out what the boundary on these generic $\Gamma$-invariant elements looks like. Keeping track of all the signs we end up with the following description of the horizontal and vertical boundaries in terms of our `$M$-basis'.

\begin{lemma}
\label{lem:fboundary}
The horizontal boundary on $f_b \cdot x$ is given by
\[
\partial^h(f_b \cdot x) = 0
\]
and the vertical boundary is given by
\[
\partial^v(f_b \cdot x) = \begin{cases}
-(-1)^{|x|}f_{b+1} \cdot xs & b \geq 0 \text{ even} \\
-(-1)^{|x|}f_{b+1} \cdot x(s+\eta) & b \geq 1 \text{ odd.}
\end{cases}
\]
\end{lemma}

\begin{lemma}
\label{lem:gboundary}
The horizontal boundary on $g_{a,b} \cdot y$ is given by
\[
\partial^h(g_{a,b} \cdot y) = 
\begin{cases}
 f_{b} \cdot y & \text{for $a=1$} \\
 - h_{a-1,b} \cdot y & \text{for $a\ge2$ even} \\
 g_{a-1,b} \cdot y\eta - h_{a-1,b} \cdot y
	& \text{for $a\ge3$ odd}
\end{cases}
\]
and the vertical boundary is given by
\[
\partial^v(g_{a,b} \cdot y) = \begin{cases}
(-1)^{|y|}g_{a,b+1} \cdot ys
	\ +\ h_{a,b+1} \cdot y & \text{for $b\ge0$ even} \\
(-1)^{|y|}g_{a,b+1} \cdot y(s+\eta)
	\ +\ h_{a,b+1} \cdot y & \text{for $b\ge1$ odd.} \\
\end{cases}
\]
\end{lemma}

\begin{lemma}
\label{lem:hboundary}
The horizontal boundary on $h_{a,b} \cdot z$ is given by
\[
\partial^h(h_{a,b} \cdot z) = \begin{cases}
h_{a-1,b} \cdot z\eta & \text{for $a \ge 2$ even} \\
0 & \text{for $a \ge 1$ odd}  \end{cases} .
\]
and the vertical boundary is given by
\[
\partial^v(h_{a,b} \cdot z) = 
\begin{cases}
- (-1)^{|z|}h_{a,b+1} \cdot z(s+\eta) & \text{for $b\ge0$ even} \\
- (-1)^{|z|}h_{a,b+1} \cdot zs & \text{for $b\ge1$ odd.}
\end{cases}
\]
\end{lemma}

We calculate the homology of $U_*$ by filtering the first tensor factor of $U_{*,*}$ and using the spectral sequence for the total complex of a bicomplex:
\[
E^1_{a,-b} = H_{-b}(\Hom_\Gamma(k,\widetilde{P}_a \otimes \Hom(P_*,M))) \Rightarrow H_{a-b}(\Hom_\Gamma(k,\widetilde{P}_* \otimes \Hom(P_*,M))) \,. 
\]
The bicomplex $(U_{*,*},\partial^h , \partial^v)$ is displayed in Figure~\ref{fig:TatecmpxforT} for the convenience of the reader. 

\begin{remark}
Let us clarify how to interpret Figure~\ref{fig:TatecmpxforT} and the matrix notation appearing in it. The horizontal and vertical boundaries are given in terms of the~`$M$-bases' $\{f_{b}\}$ for $U_{0,-b}$ and $\{g_{a,b}, h_{a,b} \}$ for~$U_{a,-b}$. We record $f_b \cdot x$ and $g_{a,b} \cdot y + h_{a,b} \cdot z$ as the column vectors
\[
\begin{bmatrix} x \end{bmatrix} \och \begin{bmatrix} y \\ z \end{bmatrix} \,,
\]
respectively. The boundaries are then indicated by multiplication with the corresponding matrices appearing in the figure. Multiplication is done, as is usual, with the matrix on the left hand side. So, to clarify: A vertical boundary $\partial^v \: U_{a,-b} \to U_{a,-b-1}$ recorded as a $2 \times 2$-matrix (with entries in $\Gamma$) and multiplied with the relevant column vector (with entries in $M$)
\[
\begin{bmatrix}
i & j \\
k & \ell
\end{bmatrix}
\begin{bmatrix}
y \\
z
\end{bmatrix} = 
\begin{bmatrix}
iy + jz \\
ky + \ell z
\end{bmatrix}
\]
indicates that this boundary is given as
\[
g_{a,b} \cdot y + h_{a,b} \cdot z \mapsto g_{a,b+1} \cdot (iy + jz) + h_{a,b+1} \cdot (ky + \ell z) \,.
\] 
Note that in this convention $\Gamma$ ends up acting on $M$ from the left, through the twist isomorphism followed by the right action. To see that the boundaries given in the matrix notation actually agree with the ones given in Lemma~\ref{lem:fboundary}, Lemma~\ref{lem:gboundary}, and Lemma~\ref{lem:hboundary}, we have to switch the position of the~$\Gamma$-values ($i$,~$j$,~$k$, and~$\ell$) and the $M$-values ($y$ and $z$), which typically introduces a sign. For example, the boundary $\partial^v \: U_{a,-b} \to U_{a,-b-1}$ for even $b$ is recorded in the figure as 
\[
\begin{bmatrix}
s & 0 \\
1 & -(s + \eta)
\end{bmatrix} \,.
\]
Left multiplication of this matrix with the column vector corresponding to $g_{a,b} \cdot y$ gives 
\[
\begin{bmatrix}
s & 0 \\
1 & -(s + \eta)
\end{bmatrix} \begin{bmatrix}
y \\
0
\end{bmatrix} = \begin{bmatrix}
sy \\
y
\end{bmatrix} \,,
\]
which tells us that this vertical boundary is given by 
\[
g_{a,b} \cdot y \mapsto g_{a,b+1} \cdot sy + h_{a,b+1} \cdot y = (-1)^{|y|} g_{a,b+1} \cdot ys + h_{a,b+1} \cdot y \,,
\]
which is indeed in agreement with Lemma~\ref{lem:gboundary}.
\end{remark}

\begin{figure}
$$
\xymatrix@R+2.5pc@C+1.5pc{
U_{0,0} \ar[d]_-{\bsm -s \esm}
& U_{1,0} \ar[l]_-{\bsm 1 & 0 \esm}
        \ar[d]_-{\bsm s & 0 \\ 1 & -(s+\eta) \esm}
& U_{2,0} \ar[l]_-{\bsm 0 & 0 \\ -1 & \eta \esm}
        \ar[d]_-{\bsm s & 0 \\ 1 & -(s+\eta) \esm}
& U_{3,0} \ar[l]_-{\bsm \eta & 0 \\ -1 & 0 \esm}
        \ar[d]_-{\bsm s & 0 \\ 1 & -(s+\eta) \esm}
& \dots \ar[l]_-{\bsm 0 & 0 \\ -1 & \eta \esm} \\
U_{0,-1} \ar[d]_-{\bsm -(s+\eta) \esm}
& U_{1,-1} \ar[l]_-{\bsm 1 & 0 \esm}
        \ar[d]_-{\bsm s+\eta & 0 \\ 1 & -s \esm}
& U_{2,-1} \ar[l]_-{\bsm 0 & 0 \\ -1 & \eta \esm}
        \ar[d]_-{\bsm s+\eta & 0 \\ 1 & -s \esm}
& U_{3,-1} \ar[l]_-{\bsm \eta & 0 \\ -1 & 0 \esm}
        \ar[d]_-{\bsm s+\eta & 0 \\ 1 & -s \esm}
& \dots \ar[l]_-{\bsm 0 & 0 \\ -1 & \eta \esm} \\
U_{0,-2} \ar[d]_-{\bsm -s \esm}
& U_{1,-2} \ar[l]_-{\bsm 1 & 0 \esm}
        \ar[d]_-{\bsm s & 0 \\ 1 & -(s+\eta) \esm}
& U_{2,-2} \ar[l]_-{\bsm 0 & 0 \\ -1 & \eta \esm}
        \ar[d]_-{\bsm s & 0 \\ 1 & -(s+\eta) \esm}
& U_{3,-2} \ar[l]_-{\bsm \eta & 0 \\ -1 & 0 \esm}
        \ar[d]_-{\bsm s & 0 \\ 1 & -(s+\eta) \esm}
& \dots \ar[l]_-{\bsm 0 & 0 \\ -1 & \eta \esm} \\
\vdots
& \vdots
& \vdots
& \vdots
&
}
$$
\caption{The bicomplex $(U_{*,*}, \partial^h, \partial^v)$
	for $\Gamma = k[s]/(s^2=\eta s)$}
\label{fig:TatecmpxforT}
\end{figure}

Before we explicitly compute the first page of the spectral sequence for the bicomplex $U_{*,*}$, we again introduce some notation.

\begin{notation}
Let $z$ be an element of $M$ and write 
\begin{align*}
u_a \cdot z &= -(-1)^{|z|} g_{a,0} \cdot z(s+\eta) - h_{a,0} \cdot z \\
&= -(-1)^{|z|} \tilde{p}_a \otimes \binom{p_0 \mapsto z(s+\eta)}{p_0 s \mapsto 0} - \tilde{p}_a s \otimes \binom{p_0 \mapsto z}{p_0 s \mapsto 0}
\end{align*}
for the specified element in $\Hom_{\Gamma}(k,\widetilde{P}_a \otimes \Hom(P_b,M))$ for $a \geq 1$ and $b \geq 0$.
\end{notation}

\begin{proposition}
\label{prop:E1pageTatecomputation}
The $E^1$-page of the bicomplex spectral sequence for $U_{*,*}$ is given as 
\[
H_{-b}(\Hom_\Gamma(k,\widetilde{P}_0 \otimes \Hom(P_*,M))) \cong \begin{cases}
f_0 \cdot \ker(s) & \text{for $b=0$,} \\
\\
{\displaystyle f_b \cdot
	\frac{\ker(s+\eta)}{\im(s)}} & \text{for $b\ge1$ odd,} \\
\\
{\displaystyle f_b \cdot
	\frac{\ker(s)}{\im(s+\eta)}} & \text{for $b\ge2$ even,} \\
\end{cases}
\]
when $a = 0$, and as
\[
H_{-b}(\Hom_\Gamma(k,\widetilde{P}_a \otimes \Hom(P_*,M))) \cong \begin{cases}
u_a \cdot M & \text{for $b=0$} \\
0 & \text{otherwise,}
\end{cases}
\]
when $a \geq 1$. 
\begin{proof}
This is essentially an exercise in linear algebra using the matrices in Figure~\ref{fig:TatecmpxforT}.  The kernels of the boundaries are computed by computing the null spaces of the corresponding matrices. Similarly, the images of boundaries are computed by computing the column spaces of the matrices. Let us give the details for the~$a \geq 1$ case, as the $a=0$ case is directly visible by inspecting the figure. The null spaces of the matrices
\[
\begin{bmatrix}
s & 0 \\
1 & -(s + \eta)
\end{bmatrix} \och \begin{bmatrix}
s + \eta & 0 \\
1 & -s
\end{bmatrix} \,
\] 
are generated by the column vectors 
\[
\begin{bmatrix}
s + \eta \\
1
\end{bmatrix} \och 
\begin{bmatrix}
s \\
1
\end{bmatrix} \,,
\]
respectively, and the column spaces are generated by the column vectors 
\[
\begin{bmatrix}
s \\
1
\end{bmatrix} \och \begin{bmatrix}
s + \eta \\
1
\end{bmatrix} \,,
\]
respectively. From this is follows that the homology is concentrated in homological degree $b=0$, in the~$a \geq 1$ case. Here it consists of elements of the form
\[
g_{a,0} \cdot (s + \eta)z + h_{a,0} \cdot z = (-1)^{|z|} g_{a,0} \cdot z(s+ \eta) + h_{a,0} \cdot z
\]
for varying $z$ in $M$. For reasons concerning the multiplicative structure, we have decided to denote the above element by $- u_a \cdot z$. 
\end{proof}
\end{proposition}

In particular, note that the above result tells us that the $E^1$-page of the spectral sequence is concentrated around the boundary of the fourth quadrant. The~$d^1$-differential $d^1 \: E^1_{a,-b} \to E^1_{a-1,-b}$ in the spectral sequence is induced by the horizontal boundary, and is by degree reasons only non-zero on the positive $a$-axis. There it is given by
\[
d^1(u_a \cdot z) = \begin{cases}
- (-1)^{|z|}f_0 \cdot z(s+\eta) & \text{for $a=1$} \\
-(-1)^{|z|}u_{a-1} \cdot zs & \text{for $a \geq 2$ even} \\
-(-1)^{|z|}u_{a-1} \cdot z(s+\eta) & \text{for $a \geq 3$ odd} 
\end{cases}
\] 
by using Lemma~\ref{lem:gboundary} and Lemma~\ref{lem:hboundary}. We conclude that the second page of the spectral sequence is concentrated along the $a$- and $b$-axes and that
\[
E_{a,-b}^2 \cong \begin{cases}
{\displaystyle f_b \cdot
	\frac{\ker(s)}{\im(s+\eta)}} & \text{for $a=0$ and $b\ge0$ even,} \\
\\
{\displaystyle f_b \cdot
	\frac{\ker(s+\eta)}{\im(s)}} & \text{for $a=0$ and $b \geq 1$ odd,} \\
\\
{\displaystyle u_{a} \cdot
	\frac{\ker(s)}{\im(s+\eta)}} & \text{for $b=0$ and $a\ge2$ even,} \\
\\
{\displaystyle u_a \cdot
	\frac{\ker(s+\eta)}{\im(s)}} & \text{for $b=0$ and $a\ge1$ odd.} \\

\end{cases}
\]
There is no room for further differentials, and the infinite cycles along the upper and left hand edges are necessarily also $d^1$-cycles in the total complex $U_*$. We conclude:

\begin{proposition}
\[
\hatExt^c_\Gamma(k,M) \cong \begin{cases}
{\displaystyle f_c \cdot
	\frac{\ker(s)}{\im(s+\eta)}} & \text{for $c\geq 0$ even,} \\
\\
{\displaystyle f_c \cdot
	\frac{\ker(s+\eta)}{\im(s)}} & \text{for $c \geq 1$ odd,} \\
\\
{\displaystyle u_{-c} \cdot
	\frac{\ker(s)}{\im(s+\eta)}} & \text{for $c \leq -2$ even,} \\
\\
{\displaystyle u_{-c} \cdot
	\frac{\ker(s+\eta)}{\im(s)}} & \text{for $c \leq -1$ odd.} \\
\end{cases}
\] 
\end{proposition}

Now all that remains is to describe the multiplicative structure. That is, given two~$\Gamma$-modules $M$ and~$N$ we want to determine the cup product
\[
\smile \: \hatExt^{c_1}_{\Gamma}(k,M) \otimes \hatExt^{c_2}_\Gamma(k,N) \longto \hatExt^{c_1 + c_2}_\Gamma(k,M \otimes N).
\]
In order to do so we need a $\Gamma$-linear chain map $\Psi \: P_* \to P_* \otimes P_*$ covering the identity of $k$, and a~$\Gamma$-linear chain map $\Phi \: \widetilde{P}_* \otimes \widetilde{P}_* \to \widetilde{P}_*$ extending the fold map, as per Section~\ref{sec:multTate}. 

\begin{lemma}
\label{lem:Psichainmap}
A $\Gamma$-linear chain map $\Psi \: P_* \to P_* \otimes P_*$ that covers the identity is given by 
\[
\Psi(p_{b}) = \sum_{b_1 + b_2 = b} p_{b_1} \otimes p_{b_2} .
\] 
By $\Gamma$-linearity we have
\[
\Psi(p_{b} s) = \sum_{b_1 + b_2 = b} p_{b_1} s \otimes p_{b_2} + p_{b_1} \otimes p_{b_2} s.
\]
This chain map is cocommutative: $\Psi = \tau \Psi$.
\begin{proof}
Note that 
\[
\partial(p_b) = p_{b-1}(s + (b-1)\eta)
\]
for $b\ge1$. To verify that $\Psi$, as specified in the statement of the lemma, is a chain map, we must show that
\begin{align*}
\partial(\Psi(p_b)) &= \sum_{b_1+b_2=b} \partial(p_{b_1}) \otimes p_{b_2}
	+ p_{b_1} \otimes \partial(p_{b_2}) \\
	&= \sum_{b_1+b_2=b} p_{b_1-1}(s + (b_1-1)\eta) \otimes p_{b_2}
        + p_{b_1} \otimes p_{b_2-1}(s + (b_2-1)\eta)
\end{align*}
is equal to
\[
\Psi(\partial(p_b)) = \Psi(p_{b-1}s) + \Psi(p_{b-1})(b-1)\eta  \,.
\]
Here
\[
\sum_{b_1+b_2=b} p_{b_1-1}s \otimes p_{b_2}
        + p_{b_1} \otimes p_{b_2-1}s
	= \Psi(p_{b-1})s \,,
\]
so it remains to check that
\[
\sum_{b_1+b_2=b} p_{b_1-1}(b_1-1)\eta  \otimes p_{b_2}
        + p_{b_1} \otimes p_{b_2-1} (b_2-1)\eta
= \Psi(p_{b-1})(b-1)\eta  \,.
\]
When $b$ is odd the terms of the left hand side cancel in pairs, and the
right hand side is zero.  When~$b$ is even only the terms with $b_1$
and $b_2$ both even contribute to the left hand side, and these add up
to~$\Psi(p_{b-1})\eta$, as required.  Finally,
\[
(\epsilon \otimes \epsilon)(\Psi(p_0)) = 1 = \epsilon(p_0) \,,
\]
so $\Psi$ is indeed a chain map covering $k \otimes k = k$.
Cocommutativity of $\Psi$ is clear from the explicit formulas.
\end{proof}
\end{lemma}
 
\begin{lemma}
\label{lem:Phichainmap}\index{fold map}
A $\Gamma$-linear chain map $\Phi \: \widetilde{P}_* \otimes \widetilde{P}_* \to \widetilde{P}_*$ that extends the fold map is given by
\begin{align*}
\Phi(\tilde p_{a_1} \otimes \tilde p_{a_2}) &= 0 \\
\Phi(\tilde p_{a_1} \otimes \tilde p_{a_2}s) &= -\tilde p_a \\
\Phi(\tilde p_{a_1}s \otimes \tilde p_{a_2}) &= -\tilde p_a \\
\Phi(\tilde p_{a_1}s \otimes \tilde p_{a_2}s) &= -\tilde p_a(s+\eta)
\end{align*}
for $a_1,a_2 \ge1$ and $a = a_1+a_2$. Furthermore,
\begin{align*}
\Phi(\tilde p_{0} \otimes \tilde p_{a_2}) &= \tilde p_{a_2} \\
\Phi(\tilde p_{0} \otimes \tilde p_{a_2} s) &= \tilde p_{a_2}s \\
\Phi(\tilde p_{a_1} \otimes \tilde p_{0}) &= \tilde p_{a_1} \\
\Phi(\tilde p_{a_1}s \otimes \tilde p_{0}) &= \tilde p_{a_1} s
\end{align*}
and $\Phi(\tilde p_0 \otimes \tilde p_0) = \tilde p_0$.
This chain map is commutative: $\Phi = \Phi \tau$.
\begin{proof}
Note that the differential in the chain complex $\widetilde{P}_*$ can be described as
\[
\tilde\partial(\tilde p_a) = \begin{cases} \tilde p_0 & \text{for $a =1$} \\ -\tilde p_{a-1}(s + a\eta) & \text{for $a \geq 2$.}\end{cases}
\]
To check that~$\Phi$, as specified in the statement of the lemma, is $\Gamma$-linear, we observe that
\begin{align*}
\Phi((\tilde p_{a_1} \otimes \tilde p_{a_2})s)
&= \Phi(\tilde p_{a_1} \otimes \tilde p_{a_2}s
	- \tilde p_{a_1} s\otimes \tilde p_{a_2}) \\
&= -\tilde p_a + \tilde p_a \\ &= 0 \\ &= \Phi(\tilde p_{a_1} \otimes \tilde p_{a_2})s
\end{align*}
and that
\begin{align*}
\Phi((\tilde p_{a_1} \otimes \tilde p_{a_2}s)s)
&= \Phi(\tilde p_{a_1}s \otimes \tilde p_{a_2}s
        + \tilde p_{a_1} \otimes \tilde p_{a_2} \eta s) \\
&= -\tilde p_a(s+\eta) \ + \tilde p_a \eta \\ &= -\tilde{p}_a s \\ &= \Phi(\tilde p_{a_1} \otimes \tilde p_{a_2}s)s \,.
\end{align*}
The check that $\Phi$ is a chain map is contained in the computations
\begin{equation*}
\begin{split}
\Phi(\partial_{\widetilde{P}_* \otimes \widetilde{P}_*}(\tilde p_{a_1} \otimes \tilde p_{a_2})) &=
\Phi(\tilde\partial(\tilde p_{a_1}) \otimes \tilde p_{a_2})
	- \Phi(\tilde p_{a_1} \otimes \tilde\partial(\tilde p_{a_2})) \\
&=
\begin{cases}
-\Phi(\tilde p_{a_1-1}(s + a_1\eta) \otimes \tilde p_{a_2}) \\ \quad \quad \quad
	+ \Phi(\tilde p_{a_1} \otimes \tilde p_{a_2-1}(s + a_2\eta))
	& \text{for $a_1$, $a_2\ge2$,} \\
 \Phi(\tilde p_0 \otimes \tilde p_{a_2}) \\ \quad \quad \quad
	+ \Phi(\tilde p_1 \otimes \tilde p_{a_2-1}(s + a_2\eta) )
	& \text{for $a_1=1$, $a_2\ge2$,} \\
- \Phi(\tilde p_{a_1-1} (s + a_1\eta) \otimes \tilde p_1) \\ \quad \quad \quad
	- \Phi(\tilde p_{a_1} \otimes \tilde p_0)
	& \text{for $a_1\ge2$, $a_2=1$}
\end{cases} \\
&=
\tilde p_{a-1} - \tilde p_{a-1} \\ &= 0 \\ &= \tilde\partial(\Phi(\tilde p_{a_1} \otimes \tilde p_{a_2}))
\end{split}
\end{equation*}
and 
\begin{equation*}
\begin{split}
\Phi(\partial_{\widetilde{P}_* \otimes \widetilde{P}_*}(\tilde p_{a_1} \otimes \tilde p_{a_2} s))
&=
\Phi(\tilde\partial(\tilde p_{a_1}) \otimes \tilde p_{a_2} s)
	- \Phi(\tilde p_{a_1} \otimes \tilde\partial(\tilde p_{a_2} s)) \\
&=
\begin{cases}
-\Phi(\tilde p_{a_1-1}(s + a_1\eta) \otimes \tilde p_{a_2}s) \\ \quad \quad \quad
	+ \Phi(\tilde p_{a_1} \otimes \tilde p_{a_2-1}(s + a_2\eta)s)
	& \text{for $a_1$, $a_2\ge2$,} \\
\Phi(\tilde p_0 \otimes \tilde p_{a_2}s)
	 \\ \quad \quad \quad + \Phi(\tilde p_1 \otimes \tilde p_{a_2-1} (s + a_2\eta)s )
	& \text{for $a_1=1$, $a_2\ge2$,} \\
-\Phi(\tilde p_{a_1-1} (s + a_1\eta) \otimes \tilde p_1s) \\ \quad \quad \quad
	- \Phi(\tilde p_{a_1} \otimes \tilde p_0s)
	& \text{for $a_1\ge2$, $a_2=1$}
\end{cases} \\
&=
\tilde p_{a-1} (s + a\eta) \\ &= -\tilde\partial(\tilde p_a) \\ &= \tilde\partial(\Phi(\tilde p_{a_1} \otimes \tilde p_{a_2}s)) \,.
\end{split}
\end{equation*}
Commutativity of $\Phi$ is clear from the explicit formulas.
\end{proof}
\end{lemma}

Now we want to use the above chain maps to compute the multiplicative structure. Recall that the cup product is induced by the composite pairing
\begin{align*}
\widetilde{P}_* \otimes \Hom(P_*, M)
	\otimes \widetilde{P}_* \otimes &\Hom(P_*, N) \\
        & \xrightarrow{1 \otimes \tau \otimes 1}
\widetilde{P}_* \otimes \widetilde P_*
	\otimes \Hom(P_*,M) \otimes \Hom(P_*, N) \\
	    &\xrightarrow{1 \otimes 1 \otimes \alpha}
\widetilde{P}_* \otimes \widetilde{P}_*
	\otimes \Hom(P_* \otimes P_*, M \otimes N) \\
        & \xrightarrow{\Phi \otimes \Psi^*}
\widetilde P_* \otimes \Hom(P_*, M \otimes N) \,.
\end{align*} 
There are two signs to be wary of here; the first one comes from twisting the factor $\Hom(P_*,M)$ past the second~$\widetilde{P}_*$-factor, and the second sign comes from using the canonical map~$\alpha$. Please refer to Equation~\ref{eq:sign}.
We note that the cocommutativity of $\Gamma$, symmetry of $\alpha$,
commutativity of $\Phi$ and cocommutativity of $\Psi$ imply that this
particular model for the cochain cup product is graded commutative.
The cup product computations, which can be found in the lemmas below, are straight-forward computations. We only include the verification of two of the lemmas, as the other two are very similar. 

\begin{lemma}
\label{lem:cupinExt}
Let $f_{b_1} \cdot m$ be a cycle with homology class in $\hatExt^{b_1}_{\Gamma}(k,M)$ and let~$f_{b_2} \cdot n$ be a cycle with homology class in $\hatExt^{b_2}_{\Gamma}(k,N)$. The cup product of these is the cycle 
\[
f_{b_1} \cdot m \smile f_{b_2} \cdot n = f_{b_1 + b_2} \cdot m \otimes n 
\]
with homology class in $\hatExt^{b_1 + b_2}_\Gamma(k,M \otimes N)$.
\end{lemma}

\begin{lemma}
Let $u_{a_1} \cdot m$ be a cycle with homology class in $\hatExt^{-a_1}_{\Gamma}(k,M)$ and let~$u_{a_2} \cdot n$ be a cycle with homology class in $\hatExt^{-a_2}_{\Gamma}(k,N)$. The cup product of these is the cycle 
\[
u_{a_1} \cdot m \smile u_{a_2} \cdot n = u_{a_1 + a_2} \cdot m \otimes n \,
\]
with homology class in $\hatExt^{-a_1 - a_2}_\Gamma(k,M \otimes N)$.
\end{lemma}

\begin{lemma}
Let $f_{0} \cdot m$ and $u_a \cdot n$ be cycles with homology classes in~$\hatExt^{0}_{\Gamma}(k,M)$ and $\hatExt^{-a}_{\Gamma}(k,N)$, respectively. The cup product of these is the cycle 
\[
f_{0} \cdot m \smile u_{a} \cdot n = u_{a} \cdot m \otimes n \,
\]
with homology class in $\hatExt^{-a}_\Gamma(k,M \otimes N)$. By graded commutativity we have
\[
u_a \cdot n \smile f_0 \cdot m = u_a \cdot n \otimes m \,.
\]
\begin{proof}
Since $f_0 \cdot m$ is assumed to be a cycle in $U_*$ we know that $m$ is an element in~$\ker(s)$. An explicit description of this cycle is then 
\[
f_0 \cdot m = \tilde{p}_0 \otimes \binom{p_0 \mapsto m}{p_0 s \mapsto 0} \,.
\]
The first map $1 \otimes \tau \otimes 1$ in the composite pairing twists the second tensor factor past the third one, so that 
\begin{multline*}
\tilde p_0 \otimes \binom{p_0 \mapsto m}{p_0 s \mapsto 0} \otimes \left(-(-1)^{|n|} \tilde p_a \otimes \binom{p_0 \mapsto n(s + \eta)}{p_0 s \mapsto 0} - \tilde{p}_a s \otimes \binom{p_0 \mapsto n}{p_0 s \mapsto 0} \right) \\ \mapsto -(-1)^{|m| + |n|} \tilde{p}_0 \otimes \tilde{p}_a \otimes \binom{p_0 \mapsto m}{p_0 s \mapsto 0} \otimes \binom{p_0 \mapsto n(s + \eta)}{p_0 s \mapsto 0} \\ - \tilde{p}_0 \otimes \tilde{p}_a s \otimes \binom{p_0 \mapsto m}{p_0 s \mapsto 0} \otimes \binom{p_0 \mapsto n}{p_0 s \mapsto 0} \,.
\end{multline*}
The second map in the composite is $1 \otimes 1 \otimes \alpha$, so that 
\begin{multline*}
-(-1)^{|m| + |n|} \tilde{p}_0 \otimes \tilde{p}_a \otimes \binom{p_0 \mapsto m}{p_0 s \mapsto 0} \otimes \binom{p_0 \mapsto n(s + \eta)}{p_0 s \mapsto 0} \\ - \tilde{p}_0 \otimes \tilde{p}_a s \otimes \binom{p_0 \mapsto m}{p_0 s \mapsto 0} \otimes \binom{p_0 \mapsto n}{p_0 s \mapsto 0} \\ \mapsto -(-1)^{|m| + |n|} \tilde{p}_0 \otimes \tilde{p}_a \otimes \begin{pmatrix} p_0 \otimes p_0 \mapsto m\otimes n(s + \eta) \\ p_0 s \otimes p_0 \mapsto 0 \\ p_0 \otimes p_0 s \mapsto 0 \\ p_0 s \otimes p_0 s \mapsto 0 \end{pmatrix} \\ - \tilde{p}_0 \otimes \tilde{p}_a s \otimes \begin{pmatrix} p_0 \otimes p_0 \mapsto m\otimes n \\ p_0 s \otimes p_0 \mapsto 0 \\ p_0 \otimes p_0 s \mapsto 0 \\ p_0 s \otimes p_0 s \mapsto 0 \end{pmatrix} \,.
\end{multline*}
Lastly, the computations of $\Psi$ and $\Phi$ given in Lemma~\ref{lem:Psichainmap} and Lemma~\ref{lem:Phichainmap} tell us that the final map in the composite is such that 
\begin{multline*}
-(-1)^{|m| + |n|} \tilde{p}_0 \otimes \tilde{p}_a \otimes \begin{pmatrix} p_0 \otimes p_0 \mapsto m\otimes n(s + \eta) \\ p_0 s \otimes p_0 \mapsto 0 \\ p_0 \otimes p_0 s \mapsto 0 \\ p_0 s \otimes p_0 s \mapsto 0 \end{pmatrix} \\ - \tilde{p}_0 \otimes \tilde{p}_a s \otimes \begin{pmatrix} p_0 \otimes p_0 \mapsto m\otimes n \\ p_0 s \otimes p_0 \mapsto 0 \\ p_0 \otimes p_0 s \mapsto 0 \\ p_0 s \otimes p_0 s \mapsto 0 \end{pmatrix} \mapsto -(-1)^{|m| + |n|} \tilde{p}_a \otimes \binom{p_0 \mapsto m \otimes n(s+\eta)}{p_0 s \mapsto 0} \\ - \tilde{p}_a s \otimes \binom{p_0 \mapsto m \otimes n}{p_0 s \mapsto 0}
\end{multline*}
where the target can be identified with $u_a \cdot (m \otimes n)$, as wanted.

By graded commutativity of $\smile$, we have 
\begin{align*}
u_a \cdot n \smile f_0 \cdot m &= (-1)^{|m||n|} f_0 \cdot m \smile u_a \cdot n \\ &= (-1)^{|m||n|} u_a \cdot m \otimes n \\ &= u_a \cdot n \otimes m \,.
\end{align*}
\end{proof}
\end{lemma}

\begin{lemma}
Let $f_1 \cdot m$ and $u_1 \cdot n$ be cycles with homology classes in~$\hatExt^1_\Gamma(k,M)$ and~$\hatExt^{-1}_\Gamma(k,N)$, respectively. Then the two cycles 
\[
f_1 \cdot m \smile u_1 \cdot n \simeq f_0 \cdot m \otimes n.
\]
are homologous in the complex $U_*$. It follows that they determine the same class in~$\hatExt^0_\Gamma(k,M \otimes N)$. 
\begin{proof}
Since $f_1 \cdot m$ and $u_1 \cdot n$ are assumed to be cycles we know that $m$ and $n$ are elements in $\ker(s+\eta)$, which directly implies that $m \otimes n$ is an element of~$\ker(s)$ since 
\begin{align*}
(m \otimes n) \cdot s &= m \otimes ns + (-1)^{|n|} ms \otimes n \\ &= m \otimes ns + m \otimes n \eta + m \otimes n \eta + (-1)^{|n|} ms \otimes n \\ &= m \otimes n(s + \eta) + (-1)^{|n|} m(s + \eta) \otimes n \\ &= 0 \,.
\end{align*}
An explicit description of the two cycles $f_1 \cdot m$ and $u_1 \cdot n$ is
\[
f_1 \cdot m = \tilde{p}_0 \otimes \binom{p_1 \mapsto m}{p_1 s \mapsto ms} \och u_1 \cdot n = -\tilde{p}_1 s \otimes \binom{p_0 \mapsto n}{p_0 s \mapsto 0} \,.
\]
The first map $1 \otimes \tau \otimes 1$ in the composite pairing twists the second tensor factor past the third one, so that
\[
- \tilde{p}_0 \otimes \binom{p_1 \mapsto m}{p_1 s \mapsto ms} \otimes \tilde{p}_1 s \otimes \binom{p_0 \mapsto n}{p_0 s \mapsto 0} \mapsto - \tilde{p}_0 \otimes \tilde{p}_{1} s \otimes \binom{p_1 \mapsto m}{p_1 s \mapsto ms} \otimes \binom{p_0 \mapsto n}{p_0 s \mapsto 0} \,.
\]
The second map in the composite is $1 \otimes 1 \otimes \alpha$, so that
\[
-\tilde{p}_0 \otimes \tilde{p}_{1} s \otimes \binom{p_1 \mapsto m}{p_1 s \mapsto ms} \otimes \binom{p_0 \mapsto n}{p_0 s \mapsto 0} \mapsto - \tilde{p}_0 \otimes \tilde{p}_{1} s \otimes \begin{pmatrix} p_1 \otimes p_0 \mapsto m\otimes n \\ p_1 s \otimes p_0 \mapsto (-1)^{|n|} ms \otimes n \\ p_1 \otimes p_0 s \mapsto 0 \\ p_1 s \otimes p_0 s \mapsto 0 \end{pmatrix} .
\] 
Lastly, the computations of $\Psi$ and $\Phi$ given in Lemma~\ref{lem:Psichainmap} and Lemma~\ref{lem:Phichainmap} tells us that the final map in the composite is such that
\[
- \tilde{p}_0 \otimes \tilde{p}_{1} s \otimes \begin{pmatrix} p_1 \otimes p_0 \mapsto m\otimes n \\ p_1 s \otimes p_0 \mapsto (-1)^{|n|} ms \otimes n \\ p_1 \otimes p_0 s \mapsto 0 \\ p_1 s \otimes p_0 s \mapsto 0 \end{pmatrix} \mapsto -\tilde{p}_1 s \otimes \binom{p_1 \mapsto m \otimes n}{p_1 s \mapsto (-1)^{|n|} ms \otimes n}  \,.
\]
The right hand term can be identified with $- h_{1,1} \cdot m \otimes n$. We conclude that
\[
f_1 \cdot m \smile u_1 \cdot n = - h_{1,1} \cdot m \otimes n \,.
\]
Note that the boundary of $g_{1,0} \cdot m \otimes n$ is 
\[
\partial(g_{1,0} \cdot m \otimes n) = f_0 \cdot m \otimes n + h_{1,1} \cdot m \otimes n \,,
\] 
which tells us that $f_0 \cdot m \otimes n$ and $-h_{1,1} \cdot m \otimes n$ are homologous in $U_*$, and hence represent the same class in~$\hatExt^0_{\Gamma}(k,M\otimes N)$.
\end{proof}
\end{lemma}

We decide to make a final change of notation.

\begin{notation}
Let $t^b \cdot m$ and $t^{-a} \cdot n$ denote the homology classes
\[
t^b \cdot m = [f_b \cdot m] \och t^{-a} \cdot n = [u_{a} \cdot n]
\] 
in $\hatExt^b_\Gamma(k,M)$ for~$b \geq 0$ and in $\hatExt^{-a}_\Gamma(k,N)$ for $a \geq 1$, respectively.
\end{notation}

Note that $t^b \cdot m$ has internal and total degrees equal to those of $f_b \cdot m$, so that 
\[
|t^b \cdot m| = |m| - b \och \|t^b \cdot m\| = |m| - 2b \,.
\]
Similarly, $t^{-a} \cdot n$ has internal and total degrees equal to that of $u_a \cdot n$, so that 
\[
|t^{-a} \cdot n| = a + |n|  \och \|t^{-a} \cdot n\| = 2a + |n|  \,.
\]
We conclude that, formally, the symbol $t$ has homological degree~$-1$, internal degree $|t| = -1$ and total degree~$\|t \| = -2$. Using this new notation we have the following theorem.

\begin{theorem}\index{complete $\Ext$!of circle group $\bT$}
\label{thm:completeExtcup}
\[
\hatExt^{c}_\Gamma(k,M) \cong \begin{cases}
{\displaystyle t^{c} \cdot
	\frac{\ker(s)}{\im(s+\eta)}} & \text{for $c$ even,} \\
\\
{\displaystyle t^{c} \cdot
	\frac{\ker(s+\eta)}{\im(s)}} & \text{for $c$ odd.} 
\end{cases}
\]
The cup product 
\[
\hatExt^{c_1}_\Gamma(k, M) \otimes \hatExt^{c_2}_{\Gamma}(k,N) \longto \hatExt^{c_1+c_2}_\Gamma(k,M \otimes N)
\]
is given by
\[
(t^{c_1} \cdot m) \smile (t^{c_2} \cdot n) = t^{c_1 + c_2} \cdot m \otimes n \,.
\]
\end{theorem}

\begin{corollary}
\label{cor:TateExtk}
\[
\hatExt^{c}_\Gamma(k,k) \cong \begin{cases}
{\displaystyle t^{c} \cdot
	\cok(\eta)} & \text{for $c$ even,} \\
{\displaystyle t^{c} \cdot
	\ker(\eta)} & \text{for $c$ odd.}  
\end{cases}
\]
In this case, the cup product
\[
\hatExt^{c_1}_\Gamma(k, k) \otimes \hatExt^{c_2}_{\Gamma}(k,k) \longto \hatExt^{c_1+c_2}_\Gamma(k,k)
\]
is given by 
\[
(t^{c_1} \cdot x) \smile (t^{c_2} \cdot y) = t^{c_1 + c_2} \cdot xy 
\]
and makes $\hatExt^{*}_\Gamma(k,k)$ into a $k$-algebra over which $\hatExt^*_\Gamma(k,M)$ is a module for any~$\Gamma$-module $M$.
\end{corollary}

\begin{remark}
\label{rmk:DGA}
Note that in Theorem~\ref{thm:completeExtcup} above, the answer is also the homology of a differential graded $\Gamma$-module 
\[
M[t,t^{-1}]
\]
with differential given by
\[
d(m) = tms \och d(t) = t^2 \eta \,,
\]
where $m$ is an element of $M$ and $t$ has homological degree~$-1$ and internal degree~$-1$. More precisely, the differential satisfies a Leibniz rule in the form
\begin{align*}
d(t^c m) &= d(t^{c})m + t^c d(m) = c t^{c+1} \eta m + t^{c+1} ms \\
& = \begin{cases} t^{c+1} m s & \text{if $c$ is even,} \\
t^{c+1} m (s+\eta) & \text{if $c$ is odd,}
\end{cases}
\end{align*}
and this gives us the same homology groups as in Theorem~\ref{thm:completeExtcup}. Note that this is also true multiplicatively: if $\mu \: M \otimes N \to L$ is a pairing of $\Gamma$-modules, then the cup product 
\[
\hatExt^{c_1}_\Gamma(k, M) \otimes \hatExt^{c_2}_{\Gamma}(k,N) \longto \hatExt^{c_1+c_2}_\Gamma(k,M \otimes N) \longto \hatExt^{c_1+c_2}_\Gamma(k,L)
\]
is precisely the one induced by the obvious pairing 
\[
M[t,t^{-1}] \otimes N[t,t^{-1}] \longto L[t,t^{-1}]
\]
on homology. 
\end{remark}

%
%
%

\chapter{Homotopy Groups of Orthogonal $G$-Spectra}
\label{sec:homotopygroups}

In this chapter we discuss some results regarding equivariant stable homotopy groups. Our chosen model for equivariant spectra is orthogonal $G$-spectra, and we recall some basic theory about these objects and their homotopy groups in Section~\ref{subsec:equivhtpygrp}. In Section~\ref{subsec:HopfalgSG} we define the main Hopf algebra that we will work with in this memoir, namely the (non-equivariant) homotopy groups of the unreduced suspension spectrum of a compact Lie group, also referred to as the spherical group ring $\bS[G]$ of that group. Since our main group of interest is the circle $\bT$, we also give an explicit description of $\pi_*(\bS[\bT])$ as an algebra over $\pi_*(\bS)$. Lastly, in Section~\ref{subsec:reshmphism} we show, under suitable projectivity assumptions, that we can sometimes describe the equivariant homotopy groups $\pi_*^G(X)$ of an orthogonal $G$-spectrum $X$ as the `$\pi_*(\bS[G])$-invariants' of the non-equivariant homotopy groups $\pi_*(X)$.  

\section{Equivariant homotopy groups}
\label{subsec:equivhtpygrp}
\index{equivariant homotopy groups|(}
\index{orthogonal $G$-spectrum|(}

Let $G$ be a compact Lie group, and let~$X$ be an orthogonal~$G$-spectrum,
as in \cite{MM02}*{\S II.2} and \cite{Sch18}*{\S3.1}. In what follows, we will always assume that~$G$ acts from the right. Recall that, in particular,~$X$ associates to each (finite-dimensional, orthogonal) $G$-representation~$V$ a based~$G$-space~$X(V)$, and to each pair $(U,V)$ of~$G$-representations a $G$-equivariant structure map $\sigma \: \Sigma^U X(V) \to X(U \oplus V)$.

We can define $G$-equivariant homotopy groups $\pi^G_*(X)$ associated to $X$. To do this, one fixes a complete~$G$-universe\index{complete $G$-universe}\footnote{A complete $G$-universe is an orthogonal representation of countably infinite dimension in which every finite dimensional $G$-representation, and their countably infinite direct sums, embeds.} $\sU$ containing a fixed copy of $\bR^\infty$. Note that the set of finite-dimensional $G$-subrepresentations of~$\sU$ is partially ordered by inclusion. For non-negative integers $q\ge0$, we define the $q$th~$G$-equivariant homotopy group of $X$ as the colimit, over this directed partially ordered set, of the sets of homotopy classes~$[f]$ of $G$-maps $f \: \Sigma^V S^q \to X(V)$:
\[
\pi^G_q(X) = \colim_V \, [\Sigma^V S^q, X(V)]^G \,.
\]
Similarly, to define the non-positive $G$-equivariant homotopy groups, we let
\[
\pi^G_{-q}(X) = \colim_V \, [\Sigma^{V - \bR^q} S^0, X(V)]^G
\]
where $V - \bR^q$ denotes the orthogonal complement of $\bR^q$ in $V$.
Here the colimit is formed over the partially ordered set of
subrepresentations~$V$ in~$\sU$ that contain~$\bR^q$.
These definitions agree for $q=0$. Each equivariant homotopy group $\pi^G_q(X)$ is naturally an abelian group.

Given a group homomorphism $H \to G$, we can view any $G$-spectrum as a $H$-spectrum. This gives rise to a restriction map\index{restriction homomorphism}
\[
\res^G_H \: \pi_*^G(X) \longto \pi_*^H(X)
\]
of graded abelian groups. In particular, any $G$-spectrum can be viewed as a non-equivariant spectrum via the inclusion homomorphism $1 \to G$, where $1$ denotes the
trivial group. In this case, we will write $\pi_*(X)$ in place of $\pi_*^1(X)$, as these are simply the ordinary non-equivariant homotopy groups of $X$ viewed as a non-equivariant orthogonal spectrum.

The category of orthogonal $G$-spectra is symmetric monoidal, with the symmetric monoidal product being denoted $\wedge$ and referred to as the \emph{smash product}\index{smash product}. The unit of this symmetric monoidal structure is the sphere spectrum $\bS$ with the trivial $G$-action. Any pairing $\phi \: X \wedge Y \to Z$ of orthogonal $G$-spectra gives rise to a pairing of the corresponding equivariant homotopy groups\index{pairing!of equivariant homotopy groups}. Consider classes $[f] \in \pi^G_p(X)$ and~$[g] \in \pi^G_q(Y)$, represented by homotopy classes of $G$-maps $f \: \Sigma^V S^p \to X(V)$ and $g \: \Sigma^W S^q \to Y(W)$, respectively. The induced pairing
\[
\phi_* \: \pi^G_p(X) \otimes \pi^G_q(Y) \longto \pi^G_{p+q}(Z)
\]
maps $[f] \otimes [g]$ to the element represented by the homotopy class
of the composite
\[
\Sigma^{V \oplus W} S^{p+q}
\overset{\cong}\longto \Sigma^V S^p \wedge \Sigma^W S^q
\overset{f\wedge g}\longto X(V) \wedge Y(W)
\overset{\phi}\longto Z(V \oplus W) \,.
\]
Similar constructions can be carried out if $p$ or $q$ is negative, although this is a bit tricky.
One approach is to use the suspension isomorphisms $E \: \pi^G_p(X)
\cong \pi^G_{p+1}(X \wedge S^1)$ and define $\phi_*$ as the composite
\begin{align*}
\pi^G_p(X) \otimes \pi^G_q(Y)
&\cong
\pi^G_{p+n}(X \wedge S^n) \otimes \pi^G_{q+m}(Y \wedge S^m) \\
&\longto
\pi^G_{p+n+q+m}(X \wedge S^n \wedge Y \wedge S^m) \\
&\xrightarrow{(-1)^{nq} \tau_*}
\pi^G_{p+q+n+m}(X \wedge Y \wedge S^n \wedge S^m) \\
&\longto
\pi^G_{p+q+n+m}(Z \wedge S^n \wedge S^m)
\cong
\pi^G_{p+q}(Z)
\end{align*}
for $n$ and~$m$ such that $p+n\ge0$ and $q+m\ge0$.
In this way, we obtain a pairing 
\[
\phi_* \: \pi^G_*(X) \otimes \pi^G_*(Y)
\longto \pi^G_*(Z)
\]
of graded abelian groups. 

\begin{remark}
More generally, given a group homomorphism $\alpha \: G \to H \times K$ and an $\alpha$-equivariant map $X \wedge Y \to Z$, where $X$, $Y$ and $Z$ are orthogonal $H$-, $K$- and~$G$-spectra, respectively, we obtain a pairing 
\[
\pi^H_*(X) \otimes \pi^K_*(Y) \longto \pi^G_*(Z) \,.
\]
Here, $\alpha$-equivariant means that the diagram 
\[
\begin{tikzcd}
X \wedge Y \wedge G_+ \arrow[r,"1 \wedge 1 \wedge \alpha"] \arrow[d] & X \wedge Y \wedge H_+ \wedge K_+ \arrow[r,"1 \wedge \tau \wedge 1"] & X \wedge H_+ \wedge Y \wedge K_+ \arrow[r] & X \wedge Y \arrow[d] \\
Z \wedge G_+ \arrow[rrr] & & & Z 
\end{tikzcd}
\]
commutes.
\end{remark}

If $R$ is a commutative (non-equivariant) orthogonal ring spectrum, with multiplication $\mu \: R \wedge R \to R$, then the induced pairing
\[
\mu_* \: \pi_*(R) \otimes \pi_*(R) \longto \pi_*(R)
\]
on (non-equivariant) homotopy groups makes $R_* = \pi_*(R)$ into a graded commutative ring. A right~$R$-module in orthogonal~$G$-spectra is an orthogonal $G$-spectrum~$X$ with an associative and unital action $\rho \: X \wedge R \to X$, defined in the category of orthogonal $G$-spectra. Here $R$ is regarded as a $G$-spectrum with trivial action. In this case, there is an induced pairing
\[
\rho_* \: \pi^G_*(X) \otimes R_* \longto \pi^G_*(X)
\]
making $\pi^G_*(X)$ into a right $R_*$-module. If $X$ and $Y$ are two~$R$-modules in orthogonal~$G$-spectra, then the canonical map $X \wedge Y \to X \wedge_R Y$ induces a pairing 
\[
\pi^G_*(X) \otimes \pi^G_*(Y) \longto \pi^G_*(X \wedge_R Y)
\]
that equalizes the two composites from $\pi^G_*(X) \otimes R_* \otimes \pi^G_*(Y)$, so that we have the induced dashed map making the diagram
\[
\begin{tikzcd}
\pi^G_*(X) \otimes R_* \otimes \pi^G_*(Y) \arrow[r,shift right] \arrow[r,shift left] & \pi^G_*(X) \otimes \pi^G_*(Y) \arrow[r] \arrow[rd] & \pi^G_*(X) \otimes_{R_*} \pi^G_*(Y) \arrow[d,dashed] \\
& & \pi^G_*(X \wedge_R Y)
\end{tikzcd}
\]
commute.
\index{orthogonal $G$-spectrum|)}
\index{equivariant homotopy groups|)}

\section{A cocommutative Hopf algebra}
\label{subsec:HopfalgSG}
\index{Hopf algebra|(}

Let us introduce the Hopf algebra that we will work with through the remainder of this memoir. The right $R$-action on $R[G] = R \wedge G_+$ is given by the composite map
\[
R \wedge G_+ \wedge R \xrightarrow{1\wedge\tau} R \wedge R \wedge G_+
	\xrightarrow{\mu\wedge1} R \wedge G_+ \,.
\]

\begin{lemma}
\label{lem:RGHopf}
If $R[G]_* = \pi_*(R \wedge G_+)$ is flat as a (right)
$R_*$-module, then $R[G]_*$ is naturally a cocommutative
Hopf algebra over~$R_* = \pi_*(R)$.
\end{lemma}

\begin{proof}
We have a pairing
\begin{multline*}
R[G]_* \otimes_{R_*} R[G]_* =
\pi_*(R \wedge G_+) \otimes_{\pi_*(R)} \pi_*(R \wedge G_+) \\
	\overset{\cdot}\longto \pi_*( (R \wedge G_+) \wedge_R (R \wedge G_+) )
	\cong \pi_*(R \wedge G_+ \wedge G_+) \,,
\end{multline*}
where the left $R_*$-action on the right hand copy of $R[G]_*$ is equal
to that obtained by twisting the right~$R_*$-action.  When $R[G]_*$
is flat as a right $R_*$-module, it follows by a well-known induction
\cite{Ada69}*{Lec.~3, Lem.~1, p.~68}
over the cells of a CW structure on the right hand copy of~$G$ that the
pairing above is an isomorphism of $R_*$-modules.

The unit inclusion $1 \to G$, group multiplication
$G \times G \to G$, collapse $G \to 1$, diagonal
$G \to G \times G$ and group inverse $G \to G$ give us
$R$-module maps
$R \to R \wedge G_+$,
$R \wedge G_+ \wedge_R R \wedge G_+ \to R \wedge G_+$,
$R \wedge G_+ \to R$,
$R \wedge G_+ \to R \wedge G_+ \wedge_R R \wedge G_+$
and $R \wedge G_+ \to R \wedge G_+$ that induce
$R_*$-module homomorphisms
\begin{align*}
\eta &\: R_* \longto R[G]_* \\
\phi &\: R[G]_* \otimes_{R_*} R[G]_* \longto R[G]_* \\
\epsilon &\: R[G]_* \longto R_* \\
\psi &\: R[G]_* \longto R[G]_* \otimes_{R_*} R[G]_* \\
\chi &\: R[G]_* \longto R[G]_*
\end{align*}
which make $R[G]_*$ a Hopf algebra over~$R_*$.  The cocommutativity of
the diagonal implies that $\psi$ is cocommutative.
\end{proof}

By the discussion in Section~\ref{sec:modulesHopf}, the category of modules over $R[G]_\ast$ is closed symmetric monoidal. Note that if $X$ is an $R$-module in orthogonal $G$-spectra, then the commuting right
$R$- and $G$-actions combine to define an action
\[
\gamma \: X \wedge_R R[G] \cong X \wedge G_+ \longto X
\]
which makes the underlying (non-equivariant) orthogonal spectrum of $X$ into a right $R[G]$-module in the category of (non-equivariant) $R$-modules. The induced pairing
\[
\gamma_* \: \pi_*(X) \otimes_{R_*} R[G]_\ast \longto \pi_*(X)
\]
gives the (non-equivariant) homotopy groups $\pi_*(X)$ the structure
of a right~$R[G]_*$-module\index{module, over Hopf algebra}. If $Y$ is a second $R$-module in orthogonal $G$-spectra, the pairing
\[
\pi_*(X) \otimes_{R_*} \pi_*(Y) \overset{\cdot}\longto \pi_*(X \wedge_R Y)
\]
is a homomorphism of $R[G]_*$-modules, where the Hopf algebra $R[G]_*$
acts diagonally on the left hand side.  Likewise,
\[
\pi_* F_R(X, Y) \longto \Hom_{R_*}(\pi_*(X), \pi_*(Y))
\]
is a homomorphism of $R[G]_*$-modules, where the Hopf algebra $R[G]_*$
acts by conjugation on the right hand side.

The case we are the most interested in is when the Lie group is the circle, so let us compute the homotopy groups of the spherical group ring of this specific group. 

\begin{proposition} \label{prop:s2isetas}
When $G = \bT = U(1)$ is the circle group,
\[
R[\bT]_* = R_*[s]/(s^2 = \eta s)
\]
with $|s| = |\eta| = 1$.  Here $s$ generates the augmentation ideal
\[
\overline{R[\bT]}_* = \ker(\epsilon \: R[\bT]_* \longto R_*)
	= R_*\{s\} \,,
\]
and $\eta$ is the image of the complex Hopf map in $\pi_1(\bS) \cong \bZ/2$. The generator $s$ is primitive, so the coproduct and involution are given by $\psi(s) = s \otimes 1 + 1 \otimes s$ and~$\chi(s) = -s$.
\end{proposition}

\begin{proof}
It suffices to prove the result for $R = \bS$. Proving this we would know that~$\bS[\bT]_\ast$ is free over~$\bS_*$, so that the case of a general ground ring spectrum~$R$ follows immediately from the isomorphism $R[\bT]_* \cong \bS[\bT]_* \otimes_{\bS_*} R_*$.

To prove the result for the sphere spectrum, we start by noting that the cofibre sequence
\[
S^0 \cong 1_+ \longto \bT_+ \longto \bT \cong S^1
\]
admits a retraction $\bT_+ \to 1_+$. Hence the induced stable
cofibre sequence
\[
\bS
	\overset{i}\longto \bS[\bT]
	\overset{p}\longto \Sigma \bS
\]
admits a retraction $c \: \bS[\bT] \to \bS$ and a section $s \: \Sigma
\bS \to \bS[\bT]$ with $ps \simeq 1$ and~$cs \simeq 0$. The maps $i$ and~$s$ represent classes in $\bS[\bT]_*$ of total (and internal) degree~$0$ and~$1$, respectively, and induce an isomorphism $\bS_* \{i,s\} \cong \bS[\bT]_*$.  Here, $i$ is the multiplicative unit and $s$
generates the augmentation ideal $\overline{\bS[\bT]}_* = \bS_*\{s\}$.
It only remains to prove that we have the relation $s^2 = \eta s$ in $\bS[\bT]_2$. This is the content of formula~(1.4.4) in \cite{Hes96}.  We give the following direct argument using the bar construction \cite{Seg68}*{\S3},~\cite{May75}*{\S7} and the bar spectral sequence\index{bar spectral sequence} \cite{Seg68}*{\S5},~
\cite{May72}*{\S11}.  We shall discuss these tools at greater length in Section~\ref{subsec:filteredEG}.

The bar construction of $\bT$ is the geometric realization $B\bT =
|B_\bullet\bT| \simeq \bC P^\infty$ of the simplicial space
\[
[q] \mapsto B_q \bT = \bT^q \,,
\]
with the usual face and degeneracy maps. There is a standard filtration of $B\bT$ by simplicial skeleta. The associated spectral sequence in (reduced) stable homotopy has $E^1$-page given as the normalised bar complex $NB_*(\bS_*,\bS[\bT]_*,\bS_*)$
\[
0 \leftarrow 0
	\overset{d^1_1}\longleftarrow \overline{\bS[\bT]}_*
	\overset{d^1_2}\longleftarrow
	\overline{\bS[\bT]}_* \otimes_{\bS_*} \overline{\bS[\bT]}_*
	\overset{d^1_3}\longleftarrow \cdots \,,
\]
which we reviewed in Construction~\ref{const:barconst}. This spectral sequence converges (strongly) to $\pi_* \Sigma^\infty(B\bT) \cong \pi_* \Sigma^\infty(\bC P^\infty)$. The part of the $E^1$-page that will be relevant to us is pictured below, with the origin in the lower left hand corner:
\[
\begin{tikzcd}
\dots & \dots & \dots & \dots & \dots \\
0 & \overline{\bS[\bT]}_2 & \overline{\bS[\bT]}_1 \otimes \overline{\bS[\bT]}_1 \arrow[l] & \dots & \dots \\
0 & \overline{\bS[\bT]}_1 & 0 & 0 & \dots \\
0 & 0 & 0 & 0 & 0 \\
\end{tikzcd}
\]
Firstly,
\[
d^1_2(x \otimes y) = \epsilon(x)y - xy + x \epsilon(y) = -xy
\]
for $x,y \in \overline{\bS[\bT]}_*$, since $x$ and~$y$ both augment to zero. With prior understanding of the stable homotopy groups of $\bC P^\infty$ we can now figure out the displayed differential. The stable class of the inclusion $S^2 \cong \bC P^1 \to \bC P^\infty$ is well-known to generate~$\pi_2 \Sigma^\infty \bC P^\infty \cong \bZ$. For degree reasons this must be detected by $\pm s$ in $E^\infty_{1,1} = E^1_{1,1} = \overline{\bS[\bT]}_1 = \bZ\{s\}$. The $4$-cell in $\bC P^2$ is attached
by the Hopf fibration\index{Hopf fibration}~$\eta$ to $\bC P^1$, so that $\pi_3 \Sigma^\infty
\bC P^\infty = 0$. This forces $\eta s \in E^1_{1,2} = \overline{\bS[\bT]}_2 =
\bZ/2\{\eta s\}$ to be a boundary in the spectral sequence.  For degree reasons, the only possibility is that 
\[
\eta s = d^1_2(s \otimes s) \,,
\]
so that $s^2 = \eta s$.

Note that the coproduct $\psi(s)$ must contain the terms $s \otimes 1$ and $1
\otimes s$ by counitality, and cannot contain other terms since
$\bS[\bT]_*$ is connected and $|s| = 1$. Hence $s$ is a primitive element of our Hopf algebra. 
\end{proof}

\begin{proposition}
When $G = \bU = Sp(1)$ is the $3$-sphere group,
\[
R[\bU]_* = R_*[t]/(t^2 = \dot\nu t)
\]
with $|t| = |\dot\nu| = 3$.  Here $t$ generates the augmentation ideal
\[
\overline{R[\bU]}_* = \ker(\epsilon \: R[\bU]_* \longto R_*)
	= R_*\{t\} \,,
\]
and $\dot\nu$ is the image of a generator of $\pi_3(\bS) \cong \bZ/24$.
The coproduct is given by $\psi(t) = t \otimes 1 + 1 \otimes t$.
\end{proposition}

\begin{proof}
Similar to the circle case.
\end{proof}

Note that we cannot assert from this line of argument that $\dot\nu$ is the image of the quaternionic Hopf map. The bar spectral sequence argument for $R = \bS$ only shows that $t^2 = \dot\nu t$ with $\dot\nu$ some generator of $\pi_3(\bS) \cong \bZ/24$. A more geometric argument might link $\dot\nu$ to the standard generator $\nu$ of $\pi_3(\bS)$, but we will not pursue this here.
\index{Hopf algebra|)}

\section{A restriction homomorphism}
\label{subsec:reshmphism}
\index{restriction homomorphism|(}
As explained earlier, the inclusion homomorphism $1 \to G$ gives rise to a map of graded abelian groups
\[
\res^G_1 \: \pi^G_*(X) \longto \pi_*(X)
\]
taking the homotopy class of a $G$-map $f \: \Sigma^V S^q \to X(V)$ to
the homotopy class of the underlying non-equivariant map, and similarly
for $G$-maps $\Sigma^{V - \bR^q} S^0 \to X(V)$. Tracing through definitions shows that $\res^G_1$ is~$R_*$-linear if $X$ is an $R$-module in orthogonal~$G$-spectra. We are interested in the following refined restriction homomorphism.

\begin{lemma} \label{lem:omegaXnat}
There is a natural $R_*$-module homomorphism
\[
\omega_X \: \pi^G_*(X) \longto \Hom_{R[G]_*}(R_*, \pi_*(X))
\]
making the triangle in the following diagram commute.
\[
\xymatrix{
\pi^G_*(X) \ar[d]_-{\omega_X} \ar[dr]^-{\res^G_1} \\
\Hom_{R[G]_*}(R_*, \pi_*(X)) \ar@{ >->}[r]
	& \pi_*(X) \ar@<1ex>[r]^-{\tilde{\gamma}} \ar@<-1ex>[r]_-{\tilde{\epsilon}}
        & \Hom_{R_*}(R[G]_*, \pi_*(X))
}
\]
Here $\tilde \epsilon$ denotes the adjoint of the trivial $R[G]_*$-action $\epsilon_*$ on $\pi_*(X)$, which equals the composite
\[
\epsilon_* \: \pi_*(X) \otimes_{R_*} R[G]_*
	\overset{1 \otimes \epsilon}\longto
	\pi_*(X) \otimes_{R_*} R_* \cong \pi_*(X) \,.
\]
Similarly $\tilde \gamma$ denotes the adjoint to the $R[G]_\ast$-action
\[
\gamma_* \: \pi_*(X) \otimes_{R_*} R[G]_* \longto \pi_*(X)
\]
on $\pi_*(X)$.
\end{lemma}

\begin{proof}
We claim that the two composite homomorphisms
\[
\xymatrix{
\pi^G_*(X) \otimes_{R_*} R[G]_* \ar[r]^-{\res^G_1 \otimes 1}
	& \pi_*(X) \otimes_{R_*} R[G]_* 
	\ar@<1ex>[r]^-{\gamma_*}
	\ar@<-1ex>[r]_-{\epsilon_*}
	& \pi_*(X)
}
\]
are equal.  This implies that the two adjoint
homomorphisms
\[
\xymatrix{
\pi^G_*(X) \ar[r]^-{\res^G_1}
        & \pi_*(X)
        \ar@<1ex>[r]^-{\tilde \gamma}
        \ar@<-1ex>[r]_-{\tilde \epsilon}
        & \Hom_{R_*}(R[G]_*, \pi_*(X))
}
\]
are equal, so that $\res^G_1$ factors uniquely
through the equalizer $\Hom_{R[G]_*}(R_*, \pi_*(X))$ of
the two right-hand arrows.

By fibrant replacement we may assume that $X$ is an $\Omega$-$G$-spectrum
\cite{MM02}*{Def.~III.3.1}, meaning that each adjoint structure map $X(V)
\to \Omega^{W-V} X(W)$ is a weak $G$-equivalence, where $V \subset W$ lie
in our fixed complete $G$-universe $\sU$. Then each element~$x$ in $\pi^G_*(X)$ is represented by the
homotopy class $[f]$ of a $G$-map $f \: S^m \to X(\bR^n)$, for suitable
non-negative integers~$m$ and~$n$.  Here $G$ acts trivially on $S^m$,
so $f$ factors through the fixed points $X(\bR^n)^G$, where the
$G$-action $\gamma$ is trivial.  It follows that $\tilde \gamma$ and
$\tilde \epsilon$ agree on $\res^G_1(x) \otimes y$ for any $y \in R[G]_*$,
as claimed.
\end{proof}

\begin{proposition} \label{prop:omegaXiso}
If $R[G]_*$ is projective as an $R_*$-module, and
$X \simeq F(G_+, Y)$ for some $R$-module $Y$ in
orthogonal $G$-spectra, then the natural homomorphism
\[
\omega_X \: \pi^G_*(X)
	\overset{\cong}\longto \Hom_{R[G]_*}(R_*, \pi_*(X))
\]
is an isomorphism.
\end{proposition}

\begin{proof}
By fibrant replacement, we may assume that $Y$ is an
$\Omega$-$G$-spectrum.  As usual we give $F(G_+, Y) \cong F_R(R[G], Y)$ the conjugate $G$-action. By naturality of~$\omega_X$ we may assume that $X = F(G_+, Y)$, in which case $X$ is also an $\Omega$-$G$-spectrum.

Let us consider the commutative diagram
\[
\begin{tikzcd}
\pi^G_*(X) \arrow[r,"\res^G_1"] \arrow[d,"\cong"]
	& \pi_*(X) \arrow[d,"\cong"] \\
\pi_*(Y) \arrow[r,"\tilde \gamma"]
	& \Hom_{R_*}(R[G]_*, \pi_*(Y)) \,.
\end{tikzcd}
\]
We make the maps involved a bit more explicit. The vertical isomorphisms are given as follows. The left hand vertical isomorphism $\pi^G_*(X) = \pi^G_*(F(G_+, Y)) \to \pi_*(Y)$ takes the homotopy class of a $G$-map $f \: S^m \to X(\bR^n)
= F(G_+, Y(\bR^n))$ bijectively to the homotopy class of $f' \: S^m \to Y(\bR^n)$ given by $f'(s) = f(s)(e)$, where $e \in G$ is the unit element of our group. The right hand vertical isomorphism is the special case $Z = R[G]$ of the natural $R[G]_*$-module homomorphism
\[
\pi_* F_R(Z, Y) \longto \Hom_{R_*}(\pi_*(Z), \pi_*(Y)) \,.
\]
Indeed, this is an isomorphism whenever $\pi_*(Z)$ is projective as an $R_*$-module. The top horizontal map is the restriction homomorphism we described at the beginning of this section, and the lower horizontal homomorphism
$\tilde \gamma$ is adjoint to the~$R[G]_*$-module action on $\pi_*(Y)$. Note that the diagram does indeed commute, since the lower and upper compositions both send the
homotopy class of the~$G$-map $f \: S^m \to F(G_+, Y(\bR^n))$ to the homomorphism $R[G]_* \to \pi_*(Y)$ induced by the left adjoint $S^m \wedge G_+ \to Y(\bR^n)$ of~$f$. 

The homomorphism $\tilde \gamma$ identifies $\pi_*(Y)$ with the
$R_*$-submodule
\[
\Hom_{R[G]_*}(R_*, \Hom_{R_*}(R[G]_*, \pi_*(Y)))
	\cong \Hom_{R[G]_*}(R[G]_*, \pi_*(Y))
\]
of its target, while $\res^G_1$ factors through $\Hom_{R[G]_*}(R_*,
\pi_*(X))$ by the previous lemma.  Hence we can apply $\Hom_{R[G]_*}(R_*,
-)$ to the right hand vertical isomorphism in the diagram above to obtain another isomorphism, and a commutative square
\[
\xymatrix{
\pi^G_*(X) \ar[r]^-{\omega_X} \ar[d]_-{\cong}
        & \Hom_{R[G]_*}(R_*, \pi_*(X)) \ar[d]^-{\cong} \\
\pi_*(Y) \ar[r]^-{\cong}
        & \Hom_{R[G]_*}(R_*, \Hom_{R_*}(R[G]_*, \pi_*(Y))) \,.
}
\]
It follows that $\omega_X$ is an isomorphism, as asserted.
\end{proof}

To handle multiplicative structure, we need the following observation.

\begin{lemma} \label{lem:omegaXmon}
The natural transformation $\omega$ is monoidal, in the sense that
the diagram
\[
\begin{tikzcd}[column sep=small]
\pi^G_*(X) \otimes_{R_*} \pi^G_*(Y) \ar[r,"\cdot"]
	\ar[d,"\omega_X \otimes \omega_Y"]
	& \pi^G_*(X \wedge_R Y) \ar[dd,"\omega_{X \wedge_R Y}"] \\
\Hom_{R[G]_*}(R_*, \pi_*(X)) \otimes_{R_*} \Hom_{R[G]_*}(R_*, \pi_*(Y)) 
\ar[d,"\alpha"] \\
\Hom_{R[G]_*}(R_*, \pi_*(X) \otimes_{R_*} \pi_*(Y)) \ar[r,"\cdot"]
	& \Hom_{R[G]_*}(R_*, \pi_*(X \wedge_R Y))
\end{tikzcd}
\]
commutes.
\end{lemma}

\begin{proof}
Since $\Hom_{R[G]_*}(R_*, \pi_*(X \wedge_R Y)) \to \pi_*(X \wedge_R Y)$
is a monomorphism it suffices to show that
\[
\xymatrix{
\pi^G_*(X) \otimes_{R_*} \pi^G_*(Y) \ar[r]^-{\cdot} \ar[d]_-{\res^G_1 \otimes \res^G_1}
	& \pi^G_*(X \wedge_R Y) \ar[d]^-{\res^G_1} \\
\pi_*(X) \otimes_{R_*} \pi_*(Y) \ar[r]^-{\cdot}
	& \pi_*(X \wedge_R Y)
}
\]
commutes, which is clear.
\end{proof}
\index{restriction homomorphism|)}

%
%
%

\chapter{Sequences of Spectra and Spectral Sequences}
\label{sec:seqspec+specseq}

In this chapter we associate a Cartan--Eilenberg system, an exact couple, and a spectral sequence to any sequence of orthogonal $G$-spectra. We identify certain well-behaved sequences called filtrations, and use these to show how pairings of sequences induce pairings of Cartan--Eilenberg systems and spectral sequences. This is essentially the content of Section~\ref{subsec:pairingsCEI} and Section~\ref{subsec:pairingsCEII}, culminating in Theorem~\ref{thm:mainmultI} and Theorem~\ref{thm:mainmultII}. We shall rely on the classical telescope construction to approximate general sequences by equivalent filtrations. Finally we discuss how pairings can be internalized in terms of the convolution product of two sequences.

\section{Cartan--Eilenberg systems}

A Cartan--Eilenberg system is an algebraic structure, introduced
in~\cite{CE56}, which determines two exact couples~\cite{Mas52} and a spectral sequence. This structure has the advantage that one can give a useful definition of a pairing of Cartan--Eilenberg systems, which determines a pairing of the corresponding spectral sequences.  These definitions were reviewed by Douady in \cite{Dou59a} and~\cite{Dou59b}. As opposed to these sources, which use cohomological indexing, we adopt homological indexing for our Cartan--Eilenberg systems, as in \cite{HR19}.

We start with some preliminary definitions. We will in particular make use of the posets $[1] = \{0 \to 1\}$ and $[2] = \{0 \to 1 \to 2\}$ regarded as categories. Note that we have three functors
\[
\delta_ 0 , \delta_1 , \delta_2 \: [1] \longto [2] \,,
\]
with subscript indicating which object of the target is skipped. In addition, we have natural transformations 
\[
i \: \delta_2 \longto \delta_1 \och p \: \delta_1 \longto \delta_0 \,.
\]

\begin{definition}[\cite{HR19}*{Def.~4.1}]
Let $(\cI, \le)$ be a linearly ordered set.
\begin{itemize}
\item Let $\cI^{[1]} = \Fun([1],\cI)$. The objects in this category are pairs $(i,j)$ in $\cI$ with $i \le j$, and we have a single morphism $(i,j) \to (i',j')$ precisely when $i \le i'$ and~$j \le j'$.
\item Let $\cI^{[2]} = \Fun([2], \cI)$. The objects in this category are triples $(i,j,k)$ in $\cI$ with $i \le j \le k$, and we have a single morphism
$(i,j,k) \to (i',j',k')$ precisely when $i \le i'$, $j \le j'$ and~$k
\le k'$.
\end{itemize}
\end{definition}

The functors $\delta_0$, $\delta_1$, and $\delta_2$ defined above induce functors 
\[
d_0 , d_1 , d_2 \: \cI^{[2]} \longto \cI^{[1]} \,.
\] 
These map $(i,j,k)$ to $(j,k)$, $(i,k)$ and~$(i,j)$, respectively. The natural transformations $i$ and $p$ induce natural transformations
\[
\iota \: d_2 \longto d_1 \och \pi \: d_1 \longto d_0
\]
with components $\iota \: (i,j) \to (i,k)$
and $\pi \: (i,k) \to (j,k)$, respectively.

Let $k$ be a graded ring and let $\cA$ be the graded abelian category of $k$-modules. The grading $\|x\|$ of a homogeneous element $x \in M$ in an object $M$ of~$\cA$ will be referred to as its \emph{total degree}\index{degree!total}.

\begin{definition}[\cite{HR19}*{Def.~4.2, Def.~6.1}]
\label{def:CEsystem}
An \emph{$\cI$-system}\index{I-system@$\cI$-system} in $\cA$ is a pair~$(H, \partial)$ where
$H \: \cI^{[1]} \to \cA$ is a functor and
$\partial \: Hd_0 \to Hd_2$ is a natural transformation of functors $\cI^{[2]} \to \cA$, such
that the triangle
\[
\xymatrix{
Hd_2 \ar[r]^-{H\iota} & Hd_1 \ar[d]^-{H\pi} \\
	& Hd_0 \ar[ul]^-{\partial}
}
\]
is exact.  We assume that $H\iota$ and $H\pi$ have total degree~$0$,
while $\partial$ has total degree~$-1$.  We generically write $\eta \:
H(i,j) \to H(i',j')$ for the total degree~$0$ morphisms in $\cA$ induced
by morphisms in~$\cI^{[1]}$.
\end{definition}

\begin{definition}
\label{def:CEsystemII}
\leavevmode
\begin{itemize}
\item A \emph{finite Cartan--Eilenberg system}\index{Cartan--Eilenberg system!finite} is a $\bZ$-system $(H, \partial)$, where $\bZ$ denotes the integers with its usual linear ordering. 
\item An \emph{extended
Cartan--Eilenberg system}\index{Cartan--Eilenberg system!extended} is a $\cI$-system $(H, \partial)$
for $\cI = \bZ \cup \{\pm\infty\}$, with the extended linear
ordering where $-\infty$ is initial and $+\infty$ is terminal.
\item An extended Cartan--Eilenberg system $(H, \partial)$ is a \emph{Cartan--Eilenberg system}\index{Cartan--Eilenberg system} if the following condition, called~(SP.5)\index{SP@(SP.5)}, is satisfied: The canonical
homomorphism
\[
\colim_j H(i,j) \overset{\cong}\longto H(i,\infty)
\]
is an isomorphism for all~$i \in \bZ$.
\end{itemize}
\end{definition}

An extended Cartan--Eilenberg system thus associates to each
pair $(i,j)$ with $-\infty \le i \le j \le \infty$ a module
$H(i,j)$, in a functorial way.  Furthermore, it associates to
each triple $(i,j,k)$ with $-\infty \le i \le j \le k \le \infty$
a long exact sequence
\[
\dots \longto H(i,j) \longto H(i,k) \longto H(j,k)
	\overset{\partial}\longto H(i,j) \longto \dots \,,
\]
where $\partial$ is a natural transformation of total degree~$-1$.
If the homomorphism in condition~(SP.5) is an isomorphism for one $-\infty
\le i < \infty$, then it is an isomorphism for every such $i$. This follows by using the 5-Lemma twice in the following map of exact sequences:
\[
\xymatrix@C-1.2pc{
\dots \ar[r]
& H(-\infty,i) \ar[r] \ar[d]^-{=}
& \colim_j H(-\infty,j) \ar[r] \ar[d]
& \colim_j H(i,j) \ar[r]^-{\partial} \ar[d]
& H(-\infty,i) \ar[r] \ar[d]^-{=}
& \dots
\\
\dots \ar[r]
& H(-\infty,i) \ar[r]
& H(-\infty, \infty) \ar[r]
& H(i, \infty) \ar[r]^-{\partial}
& H(-\infty,i) \ar[r]
& \dots \,.
}
\]
An extended Cartan--Eilenberg system determines a finite Cartan--Eilenberg
system by restriction to~$(i,j)$ with $-\infty < i \le j < \infty$.

\begin{remark}
Apart from the switch in variance, the definition given by Cartan and
Eilenberg in \cite{CE56}*{\S XV.7} corresponds to our Cartan--Eilenberg
systems.  This is also the definition recalled in \cite{McC01}*{Ex.~2.2}.
In \cite{Dou59a}*{\S~II~C}, Douady works with data defining an Adams
spectral sequence, which is concentrated in non-negative cohomological
(so: non-positive homological) filtration degrees.  He therefore
assumes that $H(i,0) = H(i,j)$ for all $i \le 0 \le j \le \infty$,
so that condition~(SP.5) is trivially satisfied.
\end{remark}

\begin{definition}[\cite{HR19}*{Def.~7.1}]
\label{def:excpl+spseq}
Let $(H, \partial)$ be a Cartan--Eilenberg system.  Let the \emph{left
couple}\index{left couple} $(A, E^1)$ be the exact couple given by 
\[A_s = H(-\infty, s) \och E^1_s = H(s-1,s)
\]
fitting together in the exact triangle
\[
\xymatrix{
A_{s-1} \ar[r] & A_s \ar[d] \\
	& E^1_s \ar[ul]^-{\partial}
}
\]
associated to the triple~$(-\infty,s-1,s)$. The \emph{abutment}\index{abutment!of exact couple} of
this exact couple is 
\[
A_\infty = \colim_s A_s \cong H(-\infty,\infty) \,.
\]
This abutment is exhaustively filtered by the images\index{image filtration}
\[
F_s A_\infty = \im(A_s \longto A_\infty) \,.
\]
\end{definition}

The Cartan--Eilenberg system and the left couple give rise to the same spectral sequence $(E^r, d^r)$. The pages of this spectral sequence\index{spectral sequence!associated to a Cartan--Eilenberg system} are given by
\[
E^r_s = Z^r_s / B^r_s
\]
and the differentials $d^r_s \: E^r_s \to E^r_{s-r}$ are of total degree~$-1$.  Here
\begin{align*}
Z^r_s &= \ker(\partial \: E^1_s \longto H(s-r,s-1)) \\
B^r_s &= \im(\partial \: H(s,s+r-1) \longto E^1_s)
\end{align*}
define the $r$-cycles\index{r-cycles@$r$-cycles} and $r$-boundaries\index{r-bounaries@$r$-boundaries} in \emph{filtration degree}\index{degree!filtration}~$s$,
respectively, and
\[
d^r_s([x]) = [\partial(\tilde x)] \,,
\]
where $x \in Z^r_s$ and $\tilde x \in H(s-r,s)$ satisfies $\eta(\tilde x) = x$.  There are preferred isomorphisms $H(E^r, d^r) \cong E^{r+1}$.  We let 
\[
Z^\infty_s =
\lim_r Z^r_s \,, \quad B^\infty_s = \colim_r B^r_s \,, \och E^\infty_s = Z^\infty_s
/ B^\infty_s \,.
\]

We refer to \cite{CE56}*{\S XV.1} or \cite{HR19}*{Prop.~4.9} for the
verification that $(E^r, d^r)$ is indeed a spectral sequence. Note in particular that it only depends on the finite part of the Cartan--Eilenberg
system~$(H, \partial)$.  The abutment and $E^\infty$-page are
related as follows.

\begin{lemma} \label{lem:monotoEinfty}
There is a natural monomorphism
\[
\beta \: \frac{F_s A_\infty}{F_{s-1} A_\infty} \longto E^\infty_s
\]
in each filtration degree~$s$.
\end{lemma}

\begin{proof}
See \cite{Boa99}*{Lem.~5.6} or \cite{HR19}*{Lem.~3.15(a)}.
\end{proof}

The main purpose of reviewing the above definitions is to let us record the following definitions of pairings of (finite and classical) Cartan--Eilenberg systems.  We assume that~$k$ is graded commutative, and write $\otimes$ in place of~$\otimes_k$.

\begin{definition} \label{def:pairfinCEsys}
\index{pairing!of finite Cartan--Eilenberg systems}
Let $(H', \partial)$, $(H'', \partial)$ and $(H,
\partial)$ be finite Cartan--Eilenberg systems in~$\cA$.  A \emph{pairing}
$\phi \: (H', H'') \to H$ of such systems is a collection of
$k$-module homomorphisms
\[
\phi_r \: H'(i-r,i) \otimes H''(j-r,j) \longto H(i+j-r,i+j)
\]
of total degree~$0$, for all $i,j \in \bZ$ and $r\ge1$.
These are required to satisfy the following two conditions:
\begin{description}
\item[SPP I]\index{SPPa@(SPP I)} Each square
\[
\xymatrix{
H'(i-r,i) \otimes H''(j-r,j) \ar[r]^-{\phi_r} \ar[d]_-{\eta \otimes \eta}
	&  H(i+j-r,i+j) \ar[d]^-{\eta} \\
H'(i'-r',i') \otimes H''(j'-r',j') \ar[r]^-{\phi_{r'}}
	&  H(i'+j'-r',i'+j')
}
\]
commutes, for all integers $i,j,i',j'$ and $r,r' \ge1$ with
$i \le i'$, $i-r \le i'-r'$, $j \le j'$ and~$j-r \le j'-r'$.
\item[SPP II]\index{SPPb@(SPP II)} In the (non-commutative) diagram
\[
\xymatrix@C-4pc{
H'(i-r,i) \otimes H''(j-r,j)
	\ar[dd]_-{\partial\otimes\eta}
	\ar[dr]^-{\phi_r}
	\ar[rr]^-{\eta\otimes\partial}
&& H'(i-1,i) \otimes H''(j-r-1,j-r) \ar[dd]^-{\phi_1} \\
& H(i+j-r,i+j) \ar[dr]^-{\partial} \\
H'(i-r-1,i-r) \otimes H''(j-1,j) \ar[rr]^-{\phi_1}
&& H(i+j-r-1,i+j-r)
}
\]
the inner composition is the sum of the two outer ones:
\[
\partial \phi_r = \phi_1 (\partial \otimes \eta)
	+ \phi_1 (\eta \otimes \partial) \,.
\]
In terms of elements, this identity in $H(i+j-r-1,i+j-r)$ can be written
\[
\partial\phi_r(x \otimes y) = \phi_1(\partial x \otimes \eta y)
	+ (-1)^{\|x\|} \phi_1(\eta x \otimes \partial y)
\]
for $x \in H'(i-r,i)$ of total degree~$\|x\|$ and $y \in H''(j-r,j)$.
\end{description}
\end{definition}

\begin{remark}
Apart from the switch in variance, this definition agrees with that
of Douady \cite{Dou59b}*{\S~II~A}, except for the fact that we ignore $r=0$, since~$\phi_0$ carries no information, and Douady omits the cases $i>0$ and $j>0$,
due to his focus on Adams spectral sequences. In the definition given in
\cite{McC01}*{Ex.~2.3}, the homomorphism $\varphi_1$ is missing from the
right hand term in his equation~(2), and the conditions $n\ge0$ and $q\ge0$
should be omitted.
\end{remark}

\begin{definition}
\label{def:sseqpairing}
Let $({}' E^r, d^r)$, $({}'' E^r, d^r)$ and $(E^r, d^r)$ be $k$-module
spectral sequences.  A \emph{pairing}\index{pairing!of spectral sequences} $\phi \: ({}' E^*, {}'' E^*)
\to E^*$ of such spectral sequences consists of a collection of
$k$-module homomorphisms
\[
\phi^r \: {}' E^r \otimes {}'' E^r \longto E^r
\]
for all $r\ge1$, such that:
\begin{enumerate}
\item The Leibniz rule\index{Leibniz rule}
\[
d^r \phi^r = \phi^r(d^r \otimes 1) + \phi^r(1 \otimes d^r)
\]
holds as an equality of homomorphisms
${}' E^r_i \otimes {}'' E^r_j \longto E^r_{i+j-r}$
for all $i,j \in \bZ$ and $r\ge1$.
\item The diagram 
\[
\begin{tikzcd}
{}' E^{r+1} \otimes {}'' E^{r+1} \arrow[r,"\phi^{r+1}"] \arrow[d] & E^{r+1} \arrow[d,"\cong"] \\
H({}' E^{r} \otimes {}'' E^{r}) \arrow[r,"H(\phi^{r})"] & H(E^{r})
\end{tikzcd}
\]
commutes for all $r \geq 1$.
\end{enumerate} 
\end{definition}

By a \emph{multiplicative spectral sequence}\index{spectral sequence!multiplicative}, we mean a spectral sequence
$(E^r, d^r)$ equipped with a pairing 
\[
\phi \: (E^*, E^*) \longto E^* \,.
\]
If $\phi^a \: E^a \otimes E^a \to E^a$ (usually with $a=1$ or $a=2$)
is associative and unital, then each pairing $\phi^r$ for~$r\ge a$
is also associative and unital, and we call $(E^r, d^r)_{r\ge a}$ an
\emph{algebra spectral sequence}\index{algebra spectral sequence}.

\begin{theorem}[\cite{Dou59b}*{Thm.~II~A}]
\label{thm:Dou59bIIA}
Let $(H', \partial)$, $(H'', \partial)$ and $(H, \partial)$ be finite
Cartan--Eilenberg systems, with associated spectral sequences referred to as $({}'
E^r, d^r)$, $({}'' E^r, d^r)$ and $(E^r, d^r)$, respectively.  Let $\phi \: (H', H'')
\to H$ be a pairing of finite Cartan--Eilenberg systems.  Then there is
a pairing $\phi \: ({}' E^*, {}'' E^*) \to E^*$ of spectral sequences,
uniquely defined by the condition $\phi^1 = \phi_1$.
\end{theorem}

\begin{proof}
Douady leaves the proof to the reader (`s'il existe'). Starting with setting
$\phi^1 \: {}' E^1_i \otimes {}'' E^1_j \to E^1_{i+j}$ equal to
\[
\phi_1 \: {}' H(i-1,i) \otimes {}'' H(j-1,j)
	\longto H(i+j-1,i+j) \,,
\]
the point is to inductively show that $d^r$ satisfies the Leibniz rule with respect to the pairing $\phi^r$ of $E^r$-pages, so that $\phi^{r+1}$ can be defined to be equal to the induced pairing in homology with respect to $d^r$. A full proof can be found in \cite{Hel17}*{Prop.~3.4.2}.
\end{proof}

We now move from finite Cartan--Eilenberg systems to classical ones. 

\begin{definition} \label{def:pairextCEsys}
\index{pairing!of Cartan--Eilenberg systems}
Let $(H', \partial)$, $(H'', \partial)$ and $(H, \partial)$ be
Cartan--Eilenberg systems.
A \emph{pairing} $\phi \: (H', H'') \to H$ of such systems consists of
a pairing $(\phi_r)_{r\ge1}$ of the restricted finite Cartan--Eilenberg
systems, together with $k$-module homomorphisms
\[
\phi_\infty \: H'(-\infty,i) \otimes H''(-\infty,j)
	\longto H(-\infty,i+j)
\]
of total degree~$0$, for all $i,j \in \bZ$.  These are required to
satisfy the following additional condition:
\begin{description}
\item[SPP III]\index{SPPc@(SPP III)} The squares
$$
\xymatrix{
H'(-\infty,i) \otimes H''(-\infty,j) \ar[r]^-{\phi_\infty}
        \ar[d]_-{\eta \otimes \eta}
& H(-\infty,i+j) \ar[d]^-{\eta} \\
H'(-\infty,i') \otimes H''(-\infty,j') \ar[r]^-{\phi_\infty}
& H(-\infty,i'+j')
}
$$
and
$$
\xymatrix{
H'(-\infty,i) \otimes H''(-\infty,j) \ar[r]^-{\phi_\infty}
        \ar[d]_-{\eta \otimes \eta}
& H(-\infty,i+j) \ar[d]^-{\eta} \\
H'(i-r,i) \otimes H''(j-r,r) \ar[r]^-{\phi_r}
& H(i+j-r,i+j)
}
$$
commute, for all integers $i \le i'$, $j \le j'$ and~$r\ge1$.
\end{description}
\end{definition}

We emphasize that the $\phi_r$ in the definition above must satisfy the conditions~(SPP~I) and~(SPP~II),
by virtue of defining a pairing of finite Cartan--Eilenberg systems.
The new condition~(SPP~III) is an analogue of~(SPP~I) for $r=\infty$.

With notation as in Definition~\ref{def:excpl+spseq}, we can rewrite $\phi_\infty$ as compatible pairings
\[
\phi_{i,j} \: A'_i \otimes A''_j \longto A_{i+j}
\]
in the corresponding left couples, for all $i, j \in \bZ$. Passing to colimits, we obtain a pairing of abutments
\[
\phi_* \: A'_\infty \otimes A''_\infty \longto A_\infty \,.
\]
This is filtration-preserving in the sense that it sends $F_i A'_\infty \otimes F_j A''_\infty$ to $F_{i+j} A_\infty$, by virtue of the commutative diagram
\[
\xymatrix{
A'_i \otimes A''_j \ar[r]^-{\phi_{i,j}} \ar@{->>}[d]
	& A_{i+j} \ar@{->>}[d] \\
F_i A'_\infty \otimes F_j A''_\infty \ar[r] \ar[d]
	& F_{i+j} A_\infty \ar@{ >->}[d] \\
A'_\infty \otimes A''_\infty \ar[r]^-{\phi_*}
	& A_\infty \rlap{\,.}
}
\]
Being filtration-preserving, the pairing $\phi_*$ then induces
pairings of filtration subquotients
\[
\bar\phi_* \: \frac{F_i A'_\infty}{F_{i-1} A'_\infty}
\otimes \frac{F_j A''_\infty}{F_{j-1} A''_\infty}
\longto \frac{F_{i+j} A_\infty}{F_{i+j-1} A_\infty}
\]
for all $i,j \in \bZ$.

In a similar way, the spectral sequence pairing $\phi = (\phi^r)_{r\ge1}$, induced by the pairing $(\phi_r)_{r \geq 1}$ per Theorem~\ref{thm:Dou59bIIA}, maps
\begin{align*}
{}' Z^r \otimes {}'' Z^r &\longto Z^r \,, \\
{}' B^r \otimes {}'' Z^r &\longto B^r \,, \\
{}' Z^r \otimes {}'' B^r &\longto B^r \,,
\end{align*}
hence also maps
\begin{align*}
{}' Z^\infty \otimes {}'' Z^\infty &\longto Z^\infty \,, \\
{}' B^\infty \otimes {}'' Z^\infty &\longto B^\infty \,, \\
{}' Z^\infty \otimes {}'' B^\infty &\longto B^\infty \,.
\end{align*}
It follows that $\phi$ also induces
$k$-module homomorphisms
\begin{equation} \label{eq:phiinfty}
\phi^\infty \: {}' E^\infty_i \otimes {}'' E^\infty_j
	\longto E^\infty_{i+j} \,,
\end{equation}
sending $[x] \otimes [y]$ to $[\phi^1(x \otimes y)]$
for any pair of infinite cycles~$x$ and~$y$.
Condition~(SPP~III) ensures that we have the following
compatibility.

\begin{proposition}
\label{prop:multabutment}
Let $\phi = (\phi_r) \: (H', H'') \to H$ be a pairing of
Cartan--Eilenberg systems, with induced pairing
$\phi = (\phi^r) \: ({}' E^*, {}'' E^*) \to E^*$
of spectral sequences, per Theorem~\ref{thm:Dou59bIIA}. Then the pairing $\phi_*$ of
filtered abutments is compatible with the pairing
$\phi^\infty$ of $E^\infty$-pages, in the sense that the
diagram
\[
\xymatrix{
\ds \frac{F_i A'_\infty}{F_{i-1} A'_\infty}
\otimes \frac{F_j A''_\infty}{F_{j-1} A''_\infty}
\ar[r]^-{\bar\phi_*} \ar[d]_-{\beta \otimes \beta}
& \ds \frac{F_{i+j} A_\infty}{F_{i+j-1} A_\infty}
	\ar@{ >->}[d]^-{\beta} \\
{}' E^\infty_i \otimes {}'' E^\infty_j \ar[r]^-{\phi^\infty}
& E^\infty_{i+j}
}
\]
commutes, for all $i,j \in \bZ$.
\end{proposition}

\begin{proof}
A detailed proof is given in \cite{Hel17}*{Prop.~3.4.4}.
\end{proof}

\begin{remark}
As a consequence of Theorem~\ref{thm:Dou59bIIA}, if $(H,\partial)$ is a multiplicative Cartan--Eilenberg system, meaning that it is equipped with a pairing $\phi \: (H,H) \to H$, then the associated spectral sequence~$(E^r,d^r)$ is also multiplicative. Moreover, Proposition~\ref{prop:multabutment} tells us that the induced pairing on the filtered abutment $A_\infty$ is compatible with the induced pairing on the $E^\infty$-page of the spectral sequence. In this situation, we say that $A_\infty$ is a \emph{multiplicative abutment}\index{abutment!multiplicative}. When $(E^r, d^r)$ \emph{converges strongly}\index{strong convergence} to~$A_\infty$, meaning that the
filtration $(F_s A_\infty)_s$ is complete Hausdorff and exhaustive, and that $\beta$ is an isomorphism,
multiplicativity of the abutment means that we can reconstruct the product
$\phi_*$ on~$A_\infty$ from the product $\phi^\infty$ on $E^\infty$,
up to the usual ambiguity created by extensions and filtration shifts.
\end{remark}

\section{Sequences}
\label{subsec:seq}

Our Cartan--Eilenberg systems will in practice be obtained from filtrations and sequences. Let us first set up some terminology, so that it is clear what we are discussing. Again, $G$ denotes a compact Lie group.

\begin{definition}
A sequential diagram $X_\star$ of orthogonal $G$-spectra and $G$-maps of the form
\[
\cdots \longto X_{i-1} \longto X_i \longto X_{i+1} \longto \cdots \,,
\]
indexed over $i \in \bZ$, is called a \emph{sequence}\index{sequence}.
\end{definition}

We can extend the sequence to be indexed over $\bZ \cup \{\pm \infty\}$ by setting 
\[
X_{-\infty} = * \och X_\infty =
\Tel(X_\star) \,,
\]
where
\[
\Tel(X_\star)
	= \bigvee_{i \in \bZ} [i,i+1]_+ \wedge X_i / {\sim}
\]
is the classical \emph{telescope construction}\index{telescope construction|(}. Here, the equivalence relation $\sim$ is given by identifying $\{i\}_+ \wedge X_{i-1}$ with $\{i\}_+ \wedge X_i$ using
the $G$-map $X_{i-1} \to X_i$. There are standard inclusions $X_i \cong \{i\}_+
\wedge X_i \subset \Tel(X_\star)$ for all $i \in \bZ$, and each diagram
\[
\xymatrix{
X_{i-1} \ar[r] \ar[dr] & X_i \ar[d] \\
& \Tel(X_\star)
}
\]
commutes up to preferred homotopy.\index{telescope construction|)}

\begin{definition}
The Cartan--Eilenberg system\index{Cartan--Eilenberg system!associated to a sequence} $(H=H(X_\star), \partial)$ associated to
the sequence $X_\star$ of orthogonal $G$-spectra is given by
\[
H(i, j) = \pi^G_*(X_i \longto X_j)
\]
for all $-\infty \le i \le j \le \infty$, and
\[
\partial \: \pi^G_*(X_j \to X_k) \longto \pi^G_{*-1}(X_i \to X_j)
\]
for all $-\infty \le i \le j \le k \le \infty$.
\end{definition}

Let us elaborate on the definition above. In the $q \geq 0$ case, the expression 
\[
H(i,j) = \pi^G_q(X\to Y) 
\]
denotes the colimit, over the partially ordered set of finite-dimensional $G$-sub-representations~$V$ of the complete $G$-universe~$\sU$, of the groups of homotopy classes $[f', f]$ of pairs $(f',f)$ of $G$-maps $f' \: \Sigma^V S^{q-1} \to X(V)$ and $f \: \Sigma^V D^q \to Y(V)$ making the square
\[
\xymatrix{
\Sigma^V S^{q-1} \ar[r] \ar[d]_-{f'} & \Sigma^V D^q \ar[d]^-{f} \\
X(V) \ar[r] & Y(V)
}
\]
commute. Similar definitions can be made for $q\le0$, but will be left to the reader. By the stability of the homotopy category of orthogonal $G$-spectra there is a natural isomorphism 
\[
\pi^G_q(X \to Y) \cong \pi^G_q(Y \cup CX) \,,
\]
where $Y \cup CX$ is the mapping cone of $X \to Y$. This isomorphism takes the homotopy class $[f',f]$ to (the image in
the colimit over~$V$ of) the homotopy class of the composite map
\[
\Sigma^V S^q \overset{\simeq}\longto \Sigma^V(D^q \cup CS^{q-1})
	\xrightarrow{f \cup Cf'} (Y \cup CX)(V) \,,
\]
where the first map is a ($V$-suspended) homotopy inverse to the collapse map $D^q \cup CS^{q-1} \to D^q/S^{q-1} \cong S^q$. 

The connecting homomorphism $\partial \: \pi^G_q(Y \to Z) \to
\pi^G_{q-1}(X \to Y)$ mentioned in the definition takes the homotopy class $[g',g]$ of a pair of $G$-maps $g' \: \Sigma^V S^{q-1} \to Y(V)$ and $g \: \Sigma^V D^q \to Z(V)$ to the homotopy class $[*,g'\pi]$ of the maps 
\[
\ast \: \Sigma^V S^{q-2} \longto * \longto X(V) \och g'\pi :\ \Sigma^V D^{q-1} \overset{\pi}\longto \Sigma^V S^{q-1} \longto Y(V) \,,
\]
where $\pi \: D^{q-1} \to S^{q-1}$ identifies $D^{q-1}/S^{q-2}$ with $S^{q-1}$. The diagram 
\[
\xymatrix{
\Sigma^V S^{q-2} \ar[r] \ar[d] & \Sigma^V D^{q-1} \ar[d]^-{\pi} \\
{*} \ar[r] \ar[d] & \Sigma^V S^{q-1} \ar[d]^-{g'} \\
X(V) \ar[r] & Y(V)
}
\]
evidently commutes. Under the isomorphism $\pi^G_{q-1}(X \to Y) \cong \pi^G_{q-1}(Y \cup CX)$ the homotopy class $[\ast ,g'\pi]$ corresponds to the homotopy class of the composite map
\[
\Sigma^V S^{q-1} \overset{g'}\longto Y(V) \longto (Y \cup CX)(V) \,.
\]
Note that the graded abelian group $\pi^G_*(X_i \to X_j)$ is functorial in $i \le j$, the homomorphism $\partial$ is natural in $i \le j \le k$, and the sequence
\[
\cdots \to \pi^G_q(X_i \to X_j) \to \pi^G_q(X_i \to X_k)
	\to \pi^G_q(X_j \to X_k) \overset{\partial}\longto
	\pi^G_{q-1}(X_i \to X_j) \to \cdots
\]
is exact for all $i \le j \le k$ and $q \in \bZ$. The canonical homomorphism
\[
\colim_j \pi^G_*(X_j) \overset{\cong}\longto \pi^G_* \Tel(X_\star)
\]
is an isomorphism, which implies that condition~(SP.5) is satisfied. Hence $(H, \partial)$ is indeed a Cartan--Eilenberg system in the sense of Definition~\ref{def:CEsystemII}.

We can extract two different exact couples~\cite{Mas52} from $(H(X_\star),\partial)$, but shall only be concerned with the `left' couple\index{left couple} of Definition~\ref{def:excpl+spseq}. Explicitly, the exact couple $(A, E^1) = (A(X_\star), E^1(X_\star))$ associated to~$X_\star$ is given by 
\[
A_s = \pi^G_*(X_s) \och E^1_s = \pi^G_*(X_{s-1} \to X_s) \,,
\] 
fitting together in the exact triangle
\[
\xymatrix{
\pi^G_*(X_{s-1}) \ar[r]
	& \pi^G_*(X_s) \ar[d] \\
& \pi^G_*(X_{s-1} \to X_s) \ar[ul]^-{\partial}
}
\]
where $\partial$ has total degree~$-1$.

Recall from~\cite{Boa99}*{Def.~5.10} that the spectral sequence associated to the unrolled exact couple~$(A, E^1)$ is said to be \emph{conditionally convergent}\index{conditional convergence} to the abutment
\begin{equation} \label{eq:AinftyTelXstar}
A_\infty = \colim_s \pi^G_*(X_s)
	\cong \pi^G_* \Tel(X_\star)
\end{equation}
if and only if
\[
A_{-\infty} = \lim_s A_s = 0 \och 
RA_{-\infty} = \Rlim_s A_s = 0 \,.
\]
Here $\Rlim = \lim^1$ denotes the (first right) derived limit of
a sequence. In view of the short exact sequence
\[
0 \longto \Rlim_s \pi^G_{*+1}(X_s) \longto \pi^G_*(\holim_s X_s)
	\longto \lim_s \pi^G_*(X_s) \longto 0
\]
this is equivalent to the condition that
\[
\pi^G_*(\holim_s X_s) = 0 \,.
\]
In particular, conditional convergence holds if $\holim_s X_s \simeq_G *$.

The \emph{spectral sequence}\index{spectral sequence!associated to a sequence} $(E^r = E^r(X_\star), d^r)_{r\ge1}$
associated to the sequence $X_\star$ (and the Cartan--Eilenberg system $(H(X_\star), \partial)$, and the exact couple~$(A(X_\star), E^1(X_\star))$) has
\[
E^1_{s,t} = \pi^G_{s+t}(X_{s-1} \to X_s)
\]
and $d^1 \: E^1_{s,t} \to E^1_{s-1,t}$ is equal to the composite
homomorphism
\[
\pi^G_{s+t}(X_{s-1} \to X_s) \overset{\partial}\longto
\pi^G_{s+t-1}(X_{s-1}) \longto
\pi^G_{s+t-1}(X_{s-2} \to X_{s-1}) \,.
\]
Here $s+t$ is the \emph{total degree}\index{degree!total}, $s$ is the \emph{filtration degree}\index{degree!filtration}, and $t$ will be called the \emph{internal degree}\index{degree!internal}.
The~$d^r$-differentials have the form
\[
d^r \: E^r_{s,t} \longto E^r_{s-r,t+r-1}
\]
and there are preferred isomorphisms $H(E^r, d^r) \cong E^{r+1}$ for all $r\ge1$.

\section{Filtrations}
\label{subsec:filt}

The category of orthogonal $G$-spectra is based topological, meaning
that it is enriched in the closed symmetric monoidal category of
compactly generated weak Hausdorff spaces with base point.

\begin{definition}
Let $I = [0,1]$, with boundary $\partial I = \{0,1\}$.
\begin{itemize}
\item A $G$-map $i \: A \to X$ of orthogonal $G$-spectra is an
\emph{$h$-cofibration}\index{h-cofibration@$h$-cofibration}\index{Hurewicz cofibration} (also called: Hurewicz cofibration) if it has the homotopy
extension property\index{homotopy extension property} with respect to any target $Z$:
\[
\xymatrix{
X \cup_A A \wedge I_+ \ar[r] \ar[d]
	& Z \\
X \wedge I_+ \ar@{-->}[ur]
	& \,.
}
\]
\item A $G$-map $p \: E \to B$ of orthogonal $G$-spectra is an
\emph{$h$-fibration}\index{h-fibration@$h$-fibration} (also called: Hurewicz fibration\index{Hurewicz fibration}) if it has the homotopy
lifting property\index{homotopy lifting property} with respect to any source~$X$:
\[
\xymatrix{
X \ar[r] \ar[d]
	& E \ar[d]^-{p} \\
X \wedge I_+ \ar[r] \ar@{-->}[ur]
	& B \rlap{\,.}
}
\]
\item Adapting \cite{SV02}*{Def.~2.4}, we say that $i \: A \to X$ is a \emph{strong $h$-cofibration}\index{strong $h$-cofibration} if the~$G$-map $X \cup_A A \wedge I_+ \to X \wedge I_+$ has the left lifting property with respect to any $h$-fibration:
\[
\xymatrix{
X \cup_A A \wedge I_+ \ar[r] \ar[d]
	& E \ar[d]^-{p} \\
X \wedge I_+ \ar[r] \ar@{-->}[ur]
	& B \rlap{\,.}
}
\]
\end{itemize}
\end{definition}

Strong $h$-cofibrations are closed under cobase change, retracts,
arbitrary sums, and sequential colimits \cite{SV02}*{Lem.~2.6}.  For each map $f \: X \to Y$ the inclusion $i_0 \: X \to Y \cup_X X \wedge I_+$ is a strong $h$-cofibration \cite{SV02}*{Rmk.~3.3(2)}.  It follows that each $q$-cofibration\index{q-cofibration@$q$-cofibration} (= Quillen cofibration\index{Quillen cofibration}, \cite{MM02}*{Def.~III.2.3}) is a strong~$h$-cofibration. Each strong~$h$-cofibration is evidently an $h$-cofibration. Our main reason for working with strong $h$-cofibrations is the availability of the following theorem.

\begin{theorem}[\cite{SV02}*{Thm.~2.7}] \label{thm:pushprodstrongcof}
If $i \: A \to X$ and $j \: B \to Y$ are strong $h$-cofibrations,
then the pushout-product map\index{pushout-product map}
\[
i \wedge 1 \cup 1 \wedge j \:
	A \wedge Y \cup_{A \wedge B} X \wedge B \longto X \wedge Y
\]
is a strong $h$-cofibration.
\end{theorem}

We can now specify well-behaved sequences, called filtrations, for which we can directly prove that pairings of sequences induce pairings of Cartan--Eilenberg systems and of spectral sequences.

\begin{definition}
Let $X_\star$ be a sequence of orthogonal $G$-spectra.  We say that~$X_\star$ is a \emph{filtration}\index{filtration} if each~$G$-map $X_{i-1} \to X_i$
for $i \in \bZ$ is a strong $h$-cofibration.
\end{definition}

In particular, if $X_\star$ is a filtration, then the $G$-maps are all $h$-cofibrations, so the canonical maps
\[
X_j \cup CX_i \longto X_j/X_i \och \Tel(X_\star) \longto \colim_i X_i = \bigcup_i X_i
\]
are $G$-equivalences, so that 
\[
H(i,j) \cong \pi^G_*(X_j/X_i) \och A_\infty
\cong \pi^G_*\left(\bigcup_i X_i\right)
\]
in the associated Cartan--Eilenberg system. 

We can always approximate a sequence $X_\star$ with an equivalent
filtration $T_\star(X)$.\index{telescope construction|(} To do this, we proceed as follows. For each integer $j$ we let
\[
T_j(X) = \{j\}_+ \wedge X_j \vee \bigvee_{i<j} [i,i+1]_+ \wedge X_i / {\sim}
\]
be the subspectrum of $\Tel(X_\star)$ with telescope coordinate in the interval $(-\infty, j]$ within the real line $(-\infty, \infty) = \bigcup_i [i,i+1]$. The sequence $T_\star(X)$ given by the inclusions
\[
\dots \longto T_{j-1}(X) \longto T_j(X) \longto T_{j+1}(X) \longto \dots
\]
is then a filtration. 

For each integer $j$ there is a deformation retraction
\[
\epsilon_j \: T_j(X) \longto X_j
\]
identifying $\{j\}_+ \wedge X_j$ with $X_j$ and mapping
$[i,i+1]_+ \wedge X_i$ to $X_j$ by the evident composition $[i,i+1]_+ \wedge X_i \to X_i \to X_j$, for each $i < j$.  The resulting diagram
\[
\xymatrix{
& \dots \ar[r]
& T_{j-1}(X) \ar[r] \ar[d]_-{\simeq_G}
& T_j(X) \ar[r] \ar[d]_-{\simeq_G}
& T_{j+1}(X) \ar[r] \ar[d]_-{\simeq_G}
& \dots \\
& \dots \ar[r]
& X_{j-1} \ar[r]
& X_j \ar[r]
& X_{j+1} \ar[r]
& \dots
}
\]
commutes, and defines a $G$-equivalence of sequences $\epsilon \: T_\star(X)\to X_\star$.  It follows that the associated maps
\[
{*} = T_{-\infty}(X) \longto X_{-\infty} = {*} \och T_\infty(X) \longto X_\infty =\Tel(X_\star)
\]
are both $G$-equivalences. Hence the induced maps of Cartan--Eilenberg systems
\[
H(T_\star(X))(i,j) = \pi^G_*(T_i(X) \longto T_j(X))
\overset{\cong}\longto
\pi^G_*(X_i \to X_j) = H(X_\star)(i,j) \,,
\]
and of spectral sequences
\[
(E^r(T_\star(X)), d^r)_{r\ge1}
	\overset{\cong}\longto (E^r(X_\star), d^r)_{r\ge1} \,,
\]
are both isomorphisms. Their common abutment for convergence
to the colimit is $A_\infty(X_\star) \cong \pi^G_* \Tel(X_\star)$,
filtered by the image subsequence\index{image filtration}
\[
F_s A_\infty(X_\star) \cong F_s \pi^G_* \Tel(X_\star)
	= \im(\pi^G_*(X_s) \to \pi^G_* \Tel(X_\star)) \,.
\]
\index{telescope construction|)}
\begin{remark}
Some form of cofibrant replacement of maps is necessary to convert general
sequences to filtrations.  We have chosen to use mapping cylinders\index{mapping cylinder} and
telescopes\index{telescope construction}, which have convenient monoidal properties.  For finite groups,
Hesselholt and Madsen~\cite{HM03}*{\S4.3} instead use a functorial $G$-CW
replacement to convert~$G$-spectra to $G$-CW spectra.  There exists a functorial $G$-CW replacement\index{functorial $G$-CW replacement} also for compact Lie groups \cite{Sey83},
but its construction is comparatively intricate, and the monoidal
properties are less clear, which may partially justify our choice.
\end{remark}

\section{Pairings of sequences}

We now turn to discussing pairings of sequences and how these behave under passage to mapping telescopes.

\begin{definition}
\label{def:pairseq}
Let $X_\star$, $Y_\star$ and $Z_\star$ be sequences of orthogonal
$G$-spectra.  A \emph{pairing}\index{pairing!of sequences} $\phi \: (X_\star, Y_\star) \to Z_\star$
is a collection of $G$-maps
\[
\phi_{i,j} \: X_i \wedge Y_j \longto Z_{i+j}
\]
for all integers~$i$ and~$j$, making the squares
\begin{equation} \label{eq:pairing}
\begin{tikzcd}[column sep = large]
X_{i-1} \wedge Y_j \arrow[r,"\phi_{i-1,j}"] \arrow[d] & Z_{i+j-1} \arrow[d]
	& X_i \wedge Y_{j-1} \arrow[l,"\phi_{i,j-1}"'] \ar[d] \\
X_i \wedge Y_j \arrow[r,"\phi_{i,j}"] & Z_{i,j}
	& X_i \wedge Y_j \arrow[l,"\phi_{i,j}"']
\end{tikzcd}
\end{equation}
commute. We say that a sequence $X_\star$ is \emph{multiplicative}\index{sequence!multiplicative} if it comes equipped with a pairing $\phi \: (X_\star, X_\star)
\to X_\star$.
\end{definition}

We note that, from \cite{MM02}*{\S II.3} and \cite{Sch18}*{\S3.5}, the smash product $X_i \wedge Y_j$ of orthogonal $G$-spectra is defined in such a way that $\phi_{i,j}$ associates to each pair of~$G$-representations $U$ and~$V$ a $G$-map of based~$G$-spaces
\[
\phi_{i,j}(U,V) \: X_i(U) \wedge Y_j(V) \longto Z_{i+j}(U \oplus V)
\,,
\]
subject to bilinearity relations for varying~$U$ and~$V$.

\begin{lemma} \label{lem:Telphi}
A pairing of sequences $\phi \: (X_\star, Y_\star) \to Z_\star$
induces a pairing of sequences $T(\phi) \: (T_\star(X),
T_\star(Y)) \to T_\star(Z)$, such that the diagram
\[
\xymatrix@C+2pc{
T_i(X) \wedge T_j(Y) \ar[r]^-{T(\phi)_{i,j}} \ar[d]_-{\simeq_G}
	& T_{i+j}(Z) \ar[d]^-{\simeq_G} \\
X_i \wedge Y_j \ar[r]^-{\phi_{i,j}}
	& Z_{i+j}
}
\]
commutes for all integers $i$ and~$j$.
\end{lemma}

\begin{proof}
Given a pairing $\phi$ of sequences, we can form $G$-maps
\begin{equation} \label{eq:XiYj-TelZstar}
[i,i+1]_+ \wedge X_i \wedge [j,j+1]_+ \wedge Y_j
	\longto [k,k+2]_+ \wedge Z_k
	\longto \Tel(Z_\star)
\end{equation}
for any integers $i$ and~$j$, with $k = i+j$. Here $[i,i+1] \times [j,j+1] \to [k,k+2]$ sends $(x, y)$ to $x+y$, while $X_i \wedge Y_j \to Z_k$ is given by $\phi_{i,j}$. The second map factors through
\[
[k,k+1]_+ \wedge Z_k \cup [k+1,k+2]_+ \wedge Z_{k+1} \,,
\]
and is given by $Z_k \to Z_{k+1}$ on $[k+1,k+2] \subset [k,k+2]$.  The
maps~\eqref{eq:XiYj-TelZstar} for varying~$i$ and~$j$ are compatible with
the identifications defining $\Tel(X_\star)$ and $\Tel(Y_\star)$, hence
combine to define a $G$-map 
\[
\Tel(\phi) \: \Tel(X_\star) \wedge \Tel(Y_\star) \longto \Tel(Z_\star) \,.
\] 
By construction, it restricts to compatible $G$-maps
\[
T(\phi)_{i,j} \: T_i(X) \wedge T_j(Y)
	\longto T_{i+j}(Z)
\]
for all integers~$i$ and~$j$, defining the pairing of sequences
$T(\phi) \: (T_\star(X), T_\star(Y)) \to T_\star(Z))$.
It is then clear that the square in the lemma commutes, and
that the vertical maps are $G$-equivariant deformation retractions.
\end{proof}

\begin{corollary} \label{cor:telxmult}\index{sequence!multiplicative}
If $(X_\star, \phi)$ is a multiplicative sequence, then so
is the sequence $(T_\star(X), T(\phi))$. Moreover, the equivalence $\epsilon \: T_\star(X) \to X_\star$ respects the multiplicative structures.
\end{corollary}

\section{Pairings of Cartan--Eilenberg systems, I}
\label{subsec:pairingsCEI}

The goal of the following two sections is to show that a pairing of sequences gives rise to a pairing of the resulting Cartan--Eilenberg systems. By Theorem~\ref{thm:Dou59bIIA} and Proposition~\ref{prop:multabutment} this is enough to guarantee that we have a pairing of the associated spectral sequences in such a way that the induced pairing on filtered abutments is compatible with the pairing on $E^\infty$-pages. Referring back to Definition~\ref{def:pairfinCEsys} and Definition~\ref{def:pairextCEsys}, we note that there are three things to check. In this section we deal with (SPP I) and (SPP III)\index{SPPa@(SPP I)}\index{SPPc@(SPP III)}.

Let $\phi \: (X_\star, Y_\star) \to Z_\star$ be a pairing of sequences. For integers $i$, $j$ and~$r$, with $r\ge1$, we define induced pairings
\[
\phi_r \: H(X_\star)(i-r,i) \otimes H(Y_\star)(j-r,j)
	\longto H(Z_\star)(i+j-r,i+j)
\]
as homomorphisms
\[
\phi_r \: \pi^G_p(X_{i-r} \to X_i) \otimes \pi^G_q(Y_{j-r} \to Y_j)
	\longto \pi^G_{p+q}(Z_{i+j-r} \to Z_{i+j}) \,.
\]
Here $p$ and~$q$ range over all integers, but for (relative)
brevity we concentrate on the case when $p\ge0$ and~$q\ge0$.
Given two pairs of vertical $G$-maps
\[
\xymatrix{
\Sigma^U S^{p-1} \ar[r] \ar[d]_-{f'} & \Sigma^U D^p \ar[d]^-{f} \\
X_{i-r}(U) \ar[r] & X_i(U)
}
\qqandqq
\xymatrix{
\Sigma^V S^{q-1} \ar[r] \ar[d]_-{g'} & \Sigma^V D^q \ar[d]^-{g} \\
Y_{j-r}(V) \ar[r] & Y_j(V)
}
\]
we first form the commutative diagram
\[
\xymatrix@C-1.4pc{
\Sigma^{U \oplus V} S^{p-1} \wedge S^{q-1} \ar[rr] \ar[dr] \ar[dd]_-{\cong}
	&& \Sigma^{U \oplus V} D^p \wedge S^{q-1} \ar[dr] \ar[dd]|\hole^(0.75){\cong} \\
	& \Sigma^{U \oplus V} S^{p-1} \wedge D^q \ar[rr] \ar[dd]_(0.25){\cong}
	&& \Sigma^{U \oplus V} D^p \wedge D^q \ar[dd]^-{\cong} \\
\Sigma^U S^{p-1} \wedge \Sigma^V S^{q-1} \ar[rr]|(0.514)\hole \ar[dr] \ar[dd]_-{f' \wedge g'}
	&& \Sigma^U D^p \wedge \Sigma^V S^{q-1} \ar[dr] \ar[dd]|\hole^(0.75){f \wedge g'} \\
	& \Sigma^U S^{p-1} \wedge \Sigma^V D^q \ar[rr] \ar[dd]_(0.25){f' \wedge g}
	&& \Sigma^U D^p \wedge \Sigma^V D^q \ar[dd]^-{f \wedge g} \\
X_{i-r}(U) \wedge Y_{j-r}(V) \ar[rr]|(0.514)\hole \ar[dr] \ar[dd]_-{\phi_{i-r,j-r}(U,V)}
	&& X_i(U) \wedge Y_{j-r}(V) \ar[dr] \ar[dd]|\hole^(0.75){\phi_{i,j-r}(U,V)} \\
	& X_{i-r}(U) \wedge Y_j(V) \ar[rr] \ar[dr]_(0.4){\phi_{i-r,j}(U,V)\quad}
	&& X_i(U) \wedge Y_j(V) \ar[d]^-{\phi_{i,j}(U,V)} \\
Z_{i+j-2r}(U \oplus V) \ar[rr]
	&& Z_{i+j-r}(U \oplus V) \ar[r]
	& Z_{i+j}(U \oplus V) \,.
}
\]
For typographical reasons, we will often suppress the stabilising
$G$-representations~$U$ and~$V$ and simply display this diagram as
\[
\xymatrix{
S^{p-1} \wedge S^{q-1} \ar[rr] \ar[dr] \ar[dd]_-{f' \wedge g'}
	&& D^p \wedge S^{q-1} \ar[dr] \ar[dd]|\hole^(0.75){f \wedge g'} \\
	& S^{p-1} \wedge D^q \ar[rr] \ar[dd]_(0.25){f' \wedge g}
	&& D^p \wedge D^q \ar[dd]^-{f \wedge g} \\
X_{i-r} \wedge Y_{j-r} \ar[rr]|(0.514)\hole \ar[dr] \ar[dd]_-{\phi_{i-r,j-r}}
	&& X_i \wedge Y_{j-r} \ar[dr] \ar[dd]|\hole^(0.75){\phi_{i,j-r}} \\
	& X_{i-r} \wedge Y_j \ar[rr] \ar[dr]_-{\phi_{i-r,j}}
	&& X_i \wedge Y_j \ar[d]^-{\phi_{i,j}} \\
Z_{i+j-2r} \ar[rr]
	&& Z_{i+j-r} \ar[r]
	& Z_{i+j} \,.
}
\]
Let
\begin{align*}
S^{p+q-1} &= S^{p-1} \wedge D^q \cup_{S^{p-1} \wedge S^{q-1}}
	D^p \wedge S^{q-1} \\
W &= X_{i-r} \wedge Y_j \cup_{X_{i-r} \wedge Y_{j-r}} X_i \wedge Y_{j-r}
\end{align*}
denote the pushouts in the squares of the upper and middle layer of the diagram, respectively. In particular, $S^{p+q-1}$ is the boundary of $D^p \wedge D^q \cong D^{p+q}$.
We then have an induced commutative diagram
\begin{equation} \label{eq:sdxywz}
\begin{tikzcd}
S^{p-1} \wedge S^{q-1} \arrow[r] \arrow[d,"f' \wedge g'"]
	& S^{p+q-1} \arrow[r] \arrow[d,"(f \wedge g)'"]
	& D^{p+q} \arrow[d,"f \wedge g"] \\
X_{i-r} \wedge Y_{j-r} \arrow[r] \arrow[d,"\phi_{i-r,j-r}"]
	& W \arrow[r] \ar[d,"\phi_W"]
	& X_i \wedge Y_j \arrow[d,"\phi_{i,j}"] \\
Z_{i+j-2r} \arrow[r]
	& Z_{i+j-r} \arrow[r]
	& Z_{i+j} \,.
\end{tikzcd}
\end{equation}
Here, $(f \wedge g)' \: S^{p+q-1} \to W$ is the induced map between the pushouts in the top and bottom square of the boxed-shaped diagram appearing above, and $\phi_W$ is the induced map in the diagram
\[
\begin{tikzcd}
X_{i-r} \wedge Y_{j-r} \arrow[r] \arrow[d] & X_i \wedge Y_{j-r} \arrow[d] \arrow[rdd , bend left,"\phi_{i,j-r}"] & \\
X_{i-r} \wedge Y_j \arrow[r] \arrow[rrd, bend right,"\phi_{i-r,j}"'] & W \arrow[rd,"\phi_W"] & \\
& & Z_{i+j-r} \,.
\end{tikzcd}
\]
We define the homomorphism 
\[
\phi_r \: \pi_p^G(X_{i-r} \to X_i) \otimes \pi_q^G(Y_{j-r} \to Y_j) \to \pi_{p+q}^G(Z_{i+j-r} \to Z_{i+j})
\]
as sending $[f', f] \otimes [g', g]$ to the homotopy class of the pair
\[
\phi_W (f \wedge g)' \: S^{p+q-1} \to Z_{i+j-r} \och \phi_{i,j} (f \wedge g) \: D^{p+q} \to Z_{i+j} \,,
\]
which is an element of $\pi^G_{p+q}(Z_{i+j-r} \to Z_{i+j})$. As a diagram, this pair is visualised as the commutative square
\[
\xymatrix{
S^{p+q-1} \ar[r] \ar[d]_-{\phi_W(f \wedge g)'}
	& D^{p+q} \ar[d]^-{\phi_{i,j}(f \wedge g)} \\
Z_{i+j-r} \ar[r]
	& Z_{i+j} \,.
}
\]
In symbols:
\[
\phi_r \: [f',f] \otimes [g',g]
	\longmapsto
	[\phi_W (f \wedge g)', \phi_{i,j} (f \wedge g)] \,.
\]
Spelled out with the stabilising $G$-representations~$U$ and~$V$, this diagram should be interpreted as the commutative diagram
\[
\xymatrix{
\Sigma^{U \oplus V} S^{p+q-1} \ar[r] \ar[d]_-{\cong}
        & \Sigma^{U \oplus V} D^{p+q} \ar[d]^-{\cong} \\
\Sigma^U S^{p-1} \wedge \Sigma^V D^q \cup \Sigma^U D^p \wedge \Sigma^V S^{q-1}
	\ar[r] \ar[d]_-{(f \wedge g)'}
        & \Sigma^U D^p \wedge \Sigma^V D^q \ar[d]^-{f \wedge g} \\
X_{i-r}(U) \wedge Y_j(V) \cup X_i(U) \wedge Y_{j-r}(V)
	\ar[r] \ar[d]_-{\phi_W(U,V)}
        & X_i(U) \wedge Y_j(V) \ar[d]^-{\phi_{i,j}(U,V)} \\
Z_{i+j-r}(U \oplus V) \ar[r]
        & Z_{i+j}(U \oplus V)
}
\]
The pushouts on the left hand side are formed along $\Sigma^U S^{p-1}
\wedge \Sigma^V S^{q-1}$ and $X_{i-r}(U) \wedge Y_{j-r}(V)$, respectively.

\begin{remark}
We note that the pushout $W$ is not generally equivalent to the
corresponding homotopy pushout, but this will hold if $X_\star$
and $Y_\star$ are filtrations.

We also note that if one only has a weak pairing, in the sense that
the squares in diagram~\eqref{eq:pairing} commute up to homotopy, then there is in
general no preferred commuting homotopy in the diagram
\[
\xymatrix{
W \ar[r] \ar[d]_-{\phi_W} & X_i \wedge Y_j \ar[d]^-{\phi_{i,j}} \\
Z_{i+j-r} \ar[r] & Z_{i+j} \rlap{\,,}
}
\]
and therefore no well-defined pairing $\phi_r$.  Any construction of
spectral sequence pairings that only assumes such compatibility at the
level of the (stable) homotopy category is therefore likely to contain
a logical gap.
\end{remark}

The pairing $\phi_r$ is evidently natural in $i$, $j$ and~$r$, in the sense that the square
\begin{equation} \label{eq:SPP-I}
\begin{tikzcd}
H(X_\star)(i-r,i) \otimes H(Y_\star)(j-r,j) \arrow[r,"\phi_r"] \arrow[d]
	& H(Z_\star)(i+j-r,i+j) \arrow[d] \\
H(X_\star)(i'-r',i') \otimes H(Y_\star)(j'-r',j') \arrow[r,"\phi_{r'}"]
	& H(Z_\star)(i'+j'-r',i'+j')
\end{tikzcd}
\end{equation}
commutes for all integers $i$, $j$, $r\ge1$, $i'$, $j'$, $r'\ge1$ with $i \le i'$, $i-r \le i'-r'$, $j \le j'$ and $j-r \le j'-r'$. As we recalled from~\cite{Dou59b}*{\S~II~A} in
Definition~\ref{def:pairfinCEsys}, this is the first (SPP~I) of two conditions for~$(\phi_r)_{r\ge1}$ to define a pairing of (finite) Cartan--Eilenberg systems.

We now check condition (SPP~III). The pairings $\phi_r$ can be extended to the case~$r = \infty$ by letting
\[
\phi_\infty \: H(X_\star)(-\infty,i) \otimes H(Y_\star)(-\infty,j)
	\longto H(Z_\star)(-\infty,i+j) \,
\]
be defined by homomorphisms
\[
\phi_\infty \: \pi^G_p(X_i) \otimes \pi^G_q(Y_j)
	\longto \pi^G_{p+q}(Z_{i+j}) \,.
\]
Given $G$-maps $f \: \Sigma^U S^p \to X_i(U)$ and $g \: \Sigma^V S^q \to Y_j(V)$, the homomorphism $\phi_\infty$ sends the homotopy classes $[f]$ and $[g]$ to the homotopy class of the composite
\[
\Sigma^{U \oplus V} S^{p+q} \cong \Sigma^U S^p \wedge \Sigma^V S^q \xrightarrow{f \wedge g} X_i(U) \wedge Y_j(V) \xrightarrow{\phi_{i,j}(U,V)} Z_{i+j}(U \oplus V) 
\]
In symbols, suppressing~$U$ and~$V$:
\[
\phi_\infty \: [f] \otimes [g] \longmapsto [\phi_{i,j}(f \wedge g)] \,.
\]
These pairings are natural
in~$i$ and~$j$, which verifies the first part of~(SPP III).

Recalling that our convention is such that $X_{-\infty} = Y_{-\infty} = Z_{-\infty} =*$, we note that the isomorphism $\pi^G_p(X_i)\cong \pi^G_p(X_{-\infty} \to X_i)$ takes the homotopy class of $f$ to the homotopy class $[*, f\pi]$ of the pair $* \: \Sigma^U S^{p-1}
\to X_{-\infty}(U)$ and $f \pi \: \Sigma^U D^p \to X_i(U)$. Here  $\pi \: D^{p} \to S^{p}$ identifies $D^{p}/S^{p-1}$ with $S^{p}$. The pairing $\phi_\infty$ then corresponds to the pairing $\phi_r$ as defined in the paragraph above, for $r=\infty$, with every reference to $X_{i-r}$, $Y_{j-r}$, $Z_{i+j-r}$ and $Z_{i+j-2r}$ replaced by~$*$. By the discussion above, it then also follows that the extended naturality condition
\begin{equation} \label{eq:SPP-III}
\begin{tikzcd}
H(X_\star)(-\infty,i) \otimes H(Y_\star)(-\infty,j) \arrow[r,"\phi_\infty"] \arrow[d]
        & H(Z_\star)(-\infty,i+j) \arrow[d] \\
H(X_\star)(i-r,i) \otimes H(Y_\star)(j-r,j) \arrow[r,"\phi_{r}"]
        & H(Z_\star)(i+j-r,i+j)
\end{tikzcd}
\end{equation}
holds for the pairings $\phi_r$ with $1 \le r \le \infty$. This is the second part of condition~(SPP~III) from Definition~\ref{def:pairextCEsys}.

\section{Pairings of Cartan--Eilenberg systems, II}
\label{subsec:pairingsCEII}

Having proved (SPP I) and (SPP III), we now turn to the second condition (SPP~II)\index{SPPb@(SPP II)} from~\cite{Dou59b}*{\S~II~A}. Recall that it says that, for
$(\phi_r)_{r\ge1}$ to define a pairing of Cartan--Eilenberg
systems, we want the Leibniz rule
\begin{equation} \label{eq:SPP-II}
\partial \phi_r = \phi_1 (\partial \otimes \eta)
	+ \phi_1 (\eta \otimes \partial)
\end{equation}
to hold. That is, we want the composite
\begin{multline*}
H(X_\star)(i-r,i) \otimes H(Y_\star)(j-r,j)
	\xrightarrow{\phi_r}
H(Z_\star)(i+j-r,i+j) \\
	\xrightarrow{\partial}
H(Z_\star)(i+j-r-1,i+j-r)
\end{multline*}
to be equal to the sum of the composite homomorphisms
\begin{multline*}
H(X_\star)(i-r,i) \otimes H(Y_\star)(j-r,j)
	\xrightarrow{\partial \otimes \eta}
H(X_\star)(i-r-1,i-r) \otimes H(Y_\star)(j-1,j) \\
	\xrightarrow{\phi_1}
H(Z_\star)(i+j-r-1,i+j-r)
\end{multline*}
and
\begin{multline*}
H(X_\star)(i-r,i) \otimes H(Y_\star)(j-r,j)
	\xrightarrow{\eta \otimes \partial}
H(X_\star)(i-1,i) \otimes H(Y_\star)(j-r-1,j-r) \\
	\xrightarrow{\phi_1}
H(Z_\star)(i+j-r-1,i+j-r) \,.
\end{multline*}
Here
\begin{align*}
\eta &\: H(X_\star)(i-r,i) \longto H(X_\star)(i-1,i) \\
\eta &\: H(Y_\star)(j-r,j) \longto H(Y_\star)(j-1,j)
\end{align*}
denote the natural maps.
Regarding signs in the Leibniz rule, we recall the convention that
\[
(\partial \otimes 1)(x \otimes y) = \partial x \otimes y \och (1
\otimes \partial)(x \otimes y) = (-1)^p x \otimes \partial y \,,
\] for $x \in \pi^G_p(X_{i-r} \to X_i)$ in total degree~$p$ of $H(X_\star)(i-r, i)$.

To verify condition~\eqref{eq:SPP-II} for a given pairing $\phi \:
(X_\star, Y_\star) \to Z_\star$, it follows from the naturality of
the boundary homomorphisms~$\partial$, and the case $i = i'$, $j = j'$,
$r \ge r' = 1$ of~\eqref{eq:SPP-I}, that it suffices to establish
the rule
\begin{equation} \label{eq:leibniz}
\partial \phi_r = \phi_r (\partial \otimes 1)
	+ \phi_r (1 \otimes \partial) \,.
\end{equation}
Here the left hand side is the composite
\begin{multline*}
H(X_\star)(i-r,i) \otimes H(Y_\star)(j-r,j)
	\xrightarrow{\phi_r}
H(Z_\star)(i+j-r,i+j) \\
	\xrightarrow{\partial}
H(Z_\star)(i+j-2r,i+j-r) \,,
\end{multline*}
and the right hand side is the sum of the two composite homomorphisms
\begin{multline*}
H(X_\star)(i-r,i) \otimes H(Y_\star)(j-r,j)
	\xrightarrow{\partial \otimes 1}
H(X_\star)(i-2r,i-r) \otimes H(Y_\star)(j-r,j) \\
	\xrightarrow{\phi_r}
H(Z_\star)(i+j-2r,i+j-r)
\end{multline*}
and
\begin{multline*}
H(X_\star)(i-r,i) \otimes H(Y_\star)(j-r,j)
	\xrightarrow{1 \otimes \partial}
H(X_\star)(i-r,i) \otimes H(Y_\star)(j-2r,j-r) \\
	\xrightarrow{\phi_r}
H(Z_\star)(i+j-2r,i+j-r) \,.
\end{multline*}
We shall now show that the identity~\eqref{eq:leibniz} holds for pairings
of filtrations of orthogonal $G$-spectra.  Thereafter we use approximation
by mapping telescopes to deduce that the identity holds for pairings of
arbitrary sequences, as well.

\begin{proposition} \label{prop:Leibpairfilt}
If $\phi \: (X_\star, Y_\star) \to Z_\star$ is a pairing of
sequences of orthogonal~$G$-spectra, and $X_\star$ and~$Y_\star$
are filtrations, then
\[
\partial \phi_r = \phi_r (\partial \otimes 1) + \phi_r (1 \otimes \partial)
\]
as homomorphisms
\[
H(X_\star)(i-r,i) \otimes H(Y_\star)(j-r,j)
	\longto H(Z_\star)(i+j-2r,i+j-r)
\]
for all integers $i,j,r$ with $r\ge1$.
\end{proposition}

\begin{proof}
In this proof we will, for the same typographical reasons as
in Section~\ref{subsec:pairingsCEI}, suppress the stabilising
representations $U$ and~$V$ implicit in the presentation of elements of
$\pi^G_p(X_{i-r} \to X_i)$ and~$\pi^G_q(Y_{j-r} \to Y_j)$ by homotopy classes of pairs~$(f',f)$ and~$(g',g)$ of $G$-maps.  The reader can reconstruct how the diagrams could be embellished with these suspensions
and shifts.

For each map $A \to B$ we have natural maps $B \to B \cup CA \to B/A$ to the homotopy cofibre and cofibre. The right hand map is
an equivalence when $A \to B$ is an $h$-cofibration. Applied to the left hand maps in diagram~\eqref{eq:sdxywz}, this gives us a commutative diagram
\[
\xymatrix@C-4.5pt{
S^{p+q-1} \ar[r] \ar[d]_-{(f \wedge g)'}
	& S^{p+q-1} \cup C(S^{p-1} \wedge S^{q-1}) \ar[r]^-{\simeq}
		\ar[d]^-{(f \wedge g)' \cup C(f' \wedge g')}
	& S^{p-1} \wedge S^q \vee S^p \wedge S^{q-1}
		\ar[d]^-{f' \wedge g'' \vee f'' \wedge g'} \\
W \ar[r] \ar[d]_-{\phi_W}
	& W \cup C(X_{i-r} \wedge Y_{j-r})
		\ar[r]^-{(\simeq)}_-{\Theta} \ar[d]^-{\Phi}
	& X_{i-r} \wedge Y_j/Y_{j-r} \vee X_i/X_{i-r} \wedge Y_{j-r} \\
Z_{i+j-r} \ar[r]
	& Z_{i+j-r} \cup CZ_{i+j-2r} \rlap{\,.}
}
\]
Here 
\[
f'' \: S^p \longto X_i/X_{i-r} \och g'' \: S^q \longto Y_j/Y_{j-r}
\]
are the quotient maps induced by $(f', f)$ and $(g', g)$, respectively, and we write 
\[
\Phi = \phi_W \cup C\phi_{i-r,j-r}
\] 
for brevity. If $X_\star$ and $Y_\star$ are filtrations, as we assume, then 
\[
X_{i-r} \wedge Y_{j-r} \longto W
\]
is a $h$-cofibration, so the collapse map 
\[
\Theta \: W \cup C(X_{i-r} \wedge Y_{j-r}) \longto X_{i-r} \wedge Y_j/Y_{j-r} \vee X_i/X_{i-r} \wedge Y_{j-r}
\]
is an equivalence.

The left hand side of Equation~\eqref{eq:leibniz} applied to
$[f',f] \otimes [g',g]$ is
\begin{equation} \label{eq:phia}
\begin{split}
\partial \phi_r([f',f] \otimes [g',g])
&= \partial [\phi_W (f \wedge g)', \phi_{i,j}(f \wedge g)] \\
& = [\ast , \phi_W(f \wedge g)'\pi]
\,.
\end{split}
\end{equation}
Under the isomorphism $\pi^G_{p+q-1}(Z_{i+j-2r} \to Z_{i+j-r}) \cong \pi^G_{p+q-1}(Z_{i+j-r} \cup CZ_{i+j-2r})$ this corresponds to the homotopy class of the composite map
\[
S^{p+q-1} \longto Z_{i+j-r} \longto Z_{i+j -r} \cup C Z_{i+j - 2r}
\]
in the diagram above. Equivalently, by the commutativity of the diagram, we can describe it as the homotopy class of the composite map
\begin{align*}
S^{p+q-1} &\longto S^{p+q-1} \cup C(S^{p-1} \wedge S^{q-1})
	\\ &\xrightarrow{(f \wedge g)' \cup C(f' \wedge g')}
	W \cup C(X_{i-r} \wedge Y_{j-r}) \\ &\xrightarrow{\Phi} Z_{i+j-r} \cup C(Z_{i+j-2r}) \,.
\end{align*}
Alternatively, we can describe it as the image $\Phi_*([a])$ of the homotopy class $[a]$ of the composite map
\begin{align*}
a \: S^{p+q-1} &\longto S^{p+q-1} \cup C(S^{p-1} \wedge S^{q-1})
	\\ &\xrightarrow{(f \wedge g)' \cup C(f' \wedge g')}
	W \cup C(X_{i-r} \wedge Y_{j-r}) 
\end{align*}
under the homomorphism
\[
\Phi_* \:
\pi^G_{p+q-1}(W \cup C(X_{i-r} \wedge Y_{j-r}))
	\longto \pi^G_{p+q-1}(Z_{i+j-r} \cup CZ_{i+j-2r}) \,.
\]
We shall confirm that Equation~\eqref{eq:leibniz} holds by writing it in the form
\[
\Phi_*([a]) = \Phi_*([b]) + (-1)^p \Phi_*([c])
\]
for some specific classes $[b]$ and $[c]$, to be defined later, and showing that 
\[
[a] = [b] + (-1)^p [c]
\] 
in $\pi^G_{p+q-1}(W \cup C(X_{i-r} \wedge Y_{j-r}))$. The latter identity will be confirmed by showing that 
\[
\Theta_*([a]) = \Theta_*([b]) + (-1)^p \Theta_*([c]) \,,
\]
where $\Theta_*$ is the isomorphism
\[
\Theta_* \: \pi^G_{p+q-1}(W \cup C(X_{i-r} \wedge Y_{j-r}))
\overset{(\cong)}\longto \pi^G_{p+q-1}(X_{i-r} \wedge Y_j/Y_{j-r}
	\vee X_i/X_{i-r} \wedge Y_{j-r})
\]
induced by the map $\Theta$ on homotopy. With this aim in mind, we note that~$\Theta_*([a])$ is the homotopy class of the composite
\begin{align*}
S^{p+q-1} &\longto S^{p+q-1} \cup C(S^{p-1} \wedge S^{q-1})
  \\ &\overset{\simeq}\longto (S^{p-1} \wedge S^q) \vee (S^p \wedge S^{q-1}) \\ &\xrightarrow{(f' \wedge g'') \vee (f'' \wedge g')}  X_{i-r} \wedge Y_j/Y_{j-r} \vee X_i/X_{i-r} \wedge Y_{j-r} \,,
\end{align*}
again by commutativity of the above diagram. Checking orientations in the boundary of $D^p \wedge D^q$, the composition of all but the last map in the displayed map has degree~$+1$ when projected to $S^{p-1} \wedge S^q$, and degree $(-1)^p$ when projected to~$S^p \wedge S^{q-1}$. Hence $\Theta_*([a])$ is the sum of the homotopy class of the composite
\begin{equation} \label{eq:psib}
\begin{split}
S^{p-1} \wedge S^q
	&\xrightarrow{f' \wedge g''} X_{i-r} \wedge Y_j/Y_{j-r}
	\\ &\xrightarrow{\In_1}
	X_{i-r} \wedge Y_j/Y_{j-r} \vee X_i/X_{i-r} \wedge Y_{j-r}
	\end{split}
\end{equation}
and $(-1)^p$ times the homotopy class of the composite
\begin{equation} \label{eq:psic}
\begin{split}
S^p \wedge S^{q-1}
	&\xrightarrow{f'' \wedge g'} X_i/X_{i-r} \wedge Y_{j-r} \\
	&\xrightarrow{\In_2}
	X_{i-r} \wedge Y_j/Y_{j-r} \vee X_i/X_{i-r} \wedge Y_{j-r} \,.
	\end{split}
\end{equation}

The first term of the right hand side of Equation~\eqref{eq:leibniz} applied to $[f',f] \otimes [g',g]$ is
\begin{equation} \label{eq:phib}
\phi_r(\partial\otimes1)([f',f] \otimes [g',g])
	= \phi_r([*, f'\pi] \otimes [g',g])
\,.
\end{equation}
Unravelling the definition of $\phi_r$, see the discussion in Section~\ref{subsec:pairingsCEI} for more details, we form the commutative diagram
\[
\xymatrix{
S^{p-2} \wedge S^{q-1} \ar[rr] \ar[dr] \ar[dd]
	&& D^{p-1} \wedge S^{q-1} \ar[dr] \ar[dd]|\hole^(0.75){f'\pi \wedge g'} \\
& S^{p-2} \wedge D^q \ar[rr] \ar[dd]
	&& D^{p-1} \wedge D^q \ar[dd]^-{f'\pi \wedge g} \\
{*} \ar[rr]|\hole \ar[dr] \ar[dd]
	&& X_{i-r} \wedge Y_{j-r} \ar[dr] \ar[dd]|\hole^(0.75){=} \\
& {*} \ar[rr] \ar[dd]
	&& X_{i-r} \wedge Y_j \ar[dd]^-{=} \\
X_{i-2r} \wedge Y_{j-r} \ar[rr]|\hole \ar[dr] \ar[dd]_-{\phi_{i-2r,j-r}}
	&& X_{i-r} \wedge Y_{j-r}
		\ar[dr] \ar[dd]|\hole^(0.75){\phi_{i-r,j-r}} \\
& X_{i-2r} \wedge Y_j \ar[rr] \ar[dr]_-{\phi_{i-2r,j}}
	&& X_{i-r} \wedge Y_j \ar[d]^-{\phi_{i-r,j}} \\
Z_{i+j-3r} \ar[rr]
	&& Z_{i+j-2r} \ar[r]
	& Z_{i+j-r} \,.
}
\]
We also introduce the pushouts
\[
S^{p+q-2} = S^{p-2} \wedge D^q \cup_{S^{p-2} \wedge S^{q-1}}
	D^{p-1} \wedge S^{q-1}
\]
and
\[
U = X_{i-2r} \wedge Y_j \cup_{X_{i-2r} \wedge Y_{j-r}}
	X_{i-r} \wedge Y_{j-r}
\]
mapping to $D^{p-1} \wedge D^q \cong D^{p+q-1}$
and $X_{i-r} \wedge Y_j$, respectively. This leads to the commutative diagram
\[
\xymatrix{
S^{p+q-2} \ar[r] \ar[d]_-{{*} \cup f'\pi \wedge g'}
	& D^{p+q-1} \ar[d]^-{f'\pi \wedge g} \\
X_{i-r} \wedge Y_{j-r} \ar[r] \ar[d] \ar@(l,l)[ddd]_-{\phi_{i-r,j-r}}
	& X_{i-r} \wedge Y_j \ar[d]^-{=} \\
U \ar[r] \ar[dd]^-{\phi_U}
	& X_{i-r} \wedge Y_j \ar[dd]^-{\phi_{i-r,j}} \ar[dr] \\
& & W \ar[dl]^-{\phi_W} \\
Z_{i+j-2r} \ar[r]
	& Z_{i+j-r} \,. &
}
\]
The class in $\pi^G_{p+q-1}(Z_{i+j-2r} \to Z_{i+j-r})$ described in Equation~\eqref{eq:phib} is represented visually as the big rectangle in the diagram, that is, by the pair of maps
\begin{align*}
\phi_{i-r,j-r} (* \cup f'\pi \wedge g') &\: S^{p+q-2} \to Z_{i+j-2r} \\ \phi_{i-r,j} (f'\pi \wedge g) &\: D^{p+q-1} \to Z_{i+j-r} \,.
\end{align*}
We can extend the diagram to the right, as follows,
\[
\xymatrix{
D^{p+q-1} \ar[r] \ar[d]_-{f'\pi \wedge g}
	& D^{p+q-1} \cup CS^{p+q-2} \ar[r]^-{\simeq}
		\ar[d]^-{f'\pi \wedge g \cup C(* \cup f'\pi \wedge g')}
	& S^{p+q-1} \ar[d]^-{f' \wedge g''} \\
X_{i-r} \wedge Y_j \ar[r] \ar[d]
	& X_{i-r} \wedge (Y_j \cup CY_{j-r})
		\ar[r]^-{(\simeq)} \ar[d]
	& X_{i-r} \wedge Y_j/Y_{j-r} \ar[d]^-{\In_1} \\
W \ar[r] \ar[d]_-{\phi_W}
	& W \cup C(X_{i-r} \wedge Y_{j-r})
		\ar[r]^-{(\simeq)}_-{\Theta} \ar[d]^-{\Phi}
	& X_{i-r} \wedge Y_j/Y_{j-r} \vee X_i/X_{i-r} \wedge Y_{j-r} \\
Z_{i+j-r} \ar[r] & Z_{i+j-r} \cup CZ_{i+j-2r} 
}
\]
where the maps marked $(\simeq)$ are equivalences by our assumption that $X_{i-r} \to X_i$ and $Y_{j-r} \to Y_j$ are $h$-cofibrations. Under the isomorphism $\pi_{p+q-1}^G(Z_{i+j-2r} \to Z_{i+j-r}) \cong \pi^G_{p+q-1}(Z_{i+j-r} \cup CZ_{i+j-2r})$ the class described in Equation~\eqref{eq:phib} is given by the composite
\begin{align*}
S^{p+q-1} &\overset{\simeq}\longto D^{p+q-1} \cup CS^{p+q-2}
	\\ &\longto X_{i-r} \wedge (Y_j \cup CY_{j-r}) \\
	&\longto W \cup C(X_{i-r} \wedge Y_{j-r})
	\\ &\overset{\Phi}\longto
		Z_{i+j-r} \cup CZ_{i+j-2r} \,,
\end{align*}
where the first map is a homotopy inverse to the collapse map.  This is the image~$\Phi_*([b])$ of the homotopy class~$[b]$ of the
composite map
\[
b: S^{p+q-1} \overset{\simeq}\longto D^{p+q-1} \cup CS^{p+q-2}
	\longto W \cup C(X_{i-r} \wedge Y_{j-r}) \,.
\]
Since the composite
\[
S^{p+q-1} \overset{\simeq}\longto D^{p+q-1} \cup CS^{p+q-2}
	\overset{\simeq}\longto S^{p+q-1} \cong S^{p-1} \wedge S^q
\]
is homotopic to the identity, $\Theta_*([b])$ is the homotopy class of the map~\eqref{eq:psib}. That is, it is the image of~$[f' \wedge g'']$ under the inclusion $(\In_1)_*$.

The second term of the right hand side of~\eqref{eq:leibniz} applied to $[f',f] \otimes [g',g]$ is
\begin{equation} \label{eq:phic}
\phi_r(1\otimes\partial)([f',f] \otimes [g',g])
= (-1)^p \phi_r([f',f] \otimes [*, g'\pi]) \,.
\end{equation}
By a similar analysis as for the first term of the sum, the class
$\phi_r([f',f] \otimes [*, g'\pi])$ in $\pi^G_{p+q-1}(Z_{i+j-2r} \to Z_{i+j-r})$ is represented by a pair of maps
\begin{align*}
\phi_{i-r,j-r}(f' \wedge g'\pi \cup *) &\: S^{p+q-2} \longto Z_{i+j-2r} \,, \\
\phi_{i,j-r}(f \wedge g'\pi) &\: D^{p+q-1} \longto Z_{i+j-r} \,,
\end{align*}
where $S^{p+q-2}$ is the boundary of $D^{p+q-1} \cong D^p \wedge D^{q-1}$. The corresponding class in
$\pi^G_{p+q-1}(Z_{i+j-r} \cup CZ_{i+j-2r})$ is the image $\Phi_*([c])$ under $\Phi_*$ of the homotopy class~$[c]$ of the composite map
\[
c \: S^{p+q-1} \overset{\simeq}\longto D^{p+q-1} \cup CS^{p+q-2}
	\longto W \cup C(X_{i-r} \wedge Y_{j-r}) \,.
\]
The class $\Theta_*([c])$ is then the homotopy class of the map~\eqref{eq:psic}. That is, it is the image of $[f'' \wedge g']$ under $(\In_2)_*$.

Summarising, we have now defined classes $[a]$, $[b]$, and $[c]$ in $\pi^G_{p+q-1}(W \cup C(X_{i-r} \wedge Y_{j-r}))$ satisfying
\[
\Theta_*([a]) = \Theta_*([b]) + (-1)^p
\Theta_*([c]) \,.
\]
Since $\Theta_*$ is an isomorphism, we deduce that 
\[
[a] = [b] + (-1)^p [c] \och \Phi_*([a]) = \Phi_*([b]) + (-1)^p \Phi_*([c]) \,.
\]
Since $\Phi_*([a])$, $\Phi_*([b])$ and $(-1)^p \Phi_*([c])$ are the three parts in Equation~\eqref{eq:leibniz} evaluated at $[f',f] \otimes [g',g]$, and $[f',f]$ and $[g',g]$ were arbitrarily chosen, it follows that Equation~\eqref{eq:leibniz}
holds whenever $X_\star$ and $Y_\star$ are filtrations.
\end{proof}

We now extend the result above to all pairings of sequences.

\begin{proposition} \label{prop:Leibpairseq}
If $\phi \: (X_\star, Y_\star) \to Z_\star$ is a pairing of
sequences of orthogonal~$G$-spectra, then
\[
\partial \phi_r = \phi_r (\partial \otimes 1) + \phi_r (1 \otimes \partial)
\]
as homomorphisms
\[
H(X_\star)(i-r,i) \otimes H(Y_\star)(j-r,j)
	\longto H(Z_\star)(i+j-2r,i+j-r)
\]
for all integers $i,j,r$ with $r\ge1$.
\end{proposition}

\begin{proof}
Let $T(\phi) \: (T_\star(X), T_\star(Y)) \to T_\star(Z)$
be the pairing of filtrations defined as in the proof of Lemma~\ref{lem:Telphi}. The equivalence $\epsilon \: T_\star(X) \to X_\star$ and its analogues for~$Y_\star$ and~$Z_\star$ are compatible with the pairings. Hence we have a commutative diagram with vertical isomorphisms
\[
\xymatrix{
H(T_\star(X))(i-r,i)
	\ar[r]^-{\partial} \ar[d]_-{\epsilon}^-{\cong}
	& H(T_\star(X))(i-2r,i-r) \ar[d]^-{\epsilon}_-{\cong} \\
H(X_\star)(i-r,i)
	\ar[r]^-{\partial}
	& H(X_\star)(i-2r,i-r) \,,
}
\]
together with its analogues for~$Y_\star$ and~$Z_\star$, and
\[
\xymatrix{
H(T_\star(X))(i-r,i) \otimes H(T_\star(Y))(j-r,j)
	\ar[r]^-{T(\phi)_r}
	\ar[d]_-{\epsilon \otimes \epsilon}^-{\cong}
        & H(T_\star(Z))(i+j-r,i+j)
	\ar[d]^-{\epsilon}_-{\cong}
\\
H(X_\star)(i-r,i) \otimes H(Y_\star)(j-r,j) \ar[r]^-{\phi_r}
        & H(Z_\star)(i+j-r,i+j) \rlap{\,,}
}
\]
for all $r\ge1$. By Proposition~\ref{prop:Leibpairfilt} applied to the pairing of filtrations $T(\phi)$ we know that
\[
\partial T(\phi)_r
  = T(\phi)_r (\partial \otimes 1) + T(\phi)_r (1 \otimes \partial)
\]
as homomorphisms
\[
H(T_\star(X))(i-r,i) \otimes H(T_\star(Y))(j-r,j)
	\longto H(T_\star(Z))(i+j-2r,i+j-r) \,.
\]
In view of the vertical isomorphisms $\epsilon$, this implies that
$\partial \phi_r = \phi_r (\partial \otimes 1) + \phi_r (1 \otimes \partial)$ as homomorphisms $H(X_\star)(i-r,i) \otimes H(Y_\star)(j-r,j)
\to H(Z_\star)(i+j-2r,i+j-r)$.
\end{proof}

This finishes the goal we set out for ourselves at the start of Section~\ref{subsec:pairingsCEI}, namely to prove that a pairing of sequences gives rise to a pairing of the resulting Cartan--Eilenberg systems\index{pairing!of sequences}\index{pairing!of Cartan--Eilenberg systems}. Let us phrase this conclusion in a theorem, so that we can refer back to it when needed.

\begin{theorem}
\label{thm:mainmultI}
A pairing $\phi \: (X_\star , Y_\star) \to Z_\star$ of sequences gives rise to a pairing $\phi \: (H(X_\star),H(Y_\star)) \to H(Z_\star)$ of the associated Cartan--Eilenberg systems, in the sense of Definition~\ref{def:pairextCEsys}.
\begin{proof}
The proof of (SPP I) and (SPP III) is the content of Section~\ref{subsec:pairingsCEI} and the proof of (SPP II) is the content of the present section.
\end{proof}
\end{theorem}

This directly gives us the following consequence for the associated spectral sequences\index{pairing!of spectral sequences}.

\begin{theorem}
\label{thm:mainmultII}
A pairing $\phi \: (X_\star , Y_{\star}) \to Z_{\star}$ of sequences of orthogonal $G$-spectra gives rise to a pairing $\phi \: (E^*(X_\star),E^*(Y_\star)) \to E^*(Z_\star)$, in the sense of Definition~\ref{def:sseqpairing}. Explicitly, we have access to a collection of homomorphisms
\[
\phi^r \: E^r(X_\star) \otimes E^r(Y_\star) \longto E^r(Z_\star)
\]
for all $r \geq 1$, such that:
\begin{enumerate}
\item The Leibniz rule\index{Leibniz rule}
\[
d^r \phi^r = \phi^r(d^r \otimes 1) + \phi^r(1 \otimes d^r)
\]
holds as an equality of homomorphisms
$E^r_i(X_\star) \otimes E^r_j(Y_\star) \longto E^r_{i+j-r}(Z_\star)$
for all $i,j \in \bZ$ and $r\ge1$.
\item The diagram 
\[
\begin{tikzcd}
E^{r+1}(X_\star) \otimes E^{r+1}(Y_\star) \arrow[r,"\phi^{r+1}"] \arrow[d] & E^{r+1}(Z_\star) \arrow[d,"\cong"] \\
H(E^{r}(X_\star) \otimes E^{r}(Y_\star)) \arrow[r,"H(\phi^r)"] & H(E^{r}(Z_\star))
\end{tikzcd}
\]
commutes for all $r \geq 1$.
\end{enumerate} 

Moreover, the induced pairing $\phi_*$ on filtered abutments is compatible with the pairing $\phi^\infty$ of $E^\infty$-pages in the sense of Proposition~\ref{prop:multabutment}. Explicitly, the diagram
\[
\xymatrix{
\ds \frac{F_i A_\infty(X_\star)}{F_{i-1} A_\infty(X_\star)}
\otimes \frac{F_j A_\infty(Y_\star)}{F_{j-1} A_\infty(Y_\star)}
\ar[r]^-{\bar\phi_*} \ar[d]_-{\beta \otimes \beta}
& \ds \frac{F_{i+j} A_\infty(Z_\star)}{F_{i+j-1} A_\infty(Z_\star)}
	\ar@{ >->}[d]^-{\beta} \\
E^\infty_i(X_\star) \otimes E^\infty_j(Y_\star) \ar[r]^-{\phi^\infty}
& E^\infty_{i+j}(Z_\star)
}
\]
commutes, for all $i,j \in \bZ$. Here the abutments are given as 
\begin{align*}
A_\infty(X_\star) &\cong \pi^G_* \Tel(X_\star) \\
A_\infty(Y_\star) &\cong \pi^G_* \Tel(Y_\star) \\
A_\infty(Z_\star) &\cong \pi^G_* \Tel(Z_\star) \,,
\end{align*}
and are filtered by the images
\begin{align*}
F_i A_\infty(X_\star) &= \im(\pi^G_*(X_i) \longto A_\infty(X_\star)) \\
F_j A_\infty(Y_\star) &= \im(\pi^G_*(Y_j) \longto A_\infty(Y_\star)) \\
F_k A_\infty(Z_\star) &= \im(\pi^G_*(Z_k) \longto A_\infty(Z_\star)) \,,
\end{align*}
respectively.

\begin{proof}
This follows from combining Theorem~\ref{thm:mainmultI} with Theorem~\ref{thm:Dou59bIIA} and Proposition~\ref{prop:multabutment}.
\end{proof}
\end{theorem}

\begin{corollary}
\label{cor:main}
If $(X_\star,\phi)$ is a multiplicative sequence of orthogonal $G$-spectra, then the associated spectral sequence $(E(X_\star),d)$ is multiplicative with multiplicative abutment.\index{sequence!multiplicative}\index{spectral sequence!multiplicative}\index{abutment!multiplicative} 
\end{corollary}

\section{The convolution product}
\label{subsec:convprod}

Given two sequences $X_\star$ and $Y_\star$ of orthogonal~$G$-spectra there is an initial pairing $\iota \: (X_\star, Y_\star) \to Z_\star$, where the sequence $Z_\star$ is given at each level by
\[
Z_k = \colim_{i+j \le k} X_i \wedge Y_j \,,
\]
with the canonical $G$-maps $Z_{k-1} \to Z_k$ between them. We call this sequence $Z_\star$ the \emph{(Day) convolution
product}\index{Day convolution}\index{convolution product} of $X_\star$ and $Y_\star$, and write 
\[
Z_\star = (X \wedge Y)_\star \,.
\]
The universal pairing $\iota \: (X_\star, Y_\star) \to
(X \wedge Y)_\star$ has components
\[
\iota_{i,j} \: X_i \wedge Y_j \longto (X \wedge Y)_{i+j} \,,
\]
each given by a structure map to the colimit. Per the discussion of Section~\ref{subsec:pairingsCEI}, the universal pairing $\iota \: (X_\star,Y_\star) \to (X \wedge Y)_\star$ induces homomorphisms
\[
\iota_r \: \pi^G_p(X_{i-r} \to X_i) \otimes \pi^G_q(Y_{j-r} \to Y_j)
	\longto \pi^G_{p+q}((X \wedge Y)_{i+j-r} \to (X \wedge Y)_{i+j})
\]
for $r\ge1$, and Theorem~\ref{thm:mainmultII} shows that the pairing~$\iota$ extends to a pairing 
\[
\iota \: E^*(X_\star) \otimes E^*(Y_\star) \to E^*((X \wedge Y)_\star)
\]
of spectral sequences, in such a way that the induced pairing on filtered abutments
\[
\bar\iota_* \:
\frac{F_i \pi^G_* \Tel(X_\star)}{F_{i-1} \pi^G_* \Tel(X_\star)}
\otimes
\frac{F_j \pi^G_* \Tel(Y_\star)}{F_{j-1} \pi^G_* \Tel(Y_\star)}
\longto
\frac{F_{i+j} \pi^G_* \Tel((X \wedge Y)_\star)}
	{F_{i+j-1} \pi^G_* \Tel((X \wedge Y)_\star)}
\]
is compatible with the induced pairing on $E^\infty$-pages. 

\begin{remark}
The colimit defining $Z_k$ can equally well be calculated over the cofinal subcategory of pairs $(i,j) \in \bZ^2$ with $k-1 \le i+j \le k$, i.e., as the colimit of the
zigzag diagram:
\[
\xymatrix{
\dots \ar[r] & X_{i-1} \wedge Y_{k-i+1} \\
& X_{i-1} \wedge Y_{k-i} \ar[r] \ar[u] & X_i \wedge Y_{k-i} \\
& & X_i \wedge Y_{k-i-1} \ar[r] \ar[u] & \dots
}
\]
\end{remark}

If $(X_\star, \phi)$ and $(Y_\star, \psi)$ are multiplicative sequences of orthogonal $G$-spectra, then the convolution product $((X \wedge Y)_\star, \phi \wedge \psi)$ is a multiplicative sequence\index{sequence!multiplicative} as well. Here, the component
\[
(\phi \wedge \psi)_{i,j} \: (X \wedge Y)_i \wedge (X \wedge Y)_j
	\longto (X \wedge Y)_{i+j}
\]
is defined as the colimit over $i_1+i_2 \le i$ and $j_1+j_2 \le j$ of the composite maps
\begin{align*}
X_{i_1} \wedge Y_{i_2} \wedge X_{j_1} \wedge Y_{j_2}
	&\xrightarrow{1 \wedge \tau \wedge 1}
X_{i_1} \wedge X_{j_1} \wedge Y_{i_2} \wedge Y_{j_2}
\xrightarrow{\phi_{i_1,j_1} \wedge \psi_{i_2,j_2}}
	X_{i_1+j_1} \wedge Y_{i_2+j_2} \\
&\xrightarrow{\iota_{i_1+j_1,i_2+j_2}} (X \wedge Y)_{i_1+j_1+i_2+j_2}
\longto (X \wedge Y)_{i+j} \,.
\end{align*}

\begin{lemma} \label{lem:iota1xyxymult}
If $(X_\star, \phi)$ and $(Y_\star, \psi)$ are multiplicative
sequences, then the homomorphism
$\iota^1 \: E^1(X_\star) \otimes E^1(Y_\star)
\to E^1((X \wedge Y)_\star)$ is multiplicative, in the sense
that the diagram
\[
\begin{tikzcd}
E^1(X_\star) \otimes E^1(Y_\star) \otimes E^1(X_\star) \otimes E^1(Y_\star)
	\arrow[r,"\iota^1 \otimes \iota^1"] \arrow[d,"1 \otimes \tau \otimes 1"',"\cong"]
& E^1((X \wedge Y)_\star) \otimes E^1((X \wedge Y)_\star)
	\arrow[dd,"(\phi\wedge\psi)^1"] \\
E^1(X_\star) \otimes E^1(X_\star) \otimes E^1(Y_\star) \otimes E^1(Y_\star)
	\ar[d,"\phi^1 \otimes \psi^1"'] \\
E^1(X_\star) \otimes E^1(Y_\star) \arrow[r,"\iota^1"]
	& E^1((X \wedge Y)_\star)
\end{tikzcd}
\]
commutes.
\end{lemma}

\begin{proof}
Let us write $\theta = \phi \wedge \psi$ for brevity. The diagrams
\[
\xymatrix@C+2pc{
X_{i_1} \wedge Y_{i_2} \wedge X_{j_1} \wedge Y_{j_2}
	\ar[r]^-{\iota_{i_1,i_2} \wedge \iota_{j_1,j_2}}
	\ar[d]_-{1 \wedge \tau \wedge 1}^-{\cong}
& (X \wedge Y)_{i_1+i_2} \wedge (X \wedge Y)_{j_1+j_2}
	\ar[dd]^-{\theta_{i_1+i_2,j_1+j_2}} \\
X_{i_1} \wedge X_{j_1} \wedge Y_{i_2} \wedge Y_{j_2}
	\ar[d]_-{\phi_{i_1,j_1} \wedge \psi_{i_2,j_2}} \\
X_{i_1+j_1} \wedge Y_{i_2+j_2}
	\ar[r]^-{\iota_{i_1+j_1,i_2+j_2}}
& (X \wedge Y)_{i_1+j_1+i_2+j_2}
}
\]
commute, and are compatible, for all $i_1$, $i_2$, $j_1$ and~$j_2$. This implies that the composite homomorphism
\begin{align*}
H(&X_\star)(i_1-r,i_1) \otimes H(Y_\star)(i_2-r,i_2)
	\otimes H(X_\star)(j_1-r,j_1) \otimes H(Y_\star)(j_2-r,j_2) \\
&\xrightarrow{\iota_r \otimes \iota_r}
H((X \wedge Y)_\star)(i_1+i_2-r,i_1+i_2)
	\otimes H((X \wedge Y)_\star)(j_1+j_2-r,j_1+j_2) \\
&\xrightarrow{\theta_r}
H((X \wedge Y)_\star)(i_1+i_2+j_1+j_2-r,i_1+i_2+j_1+j_2)
\end{align*}
is equal to the composite homomorphism
\begin{align*}
&H(X_\star)(i_1-r,i_1) \otimes H(Y_\star)(i_2-r,i_2)
	\otimes H(X_\star)(j_1-r,j_1) \otimes H(Y_\star)(j_2-r,j_2) \\
&\xrightarrow{\cong}
H(X_\star)(i_1-r,i_1) \otimes H(X_\star)(j_1-r,j_1)
	\otimes H(Y_\star)(i_2-r,i_2) \otimes H(Y_\star)(j_2-r,j_2) \\
&\xrightarrow{\phi_r \otimes \psi_r}
H(X_\star)(i_1+j_1-r,i_1+j_1)
	\otimes H(Y_\star)(i_2+j_2-r,i_2+j_2) \\
&\xrightarrow{\iota_r}
H((X \wedge Y)_\star)(i_1+j_1+i_2+j_2-r,i_1+j_1+i_2+j_2)
\end{align*}
for each $r\ge1$, where the homomorphisms $\iota_r$, $\phi_r$, $\psi_r$ and $\theta_r$ are defined as in Section~\ref{subsec:pairingsCEI}. For $r=1$, this gives the claim of the lemma.
\end{proof}

Note that for general sequences $X_\star$ and $Y_\star$ we typically have no homotopical control of their convolution product. However, if both $X_\star$ and $Y_\star$ are filtrations, then we can view each $X_i \wedge Y_j$ as a subspectrum of
\[
\colim_i X_i \wedge \colim_j Y_j = \bigcup_i X_i \wedge \bigcup_j Y_j
\]
and their colimit for $i+j \le k$ can be formed as the union
\[
(X \wedge Y)_k = \bigcup_{i+j=k} X_i \wedge Y_j \,.
\]

\begin{proposition} \label{prop:convoffiltisfilt}
If the sequences $X_\star$ and $Y_\star$ are filtrations, then their convolution product $(X \wedge Y)_\star$ is a filtration.
\end{proposition}

\begin{proof}
We must show that each map
\[
(X \wedge Y)_{k-1} \longto (X \wedge Y)_k
\]
is a strong $h$-cofibration.  This is the colimit of a sequence of maps, each of which is the cobase change of a pushout-product map
\[
X_{i-1} \wedge Y_j \cup X_i \wedge Y_{j-1}
	\longto X_i \wedge Y_j
\]
with $i+j=k$, where the pushout is formed over $X_{i-1} \wedge Y_{j-1}$. By assumption $X_{i-1} \to X_i$ and $Y_{j-1} \to Y_j$ are strong $h$-cofibrations, so the conclusion follows immediately from Theorem~\ref{thm:pushprodstrongcof}.
\end{proof}

In the special case when two arbitrary sequences $X_\star$ and $Y_\star$
are first replaced with equivalent filtrations~$T_\star(X)$ and
$T_\star(Y)$, we can give the following alternative, more explicit, argument for
why the resulting convolution product is always a filtration.

\begin{lemma} \label{lem:telxtelyfiltration}
For any two sequences~$X_\star$ and~$Y_\star$ of orthogonal $G$-spectra,
the convolution product $(T(X) \wedge T(Y))_\star$ is a filtration.
\end{lemma}

\begin{proof}
In degree~$k$,
\[
(T(X) \wedge T(Y))_k = \bigcup_{i+j=k} T_i(X) \wedge T_j(Y) \,.
\]
This is the subspectrum of $\Tel(X_\star) \wedge \Tel(Y_\star)$ 
with telescope coordinates $x$ and~$y$ satisfying $\lceil x \rceil
+ \lceil y \rceil \le k$.  Here $\lceil x \rceil$ denotes the
least integer~$i$ with $x \le i$. The inclusion
$(T(X) \wedge T(Y))_{k-1} \to (T(X) \wedge T(Y))_k$ is
then the composite of a sequence of cobase changes of maps of
the form
\[
i_0 \: A \longto B \cup_A A \wedge I_+ \,,
\]
with $A$ the double mapping cylinder of the diagram
\[
X_{i-1} \wedge Y_j \longleftarrow X_{i-1} \wedge Y_{j-1}
	\longto X_i \wedge Y_{j-1}
\]
and $B = X_i \wedge Y_j$, for $i+j=k$. Since each such map~$i_0$ is a strong $h$-cofibration, so is the structure map in $(T(X) \wedge T(Y))_\star$, as claimed.
\end{proof}

As a consequence of Proposition~\ref{prop:convoffiltisfilt}, we can write the first page of the spectral sequence associated to the convolution product of two filtrations $X_\star$ and $Y_\star$ as 
\begin{align*}
E^1_k((X \wedge Y)_\star) 
	&= \pi^G_*((X \wedge Y)_{k-1} \to (X \wedge Y)_k) \\
	&\cong \pi^G_*((X \wedge Y)_k / (X \wedge Y)_{k-1}) \\
	&\cong \bigoplus_{i+j=k} \pi^G_*(X_i/X_{i-1} \wedge Y_j/Y_{j-1})
\end{align*}
since
\[
\frac{(X \wedge Y)_k}{(X \wedge Y)_{k-1}} \cong
	\bigvee_{i+j=k} X_i/X_{i-1} \wedge Y_j/Y_{j-1} \,.
\]
Furthermore, the diagram
\begin{equation} \label{eq:E1d1forxyxy}
\begin{tikzcd}
E^1(X_\star) \otimes E^1(Y_\star) \arrow[r,"\iota^1"]
	\arrow[d,"d^1 \otimes 1 + 1 \otimes d^1"']
	& E^1((X \wedge Y)_\star) \arrow[d,"d^1"] \\
E^1(X_\star) \otimes E^1(Y_\star) \arrow[r,"\iota^1"]
	& E^1((X \wedge Y)_\star)
\end{tikzcd}
\end{equation}
commutes. To proceed we usually need more explicit control of the $d^1$-differential for $(X \wedge Y)_\star$, e.g., by use of~\eqref{eq:E1d1forxyxy} in situations where $\iota^1$ is surjective.

Suppose now that $(X_\star, \phi)$ and $(Y_\star, \psi)$ are multiplicative sequences, and also assume that the former is a filtration. This will be the situation when we filter the~$G$-Tate construction in Chapter~\ref{sec:Tatespseq}. By Corollary~\ref{cor:telxmult}, the telescopic replacement~$(T_\star(Y), T(\psi))$ is a multiplicative filtration, and by Proposition~\ref{prop:convoffiltisfilt} the convolution product $(X \wedge T(Y))_\star, \phi \wedge T(\psi))$ is then also a multiplicative filtration. Lemma~\ref{lem:iota1xyxymult} shows that $\iota^1 \: E^1(X_\star) \otimes E^1(Y_\star) \to E^1((X \wedge Y)_\star)$ is multiplicative, in the sense that the diagram
\[
\begin{tikzcd}
\txt{ $E^1(X_\star) \otimes E^1(T_\star(Y))$ \\ $\otimes$ \\ $E^1(X_\star) \otimes E^1(T_\star(Y))$ }
	\arrow[r,"\iota^1 \otimes \iota^1"] \arrow[d,"1 \otimes \tau \otimes 1"',"\cong"]
& E^1((X \wedge T(Y))_\star) \otimes E^1((X \wedge T(Y))_\star)
	\arrow[dd,"(\phi\wedge T(\psi))^1"] \\
\txt{ $E^1(X_\star) \otimes E^1(X_\star)$ \\ $\otimes$ \\ $E^1(T_\star(Y)) \otimes E^1(T_\star(Y))$ }
	\ar[d,"\phi^1 \otimes T(\psi)^1"'] \\
E^1(X_\star) \otimes E^1(T_\star(Y)) \arrow[r,"\iota^1"]
	& E^1((X \wedge T(Y))_\star)
\end{tikzcd}
\]
commutes for all $i,j \in \bZ$. In situations where $\iota^1$
is surjective, this gives us algebraic control of the product on $E^1((X
\wedge T(Y))_\star)$ in terms of the products on $E^1(X_\star)$
and $E^1(T_\star(Y)) \cong E^1(Y_\star)$.

\begin{remark}
The results of this chapter readily generalise to the case of sequences of $R$-modules in orthogonal $G$-spectra, for any fixed commutative orthogonal ring spectrum~$R$. Letting $X_\star$ denote a sequence 
\[
\cdots \longto X_{i-1} \longto X_i \longto X_{i+1} \longto \cdots
\]
of $R$-module $G$-spectra and $R$-module $G$-maps, the telescope $\Tel(X_\star)$ is an $R$-module $G$-spectrum, and the following all live in the category of $R_*$-modules: the Cartan--Eilenberg system $(H, \partial)$, the exact couple~$(A, E^1)$, the filtered abutment $A_\infty \cong \pi^G_* \Tel(X_\star)$, and the spectral sequence $(E^r, d^r)$ . The telescope filtration and equivalence
\[
\epsilon \: T_\star(X) \longto X_\star
\]
also live in the category of $R$-modules.

Given sequences $X_\star$, $Y_\star$ and~$Z_\star$ of $R$-modules in orthogonal $G$-spectra, an $R$-bilinear pairing
\[
\phi \: (X_\star, Y_\star) \longto Z_\star
\]
consists of compatible $R$-linear $G$-maps 
\[
\phi \: X_i \wedge_R Y_j \longto Z_{i+j} \,,
\]
where the usual smash product has been replaced with the smash product over~$R$. Such pairings induce~$R_*$-module
homomorphisms 
\[
\phi_r \: H(X_\star)(i-r,i) \otimes_{R_*} H(Y_\star)(j-r,j)
\longto H(Z_\star)(i+j-r,i+j) \,,
\]
where the usual tensor product has been replaced with the tensor product over~$R_*$. The Leibniz rule holds for $\phi_r$, so that $\phi$ induces an $R_*$-linear pairing of $R_*$-module spectral sequences 
\[
\phi^r \: E^r(X_\star) \otimes_{R_*}
E^r(Y_\star) \longto E^r(Z_\star) \,.
\]
The corresponding $R_*$-linear pairings $\bar\phi_*$ and $\phi^\infty$ of the filtration subquotients and~$E^\infty$-pages are compatible under the $R_*$-module monomorphism
\[
\beta \: \frac{F_i A_\infty}{F_{i-1} A_\infty} \longto E^\infty_i \,.
\]
The universal $R$-bilinear pairing $\iota \: (X_\star, Y_\star) \to Z_\star$ is given by the $R$-module convolution product $Z_\star = (X \wedge_R Y)_\star$, with 
\[
(X \wedge_R Y)_k = \colim_{i+j\le k} X_i
\wedge_R Y_j \,.
\]
If~$X_\star$ and~$Y_\star$ are $R$-module filtrations,
then $(X \wedge_R Y)_\star$ is an $R$-module filtration. For general
$R$-module sequences~$X_\star$ and~$Y_\star$ the diagram~\eqref{eq:E1d1forxyxy}
commutes after replacing~$\otimes$ and~$\wedge$ by $\otimes_{R_*}$ and
$\wedge_R$, respectively.  Finally, if~$(X_\star, \phi)$ and~$(Y_\star,
\psi)$ are multiplicative $R$-module sequences then $\iota^1 \:
E^1(X_\star) \otimes_{R_*} E^1(Y_\star) \to E^1((X \wedge_R Y)_\star)$
will be multiplicative.  This depends on the existence of $R$-module maps
\[
(23) = 1 \wedge \tau \wedge 1
	\: X_{i_1} \wedge_R Y_{i_2} \wedge_R X_{j_1} \wedge_R Y_{j_2}
\longto X_{i_1} \wedge_R X_{j_1} \wedge_R Y_{i_2} \wedge_R Y_{j_2}
\]
that are strictly compatible for varying $i_1$, $j_1$, $i_2$ and~$j_2$.
This is a point where we use the assumption that $R$ is strictly commutative, not just homotopy commutative, as an orthogonal ring spectrum.
\end{remark}

%
%
%

\chapter{The $G$-Homotopy Fixed Point Spectral Sequence}
\index{homotopy fixed point spectral sequence|(}

Given an $R$-module $X$ of orthogonal $G$-spectra we can define the $G$-homotopy fixed points as the genuine fixed points 
\[
X^{hG} = F(EG_+,X)^G = F_R(R \wedge EG_+,X)^G \,.
\]
In this chapter we construct a spectral sequence
\[
E^r_{*,*}(X) \Longrightarrow \pi_*(X^{hG})
\]
with abutment being the homotopy groups of the $G$-homotopy fixed points of $X$, for any compact Lie group $G$. This spectral sequence will be induced by the filtration, covered in Section 5.1, on the free and contractible $G$-space $EG$ coming from the simplicial bar construction. In Section 5.2, we show that this spectral sequence is multiplicative with multiplicative abutment. See Theorem~\ref{thm:pairingbarmur}. Under the assumption that $R[G]_*$ is finitely generated and projective over $R_*$ we can algebraically identify the $E^2$-page of the $G$-homotopy fixed point spectral sequence as 
\[
E^2_{*,*}(X) \cong \Ext^{-*}_{R[G]_*}(R_*,\pi_*(X)) \,,
\]
with the multiplicative structure being identified with the cup product on the right-hand side. See Theorem~\ref{thm:E2MXisExtRG}. Lastly, in Section 5.4 we discuss the relationship between the simplicial skeletal filtration on $E\bT$ and the often-used filtration coming from odd-dimensional spheres.

\section{The filtered $G$-space $EG$}
\label{subsec:filteredEG}

As always, $G$ is a compact Lie group. We let
\[
EG = B(\ast, G, G)
\]
be the free and contractible (right) $G$-space obtained by taking the geometric realization of the simplicial space
\[
[q] \mapsto B_q(*, G,G) = G^q \times G \,,
\]
with the usual face and degeneracy maps\index{bar construction} \cite{May75}*{\S7}. There is a simplicial contraction of $B_\bullet(*, G, G)$, so $EG$ is indeed contractible~\cite{May72}*{Prop.~9.8}. We let~$F_i EG$ be the image of the structure map $\Delta^i \times B_i(*, G, G) \to EG$ to the geometric realization, yielding the following exhaustive filtration \cite{May72}*{Def.~11.1}:
\[
\emptyset = F_{-1} EG \subset F_0 EG \subset
	\dots \subset F_{i-1} EG \subset F_i EG \subset \dots
	\subset EG \,.
\]
Here, the group $G$ acts freely from the right in each simplicial degree, hence also on each term in this filtration. The structure map induces a $G$-equivariant homeomorphism
\[
\Sigma^i(G^{\wedge i} \wedge G_+) \cong
\Delta^i/\partial\Delta^i\wedge G^{\wedge i} \wedge G_+
	\cong F_i EG / F_{i-1} EG
\]
for each $i\ge0$. Each smash power $G^{\wedge i} = G \wedge \dots \wedge
G$ (with $i$ copies of $G$) is formed with respect to the base point $e
\in G$ given by the unit element.

\begin{remark}
When $G$ is finite, $F_i EG$ gives the $i$-skeleton of a free $G$-CW structure on $EG$.  When $G = \bT = U(1)$ is the circle group, $F_i EG$ gives the $2i$- and $2i+1$-skeleta of a $G$-CW structure (with no odd-dimensional $G$-cells).
Similarly, when $G = \bU = Sp(1)$ is the $3$-sphere, $F_i EG$ gives the $4i$-, $4i+1$-, $4i+2$- and $4i+3$-skeleta of a $G$-CW structure. For other Lie groups the relationship is more complicated. Hence our filtration will agree with that used by Greenlees and~May in \cite{GM95}*{\S9} when $G$ is finite, be a doubly accelerated version when~$G =\bT$, and be a quadruply accelerated version when $G = \bU$. The two filtrations might be quite different for other compact Lie groups~$G$, though.
\end{remark}

We give the Cartesian product $EG \times EG$ the product filtration:
\[
F_k(EG \times EG) = \bigcup_{i+j=k} F_i EG \times F_j EG \,.
\]
Note that the diagonal $G$-map\index{diagonal map} $\Delta \: EG \to EG \times EG$, sending $x$ to $\Delta(x) = (x,x)$, is not filtration-preserving.  However, by
\cite{Seg68}*{Lem.~5.4} or~\cite{May72}*{Lem.~11.15} it is naturally
homotopic to a filtration-preserving map $D \: EG \to EG \times EG$,
which we call a \emph{diagonal approximation}\index{diagonal approximation} for~$EG$.  By inspection,
both $D$ and the natural homotopy $\Delta \simeq D$ are $G$-equivariant.
Subject to this condition, the precise choice of diagonal approximation
will not be important, only its existence.

\begin{lemma} \label{lem:Dij}
Any diagonal approximation $D$ induces a commutative diagram
of based $G$-spaces and~$G$-maps
\[
\xymatrix{
\ds \frac{F_k EG}{F_{k-1} EG} \ar[r]^-{D'_k} \ar@{ >->}[d]
	\ar@(u,u)[rr]^-{D'_{i,j}}
	& \ds \frac{F_k (EG \times EG)}{F_{k-1} (EG \times EG)}
		\ar[r]^-{\pr_{i,j}} \ar@{ >->}[d]
	& \ds \frac{F_i EG}{F_{i-1} EG} \wedge \frac{F_j EG}{F_{j-1} EG}
		\ar@{ >->}[d] \\
\ds \frac{EG}{F_{k-1} EG} \ar[r]^-{D_k}
	\ar@(d,d)[rr]_-{D_{i,j}}
	& \ds \frac{EG \times EG}{F_{k-1} (EG \times EG)} \ar[r]
	& \ds \frac{EG}{F_{i-1} EG} \wedge \frac{EG}{F_{j-1} EG}
}
\]
for all $i+j=k$.  The $G$-maps $D_{i,j}$ are compatible for
varying~$i$ and~$j$, in the sense that the squares
\[
\xymatrix{
\ds \frac{EG}{F_{i-1} EG} \wedge \frac{EG}{F_{j-1} EG} \ar[d]
& \ds \frac{EG}{F_{k-1} EG} \ar[l]_-{D_{i,j}} \ar[d] \ar[r]^-{D_{i,j}}
& \ds \frac{EG}{F_{i-1} EG} \wedge \frac{EG}{F_{j-1} EG} \ar[d] \\
\ds \frac{EG}{F_i EG} \wedge \frac{EG}{F_{j-1} EG}
& \ds \frac{EG}{F_k EG} \ar[l]_-{D_{i+1,j}} \ar[r]^-{D_{i,j+1}}
& \ds \frac{EG}{F_{i-1} EG} \wedge \frac{EG}{F_j EG}
}
\]
commute.
\end{lemma}

\begin{proof}
This follows from the inclusions
\[
D(F_{k-1} EG)
	\subset F_{k-1}(EG \times EG)
	\subset (F_{i-1} EG \times EG) \cup (EG \times F_{j-1} EG)
\]
and the splitting
\[
\frac{F_k EG}{F_{k-1} EG} \cong \bigvee_{i+j=k}
	\frac{F_i EG}{F_{i-1} EG} \wedge \frac{F_j EG}{F_{j-1} EG} \,.
\]
\end{proof}

\section{$G$-homotopy fixed points}
\label{subsec:htpyfix}

Let $X$ be an orthogonal $G$-spectrum. In this section we will construct a spectral sequence computing the homotopy groups of the~$G$-homotopy fixed points\index{homotopy fixed points} of $X$, that is, the $G$-fixed points of a fibrant replacement of the function spectrum~$F(EG_+, X)$:
\[
X^{hG} = F(EG_+,X)^{G} \,.
\]
To this end, note that the sequence of based $G$-spaces
\[
EG_+ = \frac{EG}{F_{-1} EG} \to \frac{EG}{F_0 EG} \to
	\dots \to \frac{EG}{F_{i-1} EG} \to \frac{EG}{F_i EG} \to \dots
	\to *
\]
induces a sequence 
\[
M_\star(X) = F(EG/EG_{-\star-1}, X)
\]
of orthogonal $G$-spectra. Explicitly, $M_\star(X)$ is the sequence 
\begin{multline*}
\dots \to F\left(\frac{EG}{F_i EG}, X\right)
	\to F\left(\frac{EG}{F_{i-1} EG}, X\right) \to \dots \\
	\dots \to F\left(\frac{EG}{F_0 EG}, X\right)
	\to F(EG_+, X) = F(EG_+, X) = \dots \,.
\end{multline*}

\begin{definition} \label{def:Ghfpspseq} 
The spectral sequence $(E^r(X), d^r) = (E^r(M_\star(X)), d^r)$ associated to the sequence $M_\star(X)$ above is called the \emph{$G$-homotopy fixed
point spectral sequence} of~$X$.
\end{definition}

Each map of function spectra $F(EG/F_i EG, X) \to F(EG/F_{i-1} EG, X)$ is a monomorphism of orthogonal $G$-spectra, but it is unlikely in general that these maps are $h$-cofibrations, so $M_\star(X)$ need not be a filtration. Since the sequence~$M_\star(X)$ is eventually constant, there is a natural $G$-equivalence
\[
M_\infty(X) = \Tel(M_\star(X)) \simeq_G F(EG_+, X) \,.
\]
There is also a $G$-equivalence
\[
\holim_s M_s(X) = \holim_s F\left(\frac{EG}{F_{-s-1} EG}, X\right)
	\cong F\left(\hocolim_i \frac{EG}{F_{i-1} EG}, X\right)
	\simeq_G * \,,
\]
since $\hocolim_i F_{i-1} EG \simeq_G EG$. We conclude that the $G$-homotopy fixed point spectral sequence is always conditionally convergent to the abutment
\[
A_\infty(M_\star(X)) \cong \pi^G_* \Tel(M_\star(X))
	\cong \pi^G_* F(EG_+, X) = \pi_*(X^{hG}) \,.
\]
Let us now explicitly compute the $E^1$-page of this spectral sequence.

\begin{lemma} \label{lem:E1MX}
The $E^1$-page of the $G$-homotopy fixed point spectral sequence of~$X$
is given by
\[
E^1_{-i,*}(M_\star(X)) \cong \pi^G_{-i+*} F\left(\frac{F_i EG}{F_{i-1} EG}, X\right) \cong \pi^G_{\ast} F(G_+, F(G^{\wedge i}, X))
\]
and the differential
\[
d^1_{-i,*} \: E^1_{-i,*}(X) \to E^1_{-i-1,*}(X)
\]
is contravariantly induced by the composite $G$-map
\[
\frac{F_{i+1} EG}{F_i EG} \longto
\frac{EG}{F_i EG} \simeq
\frac{EG}{F_{i-1} EG} \cup C\Bigl( \frac{F_i EG}{F_{i-1} EG} \Bigr)
\longto \Sigma \frac{F_i EG}{F_{i-1} EG} \,.
\]
\end{lemma}

\begin{proof}
The cofibre sequence
\[
\frac{F_i EG}{F_{i-1} EG} \longto \frac{EG}{F_{i-1} EG}
	\longto \frac{EG}{F_i EG}
\]
of based $G$-spaces is a homotopy cofibre sequence, hence induces a
homotopy fibre sequence
\[
M_{-i-1}(X) \longto M_{-i}(X) \longto
	F\left(\frac{F_i EG}{F_{i-1} EG}, X\right)
\]
of orthogonal $G$-spectra. It follows that
\begin{align*}
E^1_{-i,*}(X) &\cong \pi^G_{-i+*} F\left(\frac{F_i EG}{F_{i-1} EG}, X \right) \\
&\cong \pi^G_{-i+*} F(\Sigma^i(G^{\wedge i} \wedge G_+), X) \\
	&\cong \pi^G_* F(G_+, F(G^{\wedge i}, X))
\end{align*}
for $i\ge0$. The $d^1$-differential is the composite of the connecting homomorphism
\[
\pi^G_{-i+*}(M_{-i-1}(X) \to M_{-i}(X))
	\overset{\partial}\longto \pi^G_{-i-1+*}(M_{-i-1}(X))
	\cong \pi^G_{-i-1+*} F\left(\frac{EG}{F_i EG}, X\right)
\]
induced by
\[
\frac{EG}{F_i EG} \simeq \frac{EG}{F_{i-1} EG} \cup
C\left(\frac{F_i EG}{F_{i-1} EG}\right) \longto \Sigma \frac{F_i EG}{F_{i-1} EG} \,,
\]
and the homomorphism
\[
\pi^G_{-i-1+*}(M_{-i-1}(X)) \longto \pi^G_{-i-1+*}(M_{-i-2}(X) \to M_{-i-1}(X))
\]
induced by
\[
\frac{F_{i+1} EG}{F_i EG} \longto \frac{EG}{F_i EG} \,.
\]
\end{proof}

\begin{remark}
\label{rmk:fixpointPostnikov}\index{Whitehead tower}
When $G$ is finite, the spectral sequence $E^r(M_\star(X))$ agrees
with the~$G$-homotopy fixed point spectral sequence obtained from the $G$-equivariant Whitehead (or Postnikov) tower for $X$.
Greenlees and May prove this in~\cite{GM95}*{Theorem~B.8}.
When $G = \bT$ or~$\bU$ it is an accelerated version of the latter
spectral sequence. In Theorem~\ref{thm:E2MXisExtRG} we will give an
algebraic description of~$E^2(M_\star(X))$ when $X$ is an $R$-module
and $R[G]_*$ is finitely generated and projective over~$R_*$.
\end{remark}

The homotopy fixed point construction is a lax symmetric monoidal functor. To see this, let $\mu \: X \wedge Y \to Z$ be a pairing of orthogonal $G$-spectra, and recall the diagonal map $\Delta \: EG \to EG \times EG$. The associated pairing $X^{hG} \wedge Y^{hG} \to Z^{hG}$ is given by the composite
\begin{align*}
F(EG_+, X)^G \wedge F(EG_+, Y)^G
	&\overset{\alpha}\longto
F(EG_+ \wedge EG_+, X \wedge Y)^G \\
	&\overset{(\Delta_+)^*}\longto
F(EG_+, X \wedge Y)^G \\
	&\overset{\mu_*}\longto
F(EG_+, Z)^G \,.
\end{align*}
From this point of view it is hence relevant to understand how the homotopy fixed point spectral sequence interacts with multiplicative structures. First note that the maps $D_{i,j}$ from Lemma~\ref{lem:Dij} induce $G$-maps
\begin{align*}
F\left(\frac{EG}{F_{i-1} EG}, X\right) \wedge F\left(\frac{EG}{F_{j-1} EG}, Y\right)
	&\overset{\alpha}\longto
F\left(\frac{EG}{F_{i-1} EG} \wedge \frac{EG}{F_{j-1} EG} \, , \, X \wedge Y \right) \\
	&\overset{D_{i,j}^*}\longto
F\left(\frac{EG}{F_{k-1} EG} \, , \, X \wedge Y\right) \\
	&\overset{\mu_*}\longto
F\left(\frac{EG}{F_{k-1} EG} \, , \, Z\right)
\end{align*}
for $k = i+j$. These are compatible for varying $i$ and~$j$, in the sense of Definition~\ref{def:pairseq}, and so define the components $\bar \mu_{-i,-j}$ of a pairing 
\[
\bar\mu \: (M_\star(X), M_\star(Y)) \to M_\star(Z)
\]
of sequences of orthogonal $G$-spectra.

\begin{theorem} \label{thm:pairingbarmur}\index{homotopy fixed point spectral sequence!pairing}
Let $\mu \: X \wedge Y \to Z$ be a pairing of orthogonal $G$-spectra. There is then a pairing
\[
\bar\mu^r \: E^r(M_\star(X)) \otimes E^r(M_\star(Y))
        \longto E^r(M_\star(Z)) 
\]
of the associated $G$-homotopy fixed point spectral sequences, and the induced pairing~$\bar\mu_*$ on filtered abutments is compatible with the induced pairing 
\[
\bar{\mu}^\infty \: E^\infty(M_\star(X)) \otimes E^\infty(M_\star(Y)) \to E^\infty(M_\star(Z))
\]
of $E^\infty$-pages.

Moreover, the pairing~$\bar\mu^1$ of $E^1$-pages is contravariantly induced by
\[
D'_{i,j} \: \frac{F_k EG}{F_{k-1} EG}
	\longto \frac{F_i EG}{F_{i-1} EG} \wedge \frac{F_j EG}{F_{j-1} EG}
\]
under the isomorphism of Lemma~\ref{lem:E1MX}, and the pairing
\[
\bar\mu_* \: \pi_*(X^{hG}) \otimes \pi_*(Y^{hG})
        \longto \pi_*(Z^{hG})
\]
equals the pairing induced by $X^{hG} \wedge Y^{hG} \to Z^{hG}$.
\end{theorem}

\begin{proof}
In the paragraph before this theorem we noted that a map $\mu \: X \wedge Y \to Z$ of orthogonal $G$-spectra gives rise to a pairing $\bar{\mu} \: (M_\star(X),M_{\star}(Y)) \to M_{\star}(Z)$ of sequences. By Theorem~\ref{thm:mainmultII} it follows that we have an induced pairing between the associated spectral sequences, and that the induced pairing $\bar{\mu}_*$ on filtered abutments is compatible with the pairing $\bar{\mu}^\infty$ of $E^\infty$-pages.

Tracing through the definitions shows that the pairing $\bar\mu^1_{-i,-j}$ of $E^1$-pages is compatible with the pairing induced by $D'_{i,j}$ under the isomorphism
\begin{align*}
E^1_{-i,*}(M_\star(X))
&= \pi^G_{-i+*}(M_{-i-1}(X) \to M_{-i}(X)) \\
&\cong \pi^G_{-i+*} F\left(\frac{F_i EG}{F_{i-1} EG}, X\right)
\end{align*}
and its analogues for~$Y$ and~$Z$.

The abutment $A_\infty(M_\star(X)) \cong \pi^G_* F(EG_+,
X)$ is filtered by the images
\[
F_s \pi^G_* F(EG_+, X) = \im(\pi^G_* F(EG/EG_{-s-1}, X)
	\longto \pi^G_* F(EG_+, X)) \,.
\]
Note that this exhaustive filtration is constant for $s\ge0$. The pairing $\bar\mu_*$ is induced by the composite map
\begin{align*}
\bar\mu_{0,0} \: F(EG_+, X) \wedge F(EG_+, Y)
	&\xrightarrow{\alpha} F(EG_+ \wedge EG_+, X \wedge Y) \\
	&\xrightarrow{D_{0,0}^*} F(EG_+, X \wedge Y) \\
	&\xrightarrow{\mu_*} F(EG_+, Z) \,.
\end{align*}
In view of the based $G$-homotopy $\Delta_+ \simeq D_+ = D_{0,0}$,
it is also induced by the composite map
\begin{align*}
F(EG_+, X) \wedge F(EG_+, Y)
	&\xrightarrow{\alpha} F(EG_+ \wedge EG_+, X \wedge Y) \\
	&\xrightarrow{\Delta_+^*} F(EG_+, X \wedge Y) \\
	&\xrightarrow{\mu_*} F(EG_+, Z) \,,
\end{align*}
where $\Delta_+ \: EG_+ \to (EG \times EG)_+ \cong EG_+ \wedge EG_+$.
\end{proof}

\begin{corollary}\index{homotopy fixed point spectral sequence!multiplicative}
If $(X, \mu \: X \wedge X \to X)$ is a multiplicative orthogonal $G$-spectrum, then the $G$-homotopy fixed point spectral sequence $(E^r(M_\star(X)), d^r)$ is a conditionally convergent and multiplicative spectral sequence, with multiplicative abutment $\pi_*(X^{hG})$.
\end{corollary}

\section{Algebraic description of the $E^1$- and $E^2$-pages}
\label{sec:algdescriptionhfpss}

Under suitable flatness hypotheses there is an algebraic description of the first two pages of the homotopy fixed point spectral sequence. Recall from Chapter~\ref{sec:homotopygroups} that $R$ is our `ground' commutative orthogonal ring spectrum.  We write
\[
R_*(X) = \pi_*(R \wedge
X)
\]
for the associated (reduced) homology theory.  We will assume
that $R[G]_*$ is flat over~$R_*$, so that $R[G]_*$ is a cocommutative Hopf algebra over~$R_*$, per Lemma~\ref{lem:RGHopf}. Let us write
\[
E = R \wedge EG_+ \och E_i = R \wedge (F_i EG)_+ \,.
\]
Each map $E_{i-1} \to E_i$ is a $q$-cofibration, hence a strong
$h$-cofibration, so that $E_\star$ is a filtration
\[
\cdots \longto E_{i-1} \longto E_i \longto E_{i+1} \longto \cdots
\]
of $R$-modules in orthogonal $G$-spectra.  Here $E_i = *$ for $i < 0$, and 
\[
E_\infty = \Tel(E_\star) \simeq_G E \,.
\]
The $R$- and $G$-equivariant collapse map $1 \wedge c \: E =
R \wedge EG_+ \to R \wedge S^0 = R$ is a non-equivariant~$R$-equivalence, inducing an $R[G]_*$-module isomorphism $\pi_*(E) \cong R_*$.

\begin{definition} \label{def:P}
Let $(P_{*,*}, \partial) = NB_*(R_*,R[G]_*,R[G]_*)$ denote the normalised bar resolution\index{normalised bar resolution}, as defined in Construction~\ref{const:barconst}.
\end{definition}

Explicitly, the normalised bar resolution of the $R[G]_*$-module $R_*$ is a non-negative chain complex given in homological degree $n \geq 0$ as 
\[
P_{n,*} = NB_n(R_*, R[G]_*, R[G]_*)
	= \overline{R[G]}_*^{\otimes n} \otimes_{R_*} R[G]_* \,,
\]
where
\[
\overline{R[G]}_* = \cok(\eta \: R_* \to R[G]_*)
	\cong \ker(\epsilon \: R[G]_* \to R_*)
\]
denotes the augmentation (co-)ideal, and $\overline{R[G]}_*^{\otimes n}$ is its $n$-th tensor power over~$R_*$. The boundary $\partial_n \: P_{n,*} \to P_{n-1,*}$ is induced by the alternating sum of face operators
\[
\sum_{i=0}^n (-1)^i d_i
\]
for $n\ge1$, with
\[
d_i = \begin{cases}
\epsilon \otimes 1^{\otimes n}
	& \text{for $i=0$,} \\
1^{\otimes i-1} \otimes \phi \otimes 1^{\otimes n-i}
	& \text{for $0 < i \le n$.}
\end{cases}
\]
Note that the simplicial contraction \cite{May72}*{Prop.~9.8}
of $B_\bullet(R_*, R[G]_*, R[G]_*)$ shows that the augmentation
$\epsilon \: P_{0,*} = R[G]_* \to R_*$ admits an
$R_*$-linear chain homotopy inverse, so that
the augmented chain complex
\[
\dots \longto P_{q,*}
	\overset{\partial_q}\longto P_{q-1,*}
	\to \dots \to P_{1,*}
	\overset{\partial_1}\longto P_{0,*}
	\overset{\epsilon}\longto R_* \longto 0
\]
is exact.  Hence $(P_{*,*}, \partial)$ is a flat $R[G]_*$-module
resolution of $R_*$.  

\begin{lemma} \label{lem:E1isP}\index{bar spectral sequence}
If $R[G]_*$ is flat over $R_*$, then the $(E^1, d^1)$-page of the
non-equi\-variant homotopy spectral sequence
\[
E^1_{i,*} = \pi_{i+*}(E_{i-1} \to E_i)
\]
associated to $E_\star$ is isomorphic to $(P_{*,*}, \partial)$. The edge homomorphism $P_{0,*} \to \pi_*(E) \cong R_*$ is equal to
the augmentation $\epsilon \: R[G]_* \to R_*$, and makes $(P_{*,*}, \partial)$ a flat $R[G]_*$-module resolution of~$R_*$.  In particular, the spectral sequence collapses at the $E^2$-page, where is it given by
\[
E^2 = E^\infty \cong R_*
\]
concentrated in filtration degree $i=0$.
\end{lemma}

\begin{proof}
The $R$-module filtration $E_\star$ has an associated $R[G]_*$-module
spectral sequence (for non-equivariant homotopy groups) with $E^1$-page
\[
E^1_{i,*} = \pi_{i+*}(E_{i-1} \to E_i)
	\cong R_{i+*}\left(\frac{F_i EG}{F_{i-1} EG}\right)
\]
and $d^1$-differential equal to the composite
\[
R_{i+*}\left(\frac{F_i EG}{F_{i-1} EG}\right)
	\overset{\partial}\longto R_{i-1+*}(F_{i-1} EG_+)
	\longto R_{i-1+*}\left(\frac{F_{i-1} EG}{F_{i-2} EG}\right) \,.
\]
By the proof of~\cite{Seg68}*{Prop.~5.1} or~\cite{May72}*{Thm.~11.14},
$(E^1, d^1)$ is the normalized
chain complex associated to the simplicial
$R[G]_*$-module
\[
[q] \mapsto R_*(B_q(*, G, G)_+)
	= R_*((G^q \times G)_+) \,.
\]
The products
\begin{multline*}
R[G]_* \otimes_{R_*} R[G]_* \otimes_{R_*} \dots \otimes_{R_*} R[G]_* \\
\overset{\cdot}\longto
\pi_*(R \wedge G_+ \wedge_R R \wedge G_+ \wedge_R \dots \wedge_R R \wedge G_+)
\\
\cong \pi_*(R \wedge G_+ \wedge G_+ \wedge \dots \wedge G_+)
\end{multline*}
induce a homomorphism of simplicial $R[G]_*$-modules
\[
B_\bullet(R_*, R[G]_*, R[G]_*)
	\longto R_*(B_\bullet(*, G, G))_+) \,.
\]
Since $R[G]_*$ is assumed to be flat over~$R_*$ the products are isomorphisms, so that~$(E^1, d^1)$ is indeed isomorphic to the normalized chain complex associated to the simplicial $R[G]_*$-module $B_\bullet(R_*,R[G]_*,R[G]_*)$. 
\end{proof}

\begin{remark}
If $R[G]_*$ is projective over~$R_*$, then
$\overline{R[G]}_*$ is also $R_*$-projective, and each
$P_{q,*}$ is $R[G]_*$-projective by Lemma~\ref{lem:proj}. It follows that the chain complex~$(P_{*,*}, \partial)$ is a projective~$R[G]_*$-module resolution of $R_*$.
Moreover, if $R[G]_*$ is finitely generated over~$R_*$,
then so is $\overline{R[G]}_*$, and each~$P_{q,*}$ is
finitely generated as an $R[G]_*$-module. We conclude that
$(P_{*,*}, \partial)$ is a projective resolution of finite
type, in this case.
\end{remark}

To deal with the multiplicative structure of the spectral sequence we introduce the convolution product $(E \wedge_R E)_\star$. Explicitly, this is given by
\[
(E \wedge_R E)_k = R \wedge F_k(EG \times EG)_+ \,,
\]
with filtration subquotients
\[
\frac{(E \wedge_R E)_k}{(E \wedge_R E)_{k-1}} \cong
	\bigvee_{i+j=k} \frac{E_i}{E_{i-1}} \wedge_R \frac{E_j}{E_{j-1}} \,.
\]
Let $\In_{i,j}$ denote the inclusion of the $(i,j)$-th summand
in this splitting.

\begin{lemma} \label{lem:E1isPP}
The $(E^1, d^1)$-page of the homotopy spectral sequence associated to~$(E \wedge_R E)_\star$ is isomorphic to the tensor product
\[
(P_{*,*} \otimes_{R_*} P_{*,*}, \partial \otimes 1 + 1 \otimes \partial)
\,,
\]
with the same signs occurring in the boundary as specified in Section~\ref{subsec:chaincomplex}. In particular, this spectral sequence collapses at the $E^2$-page, where it is given by 
\[
E^2 = E^\infty \cong R_* \otimes_{R_*} R_* \cong R_*
\]
concentrated in filtration degree~$0$.
\end{lemma}

\begin{proof}
Theorem~\ref{thm:mainmultII} applied to the initial pairing $\iota \: (E_\star, E_\star) \to (E \wedge_R E)_\star$ gives us a pairing 
\[
\iota^r \: E^r(E_\star) \otimes_{R_*} E^r(E_\star)
	\longto E^r((E \wedge_R E)_\star)
\]
of $R[G]_*$-module spectral sequences. Since each copy of $E_\star$ is a filtration, the pairing
\[
\iota^1_{i,j} \: P_{i,*} \otimes_{R_*} P_{j,*}
	= E^1_{i,*}(E_\star) \otimes_{R_*} E^1_{j,*}(E_\star)
	\longto E^1_{k,*}((E \wedge_R E)_\star) \,,
\]
for $r=1$ and $i+j=k$, is induced by the product
\[
P_{i,*} \otimes_{R_*} P_{j,*}
	\overset{\cdot}\longto \pi_*\left(\frac{E_i}{E_{i-1}} \wedge_R \frac{E_j}{E_{j-1}}\right)
\]
and the inclusion
\[
\In_{i,j} \: \frac{E_i}{E_{i-1}} \wedge_R \frac{E_j}{E_{j-1}}
	\longto \frac{(E \wedge_R E)_k}{(E \wedge_R E)_{k-1}} \,.
\]
Since $R[G]_*$ is flat over~$R_*$, so that each $P_{i,*}$ is flat over $R_*$, the product is an isomorphism.  Adding these
together for $i+j=k$ we obtain the degree~$k$ part of an isomorphism of $R[G]_*$-module chain complexes
\[
\iota^1 \: P_{*,*} \otimes_{R_*} P_{*,*}
	\overset{\cong}\longto E^1_{*,*}((E \wedge_R E)_\star) \,.
\]
In particular, Theorem~\ref{thm:mainmultII} ensures that the tensor product boundary operator $\partial \otimes 1 + 1 \otimes \partial$ on the left hand side corresponds to the $d^1$-differential on the right hand side.  The calculation of the $E^2$-page then follows as in the proof of Proposition~\ref{prop:Pidentityextension}.
\end{proof}

\begin{lemma} \label{lem:DandPsi}
The diagonal approximation $D \: EG \to EG \times EG$
induces a map of filtrations $1 \wedge D_+ \: E_\star \to
(E \wedge_R E)_\star$ and
a chain map
\[
(1 \wedge D_+)^1 \: E^1(E_\star) \longto E^1((E \wedge_R E)_\star)
\,,
\]
which corresponds, under the isomorphisms of Lemma~\ref{lem:E1isP}
and Lemma~\ref{lem:E1isPP}, to an~$R[G]_*$-module chain map
\[
\Psi \: P_{*,*} \longto P_{*,*} \otimes_{R_*} P_{*,*} \,.
\]
In particular, the component 
\[
\Psi_{i,j} = \pr_{i,j} \circ \Psi_k \: P_{k,*} \to P_{i,*} \otimes_{R_*} P_{j,*}
\]
of $\Psi_k$, for $k=i+j$, is induced by the $G$-map $D'_{i,j}$ of Lemma~\ref{lem:Dij} and Theorem~\ref{thm:pairingbarmur}.

The chain map $\Psi$ is characterised, uniquely up to chain homotopy equivalence, by the commutative square
\[
\xymatrix{
P_{*,*} \ar[r]^-{\Psi} \ar[d]_-{\epsilon}
& P_{*,*} \otimes_{R_*} P_{*,*} \ar[d]^-{\epsilon \otimes \epsilon} \\
R_* \ar[r]^-{\cong} & R_* \otimes_{R_*} R_*
}
\]
of $R[G]_*$-module complexes.
\end{lemma}

\begin{proof}
The map of $E^1$-pages induced by the diagonal approximation
is induced by~$1 \wedge D'_k$, and the $(i,j)$-th component
in the direct sum splitting of its target can be recovered
by projecting to that summand, which is therefore induced
by $1 \wedge D'_{i,j}$.

By naturality of the edge homomorphism, we have a commutative
square of~$R[G]_*$-modules
\[
\xymatrix{
P_{0,*} \ar[r]^-{\Psi_0} \ar[d]_-{\epsilon}
  & P_{0,*} \otimes_{R_*} P_{0,*} \ar[d]^-{\epsilon \otimes \epsilon} \\
R_* \ar[r]^-{\cong} & R_* \otimes_{R_*} R_* \rlap{\,.}
}
\]
Hence the $R[G]_*$-module chain map $\Psi \: P_{*,*} \to P_{*,*}
\otimes_{R_*} P_{*,*}$ is a lift of the isomorphism $R_* \cong
R_* \otimes_{R_*} R_*$.  Since $\epsilon \: P_{*,*} \to R_*$ is an
$R[G]_*$-projective complex over~$R_*$, and $\epsilon \otimes \epsilon \: P_{*,*} \otimes_{R_*} P_{*,*} \to R_* \otimes_{R_*} R_*$ is a resolution, it follows from~\cite{Mac95}*{Thm.~III.6.1} that such a chain map $\Psi$ exists and is unique up to chain homotopy.
\end{proof}

We now suppose that $X$ is an $R$-module in orthogonal $G$-spectra. There are then compatible adjunction equivalences
\[
F_R(E/E_{i-1}, X) \cong F(EG/F_{i-1} EG, X) = M_{-i}(X)
\]
for all~$i$. The left hand side exhibits $M_\star(X)$ as a sequence of $R$-modules in orthogonal $G$-spectra, so that the $G$-homotopy fixed point spectral sequence $E^r(M_\star(X))$ is a spectral sequence of $R_*$-modules. Theorem~\ref{thm:pairingbarmur} readily generalizes: If $Y$ and $Z$ are also $R$-modules in orthogonal $G$-spectra, and $\mu \: X \wedge_R Y \to Z$ is a map in this category, then we obtain a pairing of $R_*$-module spectral sequences
\[
\bar\mu^r \: E^r(M_\star(X)) \otimes_{R_*} E^r(M_\star(Y))
	\longto E^r(M_\star(Z))
\]
such that the resulting pairing of $E^\infty$-pages is compatible
with the $R_*$-linear pairing
\[
\bar\mu_* \: \pi^G_* F(EG_+, X) \otimes_{R_*} \pi^G_* F(EG_+, Y)
	\longto \pi^G_* F(EG_+, Z)
\]
of abutments. We can now give algebraic descriptions of
the $(E^1, d^1)$-pages and the pairing $\bar\mu^1$, for
$R[G]_*$ projective over~$R_*$.

\begin{proposition} \label{prop:E1d1mu1fromPPsi}
Assume that $R[G]_*$ is projective as an $R_*$-module.
There is then a natural isomorphism
\[
E^1_{-i,*}(M_\star(X)) \cong \Hom_{R[G]_*}(P_{i,*}, \pi_*(X))
\]
of $R_*$-modules. Under this isomorphism, the $d^1$-differential
\[
d^1_{-i,*} \: E^1_{-i,*}(M_\star(X)) \longto E^1_{-i-1,*}(M_\star(X))
\]
corresponds to the boundary in the chain complex, with signs as specified in Section~\ref{subsec:chaincomplex}. The pairing
\[
\bar\mu^1 \: E^1_{-i,*}(M_\star(X)) \otimes_{R_*} E^1_{-j,*}(M_\star(Y))
	\longto E^1_{-k,*}(M_\star(Z))
\]
with $i+j=k$ is contravariantly induced by the component
\[
\Psi_{i,j} \: P_{k,*} \longto P_{i,*} \otimes_{R_*} P_{j,*}
\]
of the chain map $\Psi$.
\end{proposition}

\begin{proof}
By Lemma~\ref{lem:E1MX} and adjunction there are isomorphisms
\[
E^1_{-i,*}(M_\star(X)) \cong \pi^G_{-i+*} F(F_i EG/F_{i-1} EG, X)
\cong \pi^G_{-i+*} F_R(E_i/E_{i-1}, X) \,.
\]
Note that the spectrum appearing in the last term can be written $F_R(E_i/E_{i-1}, X) \cong F(G_+, X')$ with
\[
X' = F_R(R \wedge G^{\wedge i}, X) \cong F(G^{\wedge i}, X) \,.
\]
Under our assumption that $R[G]_*$ is projective, it follows from Proposition~\ref{prop:omegaXiso} that the natural $R_*$-module homomorphism
\[
\omega \: \pi^G_{-i+*} F_R(E_i/E_{i-1}, X)
	\overset{\cong}\longto
	\Hom_{R[G]_*}(R_*, \pi_{-i+*} F_R(E_i/E_{i-1}, X))
\]
is an isomorphism. Moreover, since $P_{i,*} = \pi_{i+*}(E_i/E_{i-1})$ is projective over $R[G]_*$ and hence also over~$R_*$, it follows that the natural $R[G]_*$-module homomorphism
\[
\pi_{-i+*} F_R(E_i/E_{i-1}, X)
	\overset{\cong}\longto \Hom_{R_*}(\pi_{i+*}(E_i/E_{i-1}), \pi_*(X))
	= \Hom_{R_*}(P_{i,*}, \pi_*(X))
\]
is an isomorphism.  Applying the functor $\Hom_{R[G]_*}(R_*, -)$
yields an isomorphism
\begin{multline*}
\Hom_{R[G]_*}(R_*, \pi_{-i+*} F_R(E_i/E_{i-1}, X)) \\
\overset{\cong}\longto
\Hom_{R[G]_*}(R_*, \Hom_{R_*}(P_{i,*}, \pi_*(X)))
\cong
\Hom_{R[G]_*}(P_{i,*}, \pi_*(X)) \,.
\end{multline*}
Composing this chain of $R_*$-module isomorphisms gives the
asserted natural isomorphism.

We now identify the $d^1$-differential. By Lemma~\ref{lem:E1MX} again, we have a commutative diagram
\[
\xymatrix{
E^1_{-i,*}(M_\star(X)) \ar[dd]_-{d^1_{-i,*}} \ar[r]^-{\cong}
& \pi^G_{-i+*} F_R(E_i/E_{i-1}, X) \ar[d] \\
& \pi^G_{-i-1+*} F_R(E/E_i, X) \ar[d] \\
E^1_{-i-1,*}(M_\star(X)) \ar[r]^-{\cong}
& \pi^G_{-i-1+*} F_R(E_{i+1}/E_i, X)
}
\]
of $R_*$-modules. By the naturality of $\omega$ in Lemma~\ref{lem:omegaXnat} the diagram
\[
\xymatrix{
\pi^G_{-i+*} F_R(E_i/E_{i-1}, X) \ar[d] \ar[r]^-{\omega}_-{\cong}
& \Hom_{R[G]_*}(R_*, \pi_{-i+*} F_R(E_i/E_{i-1}, X)) \ar[d] \\
\pi^G_{-i-1+*} F_R(E/E_i, X) \ar[d] \ar[r]^-{\omega}
& \Hom_{R[G]_*}(R_*, \pi_{-i-1+*} F_R(E/E_i, X)) \ar[d] \\
\pi^G_{-i-1+*} F_R(E_{i+1}/E_i, X) \ar[r]^-{\omega}_-{\cong}
& \Hom_{R[G]_*}(R_*, \pi_{-i-1+*} F_R(E_{i+1}/E_i, X))
}
\]
commutes. Note that these two diagrams fit together along one edge. We also have a commutative diagram of $R[G]_*$-modules
\[
\xymatrix{
\pi_{-i+*} F_R(E_i/E_{i-1}, X) \ar[d] \ar[r]^-{\cong}
& \Hom_{R_*}(P_{i,*}, \pi_*(X)) \ar[d]
	\ar@/^7pc/[dd]^(0.3){\Hom(\partial_{i+1},1)} \\
\pi_{-i-1+*} F_R(E/E_i, X) \ar[r] \ar[d]
& \Hom_{R_*}(\pi_{i+1+*}(E/E_i), \pi_*(X)) \ar[d] \\
\pi_{-i-1+*} F_R(E_{i+1}/E_i, X) \ar[r]^-{\cong}
& \Hom_{R_*}(P_{i+1,*}, \pi_*(X))
}
\]
since $\partial_{i+1} \: P_{i+1,*} \to P_{i,*}$ can be
calculated by either composite from the left to the right in
the diagram
\[
\xymatrix{
\pi_{i+1+*}(E_{i+1}/E_i) \ar[r] \ar[dr]_-{\partial}
	& \pi_{i+1+*}(E/E_i) \ar[dr]^-{\partial} \ar[d]^-{\partial} \\
& \pi_{i+*}(E_i) \ar[r]
	& \pi_{i+*}(E_i/E_{i-1}) \,.
}
\]
Applying $\Hom_{R[G]_*}(R_*, -)$ we obtain a commutative diagram
of $R_*$-modules, which fits together with the previous one.
Hence the square
\[
\xymatrix{
E^1_{-i,*}(M_\star(X)) \ar[r]^-{\cong} \ar[d]_-{d^1_{-i,*}}
	& \Hom_{R[G]_*}(P_{i,*}, \pi_*(X)) \ar[d]^-{\Hom(\partial_{i+1},1)} \\
E^1_{-i-1,*}(M_\star(X)) \ar[r]^-{\cong}
	& \Hom_{R[G]_*}(P_{i+1,*}, \pi_*(X))
}
\]
commutes, as asserted.

We now identify the multiplicative structure on the $E^1$-page. By Theorem~\ref{thm:pairingbarmur}, the diagram
\[
\xymatrix{
E^1_{-i,*}(M_\star(X)) \otimes_{R_*} E^1_{-j,*}(M_\star(Y))
	\ar[r]^-{\bar\mu^1} \ar[d]_-{\cong}
& E^1_{-k,*}(M_\star(Z))
	\ar[d]^-{\cong} \\
\pi^G_{-i+*} F_R(E_i/E_{i-1}, X)
	\otimes_{R_*} \pi^G_{-j+*} F_R(E_j/E_{j-1}, Y) \ar[r]
& \pi^G_{-k+*} F_R(E_k/E_{k-1}, Z)
}
\]
commutes, where the lower arrow is induced by
\[
1 \wedge D'_{i,j} \: E_k/E_{k-1} \to E_i/E_{i-1} \wedge_R E_j/E_{j-1} \,.
\]
Since the natural homomorphism $\omega$ is monoidal, per Lemma~\ref{lem:omegaXmon}, the composite
\begin{multline*}
\pi^G_{-i+*} F_R(E_i/E_{i-1}, X)
        \otimes_{R_*} \pi^G_{-j+*} F_R(E_j/E_{j-1}, Y)
	\longto \pi^G_{-k+*} F_R(E_k/E_{k-1}, Z) \\
	\overset{\omega}\longto
\Hom_{R[G]_*}(R_*, \pi_{-k+*} F_R(E_k/E_{k-1}, Z))
\end{multline*}
is equal to the composite
\begin{multline*}
\pi^G_{-i+*} F_R(E_i/E_{i-1}, X)
	\otimes_{R_*} \pi^G_{-j+*} F_R(E_j/E_{j-1}, Y)
	\xrightarrow{\omega \otimes \omega} \\
\Hom_{R[G]_*}(R_*, \pi_{-i+*} F_R(E_i/E_{i-1}, X)) 
	\otimes_{R_*} \Hom_{R[G]_*}(R_*, \pi_{-j+*} F_R(E_j/E_{j-1}, Y))
\\
	\overset{\alpha}\longto
\Hom_{R[G]_*}(R_*, \pi_{-i+*} F_R(E_i/E_{i-1}, X) \otimes_{R_*}
	\pi_{-j+*} F_R(E_j/E_{j-1}, Y)) \\
	\overset{\mu_*}\longto
\Hom_{R[G]_*}(R_*, \pi_{-k+*} F_R(E_k/E_{k-1}, Z)) \,.
\end{multline*}
Note that the final arrow is also induced by $1 \wedge D'_{i,j}
\: E_k/E_{k-1} \to E_i/E_{i-1} \wedge_R E_j/E_{j-1}$.
Next, we use the commutative diagram
\[
\xymatrix{
\pi_{-i+*} F_R(E_i/E_{i-1}, X) \otimes_{R_*} \pi_{-j+*} F_R(E_j/E_{j-1}, Y)
	\ar[r] \ar[d]_-{\cong}
& \pi_{-k+*} F_R(E_k/E_{k-1}, Z) \ar[d]^-{\cong} \\
\Hom_{R_*}(P_{i,*}, \pi_*(X)) \otimes_{R_*} \Hom_{R_*}(P_{j,*}, \pi_*(Y))
	\ar[r]
& \Hom_{R_*}(P_{k,*}, \pi_*(Z))
}
\]
of $R[G]_*$-modules, where the lower homomorphism is induced by $1 \wedge
D'_{i,j}$.  In view of the isomorphism $P_{i,*} \otimes_{R_*} P_{j,*}
\cong \pi_{i+j}(E_i/E_{i-1} \wedge_R E_j/E_{j-1})$ from the proof of
Lemma~\ref{lem:E1isPP}, this is the same homomorphism
as that induced by $\Psi_{i,j}$, as defined in Lemma~\ref{lem:DandPsi}.

Applying the monoidal functor $\Hom_{R[G]_*}(R_*, -)$, we obtain a
commutative square of $R_*$-modules.  Combining these results we have
a commutative square
\[
\xymatrix{
E^1_{-i,*}(M_\star(X)) \otimes_{R_*} E^1_{-j,*}(M_\star(Y))
	\ar[r]^-{\bar\mu^1} \ar[d]_-{\cong}
& E^1_{-k,*}(M_\star(Z))
	\ar[d]^-{\cong} \\
\Hom_{R[G]_*}(P_{i,*}, \pi_*(X))
	\otimes_{R_*} \Hom_{R[G]_*}(P_{j,*}, \pi_*(Y)) \ar[r]
& \Hom_{R[G]_*}(P_{k,*}, \pi_*(Z))
}
\]
where the lower homomorphism is induced by $\Psi_{i,j}$, meaning that it is equal to the composite
\begin{align*}
\Hom_{R[G]_*}&(P_{i,*}, \pi_*(X))
        \otimes_{R_*} \Hom_{R[G]_*}(P_{j,*}, \pi_*(Y)) \\
&\xrightarrow{\alpha}
\Hom_{R[G]_*}(P_{i,*} \otimes_{R_*} P_{j,*}, \pi_*(X) \otimes_{R_*} \pi_*(Y))
\\
&\xrightarrow{\Psi_{i,j}^*}
\Hom_{R[G]_*}(P_{k,*}, \pi_*(X) \otimes_{R_*} \pi_*(Y)) \\
&\xrightarrow{\mu_*}
\Hom_{R[G]_*}(P_{k,*}, \pi_*(Z)) \,.
\end{align*}
This is the same as the $(i,j)$-component of the chain map
\begin{align*}
\Hom_{R[G]_*}&(P_{*,*}, \pi_*(X))
        \otimes_{R_*} \Hom_{R[G]_*}(P_{*,*}, \pi_*(Y)) \\
&\xrightarrow{\alpha}
\Hom_{R[G]_*}(P_{*,*} \otimes_{R_*} P_{*,*}, \pi_*(X) \otimes_{R_*} \pi_*(Y))
\\
&\xrightarrow{\Psi^*}
\Hom_{R[G]_*}(P_{*,*}, \pi_*(X) \otimes_{R_*} \pi_*(Y)) \\
&\xrightarrow{\mu_*}
\Hom_{R[G]_*}(P_{*,*}, \pi_*(Z))
\end{align*}
induced by $\Psi$.
\end{proof}

As a direct consequence, we get a description of the $E^2$-page of the homotopy fixed point spectral sequence.

\begin{theorem} \label{thm:E2MXisExtRG}\index{cup product}\index{homotopy fixed point spectral sequence!pairing}
Let $G$ be a compact Lie group and let $R$ be a commutative orthogonal ring spectrum. Moreover, let $\mu \: X \wedge_R Y \to Z$ be a pairing of $R$-modules in orthogonal $G$-spectra. Assume that $R[G]_*$ is projective as an $R_*$-module. Then there is a natural isomorphism
\[
E^2_{-i,*}(M_\star(X)) \cong \Ext_{R[G]_*}^i(R_*, \pi_*(X))
\]
of $R_*$-modules, for each integer~$i$. The pairing
\[
\bar\mu^2 \: E^2_{-i,*}(M_\star(X)) \otimes_{R_*} E^2_{-j,*}(M_\star(Y))
	\longto E^2_{-i-j,*}(M_\star(Z))
\]
is given by the cup product
\[
\smile \: \Ext_{R[G]_*}^i(R_*, \pi_*(X))
	\otimes_{R_*} \Ext_{R[G]_*}^j(R_*, \pi_*(Y))
	\longto \Ext_{R[G]_*}^{i+j}(R_*, \pi_*(Z))
\]
associated to the $R[G]_*$-module pairing $\mu_* \: \pi_*(X) \otimes_{R_*}
\pi_*(Y) \to \pi_*(Z)$, in $\Ext$ over the Hopf algebra~$R[G]_*$.
\end{theorem}

\begin{proof}
By Proposition~\ref{prop:E1d1mu1fromPPsi} the first page of the spectral sequence, together with its $d^1$-differential, is identified with the chain complex $\Hom_{R[G]_*}(P_{*,*}, \pi_*(X))$ where $(P_{*,*}, \partial)$ is a projective resolution of~$R_*$. It follows that the $E^2$-page is given by the homology of this chain complex, which by definition is the graded $R_*$-module 
\[
E^2_{*,*}(X) \cong \Ext_{R[G]_*}^*(R_*,\pi_*(X)) \,.
\]
Let us now identify the multiplication on the $E^2$-page with the cup product. Let $f \: P_{i,*} \to \pi_*(X)$ and $g \: P_{j,*} \to \pi_*(Y)$ be (graded)~$R[G]_*$-module homomorphisms with $f \partial_{i+1} = 0$ and $g \partial_{j+1} = 0$.  They correspond to $i$- and $j$-cycles in $\Hom_{R[G]_*}(P_{*,*}, \pi_*(X))$ and $\Hom_{R[G]_*}(P_{*,*}, \pi_*(Y))$, respectively, with homology classes $[f] \in E^2_{-i,*}(X)$ and $[g] \in E^2_{-j,*}(Y)$. The pairing of $E^2$-pages sends $[f] \otimes [g]$ to the homology class in $E^2_{-k,*}(Z)$ of the~$k$-cycle given by the composite (graded) $R[G]_*$-module homomorphism
\[
P_{k,*} \overset{\Psi_{i,j}}\longto P_{i,*} \otimes_{R_*} P_{j,*}
	\xrightarrow{f \otimes g} \pi_*(X) \otimes_{R_*} \pi_*(Y)
	\xrightarrow{\mu_*} \pi_*(Z) \,.
\]
The verification that $\mu_* (f \otimes g) \Psi_{i,j}$ is a $k$-cycle
uses the fact that $\Psi_{i,j}$ is a component of an $R[G]_*$-module
chain map $\Psi \: P_{*,*} \to P_{*,*} \otimes_{R_*} P_{*,*}$, so that
\[
\Psi_{i,j} \partial_{k+1} = (\partial_{i+1} \otimes 1)\Psi_{i+1,j} +
(1 \otimes \partial_{j+1})\Psi_{i,j+1} \,.
\]
This is the definition of the cup product
\[
\smile \: \Ext_{R[G]_*}^*(R_*, \pi_*(X))
\otimes_{R_*} \Ext_{R[G]_*}^*(R_*, \pi_*(Y))
\longto \Ext_{R[G]_*}^*(R_*, \pi_*(Z))
\]
associated to the pairing $\mu_*$. See Section~\ref{sec:multTate}.
\end{proof}

\begin{remark}
A well-known consequence of the comparison theorem~\cite{Mac95}*{Thm.~III.6.1} is that
\[
\Ext_{R[G]_*}^i(R_*, \pi_*(X)) = H^i(\Hom_{R[G]_*}(P_{*,*}, \pi_*(X)))
\]
can be calculated with any projective $R[G]_*$-module
resolution $P_{*,*}$ of $R_*$, not necessarily the one introduced in Definition~\ref{def:P}. Likewise, by Proposition~\ref{prop:Pidentityextension}, the cup product can be calculated with any~$R[G]_*$-module chain map 
\[
\Psi \: P_{*,*} \longto P_{*,*}
\otimes_{R_*} P_{*,*}
\]
lifting $R_* \cong R_* \otimes_{R_*} R_*$, not necessarily the one induced by a given diagonal approximation~$D$.
\end{remark}

\begin{example}
When $G$ is finite,
\[
R[G]_* = R_*[G] \cong \bZ[G] \otimes_{\bZ} R_* \,,
\]
any projective $\bZ[G]$-module resolution $Q_*$ of $\bZ$ induces up to
a projective $R_*[G]$-module resolution $P_{*,*} = Q_* \otimes_{\bZ}
R_*$ of $R_*$, and any $\bZ[G]$-module diagonal approximation $\Psi \:
Q_* \to Q_* \otimes_{\bZ} Q_*$ induces up to an $R_*[G]$-module diagonal
approximation $\Psi \otimes 1 \: P_{*,*} = Q_* \otimes_{\bZ} R_* \to Q_*
\otimes_{\bZ} Q_* \otimes_{\bZ} R_* \cong P_{*,*} \otimes_{R_*} P_{*,*}$.
Hence there is a natural isomorphism
\[
\Ext_{R_*[G]}^i(R_*, \pi_*(X)) \cong \Ext_{\bZ[G]}^i(\bZ, \pi_*(X))
	= H^i(G, \pi_*(X))
\]
identifying the $E^2$-page of the $G$-homotopy fixed point spectral
sequence with the group cohomology of the $G$-module $\pi_*(X)$, and
this identification is compatible with the cup product structure on
both sides.
\end{example}

\begin{example}
When $G = \bT$ is the circle group, we showed in
Proposition~\ref{prop:s2isetas} that
\[
R[\bT]_* = R_*[s]/(s^2 = \eta s) \och \overline{R[\bT]}_* = R_* \{s\} \,.
\]
As we discussed in Definition~\ref{def:P}, the normalized bar resolution
gives a (minimal) resolution $P_{*,*} = NB_*(R_*, R[\bT]_*, R[\bT]_*)$
of $R_*$, with
$$
P_{i,*} = \overline{R[\bT]}_*^{\otimes i} \otimes_{R_*} R[\bT]_*
	\cong R[\bT]_* \{\bar p_i \} \,.
$$
Here $\bar p_i = s \otimes \dots \otimes s \otimes 1 = [s | \dots | s]1$
has homological degree~$i$, internal degree~$|\bar p_i| = i$ and total
degree~$\|p_i\| = 2i$, for $i\ge0$.  The differential is given by
$$
\partial_i(\bar p_i) = \bar p_{i-1} ((i-1) \eta + (-1)^i s)
$$
for $i\ge1$.  This means that $P_{*,*}$ is not strictly
equal to the resolution $P_*$ specified at the beginning of
Section~\ref{subsec:computation}, with $P_i = R[\bT]_* \{p_i\}$ and
$\partial_i (p_i) = p_{i-1} (s + (i-1)\eta)$, due to the sign $(-1)^i$
before the contribution from the last face operator.  However, the two
resolutions are isomorphic, by way of the chain map sending $\bar p_i$
to $(-1)^{i(i+1)/2} p_i$ for each $i\ge0$.  Even without this isomorphism,
we are free to use $P_*$ to calculate $\Ext^*_{R[\bT]_*}(R_*, \pi_*(X))$
as the homology of $\Hom_{R[\bT]_*}(P_*, \pi_*(X))$, and that
calculation was essentially done in Section~\ref{subsec:computation}.
For each $b\ge0$ the rule
\[
x \longmapsto
f_b \cdot x := \binom{p_b \mapsto x}{p_b s \mapsto xs}
\]
defines a bijection $\Sigma^{-b} \pi_*(X) \cong \Hom_{R[\bT]_*}(\Hom(P_b,
\pi_*(X)))$, and the boundary on such an element is given by
\[
\partial^v(f_b \cdot x) = \begin{cases}
-(-1)^{|x|}f_{b+1} \cdot xs & \text{for $b \geq 0$ even,} \\
-(-1)^{|x|}f_{b+1} \cdot x(s+\eta) & \text{ for $b \geq 1$ odd.}
\end{cases} 
\]
Hence we can compute the homology as
\[
\Ext_{R[\bT]_*}^b(R_*, \pi_*(X))
\cong \begin{cases}
f_0 \cdot \ker(s \: \pi_*(X) \to \pi_{*+1}(X)) & \text{for $b=0$,} \\
\\
\ds f_b \cdot \frac{\ker(s+\eta \: \pi_*(X) \to \pi_{*+1}(X))}{\im(s \: \pi_{*-1}(X) \to \pi_*(X))} & \text{for $b\ge1$ odd,} \\
\\
\ds f_b \cdot \frac{\ker(s \: \pi_*(X) \to \pi_{*+1}(X))}{\im(s+\eta \: \pi_{*-1}(X) \to \pi_*(X))} & \text{for $b\ge2$ even.}
\end{cases}
\]
Please compare with Proposition~\ref{prop:Mbasis},
Lemma~\ref{lem:fboundary}, and
Proposition~\ref{prop:E1pageTatecomputation}.

For a description of the cup product, we can use any chain map $\Psi :
P_* \to P_* \otimes_{R_*} P_*$ lifting the identity on $R_*$. Such a map
is given in Lemma~\ref{lem:Psichainmap}, so that we can compute the cup
product as
\[
f_{b_1} \cdot x \smile f_{b_2} \cdot y = f_{b_1 + b_2} \cdot x \otimes y \,.
\]
Please compare with Lemma~\ref{lem:cupinExt}. Formally writing
the class of $f_b \cdot x$ as $t^b \cdot x$, we can then express
$\Ext_{R[\bT]_*}(R_*, \pi_*(X))$ as the homology of the differential
graded~$R[\bT]_*$-module
\[
\pi_*(X) [t]
\]
with differential given by $d(x) = t x s$ and $d(t) = t^2 \eta$, for
$x \in \pi_*(X)$.  Here $t$ has homological degree~$-1$, internal
degree~$|t| = -1$ and total degree~$\|t\| = -2$.
\end{example}

\section{The odd spheres filtration}
\label{subseq:oddspheres}
\index{odd spheres filtration}
In the important case $G = \bT$, the circle action on odd-dimensional spheres provides a pleasant alternative model for $EG$.  For each $i\ge0$ let $S(i\bC) = S^{2i-1}$ be the unit sphere in $i\bC = \bC^i$, with the standard, free~$\bT$-action. We obtain an exhaustive filtration
\[
\emptyset \subset S(\bC) \subset \dots \subset
	S(i\bC) \subset S((i+1)\bC) \subset \dots \subset S(\infty\bC)
\]
of free $\bT$-spaces. Here $S((i+1)\bC)$ is obtained from $S(i\bC)$ by attaching a free~$\bT$-equivariant $2i$-cell $D^{2i} \times \bT$ along the group action map 
\[
S^{2i-1} \times \bT \cong S(i\bC) \times \bT \to S(i\bC) \,,
\]
so that $S((i+1)\bC)$ is the $2i$-skeleton in a free $\bT$-CW structure on $S(\infty\bC)$. This filtered model for a free, contractible $\bT$-CW complex was used in \cite{BR05}*{\S2} to discuss the~$\bT$-homotopy fixed point spectral sequence.

There are well-known $\bT$-equivariant homeomorphisms
\[
S((i+1)\bC) \cong \bT * \dots * \bT * \bT
\]
with $(i+1)$ copies of $\bT$, where $*$ denotes the join of spaces. These homeomorphisms are compatible for varying $i\ge0$, and $S(\infty\bC)$ is isomorphic as a filtered space to Milnor's infinite join construction\index{Milnor's infinite join construction} from \cite{Mil56}, for $G = \bT$, which we denote by 
\[
\sE G = G * G * G * \dots \,.
\]
The identifications made in the
iterated join are included among those made in geometric realization.
Hence the structure map $\Delta^i \times G^i \times G \to EG$
factors through a $G$-map
$$
q_i \: G * \dots * G * G \longto F_i EG
$$
with $(i+1)$ copies of~$G$, collapsing degenerate simplices.  These are
compatible for varying~$i$, yielding a~$G$-map $q \: \sE G \to EG$.
As explained in \cite{Seg68}*{\S3}, the Milnor join construction is a
special case $\sE G \cong E G_{\bN}$ of the two-sided bar construction
for a topological category $G_{\bN}$, and there is a continuous functor
$G_{\bN} \to G$ inducing the~$G$-maps $q_i$ and~$q$.  It follows that
the filtration-preserving diagonal approximation\index{diagonal approximation} $D_{\bN} \: E G_{\bN}
\to E G_{\bN} \times E G_{\bN}$ constructed in \cite{Seg68}*{Lem.~5.4}
is compatible with the diagonal approximation $D \: EG \to EG \times EG$
that we have used in the present memoir.  In particular, the $\bT$-map
$$
q^* \: F(E\bT_+, X) \longto F(\sE\bT_+, X) \cong F(S(\infty \bC)_+, X)
$$
maps our multiplicative sequence $M_\star(X)$ to the multiplicative tower
used in \cite{BR05}*{\S4}. Furthermore, for~$G = \bT$ the $G$-maps $q_i$ and~$q$ are equivalences, so that the two multiplicative towers of orthogonal $G$-spectra are equivalent. Hence they give isomorphic $\bT$-homotopy fixed point spectral sequences, converging to the same multiplicative filtration on the abutment.

A similar discussion applies for the $3$-sphere $G = \bU = Sp(1)$ acting on the unit spheres in $i\bH = \bH^i$, showing that $S(\infty\bH) \cong \sE \bU$ is a perfectly good alternative filtered model for $E\bU$.
\index{homotopy fixed point spectral sequence|)}

%
%
%

\chapter{The $G$-Tate Spectral Sequence} \label{sec:Tatespseq}

Given an $R$-module $X$ in orthogonal $G$-spectra we can define its $G$-Tate construction as the genuine fixed points
\[
X^{tG} = (\widetilde{EG} \wedge F(EG_+,X))^G \,,
\]
where $\widetilde{EG}$ is the mapping cone of the collapse map $EG_+ \to S^0$. In this chapter we construct an $R_*$-module spectral sequence 
\[
\hat{E}^{r}_{*,*} \Longrightarrow \pi_*(X^{tG})
\]
with abutment the $G$-equivariant homotopy groups of~$\widetilde{EG} \wedge F(EG_+,X)$, for any compact Lie group $G$. We do this by letting the filtration $E_\star$ induce a filtration $\widetilde{E}_\star$ of $R \wedge \widetilde{EG}$ and consider the so-called Hesselholt--Madsen filtration
\[
HM_\star(X) = (\widetilde{E} \wedge_R T(M(X)))_\star \,
\]
obtained by forming a convolution product. Under the assumption that $R[G]_*$ is finitely generated and projective over $R_*$ we show that the resulting spectral sequence $\hat{E}^r_{*,*}(X) = E^r_{*,*}(HM_\star(X))$ is multiplicative, as a functor of~$X$, with multiplicative abutment. With the same assumptions we also algebraically identify the~$E^2$-page as
\[
\hat E^2_{*,*}(X) \cong \hatExt^{-*}_{R[G]_*}(R_*,\pi_*(X)) \,,
\]
with the multiplicative structure given by cup product on the right-hand side. See Theorem~\ref{thm:GTatespseqismult}. To say something about the convergence of this spectral sequence we compare the Hesselholt--Madsen filtration to another filtration $GM_\star(X)$ of $\widetilde{EG} \wedge F(EG_+,X)$, dubbed the Greenlees--May filtration. While the multiplicative properties of the Greenlees--May $G$-Tate spectral sequence are less clear, it is easy to obtain convergence results for the latter spectral sequence. By the comparison we can then also obtain convergence results for the Hesselholt--Madsen $G$-Tate spectral sequence. See Section~\ref{subsec:convergence}, and in particular Theorem~\ref{thm:GTatestrongconv}.

\section{The filtered $G$-space $\wEG$}

As always, let $G$ be any compact Lie group.  Let $c \: EG_+ \to S^0$ denote the based and $G$-equivariant collapse map, and define
\[
\wEG = S^0 \cup C(EG_+)
\]
to be its reduced mapping cone, as in \cite{Car84}*{p.~198} and
\cite{GM95}*{p.~2}. Non-equivariantly, $c$ is an equivalence, so $\wEG$ is (non-equivariantly) contractible. For $i\ge0$ we let
\[
F_i \wEG = S^0 \cup C(F_{i-1} EG_+)
\]
be the mapping cone of $c$ restricted to $F_{i-1} EG_+$, where $F_{i-1} EG$
is defined as in Section~\ref{subsec:filteredEG}. For $i<0$, we set $F_i \wEG = *$. This defines an exhaustive filtration
\begin{equation} \label{eq:FiwEG}
{*} = F_{-1} \wEG \subset S^0 = F_0 \wEG \subset \dots \subset
	F_{i-1} \wEG \subset F_i \wEG \subset \dots \subset \wEG
\end{equation}
of based $G$-spaces. Each map $F_{i-1} \wEG \to F_i \wEG$ is a strong $h$-cofibration, so this is indeed a filtration, as opposed to simply a sequence. Moreover, there are homeomorphisms
\[
\frac{F_i \wEG}{F_{i-1} \wEG} \cong \Sigma \frac{F_{i-1} EG}{F_{i-2} EG}
\]
for $i\ge1$. Per Theorem~\ref{thm:pushprodstrongcof}, each pushout-product map
\[
F_{i-1} \wEG \wedge F_j \wEG \cup F_i \wEG \wedge F_{j-1} \wEG
\longto F_i \wEG \wedge F_j \wEG
\]
is a strong $h$-cofibration, with cofibre
\[
\frac{F_i \wEG}{F_{i-1} \wEG} \wedge \frac{F_j \wEG}{F_{j-1} \wEG} \,.
\]

\begin{remark}
\label{rmk:EGtildefiltcomp}
When $G$ is finite, $F_i \wEG$ gives the $i$-skeleton of a based and non-free $G$-CW structure on~$\wEG$.  When $G = \bT = U(1)$, $F_0 \wEG = S^0$ is the $0$-skeleton, while~$F_i \wEG$ for $i\ge1$ is the $2i-1$-and~$2i$-skeleton of a $G$-CW structure on $\wEG$.  Similarly, when $G = \bU = Sp(1)$, $F_i \wEG$ gives the $4i-3$-, $4i-2$-, $4i-1$- and $4i$-skeleta of a $G$-CW structure.
\end{remark}

\begin{remark}
For $G = \bT$, the $G$-equivalences $q_{i-1} \: S(i\bC) \to F_{i-1} EG$ from Section~\ref{subseq:oddspheres} induce $G$-equivalences $\tilde q_i \: S^{i\bC} \to F_i \wEG$, where we identify the one-point compactification $S^{i\bC}$ with the
mapping cone $S^0 \cup C(S(i\bC)_+)$. Hence we have a $G$-equivalence
from the exhaustive filtration
\[
{*} \to S^0 \to \dots \to S^{(i-1)\bC} \to
	S^{i\bC} \to \dots \to S^{\infty\bC}
\]
to~\eqref{eq:FiwEG}, showing that we may use $S^{\infty\bC}$ as a filtered replacement for $\wEG$, if desired.
\end{remark}

We give $\wEG \wedge \wEG$ the (convolved) smash product filtration, with
\[
F_k(\wEG \wedge \wEG) = \bigcup_{i+j=k} F_i \wEG \wedge F_j \wEG \,.
\]
The identifications $S^0 \wedge \wEG \cong \wEG \cong \wEG \wedge S^0$
agree on $S^0 \wedge S^0 \cong S^0$, hence combine to a fold map\index{fold map}
\[
\nabla \: \wEG \cup_{S^0} \wEG
	\cong \wEG \wedge S^0 \cup S^0 \wedge \wEG \longto \wEG \,.
\]
We seek a $G$-map $N \: \wEG \wedge \wEG \to \wEG$
extending $\nabla$, so that the diagram
\[
\xymatrix{
\wEG \wedge S^0 \cup S^0 \wedge \wEG \ar[r] \ar[dr]_-{\nabla}
	& \wEG \wedge \wEG \ar[d]^-{N} \\
& \wEG
}
\]
commutes.  For these pairings to induce pairing of spectral sequences, we must arrange that $N$ is filtration-preserving.
We do not know how to give a direct definition of such an extension $N \: \wEG \wedge \wEG \to \wEG$, in analogy with the explicit diagonal approximation $D \: EG \to EG \times EG$. Instead we will use obstruction theory to show that such a filtration-preserving extension $N$ of $\nabla$ exists after base change to our ground ring spectrum~$R$, assuming that $R[G]_*$ is projective over~$R_*$. See Proposition~\ref{prop:nablaextendstoN}.

\begin{definition} \label{def:tildeEstar}
Let 
\[
\widetilde E = R \wedge \wEG \och \widetilde{E}_i = R \wedge F_i \wEG \,.
\]
Each map $\widetilde E_{i-1} \to \widetilde{E}_i$ is a strong $h$-cofibration, so that $\widetilde E_\star$ is a filtration
\[
\dots \longto \widetilde E_{i-1} \longto \widetilde E_i \longto \widetilde E_{i+1}
	\longto \dots
\]
of $R$-modules in orthogonal $G$-spectra.  Here $\widetilde E_i = *$ for
$i<0$, $\widetilde E_0 = R$, and 
\[
\widetilde E_\infty = \Tel(\widetilde E_\star) \simeq_G \widetilde{E} \,.
\]
Since $\wEG$ is non-equivariantly contractible, $\pi_*(\widetilde E) = 0$.
\end{definition}

Applying non-equivariant homotopy we obtain the following unrolled exact couple
\begin{equation} \label{eq:tildeEexactcouple}
\begin{tikzcd}
\cdots \ar[r]
	& \pi_*(\widetilde{E}_{i-1}) \arrow[r,"\alpha"]
	& \pi_*(\widetilde{E}_i) \ar[r] \arrow[d,"\beta"]
	& \cdots \\
& & \pi_*(\widetilde{E}_{i-1} \to \widetilde{E}_i) \arrow[ul,"\partial"]
\end{tikzcd}
\end{equation}
with $\partial$ of total degree~$-1$. Recall the $R[G]_*$-module resolution $(P_{*,*}, \partial)$ of $R_*$, introduced in Definition~\ref{def:P}.

\begin{definition} \label{def:tildeP}
Let $(\widetilde P_{*,*}, \widetilde\partial)$ be the mapping cone of the augmentation 
\[
\epsilon \: P_{*,*} \to R_* \,,
\]
in the sense of Definition~\ref{def:mappingcone}.
\end{definition}

Explicitly, we have 
\[
\widetilde{P}_{i,*} \cong \begin{cases} R_* & \text{for $i=0$}\\ P_{i-1,*} & \text{for $i \geq 1$} \end{cases}
\] 
with boundary $\tilde{\partial} \: \widetilde{P}_{i,\ast} \to \widetilde{P}_{i-1,\ast}$ given as 
\[
\tilde{\partial} = \begin{cases}
\epsilon(x) & \text{for $i=1$} \\
- \partial(x) & \text{for $i \geq 2$.}
\end{cases}
\]
We note that $\widetilde{P}_{*,*}$ is an exact complex of flat $R_*$-modules, by our standing assumption that $R[G]_*$ is flat.  If, furthermore, $R[G]_*$ is finitely generated projective over~$R_*$, then so is each $\widetilde{P}_{i,*}$.

\begin{lemma} \label{lem:tildeE1istildeP}
If $R[G]_*$ is flat over $R_*$, then the $(E^1, d^1)$-page of the non-equi\-variant homotopy spectral sequence
\[
\widetilde E^1_{i,*} = \pi_{i+*}(\widetilde{E}_{i-1} \to \widetilde{E}_i)
\]
associated to $\widetilde{E}_\star$ is isomorphic to $(\widetilde {P}_{*,*},\tilde\partial)$. In particular, the spectral sequence collapses at the $E^2$-page, where it is given by 
\[
\widetilde E^2 = \widetilde E^\infty = 0 \,.
\]
\end{lemma}

\begin{proof}
Note that $\widetilde E$ is the mapping cone of the collapse map $1 \wedge c \: E \to R$ and can be viewed as the pushout
\[
\xymatrix{
E \ar[r] \ar[d]_-{1 \wedge c}
	& I \wedge E \ar[d] \\
R \ar[r] & \widetilde E \rlap{\,.}
}
\]
Let $I_\star$ be the filtration
\[
\ast \longto \{0,1\} \longto I \overset{=}\longto I \overset{=}\longto I \overset{=}\longto \cdots
\]
of the unit interval $I = [0,1]$, where $\partial I = \{0,1\}$ sits in filtration degree $0$. Let~$L_\star(R)$ be the non-negative filtration consisting of copies of $R$ and identity maps between them.  We then have a pushout of filtrations
\begin{equation} \label{eq:tildeEstarpushout}
\begin{tikzcd}
E_\star \arrow[r] \arrow[d,"1 \wedge c"]
	& (I \wedge E)_\star \arrow[d] \\
L_\star(R) \arrow[r] & \widetilde{E}_\star
\end{tikzcd}
\end{equation}
with colimit being the pushout square above. That this is indeed a pushout of filtrations can be checked in each filtration degree separately, noting that 
\[
(I \wedge E)_k = \partial I \wedge E_k \cup I \wedge E_{k-1} \cong E_k \cup CE_{k-1} \,.
\]
It follows as in Lemma~\ref{lem:E1isP} that we have a commutative square of associated chain complexes
\begin{equation} \label{eq:tildeE1pushout}
\xymatrix{
P_{*,*} \ar[r]^-{\cong} \ar[d]_-{\epsilon}
	& \partial I_* \otimes P_{*,*} \ar@{ >->}[r]
	& I_* \otimes P_{*,*} \ar[d] \\
R_* \ar@{ >->}[rr] && \widetilde E^1_{*,*} \rlap{\,.}
}
\end{equation}
Here $R_*$ is the chain complex consisting of $R_*$ concentrated in homological degree~0, and $I_*$ is the reduced cellular chain complex
\[
0 \longto \bZ\{i_1\} \overset{\partial_1}\longto 
	\bZ\{i_0\} \longto 0 
\]
of $I$, with $\partial_1(i_1) = i_0$. Both $i_0$
and $i_1$ have internal degree~$0$, and lie in homological degree as indicated by their subscript. The chain complex $\partial I_*$ is the subcomplex given by $\bZ\{i_0\}$ concentrated in homological degree~$0$.  Since the map
\[
P_{*,*} \cong \partial I_* \otimes P_{*,*}
 \longto I_* \otimes P_{*,*}
\]
is injective, a Mayer--Vietoris argument for the filtration subquotients of~\eqref{eq:tildeEstarpushout} shows that~\eqref{eq:tildeE1pushout} is in fact a pushout of chain
complexes. This proves that $\widetilde E^1_{*,*}$ is indeed the algebraic mapping cone of $\epsilon \: P_{*,*} \to R_*$, by the definition of the latter chain complex.
\end{proof}

\begin{lemma} \label{lem:tildeE1istildePtildeP}
The $(E^1,d^1)$-page of the non-equivariant homotopy spectral sequence
associated to $(\widetilde E \wedge_R \widetilde E)_\star$ is isomorphic to
$(\widetilde{P}_{*,*} \otimes_{R_*} \widetilde{P}_{*,*}, \tilde\partial \otimes 1 +
1 \otimes \tilde\partial)$.
\end{lemma}

\begin{proof}
This is very similar to Lemma~\ref{lem:E1isPP}.
\end{proof}

\begin{lemma} \label{lem:alphazero}
The homomorphism $ \pi_*(\widetilde{E}_{i-1}) \to \pi_*(\widetilde{E}_i)$
is zero, for each~$i$.
\end{lemma}

\begin{proof}
This follows from the exactness of $(\widetilde E^1_{*,*}, d^1) \cong (\widetilde{P}_{*,*}, \tilde\partial)$, by an induction on~$i$ in the unrolled exact couple~\eqref{eq:tildeEexactcouple}. The claim is clear for $i\le0$. Assume by induction that $\alpha \: \pi_*(\widetilde{E}_{i-1}) \to \pi_*(\widetilde{E}_i)$ is zero, for some $i\ge0$.  Then $\beta
\: \pi_{i+*}(\widetilde{E}_i) \to \widetilde{P}_{i,*}$ is injective.  Consider any
class $x \in \pi_{i+*}(\widetilde{E}_i)$.  Since $\tilde\partial_i(\beta(x))
= \beta \partial \beta(x) = 0$, exactness at $\widetilde{P}_{i,*}$ implies
that $\beta(x) = \tilde\partial_{i+1}(y) = \beta \partial(y)$ for some $y \in \widetilde P_{i+1,*}$. By injectiveness of~$\beta$ it follows that $x = \partial(y)$. Since $x$ was arbitrary, $\partial \: \widetilde{P}_{i+1,*} \to \pi_{i+*}(\widetilde{E}_i)$ is
surjective, so $\alpha \: \pi_*(\widetilde{E}_i) \to \pi_*(\widetilde{E}_{i+1})$
is zero.
\end{proof}

\begin{lemma} \label{lemma:Nunique}
There always exists an $R$-module map of orthogonal $G$-spectra
\[
N \: \widetilde{E} \wedge_R \widetilde{E} \longto \widetilde{E}
\]
extending $\nabla \: \widetilde{E} \cup_R \widetilde{E} \to \widetilde{E}$, and any two choices are homotopic.
\end{lemma}

\begin{proof}
This follows by obstruction theory, since
\[
\wEG \cup \wEG \cong \wEG \wedge S^0 \cup S^0 \wedge \wEG
	\subset \wEG \wedge \wEG
\] 
can be given the structure of a free relative $G$-CW complex,
and $\pi_*(\widetilde{E}) = 0$.
\end{proof}

The above lemma, together with the map $\Delta_+ \: EG_+ \to EG_+ \wedge EG_+$, makes sure that the Tate construction is multiplicative, in the sense that a $G$-equivariant $R$-module pairing $X \wedge_R Y \to Z$ induces an $R$-module pairing $X^{tG} \wedge_R Y^{tG} \to Z^{tG}$. See Section~\ref{subsec:GTatecon}. To arrange that the Tate spectral sequence preserves this structure we need to make sure that we can find a filtration-preserving approximation of $N$, in the same way as we could find the filtration-preserving approximation of $D$. The following proposition addresses difficulties raised in Problem~11.8 and Problem~14.8 of \cite{GM95}.

\begin{proposition} \label{prop:nablaextendstoN}\index{fold map}
Suppose that $R[G]_*$ is projective over~$R_*$.  Then there exists a filtration-preserving map
\[
N \: (\widetilde{E} \wedge_R \widetilde{E})_\star \longto \widetilde{E}_\star
\]
of $R$-modules in orthogonal $G$-spectra, extending the fold map
\[
\nabla \: \widetilde{E}_\star \cup_R \widetilde{E}_\star \cong
(\widetilde{E}_\star \wedge_R R) \cup (R \wedge_R \widetilde{E}_\star)
	\longto \widetilde{E}_\star \,.
\]
\end{proposition}

\begin{proof}
We inductively assume that $\nabla$ has been extended to a
filtration-preserving map $N_{k-1} \: (\widetilde{E} \wedge_R \widetilde{E})_{k-1} \to \widetilde{E}_{k-1}$, and show that $N_{k-1}$ can be further extended to a filtration-preserving map $N_k \: (\widetilde{E} \wedge_R \widetilde{E})_k \to \widetilde{E}_k$. It suffices to extend $N_{k-1}$ over~$\widetilde{E}_i \wedge_R \widetilde{E}_j$ for $i,j\ge1$ with $i+j=k$. In particular, there is only something to prove for $k\ge2$. Let us consider the diagram
\[
\xymatrix{
\widetilde{E}_{i-1} \wedge_R \widetilde{E}_j \cup
	\widetilde{E}_i \wedge_R \widetilde{E}_{j-1} \ar[r]^-{N_{k-1}} \ar@{ >->}[d]
	& \widetilde{E}_{k-1} \ar[d]^-{\alpha} \\
\widetilde{E}_i \wedge_R \widetilde{E}_j \ar@{-->}[r]^-{N_{i,j}} \ar[d]
	&  \widetilde{E}_k \\
\widetilde{E}_i/\widetilde{E}_{i-1} \wedge_R \widetilde{E}_j/\widetilde{E}_{j-1}
}
\]
where the left hand column is a (Hurewicz) cofibre sequence. By the homotopy extension property, in order to find a dashed map $N_{i,j}$ making the diagram commute, it suffices to find an extension up to homotopy of $\alpha \circ N_{k-1}$. Let
\[
W = E_{i-1}/E_{i-2} \wedge_R E_{j-1}/E_{j-2}
\cong R \wedge G^{\wedge i-1} \wedge G_+ \wedge G^{\wedge j-1} \wedge G_+
\]
so that $\Sigma^2 W \cong \widetilde{E}_i/\widetilde{E}_{i-1} \wedge_R \widetilde{E}_j/\widetilde{E}_{j-1}$. There is then a (stably defined) homotopy cofibre sequence
\[
\Sigma W
	\overset{\partial}\longto \widetilde{E}_{i-1} \wedge_R \widetilde{E}_j
		\cup \widetilde{E}_i \wedge_R \widetilde{E}_{j-1}
	\longto \widetilde{E}_i \wedge_R \widetilde{E}_j
	\longto \Sigma^2 W
\]
and it suffices to prove that $\alpha \circ N_{k-1} \circ \partial \:
\Sigma W \to \widetilde{E}_k$ is null-homotopic.  We confirm this by showing
that $\alpha$ induces the trivial homomorphism
\[
\alpha_* \: [\Sigma W, \widetilde{E}_{k-1}]^G_R
	\longto [\Sigma W, \widetilde{E}_k]^G_R \,,
\]
where $[-,-]^G_R$ denotes homotopy classes of $G$-maps of $R$-modules in
orthogonal $G$-spectra.  Note that $G$ acts
diagonally on the two copies of $G_+$ in $W$, so that there is
an untwisting isomorphism $W \cong  V \wedge G_+$ where
\[
V = R \wedge G^{\wedge i-1} \wedge G^{\wedge j-1} \wedge G_+
\]
has trivial $G$-action. By adjunction we can therefore rewrite
the homomorphism above as
\[
\alpha_* \: [\Sigma V, \widetilde{E}_{k-1}]_R
	\longto [\Sigma V, \widetilde{E}_k]_R
\]
where $[-,-]_R$ denotes homotopy classes of maps of (non-equivariant)
$R$-modules.  By our assumption that $R[G]_*$ is $R_*$-projective,
it follows that
\[
\pi_*(V) \cong \overline{R[G]}_*^{\otimes i-1}
	\otimes_{R_*} \overline{R[G]}_*^{\otimes j-1} \otimes_{R_*} R[G]_*
\]
is $R_*$-projective. Hence we can rewrite $\alpha_*$ as the homomorphism
\[
\Hom_{R_*}(\Sigma \pi_*(V), \pi_*(\widetilde{E}_{k-1}))
\longto
\Hom_{R_*}(\Sigma \pi_*(V), \pi_*(\widetilde{E}_k))
\]
given by composition with $\alpha \: \pi_*(\widetilde{E}_{k-1})
\to \pi_*(\widetilde{E}_k)$.  By Lemma~\ref{lem:alphazero} that
homomorphism is zero, which completes the proof.
\end{proof}

\begin{definition} \label{def:Phi}
Suppose that $R[G]_*$ is projective over~$R_*$.  Let
\[
\Phi \: \widetilde{P}_{*,*} \otimes_{R_*} \widetilde{P}_{*,*}
	\longto \widetilde{P}_{*,*}
\]
be the $R[G]_*$-module chain map that corresponds, under the isomorphisms of Lemma \ref{lem:tildeE1istildeP} and Lemma~\ref{lem:tildeE1istildePtildeP}, to the pairing $N^1$ of $(E^1, d^1)$-pages induced by the filtration-preserving map $N \: (\widetilde{E} \wedge_R \widetilde{E})_\star \to \widetilde{E}_\star$ of Proposition~\ref{prop:nablaextendstoN}.
\end{definition}

\begin{lemma}
Suppose that $R[G]_*$ is projective over~$R_*$. Then the map
\[
\Phi \: \widetilde{P}_{*,*} \otimes_{R_*} \widetilde{P}_{*,*}
\to \widetilde{P}_{*,*}
\]
is uniquely characterized, up to $R[G]_*$-module chain homotopy, by being an~$R[G]_*$-module chain map that extends the fold map $\nabla$.
\end{lemma}

\begin{proof}
By construction, $\Phi$ extends the fold map, and it follows that this map is unique up chain homotopy equivalence by Proposition~\ref{prop:tildePfoldextension}.
\end{proof}

\section{The $G$-Tate construction}
\label{subsec:GTatecon}

Let $X$ be an $R$-module in orthogonal $G$-spectra\footnote{If $X$ is an orthogonal $G$-spectrum without $R$-action, the discussion in this section applies to $R \wedge X$ in place of $X$.}. In this section, we discuss the Tate construction and its multiplicative properties.

\begin{definition}\index{Tate construction}
The \emph{$G$-Tate construction} $X^{tG}$ is the $G$-fixed
point spectrum of (a fibrant replacement of) $\wEG \wedge F(EG_+, X)$:
\[
X^{tG} = \left( \wEG \wedge F(EG_+, X) \right)^G
\]
\end{definition}

Note that the homotopy groups
\[
\pi_*(X^{tG}) \cong \pi^G_*(\wEG \wedge F(EG_+, X))
\]
naturally form an $R_*$-module, and that we can write
\[
\wEG \wedge F(EG_+, X) \cong \widetilde{E} \wedge_R F_R(E, X) \,.
\]
The inclusion $S^0 \to \wEG$ induces a $G$-map
\[
F(EG_+, X) \cong S^0 \wedge F(EG_+, X) \longto \wEG \wedge F(EG_+, X)
\]
and a map of their corresponding $G$-fixed points $X^{hG} \longto X^{tG}$. We can write these as maps of $R$-modules, using the inclusion $R \to \widetilde E$, to obtain a $G$-map
\[
F_R(E, X) \cong R \wedge_R F_R(E, X) \longto \widetilde E \wedge_R F_R(E, X)
\]
and a canonical map
\[
\gamma \: X^{hG} = F_R(E, X)^G
	\longto \left(\widetilde E \wedge_R F_R(E, X)\right)^G = X^{tG} \,,
\]
inducing a homomorphism $\pi_*(X^{hG}) \to \pi_*(X^{tG})$ of $R_*$-modules.

The Tate construction interacts well with the multiplicative structure on homotopy fixed points we described in the paragraph preceding Remark~\ref{rmk:fixpointPostnikov}. Note first that given a pairing $\mu \: X \wedge_R Y \to Z$ of $R$-modules in orthogonal $G$-spectra, the $R$-module pairing $X^{hG} \wedge_R Y^{hG} \to Z^{hG}$ extends to $R$-module pairings $X^{tG} \wedge_R Y^{hG} \to Z^{tG}$ and $X^{hG} \wedge_R Y^{tG} \to Z^{tG}$. The first is given by a composite
\begin{align*}
\left(\widetilde E \wedge_R F_R(E, X)\right)^G \wedge_R F_R(E, Y)^G
	&\xrightarrow{\wedge}
\left(\widetilde E \wedge_R F_R(E, X) \wedge_R F_R(E, Y)\right)^G \\
	&\xrightarrow{1 \wedge \alpha}
\left(\widetilde E \wedge_R F_R(E \wedge_R E, X \wedge_R Y)\right)^G \\
	&\xrightarrow{1 \wedge (1 \wedge \Delta_+)^*}
\left(\widetilde E \wedge_R F_R(E, X \wedge_R Y)\right)^G \\
	&\xrightarrow{1 \wedge \mu_*}
\left(\widetilde E \wedge_R F_R(E, Z)\right)^G \,.
\end{align*}
The second one is similar, and left to the reader. The two pairings induce $R_*$-module pairings $\pi_*(X^{tG}) \otimes_{R_*} \pi_*(Y^{hG}) \to \pi_*(Z^{tG})$ and $\pi_*(X^{hG}) \otimes_{R_*} \pi_*(Y^{tG}) \to \pi_*(Z^{tG})$. These pairings are all compatible via the canonical map, meaning that the $R$-module diagram
\[
\xymatrix{
X^{tG} \wedge_R Y^{hG} \ar[d]
	& X^{hG} \wedge_R Y^{hG} \ar[l]_-{\gamma\wedge1}
		\ar[d] \ar[r]^-{1\wedge\gamma}
	& X^{hG} \wedge_R Y^{tG} \ar[d] \\
Z^{tG}
	& Z^{hG} \ar[l]_-{\gamma} \ar[r]^-{\gamma}
	& Z^{tG} \,,
}
\]
and the induced $R_*$-module diagram both commute. Per Lemma~\ref{lemma:Nunique} we can choose a unique (up to homotopy) extension $N \: \widetilde E \wedge_R \widetilde E \to \widetilde E$ of the fold map $\nabla \: \widetilde E \cup_R \widetilde E \to \widetilde E$, in the category of~$R$-modules in orthogonal $G$-spectra. We can then promote the two $R$-module pairings to an $R$-module pairing $X^{tG} \wedge_R Y^{tG} \to Z^{tG}$, given by the composite
\begin{align*}
\left(\widetilde E \wedge_R F_R(E, X)\right)^G &\wedge_R \left(\widetilde E \wedge_R F_R(E, Y)\right)^G \\
	&\xrightarrow{\wedge}
\left(\widetilde E \wedge_R F_R(E, X) \wedge_R \widetilde E \wedge_R F_R(E, Y)\right)^G \\
	&\xrightarrow{1 \wedge \tau \wedge 1}
\left(\widetilde E \wedge_R \widetilde E \wedge_R F_R(E, X) \wedge_R F_R(E, Y)\right)^G \\
	&\xrightarrow{1 \wedge 1 \wedge \alpha}
\left(\widetilde E \wedge_R \widetilde E \wedge_R F_R(E \wedge_R E, X \wedge_R Y)\right)^G \\
	&\xrightarrow{N \wedge (1 \wedge \Delta_+)^*}
\left(\widetilde E \wedge_R F_R(E, X \wedge_R Y)\right)^G \\
	&\xrightarrow{1 \wedge \mu_*}
\left(\widetilde E \wedge_R F_R(E, Z)\right)^G \,.
\end{align*}
These pairings are also compatible via the canonical map, meaning that
the $R$-module diagram
\[
\xymatrix{
X^{tG} \wedge_R Y^{hG} \ar[r]^-{1\wedge\gamma} \ar[dr]
	& X^{tG} \wedge_R Y^{tG}
		\ar[d]
	& X^{hG} \wedge_R Y^{tG} \ar[l]_-{\gamma\wedge1} \ar[dl] \\
	& Z^{tG}
	& 
}
\]
and the induced $R_*$-module diagram both commute.  Taken together,
these diagrams show that 
\[
\gamma \: X^{hG} \to X^{tG} \och \gamma_* \:
\pi_*(X^{hG}) \to \pi_*(X^{tG})
\]
are multiplicative. We would now like to access $\pi_*(X^{tG})$ and the pairings above through filtrations and their associated spectral sequences.

\section{The Hesselholt--Madsen filtration}
\label{subsec:HMfilt}
\index{Hesselholt--Madsen filtration|(}
We can now generalize the filtration of $X^{tG}$ from \cite{HM03}*{\S4.3} to the case of compact Lie groups~$G$.

\begin{definition} 
Let
\[
HM_\star(X) = (\widetilde E \wedge_R T(M(X)))_\star
\]
be the filtration
\[
\dots \to HM_{k-1}(X) \to HM_k(X) \to HM_{k+1}(X) \to \dots
\]
of $R$-modules in orthogonal $G$-spectra given by the Day convolution product of the filtrations $\widetilde{E}_\star$ and~$T_\star(M(X))$.
\end{definition}

Recall: we introduced the filtration $\widetilde{E}_\star$ in Definition~\ref{def:tildeEstar}, the sequence~$M_\star(X)$ in Section~\ref{subsec:htpyfix}, and its telescopic approximation~$T_\star(M(X))$ in Section~\ref{subsec:filt}. The convolution product of $\widetilde{E}_\star$ and $T_\star(M(X))$ was defined in Section~\ref{subsec:convprod}, and is a filtration by Proposition~\ref{prop:convoffiltisfilt}.
We can realize
\[
HM_k(X) = \bigcup_{i+j=k} \widetilde{E}_i \wedge_R T_j(M(X))
\]
as a subspectrum of $\widetilde E \wedge_R \Tel(M_\star(X))$. The
structure maps $HM_{k-1}(X) \to HM_k(X)$ are then inclusions of
subspectra. These are (strong) $h$-cofibrations, so the canonical map
\[
\Tel(HM_\star(X)) \longto \colim_k HM_k(X)
	= \widetilde E \wedge_R \Tel(M_\star(X))
\]
is an equivalence.  Since $M_j(X) = F_R(E, X)$ for all $j\ge0$
there is a deformation retraction
\[
\Tel(M_\star(X)) \overset{\simeq_G}\longto F_R(E, X)
\]
and a further equivalence
\[
\widetilde E \wedge_R \Tel(M_\star(X))
	\overset{\simeq_G}\longto \widetilde E \wedge_R F_R(E, X)
	\cong \wEG \wedge F(EG_+, X) \,.
\]

\begin{definition} \label{def:hatErX}\index{Hesselholt--Madsen--Tate spectral sequence}
Let $X$ be an $R$-module in orthogonal $G$-spectra. We define the~\emph{$G$-Tate spectral sequence} for $X$ to be the $R_*$-module spectral sequence $(\hat E^r(X), d^r)$ associated to the filtration $HM_\star(X)$ with
\[
\hat E^r(X) = E^r(HM_\star(X))
\]
for each $r\ge1$.
\end{definition}

The abutment of the $G$-Tate spectral sequence for $X$ is the colimit
\[
A_\infty(HM_\star(X)) \cong \pi^G_* \Tel(HM_\star(X))
	\cong \pi^G_*(\widetilde E \wedge_R F_R(E, X))
	\cong \pi_*(X^{tG}) \,,
\]
filtered by the image submodules
\[
F_k \pi_*(X^{tG})
	= \im\big(\pi^G_*(HM_k(X)) \to \pi^G_* \Tel(HM_\star(X))
		\cong \pi_*(X^{tG})\big) \,.
\]

\begin{remark} \label{rmk:noweakconvergence}
In general, we do not claim that the $G$-Tate spectral sequence converges to the stated abutment, neither in the conditional nor in the weak sense. As we recalled in Section~\ref{subsec:seq}, conditional convergence to the colimit holds if $\holim_k HM_k(X) \simeq_G *$.  The latter condition would follow from an interchange of homotopy colimits and homotopy limits. More precisely, for each $a\ge0$ and integer~$k$, consider the subspectrum
\[
S_{a,k} = \bigcup_{\substack{i+j=k\\i\le a}}
	\widetilde{E}_i \wedge_R T_j(M(X))
\]
of $\widetilde E \wedge_R \Tel(M_\star(X))$. Then $\hocolim_a S_{a,k}
\simeq_G HM_k(X)$, and the sufficient condition $\holim_k HM_k(X)
\simeq_G *$ for conditional convergence is equivalent to
\begin{equation} \label{eq:limcolimS}
\holim_k \hocolim_a S_{a,k} \simeq_G * \,.
\end{equation}
On the other hand, $\holim_j \widetilde{E}_i \wedge_R T_j(M(X)) \simeq_G *$
for each~$i$, since $F_i \wEG$ is a finite~$G$-CW space.  It follows
by induction that $\holim_k S_{a,k} \simeq_G *$ for each finite~$a$, which implies that
\begin{equation} \label{eq:colimlimS}
\hocolim_a \holim_k S_{a,k} \simeq_G * \,.
\end{equation}
Without further hypotheses we do not see how to deduce~\eqref{eq:limcolimS}
from~\eqref{eq:colimlimS}.
(However, for $G = \bT$ see~\cite{BM17}*{Lemma~3.16}.)
\end{remark}
\index{Hesselholt--Madsen filtration|)}

\section{Algebraic description of $\hat E^1$ and $\hat E^2$}

Under the assumption that~$R[G]_*$ is finitely generated projective over $R_*$, we can algebraically describe the $E^1$- and~$E^2$-pages of the $G$-Tate spectral sequence, in the same way as we did for the $G$-homotopy fixed point spectral sequence in Section~\ref{sec:algdescriptionhfpss}.

\begin{proposition} \label{prop:TateE1isTatecx}
Suppose that $R[G]_*$ is $R_*$-projective.  There is then a natural isomorphism of $R_*$-module chain complexes
\[
E^1_{*,*}(HM_\star(X)) \cong \Hom_{R[G]_*}(R_*,
	\widetilde{P}_{*,*} \otimes_{R_*} \Hom_{R_*}(P_{*,*}, \pi_*(X)) ) \,.
\]
In the notation of Definition~\ref{def:Tatecx} and Definition~\ref{def:hatErX}, we have
\[
\hat E^1_{*,*}(X) \cong \Hom_{R[G]_*}(R_*, \hm_*(\pi_*(X))) \,,
\]
where the $d^1$-differential on the left hand side corresponds to $\Hom(1,\partial_{\hm})$ on the right hand side.
\end{proposition}

\begin{proof}
We first check that the natural restriction homomorphism
\begin{equation} \label{eq:omegaforHMX}
\omega \: E^1_{*,*}(HM_\star(X))
	\overset{\cong}\longto \Hom_{R[G]_*}(R_*,
	E^1_{*,*}((\widetilde E \wedge_R T(M(X)))_\star))
\end{equation}
from Lemma~\ref{lem:omegaXnat}
is an isomorphism of $R_*$-module chain complexes, where $HM_\star(X)$
at the left hand side is treated as an~$R$-module filtration in
orthogonal $G$-spectra, while~$(\widetilde E \wedge_R T(M(X)))_\star$ at the right hand side refers to the underlying $R$-module filtration in non-equivariant orthogonal spectra, with the residual $R[G]$-module action. We first note that we have
\begin{align*}
E^1_{k,*}(HM_\star(X))
	&= \pi^G_{k+*}(HM_{k-1}(X) \to HM_k(X)) \\
	&\cong \pi^G_{k+*}(HM_k(X) / HM_{k-1}(X))
\end{align*}
while
\begin{align*}
E^1_{k,*}((\widetilde E \wedge_R T(M(X)))_\star)
	&= \pi_{k+*}(HM_{k-1}(X) \to HM_k(X)) \\
        &\cong \pi_{k+*}(HM_k(X) / HM_{k-1}(X)) \,.
\end{align*}
Secondly, we note that
\begin{align*}
HM_k(X) / HM_{k-1}(X)
	&\cong \bigvee_{i+j=k} \widetilde{E}_i / \widetilde{E}_{i-1}
	\wedge_R T_j(M(X)) / T_{j-1}(M(X)) \\
	&\cong \bigvee_{i+j=k} \widetilde{E}_i / \widetilde{E}_{i-1}
	\wedge_R (M_j(X) \cup CM_{j-1}(X)) \\
	&\simeq_G \bigvee_{i+j=k} \widetilde{E}_i / \widetilde{E}_{i-1}
	\wedge_R F_R(E_{-j}/E_{-j-1}, X) \,,
\end{align*}
which is moreover $G$-equivalent to $F(G_+, X')$ for
\[
X' = \bigvee_{i+j=k} \widetilde{E}_i / \widetilde{E}_{i-1}
        \wedge \Sigma^j F(G^{\wedge -j}, X) \,.
\]
This uses that $R[G]$ is dualisable (see Definition~\ref{def:dual}). Proposition~\ref{prop:omegaXiso} therefore implies that the natural restriction homomorphism~(\ref{eq:omegaforHMX}) is an isomorphism in every homological degree. We check that it is also an isomorphism of chain complexes. The $d^1$-differential in the spectral sequence appearing at the left hand side corresponds to the composition
\begin{align*}
\pi^G_{k+*}(HM_k(X) / HM_{k-1}(X))
	&\overset{\partial}\longto \pi^G_{k-1+*}(HM_{k-1}(X)) \\
	&\longto \pi^G_{k-1+*}(HM_{k-1} / HM_{k-2}(X))) \,,
\end{align*}
which by the naturality of $\omega$ corresponds to $\Hom(1, d^1_{k,*})$, where $d^1_{k,\ast}$ is the $d^1$-differential in the spectral sequence appearing at the right hand side. This is given by the composite
\begin{align*}
\pi_{k+*}(HM_k(X) / HM_{k-1}(X))
	&\overset{\partial}\longto \pi_{k-1+*}(HM_{k-1}(X)) \\
	&\longto \pi_{k-1+*}(HM_{k-1} / HM_{k-2}(X))) \,.
\end{align*}
Hence~\eqref{eq:omegaforHMX} is indeed an isomorphism of chain complexes.

We now want to identify $E^1_{*,*}((\widetilde E \wedge_R T(M(X)))_\star)$ with the complex $\hm_*(\pi_*(X))$. For this aim, we can use the canonical pairing
\[
\iota \: (\widetilde{E}_\star, T_\star(M(X)))
	\longto (\widetilde E \wedge_R T(M(X)))_\star
\]
of $R$-module filtrations to obtain an $R[G]_*$-module chain map of the associated~$E^1$-pages
\[
\iota^1 \: E^1(\widetilde{E}_\star) \otimes_{R_*} E^1(T_\star(M(X)))
	\longto E^1((\widetilde E \wedge_R T(M(X)))_\star)
\]
as in Theorem~\ref{thm:mainmultII}, but in the non-equivariant setting. This map is the direct sum of the maps
\begin{multline*}
\pi_{i+*}(\widetilde{E}_i/\widetilde{E}_{i-1})
\otimes_{R_*}
\pi_{j+*}(T_j(M(X))/T_{j-1}(M(X))) \\
\overset{\cdot}\longto
\pi_{k+*}(\widetilde{E}_i/\widetilde{E}_{i-1} \wedge_R T_j(M(X))/T_{j-1}(M(X)))
\end{multline*}
for $i+j=k$. Note that each one of these maps is an isomorphism, because $\widetilde{P}_{i,*} = \pi_{i+*}(\widetilde{E}_i/\widetilde{E}_{i-1})$ is projective, hence flat, over~$R_*$. We conclude that $\iota^1$ is an isomorphism of $R[G]_*$-chain complexes, and thus induces an isomorphism
\begin{multline}
\Hom(1, \iota^1) \:
\Hom_{R[G]_*}(R_*, E^1(\widetilde{E}_\star) \otimes_{R_*} E^1(T_\star(M(X))))
\\
	\overset{\cong}\longto
\Hom_{R[G]_*}(R_*, E^1((\widetilde E \wedge_R T(M(X)))_\star))
\end{multline}
of $R_*$-module complexes.

The equivalence $\epsilon \: T_\star(M(X)) \to M_\star(X)$ induces an isomorphism
$$
\epsilon \: E^1(T_\star(M(X))) \overset{\cong}\longto E^1(M_\star(X))
$$
of $R[G]_*$-module chain complexes, which in turn induces an isomorphism
\begin{multline}
\Hom(1, 1 \otimes \epsilon) \:
\Hom_{R[G]_*}(R_*, E^1(\widetilde{E}_\star) \otimes_{R_*} E^1(T_\star(M(X))))
\\
	\overset{\cong}\longto
\Hom_{R[G]_*}(R_*, E^1(\widetilde{E}_\star) \otimes_{R_*} E^1(M_\star(X))) \,,
\end{multline}
of $R_*$-module chain complexes. Finally, we have
\[
E^1_{j,*}(M_\star(X)) \cong \pi_{j+*}(F_R(E_{-j} / E_{-j-1}, X))
\cong \Hom_{R_*}(P_{-j,*}, \pi_*(X))
\]
as $R[G]_*$-modules, because $\pi_{-j+*}(E_{-j} / E_{-j-1}) \cong
P_{-j,*}$ is $R_*$-projective. The $d^1$-differentials correspond
to $\tilde\partial$ and $\Hom(\partial, 1)$ by the argument in the proof of Proposition~\ref{prop:E1d1mu1fromPPsi}. Hence we have an isomorphism
\begin{multline} \label{eq:Tatecomplex}
\Hom_{R[G]_*}(R_*, E^1_{*,*}(\widetilde{E}_\star)
		\otimes_{R_*} E^1_{*,*}(M_\star(X))) \\
	\cong \Hom_{R[G]_*}(R_*,
	\widetilde{P}_{*,*} \otimes_{R_*} \Hom_{R_*}(P_{*,*}, \pi_*(X)))
\end{multline}
of $R_*$-module complexes.

When strung together, the numbered isomorphisms~\eqref{eq:omegaforHMX} through~\eqref{eq:Tatecomplex} establish the asserted identification of the $G$-Tate spectral sequence $(E^1, d^1)$-page for the orthogonal $G$-spectrum~$R$-module~$X$ with the Tate complex for the $R[G]_*$-module~$\pi_*(X)$.
\end{proof}

\begin{theorem} \label{thm:TateE2isTateExt}\index{complete $\Ext$}
Let $X$ be an $R$-module in orthogonal $G$-spectra, and
suppose that~$R[G]_*$ is $R_*$-projective. Then there is a
natural $R_*$-module isomorphism
\[
\hat E^2_{i,*}(X) = E^2_{i,*}(HM_\star(X)) \cong
	\hatExt_{R[G]_*}^{-i}(R_*, \pi_*(X)) \,,
\]
for each integer~$i$.
\end{theorem}

\begin{proof}
This is immediate by passage to homology from
Proposition~\ref{prop:TateE1isTatecx}. Here we are using the definition of Hopf
algebra Tate cohomology given in Definition~\ref{def:completeExt}.
\end{proof}

We now go on to discuss the multiplicative structure of the Tate spectral sequence. Let $\mu \: X \wedge_R Y \to Z$ be a pairing of $R$-modules in orthogonal $G$-spectra. As discussed in the paragraph before Theorem~\ref{thm:pairingbarmur}, the diagonal approximation~$D$ and~$\mu$ combine to define a pairing $\bar\mu \: (M_\star(X), M_\star(Y)) \to M_\star(Z)$ of sequences of~$R$-modules in orthogonal $G$-spectra. By Lemma~\ref{lem:Telphi} we
have an induced pairing
\[
T(\bar\mu) \: (T_\star(M(X)), T_\star(M(Y)))
	\longto T_\star(M(Z))
\]
of filtrations. By Proposition~\ref{prop:nablaextendstoN} there is also a pairing of filtrations
\[
N \: (\widetilde{E}_\star, \widetilde{E}_\star) \longto \widetilde{E}_\star \,
\]
which extends the fold map. Hence $(\widetilde{E}_\star, N)$ is a multiplicative $R$-module filtration in orthogonal $G$-spectra. We can now form the induced pairing of convolution filtrations
\[
\theta = N \wedge T(\bar\mu) \: (HM_\star(X), HM_\star(Y))
	\longto HM_\star(Z) \,.
\]
This has components
\[
\theta_{i,j} \: HM_i(X) \wedge_R HM_j(Y) \longto HM_{i+j}(Z)
\]
given by the union over $i_1+i_2 = i$ and $j_1 + j_2 = j$ of the composite maps
\begin{align*}
\widetilde{E}_{i_1} \wedge_R T_{i_2}(M(X))
	\wedge_R \widetilde{E}_{j_1} &\wedge_R T_{j_2}(M(Y)) \\
&\xrightarrow{1 \wedge \tau \wedge 1}
\widetilde{E}_{i_1} \wedge_R \widetilde{E}_{j_1}
	\wedge_R T_{i_2}(M(X)) \wedge_R T_{j_2}(M(Y)) \\
&\xrightarrow{N_{i_1,j_1} \wedge T(\bar\mu)_{i_2,j_2}}
\widetilde{E}_{i_1+j_1} \wedge_R T_{i_2+j_2}(M(Z)) \\
&\xrightarrow{\iota_{i_1+j_1,i_2+j_2}}
HM_{i+j}(Z) \,.
\end{align*}
Viewing $HM_i(X)$ as a subspectrum of $\widetilde E \wedge_R \Tel(M_\star(X))$ (and similarly for $Y$ and~$Z$ in place of $X$) the maps $\theta_{i,j}$ are compatible with the composite map
\begin{align*}
\widetilde E \wedge_R F_R(E, X) \wedge_R \widetilde E \wedge_R F_R(E, Y)
	&\xrightarrow{1 \wedge \tau \wedge 1}
\widetilde E \wedge_R \widetilde E \wedge_R F_R(E, X) \wedge_R F_R(E, Y) \\
	&\xrightarrow{1 \wedge \alpha }
\widetilde E \wedge_R \widetilde E \wedge_R F_R(E \wedge_R E, X \wedge_R Y) \\
	&\xrightarrow{N \wedge 1}
\widetilde E \wedge_R F_R(E \wedge_R E, X \wedge_R Y) \\
	&\xrightarrow{(1 \wedge D_+)^*}
\widetilde E \wedge_R F_R(E, X \wedge_R Y) \\
	&\xrightarrow{1 \wedge \mu_*}
\widetilde E \wedge_R F_R(E, Z) \,.
\end{align*}
This is $G$-homotopic to the corresponding map with $\Delta$ in
place of $D$, which defines the product
\[
\theta_* \: \pi_*(X^{tG}) \otimes_{R_*} \pi_*(Y^{tG}) \to \pi_*(Z^{tG})
\]
that we introduced in Section~\ref{subsec:GTatecon}.  Hence this product is filtration-preserving, taking $F_i \pi_*(X^{tG}) \otimes_{R_*} F_j \pi_*(Y^{tG})$ to $F_{i+j} \pi_*(Z^{tG})$ for all $i$ and~$j$. We write $\bar\theta_*$ for the induced pairing of filtration subquotients.

\begin{theorem} \label{thm:GTatespseqismult} 
Let $\mu \: X \wedge_R Y \to Z$ be a pairing of $R$-modules in orthogonal~$G$-spectra, and assume that $R[G]_*$ is projective over $R_*$. The pairing 
\[
\theta = N \wedge T(\bar\mu) \: (HM_\star(X), HM_\star(Y)) \to HM_\star(Z)
\]
of filtrations induces a pairing of $G$-Tate spectral sequences
\[
\theta \: \hat E^*(X) \otimes_{R_*} \hat E^*(Y)
	\longto \hat E^*(Z)
\]
in the sense of Definition~\ref{def:sseqpairing}. Moreover, the induced pairing $\theta^\infty$ of $E^\infty$-pages is compatible with the pairing $\bar\theta_*$ of filtration subquotients, in the sense of Proposition~\ref{prop:multabutment}.
\end{theorem}

\begin{proof}
This is a direct consequence of Theorem~\ref{thm:mainmultII}.
\end{proof}

\begin{corollary} 
If $(X, \mu \: X \wedge_R X \to X)$ is a multiplicative $R$-module in orthogonal $G$-spectra, then $(HM_\star(X), N \wedge T(\bar\mu))$ is a multiplicative filtration, and the $G$-Tate spectral sequence
\[
(\hat E^r(X), d^r) = (E^r(HM_\star(X)), d^r)
\]
is a multiplicative spectral sequence with multiplicative
abutment $\pi_*(X^{tG})$.\index{Hesselholt--Madsen--Tate spectral sequence!multiplicative}
\begin{proof}
This follows from Corollary~\ref{cor:main}.
\end{proof}
\end{corollary}

\begin{proposition} \label{prop:theta1fromPhiPsi}
Let $\mu \: X \wedge_R Y \to Z$ be a pairing of $R$-modules in orthogonal $G$-spectra and assume that $R[G]_*$ is $R_*$-projective. Under the isomorphism of Proposition~\ref{prop:TateE1isTatecx}, the pairing
\[
\theta^1 \: \hat E^1(X) \otimes_{R_*} \hat E^1(Y)
	\longto \hat E^1(Z)
\]
corresponds to the pairing covariantly induced by $\Phi \: \widetilde{P}_{*,*} \otimes_{R_*} \widetilde{P}_{*,*} \longto \widetilde{P}_{*,*}$ and contravariantly induced by $\Psi \: P_{*,*} \longto P_{*,*} \otimes_{R_*} P_{*,*}$, as in Section~\ref{sec:multTate}.
\end{proposition}

\begin{proof}
For typographical reasons we will use the abbreviation
\[
M^{R[G]_*} = \Hom_{R[G]_*}(R_*, M)
\]
in what follows, for various $R[G]_*$-modules~$M$. In the same way as in the proof of Proposition~\ref{prop:TateE1isTatecx}, the notation $E^1(HM_\star(X))$ will refer to the $E^1$-page of the associated spectral sequence on equivariant homotopy groups, while $E^1((\widetilde E \wedge_R T(M(X)))_\star)$ will refer to the $E^1$-page of the associated non-equivariant spectral sequence. 

We first note some results regarding multiplicative compatibility. Firstly, the natural homomorphism~$\omega$ is monoidal by Lemma~\ref{lem:omegaXmon}, so the diagram
\[
\xymatrix{
E^1(HM_\star(X)) \otimes_{R_*} E^1(HM_\star(Y))
	\ar[rr]^-{\theta^1} \ar[d]_-{\omega \otimes \omega}^-{\cong}
&& E^1(HM_\star(Z)) \ar[dd]^-{\omega}_-{\cong} \\
\txt{$E^1((\widetilde E \wedge_R T(M(X)))_\star)^{R[G]_*}$ \\
$\otimes_{R_*}$ \\
$E^1((\widetilde E \wedge_R T(M(Y)))_\star)^{R[G]_*}$}
	\ar[d]_{\alpha} \\
\left(\txt{$E^1((\widetilde E \wedge_R T(M(X)))_\star)$ \\
$\otimes_{R_*}$ \\
$E^1((\widetilde E \wedge_R T(M(Y)))_\star)$}\right)^{R[G]_*}
	\ar[rr]^-{(\theta^1)^{R[G]_*}}
&& E^1((\widetilde E \wedge_R T(M(Z)))_\star)^{R[G]_*}
}
\]
commutes. Secondly, by a slight generalization of Lemma~\ref{lem:iota1xyxymult}, the pairing
\[
\iota^1 \: E^1(\widetilde{E}_\star) \otimes_{R_*} E^1(T_\star(M(X)))
	\longto E^1((\tilde E \wedge_R T(M(X)))_\star)
\]
and its variants for $Y$ and~$Z$ in place of $X$ are multiplicatively compatible in the sense that the diagram 
\[
\xymatrix{
\txt{$E^1(\widetilde{E}_\star) \otimes_{R_*} E^1(T_\star(M(X)))$ \\
$\otimes_{R_*}$ \\
$E^1(\widetilde{E}_\star) \otimes_{R_*} E^1(T_\star(M(Y)))$}
\ar[rr]^-{1 \otimes \tau \otimes 1}
\ar[dd]_-{\iota^1 \otimes \iota^1}^-{\cong}
&&
\txt{$E^1(\widetilde{E}_\star) \otimes_{R_*} E^1(\widetilde{E}_\star)$ \\
$\otimes_{R_*}$ \\
$E^1(T_\star(M(X))) \otimes_{R_*} E^1(T_\star(M(Y)))$}
\ar[d]^-{N^1 \otimes T(\bar\mu)^1}
\\
&& E^1(\widetilde{E}_\star) \otimes_{R_*} E^1(T_\star(M(Z)))
	\ar[d]^-{\iota^1}_-{\cong} \\
\txt{$E^1((\widetilde{E} \wedge_R T(M(X)))_\star)$ \\
$\otimes_{R_*}$ \\
$E^1((\widetilde{E} \wedge_R T(M(Y)))_\star)$}
\ar[rr]^-{\theta^1}
&& E^1((\widetilde{E} \wedge_R T(M(Z)))_\star)
}
\]
commutes. Note that, by Definition~\ref{def:Phi}, the map 
\[
N^1 \: E^1(\widetilde{E}_\star) \otimes_{R_*}
E^1(\widetilde{E}_\star) \longto E^1(\widetilde{E}_\star)
\]
corresponds to $\Phi \: \widetilde{P}_{*,*} \otimes_{R_*} \widetilde{P}_{*,*} \to \widetilde{P}_{*,*}$ under the isomorphism $E^1_{*,*}(\widetilde{E}_\star) \cong \widetilde{ P}_{*,*}$. As discussed in the proofs of Proposition~\ref{prop:Leibpairseq},
Theorem~\ref{thm:pairingbarmur} and (the non-equivariant version of) Proposition~\ref{prop:E1d1mu1fromPPsi},
\[
T(\bar\mu)^1 \:
	E^1(T_\star(M(X))) \otimes_{R_*} E^1(T_\star(M(Y)))
	\longto E^1(T_\star(M(Z)))
\]
corresponds to the composite homomorphism
\begin{align*}
\Hom_{R_*}(P_{*,*}, \pi_*(X)) &\otimes_{R_*} \Hom_{R_*}(P_{*,*}, \pi_*(Y)) \\
&\overset{\alpha}\longto
\Hom_{R_*}(P_{*,*} \otimes_{R_*} P_{*,*}, \pi_*(X) \otimes_{R_*} \pi_*(Y)) \\
&\overset{\Psi^*}\longto
\Hom_{R_*}(P_{*,*}, \pi_*(X) \otimes_{R_*} \pi_*(Y)) \\
&\overset{\mu_*}\longto
\Hom_{R_*}(P_{*,*}, \pi_*(Z))
\end{align*} 
under the isomorphisms
\[
E^1(T_\star(M(X)))
\cong
E^1(M_\star(X))
\cong
\Hom_{R_*}(P_{*,*}, \pi_*(X))
\]
and their variants with $Y$ and~$Z$ in place of $X$.

Combining all of these results, we have shown that $\theta^1$ corresponds to the composite
\begin{align*}
&\hm(\pi_*(X))^{R[G]_*} \otimes_{R_*} \hm(\pi_*(Y))^{R[G]_*}
\xrightarrow{\alpha}
(\hm_*(\pi_*(X)) \otimes_{R_*} \hm(\pi_*(Y)))^{R[G]_*} \\
&\xrightarrow{1 \otimes \tau \otimes 1}
(\widetilde{P}_{*,*} \otimes_{R_*} \widetilde{P}_{*,*} \otimes_{R_*}
	\Hom_{R_*}(P_{*,*}, \pi_*(X)) \otimes_{R_*}
	\Hom_{R_*}(P_{*,*}, \pi_*(Y)))^{R[G]_*} \\
&\xrightarrow{1 \otimes 1 \otimes \alpha}
(\widetilde{P}_{*,*} \otimes_{R_*} \widetilde{P}_{*,*} \otimes_{R_*}
	\Hom_{R_*}(P_{*,*} \otimes_{R_*} P_{*,*},
	\pi_*(X) \otimes_{R_*} \pi_*(Y)))^{R[G]_*} \\
&\xrightarrow{\Phi \otimes 1}
(\widetilde{P}_{*,*} \otimes_{R_*}
	\Hom_{R_*}(P_{*,*} \otimes_{R_*} P_{*,*},
	\pi_*(X) \otimes_{R_*} \pi_*(Y)))^{R[G]_*} \\
&\xrightarrow{1 \otimes \Psi^*}
\hm_*(\pi_*(X) \otimes_{R_*} \pi_*(Y))^{R[G]_*} \\
&\xrightarrow{1 \otimes \mu_*}
\hm_*(\pi_*(Z))^{R[G]_*} \,,
\end{align*}
where we have abbreviated 
\[
\hm_*(\pi_*(X)) = \widetilde{P}_{*,*} \otimes_{R_*} \Hom_{R_*}(P_{*,*},\pi_*(X)) \,.
\]
Note that this is the pairing that induces the cup product, as in Section~\ref{sec:multTate}.
\end{proof}

\begin{theorem} \label{thm:E2Tatepairingiscup}\index{cup product!on complete $\Ext$}
Let $\mu \: X \wedge_R Y \to Z$ be a pairing of $R$-modules in orthogonal~$G$-spectra, and assume that $R[G]_*$ is $R_*$-projective. Then the pairing
\[
\theta^2 \: E^2_{i,*}(HM_\star(X)) \otimes_{R_*} E^2_{j,*}(HM_\star(Y))
	\longto E^2_{i+j,*}(HM_\star(Z))
\]
of $G$-Tate spectral sequence $E^2$-pages corresponds, under
the isomorphism of Theorem~\ref{thm:TateE2isTateExt}, 
to the cup product
\[
\smile \: \hatExt_{R[G]_*}^{-i,*}(R_*, \pi_*(X))
\otimes_{R_*} \hatExt_{R[G]_*}^{-j,*}(R_*, \pi_*(Y))
\longto \hatExt_{R[G]_*}^{-i-j,*}(R_*, \pi_*(Z))
\]
associated to $\mu_* \: \pi_*(X) \otimes_{R_*} \pi_*(Y) \to \pi_*(Z)$.
\end{theorem}

\begin{proof}
This is immediate by passage to homology from
Proposition~\ref{prop:theta1fromPhiPsi}. See Section~\ref{sec:multTate} for the definition of the cup product in Hopf algebra Tate cohomology.
\end{proof}

\begin{corollary}
If $(X, \mu)$ is a multiplicative $R$-module in orthogonal $G$-spectra, then the product in $\hat E^2(X) = E^2(HM_\star(X))$ corresponds to the cup product in
$\hatExt_{R[G]_*}^*(R_*, \pi_*(X))$ that is associated to the product~$\mu_*$ in $\pi_*(X)$.
\end{corollary}

Note independence of the particular choices of maps $D$ and $N$,
since the resulting chain maps $\Psi$ and $\Phi$ are unique up to homotopy, per Proposition~\ref{prop:Pidentityextension} and Proposition~\ref{prop:tildePfoldextension}.

\section{The Greenlees--May filtration}
\label{subsec:GMfilt}
\index{Greenlees--May filtration|(}
Greenlees~\cite{Gre87}*{\S1} spliced the filtration $F_\star \wEG$ with its Spanier--Whitehead dual to obtain a sequence of $G$-spectra
\[
\cdots \longto D(F_2 \wEG) \longto D(F_1 \wEG) \longto \bS \longto \Sigma^\infty F_1 \wEG \longto
  \Sigma^\infty F_2 \wEG \longto \cdots
\]
with mapping telescope equivalent to $\wEG$. The induced sequence
\[
\dots \to D(F_1 \wEG) \wedge F(EG_+, X)
\to \Sigma^\infty F(EG_+, X)
\to F_1 \wEG \wedge F(EG_+, X)
\to \cdots
\]
was used in \cite{GM95}*{(9.5), Thm.~10.3} to define a spectral sequence with abutment being the homotopy groups of the $G$-Tate construction on $X$. In this section, we will define a spliced filtration $GM_\star(X)$ with a  map to the Hesselholt--Madsen filtration~$HM_\star(X)$, and show that the induced map of $G$-homotopy spectral sequences
\[
\check E^r(X) = E^r(GM_\star(X))
	\longto E^r(HM_\star(X)) = \hat E^r(X)
\]
is an isomorphism for $r\ge2$. Thereafter we show that $GM_\star(X)$ is equivalent to the spliced sequence of Greenlees and May, at least for finite groups~$G$. For other compact Lie groups the sequences will differ in the same way that our filtration~$F_\star \wEG$ differs from the $G$-CW skeletal filtration. See Remark~\ref{rmk:EGtildefiltcomp}.

\begin{definition}
Recall the filtration $\widetilde{E}_\star$ from Definition~\ref{def:tildeEstar} and let~$GM_\star(X)$ be the filtration of orthogonal~$G$-spectra defined as 
\[
GM_k(X) = \begin{cases}
\widetilde{E}_k \wedge_R T_0(M(X)) & \text{for $k\ge0$,} \\
\widetilde{E}_0 \wedge_R T_k(M(X)) & \text{for $k\le0$.}
\end{cases}
\]
The structure maps $GM_{k-1}(X) \to GM_k(X)$ for $k\ge1$ are
induced by the maps $\widetilde{E}_{k-1} \to \widetilde{E}_k$ in the filtration $\widetilde{E}_\star$, while the maps for $k\le 0$ are those of $T_\star(M(X))$. We refer to the filtration $GM_\star(X)$ as the \emph{Greenlees--May filtration}.
\end{definition}

\begin{notation}
\index{Greenlees--May--Tate spectral sequence}
Let
\[
\check E^r(X) = E^r(GM_\star(X))
\]
denote the $G$-homotopy spectral sequence associated to the filtration $GM_\star(X)$.  
\end{notation}

We now discuss the map of filtrations between the Greenlees--May filtration and the Hesselholt--Madsen filtration.

\begin{lemma} \label{lem:alphamapfilt}
The inclusions
$\widetilde{E}_k \wedge_R T_0(M(X)) \to (\widetilde{E} \wedge_R T(M(X)))_k$
for $k\ge0$ and
$\widetilde{E}_0 \wedge_R T_k(M(X)) \to (\widetilde{E} \wedge_R T(M(X)))_k$
for $k\le0$ define a map of filtrations
\[
\alpha \: GM_\star(X) \longto HM_\star(X)
\]
of $R$-modules in orthogonal $G$-spectra.
The induced maps of mapping telescopes and colimits
\[
\xymatrix{
\Tel(GM_\star(X)) \ar[r]^-{\simeq_G} \ar[d]_-{\simeq_G}
& \Tel(HM_\star(X)) \ar[d]^-{\simeq_G} \\
\widetilde E \wedge_R T_0(M(X)) \ar[r]^-{\simeq_G}
	& \widetilde E \wedge_R \Tel(M_\star(X))
}
\]
are all equivalences.
\end{lemma}

\begin{proof}
Recall from Section~\ref{subsec:HMfilt} that
\[
HM_k(X) = \bigcup_{i+j=k} \widetilde{E}_i \wedge_R T_j(M(X))
\]
as a subspectrum of $\widetilde E \wedge_R \Tel(M(X))$. The existence of the filtered map $\alpha$ is then clear. The vertical maps from mapping telescopes to colimits are equivalences, since~$GM_\star(X)$ and $HM_\star(X)$ are both filtrations.  The lower horizontal map is also an equivalence, since the sequence $M_\star(X)$ is constant for $\star \ge 0$.
\end{proof}

As a consequence of the above lemma, there is a map 
\[
\alpha \: \check E^*(X) \to \hat
E^*(X)
\]
of $R_*$-module spectral sequences.

\begin{remark}\index{Greenlees--May filtration!multiplicative}
Recall from Proposition~\ref{prop:nablaextendstoN} that, under the assumption that the Hopf algebra~$R[G]_*$ is~$R_*$-projective, we have a filtration-preserving pairing $N \: (\widetilde{E}_\star,\widetilde{E}_\star) \to \widetilde{E}_\star$. However, when $(X, \mu)$ is multiplicative, the induced pairing 
\[
N \wedge T(\bar\mu) \: (HM_\star(X), HM_\star(X)) \longto HM_\star(X)
\]
does usually not restrict to a multiplication on $GM_\star(X)$.  For instance, $GM_a(X) \wedge_R GM_{-b}(X)$ with~$a>0$ and $b>0$ maps to $HM_a(X) \wedge_R HM_{-b}(X)$ and $\widetilde{E}_a \wedge_R T_{-b}(M(X))$ in $HM_{a-b}(X)$, which is hardly ever in $GM_{a-b}(X)$. Hence $GM_\star(X)$ is not a multiplicative filtration, and $\check E^r(X)$ is not evidently
a multiplicative spectral sequence. Nonetheless, we will show that $\check E^r(X)$ is isomorphic to the $G$-Tate spectral sequence $\hat E^r(X)$ for $r\ge2$, which we showed to be multiplicative in Theorem~\ref{thm:GTatespseqismult}. This will then show that $(\check E^r(X), d^r)$ is also multiplicative, at least for $r\ge2$.
\end{remark}

Thinking only about the additive properties of the spectral sequence
$\check E^r(X)$, we can safely replace the filtration $GM_\star(X)$
with a simpler, but equivalent, sequence.

\begin{lemma} \label{lem:GMeqGM'}
There is an equivalence from $GM_\star(X)$ to the sequence
$GM'_\star(X)$ with
$$
GM'_k(X) =
\begin{cases}
\widetilde{E}_k \wedge_R F_R(E, X) & \text{for $k\ge0$,} \\
F_R(E/E_{-k-1}, X) & \text{for $k\le0$.}
\end{cases}
$$
\end{lemma}

\begin{proof}
The equivalences $\epsilon \: T_k(M(X)) \to M_k(X)$ induce
the following commutative diagram.
$$
\xymatrix@C-1pc{
\dots \ar[r]
& \widetilde{E}_0 \wedge_R T_{-1}(M(X)) \ar[r] \ar[d]_-{\epsilon}^-{\simeq}
& \widetilde{E}_0 \wedge_R T_0(M(X)) \ar[r] \ar[d]_-{\epsilon}^-{\simeq}
& \widetilde{E}_1 \wedge_R T_0(M(X)) \ar[r] \ar[d]_-{\epsilon}^-{\simeq}
& \dots \\
\dots \ar[r]
& M_{-1}(X) \ar[r]
& F_R(E, X) \ar[r]
& \widetilde{E}_1 \wedge_R F_R(E, X) \ar[r]
& \dots
}
$$
Here $M_k(X) = F_R(E/E_{-k-1}, X)$ for $k\le0$.
\end{proof}

We refer to \cite{LMS86}*{\S III.1} for the basic Spanier--Whitehead
duality theory in a closed symmetric monoidal category.  In the case of (the homotopy category of)~$R$-modules in orthogonal $G$-spectra, we refer to the objects called `finite' by Lewis and~May as `dualisable'.

\begin{definition} \label{def:dual}\index{Spanier--Whitehead dual}
For an $R$-module $X$ in orthogonal $G$-spectra, let
\[
D(X) = F_R(X, R)
\]
be its \emph{functional dual}\index{functional dual!of spectrum}.  For dualisable $X$ we refer to $D(X)$ as the \emph{Spanier--Whitehead dual}\index{Spanier--Whitehead dual} of $X$. There are natural maps
\[
\rho \: X \to D(D(X)) \och \nu \: DX \wedge_R Y \to F_R(X, Y) \,,
\]
which are equivalences when $X$ is dualisable. 
\end{definition}

\begin{lemma} \label{lem:tildeEidualizable}
Each term in the filtration $\widetilde{E}_\star$ is dualisable.
\end{lemma}

\begin{proof}
We can give $G$ a finite CW structure, with $e$ as a $0$-cell.
It follows that the bar construction is a finite $G$-CW space in each simplicial degree $B_q(*, G, G) = G^q \times G$, so that $G^{\wedge q} \wedge G_+$ is a finite $G$-CW space and $R \wedge G^{\wedge q} \wedge G_+$ is a dualisable $R$-module in orthogonal $G$-spectra. By induction, this implies that $E_{i-1}$ is dualisable, and therefore the mapping cone $\widetilde{E}_i$ is also dualisable, for each $i\ge0$.
\end{proof}

\begin{lemma} \label{lem:checkE1isgmcx}
The $E^1$-page of the spectral sequence associated to~$GM'_\star(X)$ is the $R_*$-module chain complex with
\[
\check E^1_{\ast,\ast}(X) \cong \Hom_{R[G]_*}(R_*,\gm_*(\pi_*(X))) \,,
\]
where we use the notation of Definition~\ref{def:gmcx}.
\end{lemma}

\begin{proof}
For $\star \le 0$ the sequence $GM'_\star(X)$ agrees with the sequence
$M_\star(X)$ from Section~\ref{subsec:htpyfix}, so $(\check E^1_*(X),
d^1)$ for $*\le0$ agrees with $\Hom_{R[G]_*}(P_{-*,*}, \pi_*(X))$
by Proposition~\ref{prop:E1d1mu1fromPPsi}.

For $\star \ge 0$ the sequence $GM'_\star(X)$ agrees with the
filtration $\widetilde{E}_\star \wedge_R F_R(E, X)$.  Its subquotients
for $i\ge1$ are of the form
$(\widetilde{E}_i / \widetilde{E}_{i-1}) \wedge_R F_R(E, X)$
with
$$
\widetilde{E}_i / \widetilde{E}_{i-1} \cong \Sigma(E_{i-1} / E_{i-2})
\cong R \wedge \Sigma^i(G^{\wedge i-1} \wedge G_+) \,.
$$
Let $d$ be the dimension of~$G$.  Since $G_+$ is stably dualisable,
with Spanier--Whitehead dual $D(G_+) \simeq_G \Sigma^{-d} G_+$, each
subquotient above is equivalent to $F(G_+, X')$ for some $R$-module $X'$
in orthogonal $G$-spectra.
It follows from Proposition~\ref{prop:omegaXiso} that
$$
\check E^1_*(X) \cong \Hom_{R[G]_*}(R_*,
	E^1_*(\widetilde{E}_\star \wedge_R F_R(E, X)))
$$
for $*\ge1$.  Here
\begin{align*}
E^1_i(\widetilde{E}_\star \wedge_R F_R(E, X))
&= \pi_{i+*}(\widetilde{E}_{i-1} \wedge_R F_R(E, X)
	\to \widetilde{E}_i \wedge_R F_R(E, X))  \\
&\cong \pi_{i+*}((\widetilde{E}_i/\widetilde{E}_{i-1} \wedge_R F_R(E, X))) \\
&\cong \pi_{i+*}(\widetilde{E}_i/\widetilde{E}_{i-1})
	\otimes_{R_*} \pi_* F_R(E, X) \\
&\cong \widetilde{P}_{i,*} \otimes_{R_*} \pi_*(X)
\end{align*}
for $i\ge1$, since $\widetilde{P}_{i,*} = \pi_{i+*}(\widetilde {E}_i/\widetilde{E}_{i-1})$
is projective, hence flat, over $R_*$, and $c \: E \to R$ induces an
isomorphism $\pi_*(X) \cong \pi_* F_R(E, X)$ of $R[G]_*$-modules.
This shows that $(\check E^1_*(X), d^1)$ for $*\ge1$ agrees with
$\Hom_{R[G]_*}(R_*, \widetilde{P}_{*,*} \otimes_{R_*} \pi_*(X))$.

It remains to verify that $d^1 \: \check E^1_1(X) \to \check E^1_0(X)$
is as asserted.  By definition, it is given by the left-to-right composite
in the following diagram.
\[
\xymatrix{
\pi^G_{1+*}((\widetilde{E}_1/\widetilde{E}_0) \wedge_R X) \ar[r]^-{\partial}
	\ar[d]_-{1 \wedge c^*}^-{\cong}
& \pi^G_*(X) \ar[d]_-{c^*} \\
\pi^G_{1+*}((\widetilde{E}_1/\widetilde{E}_0) \wedge_R F_R(E, X)) \ar[r]^-{\partial}
& \pi^G_*(F_R(E, X)) \ar[r]
& \pi^G_*(F_R(E_0, X))
}
\]
By naturality of $\omega$, as in Lemma~\ref{lem:omegaXnat}, this is obtained
from the left-to-right composite
\[
\xymatrix{
\pi_{1+*}((\widetilde{E}_1/\widetilde{E}_0) \wedge_R X) \ar[r]^-{\partial}
& \pi_*(X) \ar[d]_-{c^*}^-{\cong} \\
& \pi_*(F_R(E, X)) \ar[r]
& \pi_*(F_R(E_0, X))
}
\]
by applying $\Hom_{R[G]_*}(R_*, -)$.  Under the isomorphisms above,
this is the composition
\begin{align*}
\widetilde{P}_{1,*} \otimes_{R_*} \pi_*(X)
& \overset{\tilde\partial_1}\longto
	\widetilde{P}_{0,*} \otimes_{R_*} \pi_*(X) \\
& \cong \pi_*(X) \cong \Hom_{R_*}(R_*, \pi_*(X))
	\overset{\epsilon^*}\longto \Hom_{R_*}(P_0, \pi_*(X)) \,.
\end{align*}
As we made explicit in Proposition~\ref{prop:gmMexplicit}, this equals
the boundary $\gm_1(\pi_*(X)) \to \gm_0(\pi_*(X))$.
\end{proof}

\begin{proposition}
\label{prop:GMHMfiltrationiso}
The filtration-preserving map $\alpha \: GM_\star(X) \to HM_\star(X)$
induces an isomorphism of spectral sequences
\[
\alpha^r \: \check E^r(X) \overset{\cong}\longto \hat E^r(X)
\]
for $r\ge2$.
\end{proposition}

\begin{proof}
Comparing Proposition~\ref{prop:TateE1isTatecx} and
Lemma~\ref{lem:checkE1isgmcx} shows that 
\[
\alpha^1 \: \check E^1(X)
\to \hat E^1(X)
\] is the chain map $\Hom(1, \alpha)$ shown to be a
quasi-isomorphism in Proposition~\ref{prop:gmtohmisqi}.  Hence
$\alpha^2 = H(\alpha^1, d^1)$ is an isomorphism, which implies
that $\alpha^r$ is an isomorphism for each $r\ge2$.
\end{proof}

Following \cite{Gre87}*{\S1}, we can splice the filtration
\[
R \cong \widetilde{E}_0 \longto \widetilde{E}_1 \longto \widetilde{E}_2 \longto \dots
\]
with the Spanier--Whitehead dual sequence
\[
\dots \longto D(\widetilde{E}_2) \longto D(\widetilde{E}_1)
	\longto D(\widetilde{E}_0) \cong R
\]
to obtain a bi-infinite sequence
\begin{equation} \label{eq:biinfinite}
\dots \longto D(\widetilde{E}_2) \longto D(\widetilde{E}_1)
	\longto R \longto \widetilde{E}_1 \longto \widetilde{E}_2 \longto \dots
\end{equation}
of dualisable $R$-modules in orthogonal $G$-spectra.  This is the
sequence~\cite{GM95}*{(9.5)} used by Greenlees and May to define their
Tate spectral sequence, at least for finite~$G$. For $G = \bT = U(1)$
or $\bU = Sp(1)$ they instead repeat each term in this sequence two or
four times, respectively.  For other compact Lie groups, the connection
is less direct.

\begin{proposition}
There is a zigzag of equivalences from $GM'_\star(X)$ to the sequence~$GM''_\star(X)$ with
\[
GM''_k(X) =
\begin{cases}
\widetilde{E}_k \wedge_R F_R(E, X) & \text{for $k\ge0$,} \\
D(\widetilde{E}_{-k}) \wedge_R F_R(E, X) & \text{for $k\le0$.}
\end{cases}
\]
Hence the spectral sequence $\check E^r(X)$ is isomorphic to the
Greenlees--May Tate spectral sequence~\cite{GM95}*{Thm.~10.3} for
$\pi^G_*$ applied to the sequence $GM''_\star(X)$.
\end{proposition}

\begin{proof}
The zigzag of equivalences connecting $GM'_\star(X)$ to
$GM''_\star(X)$ consists of identity maps for $\star\ge0$.
For $\star\le0$ it takes the following form:
$$
\xymatrix@C-1pc{
\dots \ar[r]
  & F_R(E/E_1, X) \ar[r] \ar[d]^-{\simeq_G}
  & F_R(E/E_0, X) \ar[r] \ar[d]^-{\simeq_G}
  & F_R(E, X) \ar[d]^{=} \\
\dots \ar[r]
  & F_R(E \cup CE_1, X) \ar[r]
  & F_R(E \cup CE_0, X) \ar[r]
  & F_R(E, X) \\
\dots \ar[r]
  & F_R(\widetilde{E}_2 \wedge_R E, X) \ar[u]^-{\tilde\Delta^*}_-{\simeq_G} \ar[r]
  & F_R(\widetilde{E}_1 \wedge_R E, X) \ar[u]^-{\tilde\Delta^*}_-{\simeq_G} \ar[r]
  & F_R(\widetilde{E}_0 \wedge_R E, X) \ar[u]^-{\tilde\Delta^*}_-{\cong} \\
\dots \ar[r]
  & D(\widetilde{E}_2) \wedge_R F_R(E, X) \ar[u]^-{\nu}_-{\simeq_G} \ar[r]
  & D(\widetilde{E}_1) \wedge_R F_R(E, X) \ar[u]^-{\nu}_-{\simeq_G} \ar[r]
  & D(\widetilde{E}_0) \wedge_R F_R(E, X) \ar[u]^-{\nu}_-{\cong} \\
}
$$
The two top rows are equivalent because each quotient map
$E \cup CE_{i-1} \to E/E_{i-1}$ is an equivalence, since
$E_{i-1} \to E$ is a (strong) $h$-cofibration.
The equivalence between the middle two rows is induced by the
map $\widetilde\Delta$ of mapping cones associated to the
diagonal equivalence $\Delta \: E_{i-1} \to E_{i-1} \wedge_R E$:
\[
\xymatrix{
E_{i-1} \ar[r] \ar[d]_-{\Delta}^-{\simeq_G}
	& E \ar[r] \ar[d]^-{=}
	& E \cup CE_{i-1} \ar[d]_-{\widetilde\Delta}^-{\simeq_G} \\
E_{i-1} \wedge_R E \ar[r]^-{c \wedge 1}
	& R \wedge_R E \ar[r]
	& \widetilde{E}_i \wedge_R E \,.
}
\]
The lower two rows are equivalent because each $\widetilde{E}_i$
is dualisable by Lemma~\ref{lem:tildeEidualizable}.
\end{proof}

\begin{remark}
Our comparison of the Hesselholt--Madsen Tate spectral sequence
\[
\hat E^r(X) = E^r(HM_\star(X))
\]
and the Greenlees--May Tate spectral sequence
\[
\check E^r(X) = E^r(GM_\star(X))
	\cong E^r(GM'_\star(X))
	\cong E^r(GM''_\star(X))
\]
is a little different from that of \cite{HM03}*{Rmk.~4.3.6}, since we
obtain the Greenlees sequence $GM''_\star(X)$ by splicing the two
perpendicular edges $(i=0, j\le0)$ and $(i\ge0, j=0)$ of the bifiltration
\[
HM_{i,j}(X) = 
\widetilde{E}_i \wedge_R T_j(M(X))
\simeq_G
\widetilde{E}_i \wedge_R F_R(E/E_{-j-1}, X) \,,
\]
while Hesselholt and Madsen first invert a quasi-isomorphism, so as to
position both halves of the Greenlees sequence on the line $j=0$.
\end{remark}

\begin{remark} \label{rem:GM-vs-Whitehead}\index{odd spheres filtration}\index{Whitehead tower}
In the case of the circle group $G = \bT$ we can work over $R = \bS$
and use the odd spheres $S((k+1)\bC)$ to filter $E\bT =
S(\infty\bC)$, so that $\widetilde{E}_k = S^{k\bC}$ equals (the suspension
spectrum of) a representation sphere.  Then $D(\widetilde{E}_{-k})
= D(S^{-k\bC}) = S^{k\bC}$ is a virtual representation sphere for
each $k<0$.  For brevity, let us also write~$\widetilde{E}_k$ for the
latter $\bT$-spectra, so that $\{\widetilde{E}_k\}_{k \in \bZ}$ is the
bi-infinite sequence~\eqref{eq:biinfinite}, with cofibre sequences
$\widetilde{E}_{k-1} \to \widetilde{E}_k \to \Sigma^{2k-1} \bT_+$.
The Greenlees--May spectral sequence associated to the sequence
$$
\widetilde{E}_\star \wedge F(E\bT_+, X)
$$
of
$\bT$-spectra has $E^1$-page
$$
\check E^1_{k,*}(X)
  = \pi_{k+*}^{\bT}(\Sigma^{2k-1} \bT_+ \wedge F(E\bT_+, X))
  \cong \pi_{*-k}(X)
$$
for each $k \in \bZ$, which we can formally write as $t^{-k} \cdot
\pi_*(X)$ with $t$ in bidegree $(-1,-1)$.  As remarked earlier, we may
reindex the filtration and spectral sequence so as to put $t$ in bidegree
$(-2,0)$, in which case $E^2_{2k,*}(X) \cong t^{-k} \cdot \pi_*(X)$
and $E^2_{2k-1,*}(X) = 0$.
Maunder~\cite{Mau63}*{Thm.~3.3}, as well as Greenlees and May \cite{GM95}*{Thm.~B.8}, prove that the latter spectral
sequence is isomorphic to one associated to the tower of $\bT$-spectra
$$
\dots \longto \widetilde{E\bT} \wedge F(E\bT_+, X^{\ell+1})
	\longto \widetilde{E\bT} \wedge F(E\bT_+, X^\ell) \longto \dots \,.
$$
Here $\{X^\ell\}_\ell$ denotes a $\bT$-equivariant Whitehead tower
for~$X$, with homotopy fibre sequences $X^{\ell+1} \to X^\ell \to
\Sigma^\ell H\pi_\ell(X)$.  The latter spectral sequence is indexed
so that
$$
E^2_{*,\ell}(X)
  = \pi_{*+\ell}^{\bT}(\widetilde{E\bT}
	\wedge F(E\bT_+, \Sigma^\ell H\pi_\ell(X)))
	\cong \pi_*(H\pi_\ell(X)^{t\bT})
$$
for each integer $\ell$.  In particular 
\[
\pi_{2k}(H \pi_\ell(X)^{t\bT})
\cong t^{-k} \cdot \pi_\ell(X) \och \pi_{2k-1}(H \pi_\ell(X)^{t\bT}) =
0 \,,
\]
so that, formally, $\pi_* (H \pi_\ell(X)^{t\bT}) \cong \pi_\ell(X)
[t, t^{-1}]$.  Furthermore, Greenlees and May argue that the latter
spectral sequence is multiplicative, with respect to some topologically
defined pairings of the form
$$
\pi_* (H \pi_i(X)^{t\bT}) \otimes \pi_* (H \pi_j(Y)^{t\bT})
	\longto \pi_* (H \pi_{i+j}(Z)^{t\bT}) \,.
$$
However, as is implicit in \cite{GM95}*{Prob.~14.8}, they do not
establish that these topological pairings agree with the evident
algebraic pairings
$$
\pi_i(X) [t, t^{-1}] \otimes \pi_j(Y) [t, t^{-1}]
	\longto \pi_{i+j}(Z) [t, t^{-1}] \,.
$$
Hence they do not assert that the isomorphism $E^2_{*,*}(X) \cong \pi_*(X)
[t, t^{-1}]$ takes the topological product to the algebraic product.
In particular, the higher differentials in this spectral sequence are
known to obey a Leibniz rule, but conceivably not with respect to the
most evident algebraic product.

Nonetheless, we can confirm directly that the first differential in
each of these spectral sequences is a derivation with respect to the
algebraic product.  To express this, we return to the indexing used
elsewhere in the memoir, i.e., to the Greenlees--May spectral sequence
$\check E^r_{*,*}(X)$.  Up to the technical issue we have pointed out about
compatibility of product structures, the following result is due to
Hesselholt \cite{Hes96}*{Lem.~1.4.2}.
\end{remark}

\begin{proposition} \label{prop:d1-in-GM-sp-seq}
Let $X$ be any orthogonal $\bT$-spectrum, so that $\pi_*(X)$
is a right $\bS[\bT]_*$-module.  There is a natural isomorphism
$$
\check E^1_{*,*}(X) \cong \pi_*(X) [t, t^{-1}]
$$
with $t$ in bidegree $(-1,-1)$, such that $d^1 \: \check E^1_{k,*}(X)
\to \check E^1_{k-1,*}(X)$ corresponds to the differential
$d \: t^{-k} \cdot \pi_*(X) \to t^{-k+1} \cdot \pi_*(X)$ given
by
$$
d(t^{-k} \cdot x) = \begin{cases}
t^{-k+1} \cdot x s & \text{for $k$ even,} \\
t^{-k+1} \cdot x (s + \eta) & \text{for $k$ odd.}
\end{cases}
$$
\end{proposition}

\begin{proof}
By naturality of the Greenlees--May spectral sequence with respect to
$\bT$-maps $x \: \Sigma^\ell \bS[\bT] \to X$, corresponding to homotopy
classes $x \in \pi_\ell(X)$, it suffices to prove the result in the case $X = \bS[\bT]$ and $x = 1 \in \pi_0(\bS[\bT])$. 

Consider the case $X = H\bZ[\bT]$. We have $\pi_*(H\bZ[\bT]) = \bZ\{1,\sigma\}$ and $H\bZ[\bT]^{t\bT} \simeq *$ since $H\bZ$, as a $\bT$-spectrum, is induced up from $H\bZ$. For bidegree reasons the $\bT$-Tate spectral
sequence must collapse to zero at the $E^2$-page, which forces 
\[
d(t^{-k} \cdot 1) = \pm t^{-k+1} \cdot \sigma \,.
\]
Here, we can iteratively fix the sign of $t^{-k}$ implicit in the identification $\check E^1_{k,*}(X) \cong t^{-k} \cdot \pi_*(X)$ so that each of these signs is
a plus.  By naturality with respect to the Hurewicz homomorphism $\bS[\bT] \to H\bZ[\bT]$ it follows
that 
\[
d(t^{-k} \cdot 1) \equiv t^{-k+1} \cdot s \mod t^{-k+1} \cdot \eta
\]
in the $\bT$-Tate spectral sequence for $\bS[\bT]$, since $\pi_1(\bS[\bT]) = \bZ\{s\} \oplus \bZ/2\{\eta\}$
with the Hurewicz homomorphism mapping $s$ to $\sigma$.

Now consider the case $X = \bS$ with trivial $\bT$-action. The part $k\ge1$ of the Greenlees--May spectral sequence maps to the Atiyah--Hirzebruch spectral sequence for~$\Sigma^2
\bC P^\infty_+$. Since the $2k$-cell in $\Sigma^2 \bC P^\infty$ is
stably attached to the $2k-2$-cell by~$k \eta$, it follows that 
\[
d(t^{-k} \cdot 1) = t^{-k+1} \cdot k \eta
\]
for $k \ge 2$.  Similarly, the part
$k\le0$ receives a map from the Atiyah--Hirzebruch spectral sequence
for $D(\bC P^\infty_+)$, where the $-2k$-cell is attached
to the $-2k-2$-cell by $k\eta$, see~\cite{Mos68}*{Prop.~5.1}, so that 
\[
d(t^{-k} \cdot 1) = t^{-k+1} \cdot k \eta
\]
for $k \le 0$, as well. Finally, for $k=1$ the differential is induced by the
composite $\bT$-map
\[
\Sigma^{-1} \widetilde{E}_1/\widetilde{E}_0 \wedge F(E\bT_+, \bS)
\longto
\widetilde{E}_0 \wedge F(E\bT_+, \bS)
\longto
\widetilde{E}_0/\widetilde{E}_{-1} \wedge F(E\bT_+, \bS) \,,
\]
which we can rewrite in terms of the counit $\epsilon \: \bT_+ \to S^0$
and its Spanier--Whitehead dual $D(\epsilon) \: \bS \to D(\bT_+)$ as
\[
F(E\bT_+, \bT_+)
\xrightarrow{\epsilon_*}
F(E\bT_+, \bS)
\xrightarrow{D(\epsilon)_*}
F(E\bT_+, D(\bT_+)) \,.
\]
Passing to $\bT$-fixed points, this is a composite
$$
\Sigma \bS \simeq \bS[\bT]^{h\bT} \longto \bS^{h\bT} \longto
D(\bT_+)^{h\bT} \simeq \bS \,,
$$
which we claim equals $\eta \in \pi_1(\bS)$.
This can be seen using the Pontryagin--Thom collapse $t \: S^{\bC} \to
S^{\bC}/S^0 \cong \Sigma(\bT_+)$ associated to the embedding $\bT
\subset \bC$, and the untwisting isomorphism $\zeta \: \Sigma^2(\bT_+)
\cong \Sigma^{\bC}(\bT_+)$.  Then $\zeta(1\wedge t) \: \Sigma S^{\bC}
\to \Sigma^{\bC}(\bT_+)$ defines a stable $\bT$-map $S^1 \to \bT_+$.
The composite $\epsilon \zeta (1 \wedge t) \: \Sigma S^{\bC} \to S^{\bC}$
has (non-equivariant) Hopf invariant~$\pm 1$, since the preimages of
two generic points in~$S^{\bC}$ are circles in $\Sigma S^{\bC}$ with
that linking number.

This proves $d(t^{-1} \cdot 1) = t^0 \cdot \eta$ in the $\bT$-Tate spectral sequence for $\bS$ and, by naturality with
respect to $\bS[\bT] \to \bS$, the asserted formulas follow.
\end{proof}
\index{Greenlees--May filtration|)}

\section{Convergence}
\label{subsec:convergence}

In this section we use Proposition~\ref{prop:GMHMfiltrationiso} to deduce convergence results for the Hesselholt--Madsen spectral sequence $\hat E^r(X)$ from corresponding results for the Greenlees--May spectral sequence~$\check E^r(X)$.

\begin{lemma} \label{lem:alphainftyiso}
The map of abutments
\[
\alpha_\infty \: A_\infty(GM_\star(X))
	\overset{\cong}\longto
	A_\infty(HM_\star(X)) \cong \pi_*(X^{tG})
\]
is an isomorphism. Hence the homomorphism
\[
\alpha_s \: F_s A_\infty(GM_\star(X))
	\longto F_s A_\infty(HM_\star(X))
\]
is injective, for each $s \in \bZ$.
\end{lemma}

\begin{proof}
The first assertion follows from Lemma~\ref{lem:alphamapfilt}. The commutative diagram
\[
\xymatrix{
A_s(GM_\star(X)) \ar@{->>}[r] \ar[d]
& F_s A_\infty(GM_\star(X)) \ar@{ >->}[r] \ar[d]_-{\alpha_s}
& A_\infty(GM_\star(X)) \ar[d]_-{\alpha_\infty}^-{\cong} \\
A_s(HM_\star(X)) \ar@{->>}[r]
& F_s A_\infty(HM_\star(X)) \ar@{ >->}[r]
& A_\infty(HM_\star(X))
}
\]
implies the injectivity assertion.
\end{proof}

\begin{lemma} \label{lem:checkEXcondconv}
\index{Greenlees--May--Tate spectral sequence!conditional convergence}
The spectral sequence
\[
\check E^r(X) = E^r(GM_\star(X))
	\Longrightarrow \pi_*(X^{tG})
\]
is conditionally convergent.
\end{lemma}

\begin{proof}
The spectral sequence associated to $GM_\star(X)$ is conditionally
convergent to the colimit, in the sense of \cite{Boa99}*{Def.~5.10},
whenever $\holim_s GM_s(X) \simeq_G *$.  Since the
sequences $GM_\star(X)$ and $GM'_\star(X)$ are equivalent by
Lemma~\ref{lem:GMeqGM'}, we may equally well verify that $\holim_s
GM'_s(X) \simeq_G *$. But 
\[
GM'_s(X) = F_R(E/E_{s-1}, X) = M_s(X)
\] 
for $s\le 0$, and we saw in Section~\ref{subsec:htpyfix} that $\holim_s M_s(X) \simeq_G *$. Hence the Greenlees--May $G$-Tate spectral sequence $\check E^r(X)$ is always conditionally convergent.
\end{proof}

\begin{lemma}
The maps of $E^\infty$- and $RE^\infty$-pages
\begin{align*}
E^\infty(GM_\star(X)) &\overset{\cong}\longto E^\infty(HM_\star(X)) \\
RE^\infty(GM_\star(X)) &\overset{\cong}\longto RE^\infty(HM_\star(X))
\end{align*}
are isomorphisms.
\end{lemma}

\begin{proof}
Recall from \cite{Boa99}*{(5.1)} that for each spectral sequence $(E^r,
d^r)$ there are filtrations
\[
0 = B^1 \subset B^2 \subset B^3 \subset \dots
	\subset Z^3 \subset Z^2 \subset Z^1 = E^1
\]
with $E^r \cong Z^r/B^r$ for $r\ge1$. We set 
\[
B^\infty = \colim_r B^r \,, \quad Z^\infty = \lim_r Z_r \,, \quad E^\infty  = Z^\infty/B^\infty \,, \quad RE^\infty = \Rlim_r Z_r \,.
\]
Letting $\bar B^r = B^r/B^2$ and $\bar Z^r = Z^r/B^2$ for $r\ge2$, we obtain a filtration
\[
0 = \bar B^2 \subset \bar B^3 \subset \dots
	\subset \bar Z^3 \subset \bar Z^2 \cong E^2
\]
with $E^r \cong \bar Z^r/\bar B^r$ for $r\ge2$. Let 
\[
\bar B^\infty = \colim_r \bar B^r \och \bar Z^\infty = \lim_r \bar Z_r \,.
\]
Then $\bar B^\infty \cong B^\infty/B^2$, while $\bar Z^\infty \cong Z^\infty/B^2$ and $\Rlim_r Z^r \cong \Rlim_r \bar Z^r$ by the $\lim$-$\Rlim$ exact sequence. Hence $E^\infty \cong \bar Z^\infty / \bar B^\infty$ by the
Noether isomorphism, and $RE^\infty \cong \Rlim_r \bar Z^r$.

A map of spectral sequences inducing an isomorphism of $E^2$-pages will
by induction induce isomorphisms of $\bar B^r$- and $\bar Z^r$-pages for
all $r\ge2$, and therefore also of $\bar B^\infty$-,
$\bar Z^\infty$-, $E^\infty$- and $RE^\infty$-pages.
\end{proof}

When $X$ is bounded below, and $G$ is finite or equal to $\bT = U(1)$ or~$\bU = Sp(1)$, the $E^1$-pages $\check E^1(X)$ and $\hat E^1(X)$
are both concentrated in half-planes with entering differentials
\cite{Boa99}*{\S7}. However, for more general groups~$G$
(such as $\bT \times \bU$)
the $E^1$-page
$\check E^1(X)$ occupies a region that is only bounded by a broken line,
and~$\hat E^1(X)$ may not be bounded in any ordinary sense. We therefore need to discuss convergence for the spectral sequences $\check E^r(X)$ and~$\hat E^r(X)$ in the generality of whole-plane spectral sequences
\cite{Boa99}*{\S8}.

\begin{definition}
\index{Boardman's whole plane obstruction group}
Let $(A, E^1)$ be the exact couple associated to a
Cartan--Eilenberg system $(H, \partial)$.
Boardman's \emph{whole-plane obstruction group}~$W$
is defined in \cite{Boa99}*{Lem.~8.5} by an expression
$$
W = \colim_s \Rlim_r K_\infty \im^r A_s \,.
$$
We refer to Boardman's paper for an explanation of the notation.
By \cite{HR19}*{Thm.~7.5} there is an isomorphism
$$
W \cong \ker(\kappa)
$$
where
$$
\kappa \: \colim_j \lim_i H(i,j) \longto \lim_i \colim_j H(i,j)
$$
is the interchange morphism, which is always surjective.
\end{definition}

While $W$ is defined in terms of the underlying exact couple (or
Cartan--Eilenberg system), Boardman gives the following criterion for
the vanishing of~$W$, which is internal to the spectral sequence.

\begin{lemma}[\cite{Boa99}*{Lem.~8.1}] \label{lem:Wzero}
Suppose that for each~$m$, there exist numbers~$u(m)$ and~$v(m)$
such that for all $u \ge u(m)$ and $v \ge v(m)$, the differential
\[
d^{u+v} \: E^{u+v}_{u,m-u} \longto E^{u+v}_{-v,m+v-1}
\]
vanishes.  Then $W = 0$.
\end{lemma}

\begin{remark}
If for some fixed~$r$ the $E^r$-page of the spectral sequence is bounded
from the side of entering differentials, in the sense that for each~$m$
there is a number~$u(m)$ such that $E^r_{u,m-u} = 0$ for all $u \ge u(m)$,
then Boardman's vanishing criterion is satisfied with $v(m) = r - u(m)$.
Hence $W = 0$ in these cases.

Alternatively, if the spectral sequence collapses at the $E^r$-page,
so that $d^r$ and all later differentials are zero,
then Boardman's vanishing criterion is satisfied with~$u(m) = r$ and
$v(m) = 0$.  Thus $W = 0$ also in these cases.
\end{remark}

\begin{theorem} \label{thm:checkErstrongiffREWzero}\index{Greenlees--May--Tate spectral sequence!strong convergence}
The spectral sequence 
\[
\check E^r(X) = E^r(GM_\star(X))
\]
converges
strongly to $A_\infty(GM_\star(X)) \cong \pi_*(X^{tG})$ if and only if
$RE^\infty = 0$ and $W = 0$ for this spectral sequence.
\end{theorem}

\begin{proof}
We saw that $\check E^r(X)$ is conditionally convergent in
Lemma~\ref{lem:checkEXcondconv}.  Hence the statements `$RE^\infty =
0$ and $W = 0$' and `the spectral sequence is strongly convergent'
are equivalent by \cite{Boa99}*{Thm.~8.10}.
\end{proof}

\begin{theorem} \label{thm:strongimpliesstrong}\index{Hesselholt--Madsen--Tate spectral sequence!strong convergence}
If the Greenlees--May spectral sequence
\[
\check E^r(X) = E^r(GM_\star(X))
	\Longrightarrow \pi_*(X^{tG})
\]
is strongly convergent, then the Hesselholt--Madsen spectral sequence 
\[
\hat E^r(X) = E^r(HM_\star(X))
	\Longrightarrow \pi_*(X^{tG})
\]
is strongly convergent, as well. Moreover, $F_s A_\infty(GM_\star(X)) = F_s
A_\infty(HM_\star(X))$ for all integers $s$.
\end{theorem}

\begin{proof}
We assume $\check E^r(X)$ is strongly convergent. Explicitly, this means that the exhaustive filtration $(F_s A_\infty(GM_\star(X)))_s$ of
$A_\infty(GM_\star(X)) \cong \pi_*(X^{tG})$ is
complete Hausdorff, and the left hand monomorphism~$\beta$
in the commutative square 
\[
\xymatrix{
\ds \frac{F_s A_\infty(GM_\star(X))}{F_{s-1} A_\infty(GM_\star(X))}
\ar[r]^-{\bar\alpha_s} \ar[d]_-{\beta}^-{\cong}
& \ds \frac{F_s A_\infty(HM_\star(X))}{F_{s-1} A_\infty(HM_\star(X))}
\ar@{ >->}[d]_-{\beta} \\
E^\infty_s(GM_\star(X)) \ar[r]^-{\alpha^\infty_s}_-{\cong}
& E^\infty_s(HM_\star(X))
}
\]
is an isomorphism. It follows that the right hand monomorphism $\beta$ is also an isomorphism. Since the filtration $(F_s A_\infty(HM_\star(X)))_s$
is exhaustive, this means that~$\hat E^r(X)$ converges weakly to
$\pi_*(X^{tG})$. It also follows that the upper homomorphism $\bar\alpha_s$
is an isomorphism. By induction, this implies that the map of filtration quotients
\[
\frac{F_t A_\infty(GM_\star(X))}{F_s A_\infty(GM_\star(X))}
\overset{\cong}\longto
\frac{F_t A_\infty(HM_\star(X))}{F_s A_\infty(HM_\star(X))}
\]
is an isomorphism for all integers $s \le t$.

Passing to colimits over~$t$, and using the fact that
\[
\alpha_\infty \: A_\infty(GM_\star(X))
\overset{\cong}\longto
A_\infty(HM_\star(X))
\]
is an isomorphism by Lemma~\ref{lem:alphainftyiso}, we deduce that
\[
\alpha_s \: F_s A_\infty(GM_\star(X)) \to F_s A_\infty(HM_\star(X))
\]
is an isomorphism, for each $s \in \bZ$.  The filtration $(F_s
A_\infty(HM_\star(X)))_s$ is therefore complete and Hausdorff, meaning
that $\hat E^r(X)$ converges strongly to~$\pi_*(X^{tG})$.
\end{proof}

Combining these results we obtain the following theorem, which often
compensates for the problem that we do not \emph{a priori} know
when $\hat E^r(X) = E^r(HM_\star(X))$ is conditionally convergent,
cf.~Remark~\ref{rmk:noweakconvergence}.

\begin{theorem} \label{thm:GTatestrongconv}
If $RE^\infty = 0$ and Boardman's vanishing criterion for~$W$ from
Lemma~\ref{lem:Wzero} is satisfied for the $G$-Tate spectral sequence
$\hat E^r(X) = E^r(HM_\star(X))$, then this spectral sequence converges
strongly and conditionally to $A_\infty(HM_\star(X)) \cong \pi_*(X^{tG})$.
\end{theorem}

Note that we are not just assuming that $W = 0$ for $\hat E^r(X)$,
but that this group vanishes for the reason given by Boardman's criterion.

\begin{proof}
Since $\check E^r(X) \cong \hat E^r(X)$ for $r\ge2$, the vanishing of
$RE^\infty$ for $\hat E^r(X)$ implies the vanishing of~$RE^\infty$
for $\check E^r(X)$.  Furthermore, the hypothesis of Boardman's
criterion for $\hat E^r(X)$ implies the same hypothesis for
$\check E^r(X)$.  Hence $\check E^r(X)$ converges strongly by
Theorem~\ref{thm:checkErstrongiffREWzero}, which implies that $\hat
E^r(X)$ converges strongly by Theorem~\ref{thm:strongimpliesstrong}.
By \cite{Boa99}*{Thm.~8.10} strong convergence and the vanishing
of~$RE^\infty$ and~$W$ imply conditional convergence.
\end{proof}

\section{Summary: The $\bT$-Tate spectral sequence}
\label{subsec:TTatesummary}

The main example we had in mind when writing this memoir was $G = \bT$. Note that when discussing \emph{the $\bT$-Tate spectral sequence for a $\bT$-spectrum $X$} one could really refer to at least\footnote{There is also at least one more Tate spectral sequence, namely the one arising from a Postnikov or Whitehead tower of $X$; see Remark~\ref{rmk:fixpointPostnikov} and~Remark~\ref{rem:GM-vs-Whitehead}.} two different spectral sequences: one arising from the Greenlees--May filtration and one from the Hesselholt--Madsen filtration. The first has better convergence properties, while the latter has better multiplicative properties. Fortunately there are quite good comparison results between the two, as covered in Section~\ref{subsec:convergence}.

Let us start by summarizing the additive results regarding the
Greenlees--May and Hesselholt--Madsen versions of the $\bT$-Tate
spectral sequence.  We work over $R = \bS$, and write $\otimes$
for $\otimes_{\bS_*}$.  We first note that by virtue of $X$ being a
$\bT$-spectrum, there is an action
\[
\gamma \: X \wedge \bT_+ \cong X \wedge \bS[\bT] \longto X
\]
which makes $X$ into a right module over the spherical group ring $\bS[\bT]$. The induced pairing 
\[
\gamma_* \: \pi_*(X) \otimes \bS[\bT]_* \longto \pi_*(X)
\]
on homotopy groups then gives $\pi_*(X)$ the structure of a right module over the Hopf algebra
\[
\bS[\bT]_* \cong \bS_*[s]/(s^2=\eta s)  \,, \quad |s| = 1 \,.
\]
Here $\eta$ is the image of the complex Hopf map in $\pi_1(\bS) \cong \bZ/2$. Note that this Hopf algebra is finitely generated and projective over $\bS_*$. We denote the image $\gamma_*(x \otimes s)$ by $x s$. 

There is a minimal projective $\bS[\bT]_*$-module resolution $P_*$
of $\bS_*$, with $P_k = \bS[\bT]\{p_k\}$ and $\partial(p_k) = p_{k-1}
(s + (k-1) \eta)$.  Let $\widetilde{P}_*$ be the mapping cone of the
augmentation $\epsilon \: P_* \to \bS_*$. Let the complete resolution
$\hat P_*$ be the fibre product $P_* \times_{\bS_*} D(\widetilde{P}_*)$,
which is obtained by splicing $P_*$ with its dual.

\begin{theorem}[Greenlees--May--Tate spectral sequence]\index{Greenlees--May--Tate spectral sequence!for the circle group $\bT$}
Given an orthogonal $\bT$-spectrum $X$, there is a filtration $GM_\star(X)$ of orthogonal $\bT$-spectra, and an associated~$\bS_*$-module spectral sequence
\[
\check{E}^r(X) = E^r(GM_\star(X))
\]
with abutment
\[
A_\infty(GM_\star(X)) \cong \pi_*(X^{t\bT})
\]
filtered by the images $\im(\pi_* GM_s(X) \to \pi_*(X^{t\bT}))$.
We refer to this spectral sequence as the Greenlees--May $\bT$-Tate spectral sequence for $X$. The following hold:
\begin{description}
\item[$E^1$-page] The $E^1$-page of the Greenlees--May $\bT$-Tate spectral sequence can be written 
\[
E^1_{*,*}(GM_\star(X)) \cong \Hom_{\bS[\bT]_*}(\hat{P}_*,\pi_*(X))
\]
where $\hat{P}_*$ is a complete resolution of $\bS_*$ as a trivial $\bS[\bT]_*$-module.  For the minimal such resolution we can write
$$
E^1_{*,*}(GM_\star(X)) \cong \pi_*(X) [t, t^{-1}]
$$
with $t$ in bidegree~$(-1,-1)$, and then $d^1(t^c \cdot x) = t^{c+1} \cdot x (s + c \eta)$ for all $c \in \bZ$ and $x \in \pi_*(X)$.
\item[Convergence] The Greenlees--May spectral sequence converges conditionally to the abutment. It converges strongly to the abutment if and only if the derived $E^\infty$-page $RE^\infty$ and Boardman's whole plane obstruction group $W$ are both trivial.
\end{description}
\begin{proof}
The first statement is Lemma~\ref{lem:checkE1isgmcx} combined with Proposition~\ref{prop:gmishomhatP} and~Proposition~\ref{prop:d1-in-GM-sp-seq}. The second statement is Lemma~\ref{lem:checkEXcondconv} combined with Theorem~\ref{thm:checkErstrongiffREWzero}.
\end{proof}
\end{theorem}

\begin{theorem}[Hesselholt--Madsen--Tate spectral sequence]\index{Hesselholt--Madsen--Tate spectral sequence!for the circle group $\bT$}
\label{thm:main3}
Given an orthogonal $\bT$-spectrum $X$, there is a filtration $HM_\star(X)$ of orthogonal $\bT$-spectra, and an associated $\bS_*$-module spectral sequence
\[
\hat{E}^r(X) = E^r(HM_\star(X))
\]
with abutment
\[
A_\infty(HM_\star(X)) \cong \pi_*(X^{t\bT})
\]
filtered by the images $\im(\pi_* HM_\star (X) \to \pi_*(X^{t\bT}))$.
We refer to this spectral sequence as the Hesselholt--Madsen $\bT$-Tate spectral sequence for $X$. The following hold:
\begin{description}
\item[$E^1$-page] The $E^1$-page of the Hesselholt--Madsen $\bT$-Tate spectral sequence can be written  
\[
E^1_{*,*}(HM_\star(X)) \cong \Hom_{\bS[\bT]_*}(\bS_*,\widetilde{P}_* \otimes \Hom(P_{*},\pi_*(X)))
\]
where $P_*$ is a projective resolution of $\bS_*$ as a trivial $\bS[\bT]_*$-module and $\widetilde{P}_*$ denotes the mapping cone of the augmentation $\epsilon \: P_* \to \bS_*$.
\item[$E^2$-page] The $E^2$-page is given in terms of Hopf algebra
Tate cohomology, alias complete $\Ext$, as
$$
E^2_{*,*}(HM_\star(X)) \cong \hatExt^{-*}_{\bS[\bT]_*}(\bS_*, \pi_*(X)) \,.
$$
\item[Convergence] If the Greenlees--May $\bT$-Tate spectral sequence for~$X$ converges strongly, then the Hesselholt--Madsen $\bT$-Tate spectral sequence for the same spectrum also converges strongly. Moreover, the two associated filtrations of $\pi_*(X^{t\bT})$ agree. 
\end{description}
\end{theorem}

\begin{proof}
The first statement is Proposition~\ref{prop:TateE1isTatecx}, the second is Theorem~\ref{thm:TateE2isTateExt}, and the third statement is Theorem~\ref{thm:strongimpliesstrong}.
\end{proof}

Worth pointing out is that the Greenlees--May and the Hesselholt--Madsen versions of the $\bT$-Tate spectral sequence are isomorphic from the $E^2$-page and on, per Proposition~\ref{prop:GMHMfiltrationiso}. In particular, the $E^2$-page of both spectral sequences is given by
\begin{align*}
\hat E^2_{-c,*}(X) \cong \check E^2_{-c,*}(X)
&\cong \hatExt_{\bS[\bT]_*}^c(\bS_*,\pi_*(X)) \\
&\cong \begin{cases}
{\displaystyle
	\frac{\ker(s \: \pi_*(X) \to \pi_{*+1}(X))}{\im(s+\eta \: \pi_{*-1}(X) \to \pi_*(X))}} & \text{for $c$ even,} \\
\\
{\displaystyle
	\frac{\ker(s+\eta \: \pi_*(X) \to \pi_{*+1}(X))}{\im(s \: \pi_{*-1}(X) \to \pi_*(X))}} & \text{for $c$ odd,} 
\end{cases}
\end{align*}
where the last isomorphism is the result of the computation of Section~\ref{subsec:computation}.

Regarding convergence, we note that Lemma~\ref{lem:Wzero} gives a criterion, internal to the spectral sequence itself, for when Boardman's whole-plane obstruction vanishes. In particular, if $X$ is bounded below, either version of the Tate spectral sequence is a half-plane spectral sequence with entering differentials
(at least from the $E^2$-page), which guarantees this. In the applications we have in mind, we are in this situation if we consider topological Hochschild homology\index{topological Hochschild homology} $X = \THH(B)$ for some connective orthogonal ring spectrum $B$.  

Let us now summarise the multiplicative structure of the two spectral sequences discussed.

\begin{theorem}
\leavevmode
\begin{description}
\item[Multiplicativity] The Hesselholt--Madsen $\bT$-Tate spectral sequence is multiplicative, in the sense that any pairing $\phi \: X \wedge Y \to Z$ of orthogonal $\bT$-spectra gives rise to a pairing $\phi \: (\hat{E}^*(X),\hat{E}^*(Y)) \to \hat{E}^*(Z)$ of the associated spectral sequences. Explicitly, $\phi$ gives rise to homomorphisms
\[
\phi^r \: \hat E^r(X_\star) \otimes \hat E^r(Y_\star) \longto \hat E^r(Z_\star)
\]
for all $r \geq 1$, such that:
\begin{enumerate}
\item The Leibniz rule
\[
d^r \phi^r = \phi^r(d^r \otimes 1) + \phi^r(1 \otimes d^r)
\]
holds as an equality of homomorphisms
\[
\hat E^r_i(X) \otimes \hat E^r_j(Y) \longto \hat E^r_{i+j-r}(Z)
\]
for all $i,j \in \bZ$ and $r\ge1$.
\item The diagram 
\[
\begin{tikzcd}
\hat E^{r+1}(X) \otimes \hat E^{r+1}(Y) \arrow[r,"\phi^{r+1}"] \arrow[d] & \hat E^{r+1}(Z) \arrow[d,"\cong"] \\
H(\hat E^{r}(X) \otimes \hat E^{r}(Y)) \arrow[r,"H(\phi^{r})"] & H(\hat E^{r}(Z))
\end{tikzcd}
\]
commutes for all $r \geq 1$.
\item
The pairing
\[
\phi^2 \: \hat E^2(X_\star) \otimes \hat E^2(Y_\star) \longto \hat E^2(Z_\star)
\]
agrees with the cup product
\[
\smile \: \hatExt_{\bS[\bT]_*}^{-i,*}(\bS_*, \pi_*(X))
\otimes \hatExt_{\bS[\bT]_*}^{-j,*}(\bS_*, \pi_*(Y))
\longto \hatExt_{\bS[\bT]_*}^{-i-j,*}(\bS_*, \pi_*(Z)) \,.
\]
\end{enumerate} 
For $r\ge2$ the same statements hold for the Greenlees--May spectral sequence, with $\check E^r$ in place of $\hat E^r$.

\item[Multiplicative abutment] We have an induced pairing 
\[
\phi_* \: \pi_*(X^{t\bT}) \otimes \pi_*(Y^{t\bT}) \longto \pi_*(Z^{t\bT})
\]
of abutments with the Hesselholt--Madsen filtrations, which is compatible with the pairing $\phi^\infty$ of $E^\infty$-pages. Explicitly, the diagram
\[
\xymatrix{
\ds \frac{F_i \pi_*(X^{t\bT})}{F_{i-1} \pi_*(X^{t\bT})}
\otimes \frac{F_j \pi_*(Y^{t\bT})}{F_{j-1} \pi_*(Y^{t\bT})}
\ar[r]^-{\bar\phi_*} \ar[d]_-{\beta \otimes \beta}
& \ds \frac{F_{i+j} \pi_*(Z^{t\bT})}{F_{i+j-1} \pi_*(Z^{t\bT})}
	\ar@{ >->}[d]^-{\beta} \\
\hat E^\infty_i(X) \otimes \hat E^\infty_j(Y) \ar[r]^-{\phi^\infty}
& \hat E^\infty_{i+j}(Z)
}
\]
commutes, for all $i,j \in \bZ$. 

If the Greenlees--May spectral sequence is strongly convergent, then the same statements hold for the Greenlees--May filtrations and $\check E^\infty_{*,*}$.
\end{description}
\begin{proof}
For the Hesselholt--Madsen $\bT$-Tate spectral sequence this is Theorem~\ref{thm:GTatespseqismult} and Theorem~\ref{thm:E2Tatepairingiscup}.  The statements about multiplicativity of the $E^r$-pages and $d^r$-differentials can be transported to the Greenlees--May spectral sequence for $r\ge2$ by way of the isomorphism of Proposition~\ref{prop:GMHMfiltrationiso}.  The statements about multiplicativity of filtered abutments carry over to the Greenlees--May spectral sequence when~$GM_\star(X)$ and $HM_\star(X)$ induce the same filtration on $\pi_*(X^{t\bT})$, which holds under the hypothesis of strong convergence by Theorem~\ref{thm:strongimpliesstrong}.
\end{proof}
\end{theorem}

Recalling the discussion of Remark~\ref{rmk:DGA}, in the context of the circle group, the Hopf algebra Tate cohomology can also be described as the homology of the differential graded $\bS[\bT]_*$-module
\[
\pi_*(X)[t, t^{-1}] \,, \quad |t| = -1
\]
with differential characterised by 
\[
d(x) = txs \och d(t) = t^2 \eta \,.
\]
Moreover, given a pairing $X \wedge Y \to Z$ we have an induced pairing $\pi_*(X) \otimes \pi_*(Y) \to \pi_*(Z)$ on homotopy groups, and the cup product on Tate cohomology is precisely the one induced by the obvious map
\[
\pi_*(X)[t,t^{-1}] \otimes \pi_*(Y)[t,t^{-1}] \longto \pi_*(Z)[t,t^{-1}]
\]
on homology. By Theorem~\ref{thm:E2Tatepairingiscup}, the
multiplicative structure on the second page of (both versions of)
the $\bT$-Tate spectral sequence corresponds to this cup product.
By Proposition~\ref{prop:d1-in-GM-sp-seq} we can formally impose this
algebra structure on the Greenlees--May $E^1$-page, in which case
the $d^1$-differential is a derivation, and this lets us extend the
multiplicativity statement for the Greenlees--May $\bT$-Tate spectral
sequence to the range $r\ge1$, in place of $r\ge2$, even if there is no
underlying topological source of the pairing of $E^1$-pages.

\begin{remark}
\label{rmk:Blumberg-Mandell}
Blumberg and Mandell also set up $\bT$-Tate spectral sequences in line with Greenlees--May and Hesselholt--Madsen in~\cite{BM17}; let us elaborate on how our spectral sequences compare to theirs. While the homotopy theoretical technicalities are resolved in a different manner to what we have done in this memoir, the main ideas are the same: their Hesselholt--Madsen filtration~\cite{BM17}*{Section 12} is essentially the same as ours. It is worth noting that Blumberg and Mandell prove that conditional convergence for the Hesselholt--Madsen $\bT$-Tate spectral sequence always holds~\cite{BM17}*{Lemma 3.16, p. 38}, something that we have not proved in this memoir. We believe that a similar argument to theirs works to show that the Hesselholt--Madsen $G$-Tate spectral sequence is conditionally convergent for all compact Lie groups $G$, but refrain from making any definite statement before the details are checked.\index{Hesselholt--Madsen--Tate spectral sequence!conditional convergence} 

As mentioned in Section~\ref{subsec:filt}, some form of cofibrant replacement of maps is necessary to solve homotopy theoretical difficulties when dealing with sequences of spectra. Blumberg--Mandell deal with these by referring to model structures on such categories, with a focus on $h$-cofibrations~\cite{BM17}*{Section 6}\index{h-cofibration@$h$-cofibration}. While such model structures do exist, and we also started out trying to solve technicalities in such a manner, we were ultimately unable to locate a reference for why such a model structure is monoidal, which is needed for multiplicative structures on spectral sequences. This was one of the main reasons that we ended up choosing to work with explicit models for our homotopy colimits (with convenient monoidal properties), as well as with strong $h$-cofibrations (so that we can use Theorem~\ref{thm:pushprodstrongcof}), rather than referring to the framework of model categories. 
\end{remark}

\appendix

\backmatter
\bibliographystyle{amsalpha}

\printindex

\end{document}